\documentclass[leqno]{article}

\def\diam{\mathop{\rm diam}}
\def\re{\mathop{\rm Re}}
\def\im{\mathop{\rm Im}}

\newtheorem{theorem}{Theorem}
\newtheorem{lemma}[theorem]{Lemma}
\newtheorem{proposition}[theorem]{Proposition}
\newtheorem{definition}[theorem]{Definition}
\newtheorem{corollary}[theorem]{Corollary}

\newcommand{\begintheorem}{\addtocounter{equation}{1}\begin{theorem}}
\newcommand{\beginlemma}{\addtocounter{equation}{1}\begin{lemma}}
\newcommand{\beginproposition}{\addtocounter{equation}{1}\begin{proposition}}
\newcommand{\begindefinition}{\addtocounter{equation}{1}\begin{definition}}
\newcommand{\begincorollary}{\addtocounter{equation}{1}\begin{corollary}}

\begin{document}

\title{Sums, rearrangements, and norms}

\author{Stephen Semmes\\
        Rice University}

\date{}

\maketitle

\begin{abstract}
These informal notes deal with a number of questions related to sums
and integrals in analysis.
\end{abstract}

\tableofcontents

\part{Basic notions}

\section{Real and complex numbers}
\label{real and complex numbers}
\setcounter{equation}{0}

        Of course, the real numbers ${\bf R}$ are contained in the
complex numbers ${\bf C}$, and every $z \in {\bf C}$ can be expressed
as $z = x + y \, i$, where $x, y \in {\bf R}$ and $i^2 = -1$.  In this
case, $x$ and $y$ are called the real and imaginary parts of $z$,
respectively.  The \emph{complex conjugate} $\overline{z}$ of $z$ is
given by
\begin{equation}
        \overline{z} = x - y \, i.
\end{equation}
It is easy to see that
\begin{equation}
        \overline{z + w} = \overline{z} + \overline{w}
\end{equation}
and
\begin{equation}
        \overline{z \, w} = \overline{z} \, \overline{w}
\end{equation}
for every $z, w \in {\bf C}$.  The \emph{modulus} $|z|$ of $z$ is given by
\begin{equation}
        |z| = (x^2 + y^2)^{1/2}.
\end{equation}
Thus
\begin{equation}
        |z|^2 = z \, \overline{z}.
\end{equation}
This implies that
\begin{equation}
        |z \, w|^2 = (z \, w) \, \overline{z \, w}
                   = z \, w \, \overline{z} \, \overline{w} = |z|^2 \, |w|^2
\end{equation}
for every $z, w \in {\bf C}$, and hence
\begin{equation}
        |z \, w| = |z| \, |w|.
\end{equation}
Note that the modulus of a real number is the same as its absolute
value, and that the modulus of $z = x + y \, i \in {\bf C}$ is the
same as the Euclidean norm of $(x, y) \in {\bf R}^2$.

\section{Rearrangements}
\label{rearrangements}
\setcounter{equation}{0}

        Let $\sum_{j = 1}^\infty a_j$ be an infinite series of real or
complex numbers.  If $\pi$ is a one-to-one mapping from the set ${\bf
Z}_+$ of positive integers onto itself, then the series
\begin{equation}
\label{sum_{j = 1}^infty a_{pi(j)}}
        \sum_{j = 1}^\infty a_{\pi(j)}
\end{equation}
is said to be a \emph{rearrangement} of $\sum_{j = 1}^\infty a_j$.

        Remember that $\sum_{j = 1}^\infty a_j$ converges if the
sequence of partial sums $\sum_{j = 1}^n a_j$ converges as $n \to
\infty$.  If $a_j$ is a nonnegative real number for each $j$, then the
partial sums are monotone increasing, and convergence is equivalent to
boundedness of the partial sums.  In this case, convergence of
$\sum_{j = 1}^\infty a_j$ implies the convergence of every
rearrangement (\ref{sum_{j = 1}^infty a_{pi(j)}}), and the values of
these sums are the same.  More precisely,
\begin{equation}
        \sum_{j = 1}^n a_{\pi(j)} \le \sum_{j = 1}^N a_j
\end{equation}
when $\pi(1), \ldots, \pi(n) \le N$, so that the boundedness of the
partial sums of $\sum_{j = 1}^\infty a_j$ implies the boundedness
of the partial sums of (\ref{sum_{j = 1}^infty a_{pi(j)}}).  Similarly,
\begin{equation}
        \sum_{j = 1}^n a_j \le \sum_{j = 1}^N a_{\pi(j)}
\end{equation}
when $\pi^{-1}(1), \ldots, \pi^{-1}(n) \le N$, and these two simple
extimates imply that the suprema of the partial sums of $\sum_{j
=1}^\infty a_j$ and (\ref{sum_{j = 1}^infty a_{pi(j)}}) are the same.

        An infinite series $\sum_{j = 1}^\infty a_j$ of real or
complex numbers is said to converge absolutely if $\sum_{j = 1}^\infty
|a_j|$ converges.  It is well known that absolute convergence implies
convergence, by the Cauchy criterion.  If $\sum_{j = 1}^\infty a_j$
converges absolutely, then the preceding discussion implies that
(\ref{sum_{j = 1}^infty a_{pi(j)}}) also converges absolutely, and one
can show that the two sums have the same value.  This is trivial when
$a_j = 0$ for all but finitely many $j$, and otherwise $\sum_{j =
1}^\infty a_j$ can be approximated by series with this property.
Alternatively, $\sum_{j = 1}^\infty a_j$ may be expressed as a linear
combination of convergent series whose terms are nonnegative real
numbers, so that the equality of the sums reduces to the previous
case.

\section{Generalized convergence}
\label{generalized convergence}
\setcounter{equation}{0}

        Let $E$ be a nonempty set, and let $f(x)$ be a real or
complex-valued function on $E$.  Let us say that $\sum_{x \in E} f(x)$
converges in the generalized sense if there is a $\lambda \in {\bf R}$
or ${\bf C}$, as appropriate, such that for each $\epsilon > 0$
there is a finite set $A_\epsilon \subseteq E$ for which
\begin{equation}
\label{|sum_{x in B} f(x) - lambda| < epsilon}
        \biggl|\sum_{x \in B} f(x) - \lambda\biggr| < \epsilon
\end{equation}
whenever $B \subseteq E$ is a finite set that satisfies $A_\epsilon
\subseteq B$.  It is easy to see that such a $\lambda$ is unique when
it exists, in which case $\sum_{x \in E} f(x)$ is defined to be
$\lambda$.

        If $f(x)$ has this property and $\pi$ is a one-to-one mapping
of $E$ onto itself, then $f(\pi(x))$ has the same property, and
\begin{equation}
        \sum_{x \in E} f(\pi(x)) = \sum_{x \in E} f(x).
\end{equation}
This follows from the fact that
\begin{equation}
        \sum_{x \in A} f(\pi(x)) = \sum_{x \in \pi(A)} f(x)
\end{equation}
for every finite set $A \subseteq E$.  Thus this definition of
$\sum_{x \in E} f(x)$ is automatically invariant under rearrangements.

        Suppose that $f(x)$ is a nonnegative real number for each $x
\in E$.  If the partial sums $\sum_{x \in A} f(x)$ over finite subsets
$A$ of $E$ are uniformly bounded, then $\sum_{x \in E} f(x)$ converges
in the generalized sense, and
\begin{equation}
 \quad  \sum_{x \in E} f(x) = \sup \bigg\{\sum_{x \in A} f(x) : A \subseteq E
                              \hbox{ has only finitely many elements}\bigg\}.
\end{equation}
If $f(x)$ is a real or complex-valued function on $E$ such that the
sums $\sum_{x \in A} |f(x)|$ over finite sets $A \subseteq E$ are
bounded, then $\sum_{x \in E} f(x)$ also converges in the generalized
sense.  This follows by expressing $f(x)$ as a linear combination of
nonnegative real-valued functions for which the partial sums over
finite subsets of $E$ are bounded.

        Conversely, if $f(x)$ is a real or complex-valued function on
$E$ such that $\sum_{x \in E} f(x)$ converges in the generalized
sense, then the sums $\sum_{x \in A} |f(x)|$ over finite subsets $A$
of $E$ are uniformly bounded.  To see this, one can take $\epsilon =
1$ in the definition of convergence to get a finite set $A_1 \subseteq
E$ for which the partial sums $\sum_{x \in B} f(x)$ over finite
subsets $B$ of $E$ with $A_1 \subseteq B$ are uniformly bounded.  This
implies that the partial sums $\sum_{x \in A} f(x)$ over arbitrary
finite sets $A \subseteq E$ are bounded, by taking $B = A \cup A_1$,
and using the fact that the sums over subsets of $A_1$ are bounded.
The boundedness of the partial sums of $|f(x)|$ can then be obtained
by applying this to finite sets $A \subseteq E$ on which $f(x)$ is
positive or negative in the real case, or on which the real or
imaginary parts of $f(x)$ are positive or negative in the complex
case.

\section{Nets}
\label{nets}
\setcounter{equation}{0}

        A partially ordered set $(A, \prec)$ is said to be a
\emph{directed system} if for every $a, b \in A$ there is a $c \in A$
such that $a, b \prec c$.  A \emph{net} $\{x_a\}_{a \in A}$ indexed by
$A$ assigns to each $a \in A$ an element $x_a$ of a set $X$.  If $X$
is a topological space, then the net $\{x_a\}_{a \in A}$ converges to
$x \in X$ if for every open set $U \subseteq X$ with $x \in U$ there
is an $a \in A$ such that $x_b \in U$ when $b \in A$ and $a \prec b$.
This reduces to the usual definition of convergence of a sequence when
$A$ is the set of positive integers with the standard ordering.  Now
let $E$ be a nonempty set, and let $f(x)$ be a real or complex-valued
function on $E$.  The collection of nonempty finite subsets of $E$ is
partially ordered by inclusion, and defines a directed system.  More
precisely, any two finite subsets of $E$ is contained in their union,
which is also a finite subset of $E$.  Consider the net associated to
this directed system that assigns to each nonempty finite set $B
\subseteq E$ the real or complex number $\sum_{x \in B} f(x)$.  It is
easy to see that convergence of this net in ${\bf R}$ or ${\bf C}$, as
appropriate, is the same as convergence of $\sum_{x \in E} f(x)$ in
the sense described in the previous section.

\section{Norms on vector spaces}
\label{norms on vector spaces}
\setcounter{equation}{0}

        Let $V$ be a vector space over the real or complex numbers.
A \emph{norm} on $V$ is a nonnegative real-valued function $\|v\|$
defined for $v \in V$ such that $\|v\| = 0$ if and only if $v = 0$,
\begin{equation}
\label{||t v|| = |t| ||v||}
        \|t \, v\| = |t| \, \|v\|
\end{equation}
for every $v \in V$ and $t \in {\bf R}$ or ${\bf C}$, as appropriate, and
\begin{equation}
\label{||v + w|| le ||v|| + ||w||}
        \|v + w\| \le \|v\| + \|w\|
\end{equation}
for every $v, w \in V$.

        A set $E \subseteq V$ is said to be \emph{convex} if for every
$v, w \in E$ and $t \in {\bf R}$ with $0 < t < 1$,
\begin{equation}
        t \, v + (1 - t) \, w \in E.
\end{equation}
If $\|v\|$ is a norm on $V$ and
\begin{equation}
\label{B_1 = {v in V : ||v|| le 1}}
        B_1 = \{v \in V : \|v\| \le 1\}
\end{equation}
is the corresponding closed unit ball, then it is easy to see that
$B_1$ is a convex set in $V$.

        Conversely, suppose that $\|v\|$ is a nonnegative real-valued
function on $V$ that satisfies the positivity condition $\|v\| > 0$
when $v \ne 0$ and the homogeneity condition (\ref{||t v|| = |t|
||v||}).  If $B_1$ is convex, then one can show that $\|v\|$ satisfies
the triangle inequality (\ref{||v + w|| le ||v|| + ||w||}), and hence
that $\|v\|$ is a norm.  To see this, let $v, w \in V$ be given, with
$v, w \ne 0$, since otherwise (\ref{||v + w|| le ||v|| + ||w||}) is
trivial.  Put
\begin{equation}
        v' = \frac{v}{\|v\|}, \quad w' = \frac{w}{\|w\|},
\end{equation}
so that $\|v'\| = \|w'\| = 1$.  Thus $v', w' \in B_1$, and hence
\begin{equation}
        \|t \, v' + (1 - t) \, w'\| \le 1
\end{equation}
when $t \in {\bf R}$ and $0 < t < 1$, by hypothesis.  If $t = \|v\| /
(\|v\| + \|w\|)$, then $1 - t = \|w\| / (\|v\| + \|w\|)$, and
\begin{equation}
        t \, v' + (1 - t) \, w' = \frac{v + w}{\|v\| + \|w\|}.
\end{equation}
Therefore
\begin{equation}
        \biggl\|\frac{v + w}{\|v\| + \|w\|}\biggr\| \le 1,
\end{equation}
which implies (\ref{||v + w|| le ||v|| + ||w||}), as desired.

\section{Bounded functions}
\label{bounded functions}
\setcounter{equation}{0}

        Let $E$ be a nonempty set, and consider the spaces
$\ell^\infty(E, {\bf R})$, $\ell^\infty(E, {\bf C})$ of real or
complex-valued functions on $E$ that are bounded.  It is sometimes
convenient to use the notation $\ell^\infty(E)$ to refer to either of
these spaces, which are vector spaces with respect to pointwise
addition and scalar multiplication.  The supremum or $\ell^\infty$
norm is defined as usual by
\begin{equation}
        \|f\|_\infty = \sup \{|f(x)| : x \in E\}.
\end{equation}
It is easy to see that this is a norm on $\ell^\infty(E)$, because of
the triangle inequality for the ordinary absolute value on ${\bf R}$
or modulus on ${\bf C}$.

\section{Summable functions}
\label{summable functions}
\setcounter{equation}{0}

        A real or complex-valued function $f(x)$ on a nonempty set $E$
is said to be \emph{summable} if the partial sums $\sum_{x \in A}
|f(x)|$ over nonempty finite subsets $A$ of $E$ are uniformly bounded.
This is equivalent to the convergence of $\sum_{x \in E} |f(x)|$ in
the sense of Section \ref{generalized convergence}, whose value is
equal to the supremum of $\sum_{x \in A} |f(x)|$ over all nonempty
finite sets $A \subseteq E$.  Let $\ell^1(E, {\bf R})$, $\ell^1(E,
{\bf C})$ be the spaces of summable real or complex-valued functions
on $E$, respectively, which may also be denoted by $\ell^1(E)$ to
include both cases at the same time.  It is easy to see that these are
vector spaces with respect to pointwise addition and scalar
multiplication, and that
\begin{equation}
        \|f\|_1 = \sum_{x \in E} |f(x)|
\end{equation}
defines a norm on these spaces.

\section{$p$-Summable functions}
\label{p-summable functions}
\setcounter{equation}{0}

        Let $f(x)$ be a real or complex-valued function on a nonempty
set $E$, and let $p$ be a positive real number.  If $|f(x)|^p$ is a
summable function on $E$, then we say that $f(x)$ is
\emph{$p$-summable} on $E$.  The spaces of real or complex-valued
$p$-summable functions on $E$ are denoted $\ell^p(E, {\bf R})$,
$\ell^p(E, {\bf C})$, respectively, or simply $\ell^p(E)$ to include
both cases at the same time.  One can check that these are vector
spaces over the real or complex numbers, as appropriate, with respect
to pointwise addition and scalar multiplication of functions.

        If $f$ is a $p$-summable function on $E$, then put
\begin{equation}
        \|f\|_p = \Big(\sum_{x \in E} |f(x)|^p\Big)^{1/p}.
\end{equation}
This satisfies the positivity and homogeneity properties of a norm on
$\ell^p(E)$ for every $p > 0$.  Let us check that this is a norm on
$\ell^p(E)$ when $p \ge 1$.  As in Section \ref{norms on vector
spaces}, it suffices to show that the closed unit ball in $\ell^p(E)$
associated to $\|f\|_p$ is convex when $p \ge 1$.  Equivalently, if
$f$, $g$ are $p$-summable functions on $E$ such that $\|f\|_p, \|g\|_p
\le 1$, then we would like to check that
\begin{equation}
        \|t \, f + (1 - t) \, g\|_p \le 1
\end{equation}
when $t \in {\bf R}$ and $0 < t < 1$.  The main point is that
\begin{eqnarray}
 |t \, f(x) + (1 - t) \, g(x)|^p & \le & (t \, |f(x)| + (1 - t) \, |g(x)|)^p \\
                         & \le & t \, |f(x)|^p + (1 - t) \, |g(x)|^p \nonumber
\end{eqnarray}
for every $x \in E$, because of the convexity of the function
$\phi_p(r) = r^p$ on the nonnegative real numbers when $p \ge 1$.  Hence
\begin{equation}
        \sum_{x \in E} |t \, f(x) + (1 - t) \, g(x)|^p
 \le t \sum_{x \in E} |f(x)|^p + (1 - t) \sum_{x \in E} |g(x)|^p \le 1.
\end{equation}

\section{Monotonicity}
\label{monotonicity}
\setcounter{equation}{0}

        Let $p$ be a positive real number, and let $f$ be a real or
complex-valued $p$-summable function on a nonempty set $E$.  Clearly
\begin{equation}
        |f(x)| \le \|f\|_p
\end{equation}
for every $x \in E$, which implies that $f$ is bounded and satisfies
\begin{equation}
        \|f\|_\infty \le \|f\|_p.
\end{equation}
If $q \ge p$, then $f$ is also $q$-summable, because
\begin{equation}
 |f(x)|^q \le \|f\|_\infty^{q - p} \, |f(x)|^p \le \|f\|_p^{q - p} \, |f(x)|^p
\end{equation}
for every $x \in E$.  Moreover,
\begin{equation}
 \|f\|_q^q = \sum_{x \in E} |f(x)|^q
            \le \|f\|_p^{q - p} \sum_{x \in E} |f(x)|^p = \|f\|_p^q,
\end{equation}
and hence
\begin{equation}
        \|f\|_q \le \|f\|_p.
\end{equation}

        If $q = 1$, then we get that
\begin{equation}
        \Big(\sum_{x \in E} |f(x)|\Big)^p \le \sum_{x \in E} |f(x)|^p
\end{equation}
when $f$ is $p$-summable and $0 < p \le 1$.  In particular,
\begin{equation}
\label{(a + b)^p le a^p + b^p}
        (a + b)^p \le a^p + b^p
\end{equation}
for every pair of nonnegative real numbers $a$, $b$ when $0 < p \le 1$,
by applying the previous inequality to a set $E$ with exactly two elements.
Conversely, one can apply (\ref{(a + b)^p le a^p + b^p}) repeatedly to get
\begin{equation}
        \Big(\sum_{j = 1}^n a_j\Big)^p \le \sum_{j = 1}^n a_j^p
\end{equation}
for any positive integer $n$ and nonnegative real numbers $a_1, \ldots, a_n$,
which implies the analogous inequality for arbitrary sums by passing to a
suitable limit.

\section{$p$-Norms, $0 < p \le 1$}
\label{p-norms, 0 < p le 1}
\setcounter{equation}{0}

        Let $V$ be a vector space over the real or complex numbers,
and let $\|v\|$ be a nonnegative real-valued function on $V$ such that
$\|v\| > 0$ when $v \ne 0$ and
\begin{equation}
        \|t \, v \| = |t| \, \|v\|
\end{equation}
for every $v \in V$ and $t \in {\bf R}$ or ${\bf C}$, as appropriate.
We say that $\|v\|$ is a \emph{$p$-norm}, $0 < p \le 1$, if in addition
\begin{equation}
        \|v + w\|^p \le \|v\|^p + \|w\|^p
\end{equation}
for every $v, w \in V$.  This reduces to the ordinary triangle
inequality (\ref{||v + w|| le ||v|| + ||w||}) when $p = 1$, so that a
$1$-norm is the same as a norm.  For example, $\|f\|_p$ defines a
$p$-norm on $\ell^p(E)$ for any nonempty set $E$ when $0 < p \le 1$,
because of (\ref{(a + b)^p le a^p + b^p}).

        Equivalently, $\|v\|$ is a $p$-norm when
\begin{equation}
        \|v + w\| \le (\|v\|^p + \|w\|^p)^{1/p}
\end{equation}
for every $v, w \in V$.  As in the previous section, the right side of
this inequality is monotone decreasing in $p$.  Hence a $p$-norm is
also a $\widetilde{p}$-norm when $0 < \widetilde{p} \le p \le 1$.

        Let $B_1$ be the closed unit ball associated to $\|v\|$,
as in (\ref{B_1 = {v in V : ||v|| le 1}}).  If $\|v\|$ is a $p$-norm,
then
\begin{equation}
        a \, v + b \, w \in B_1
\end{equation}
whenever $v, w \in B_1$ and $a$, $b$ are nonnegative real numbers such
that $a^p + b^p \le 1$.  Conversely, let us check that this property
implies that $\|v\|$ is a $p$-norm, as in Section \ref{norms on vector
spaces}.  Let $v$, $w$ be nonzero vectors in $V$, and put $v' =
v/\|v\|$, $w' = w/\|w\|$, as before.  Also put
\begin{equation}
        a = \frac{\|v\|}{(\|v\|^p + \|w\|^p)^{1/p}}, \quad 
         b = \frac{\|w\|}{(\|v\|^p + \|w\|^p)^{1/p}}.
\end{equation}
Thus
\begin{equation}
        a^p + b^p =
    \frac{\|v\|^p}{\|v\|^p + \|w\|^p} + \frac{\|w\|^p}{\|v\|^p + \|w\|^p} = 1,
\end{equation}
and hence
\begin{equation}
        a \, v' + b \, w' = \frac{v + w}{(\|v\|^p + \|w\|^p)^{1/p}} \in B_1.
\end{equation}
This implies the $p$-norm version of the triangle inequality when $v, w \ne 0$,
and of course it is trivial when $v$ or $w$ is equal to $0$.

\section{Metric spaces}
\label{metric spaces}
\setcounter{equation}{0}

        Remember that a metric space is a set $M$ with a nonnegative
real-valued function $d(x, y)$ defined for $x, y \in M$ such that
$d(x, y) = 0$ if and only if $x = y$,
\begin{equation}
        d(y, x) = d(x, y)
\end{equation}
for every $x, y \in M$, and
\begin{equation}
        d(x, z) \le d(x, y) + d(y, z)
\end{equation}
for every $x, y, z \in M$.  If $V$ is a real or complex vector space
equipped with a norm $\|v\|$, then
\begin{equation}
\label{d(v, w) = ||v - w||}
        d(v, w) = \|v - w\|
\end{equation}
is a metric on $V$.  Similarly, if $\|v\|$ is a $p$-norm on $V$ for
some $p$, $0 < p \le 1$, then
\begin{equation}
\label{d(v, w) = ||v - w||^p}
        d(v, w) = \|v - w\|^p
\end{equation}
is a metric on $V$.

        Let $(M, d(x, y))$ be a metric space.  A sequence $\{x_j\}_{j
= 1}^\infty$ of elements of $M$ is said to \emph{converge} to $x \in
M$ if for every $\epsilon > 0$ there is an $L \ge 1$ such that
\begin{equation}
        d(x_j, x) < \epsilon
\end{equation}
for every $j \ge L$.  We say that $\{x_j\}_{j = 1}^\infty$ is a
\emph{Cauchy sequence} if for every $\epsilon > 0$ there is an $L \ge
1$ such that
\begin{equation}
        d(x_j, x_l) < \epsilon
\end{equation}
for every $j, l \ge L$.  It is easy to check that every convergent
sequence is a Cauchy sequence, and a metric space is said to be
\emph{complete} if every Cauchy sequence converges to an element of
the space.  For example, it is well known that the real and complex
numbers are complete with respect to their standard metrics.

        If $\{x_j\}_{j = 1}^\infty$ is a sequence of elements of $M$
with the property that
\begin{equation}
\label{sum_{j = 1}^infty d(x_j, x_{j + 1})}
        \sum_{j = 1}^\infty d(x_j, x_{j + 1})
\end{equation}
converges, then $\{x_j\}_{j = 1}^\infty$ is a Cauchy sequence in $M$.
This uses the triangle inequality to get that
\begin{equation}
        d(x_k, x_l) \le \sum_{j = k}^{l - 1} d(x_j, x_{j + 1})
\end{equation}
when $k < l$.  If $M$ is complete, then it follows that $\{x_j\}_{j =
1}^\infty$ converges in $M$.  Converesely, if $\{x_j\}_{j = 1}^\infty$
is a Cauchy sequence in $M$, then there is a subsequence
$\{x_{j_n}\}_{n = 1}^\infty$ of $\{x_j\}_{j = 1}^\infty$ such that
\begin{equation}
        d(x_{j_n}, x_{j_{n + 1}}) \le 2^{-n}
\end{equation}
for each $n$, which implies that
\begin{equation}
        \sum_{n = 1}^\infty d(x_{j_n}, x_{j_{n + 1}})
\end{equation}
converges.  If this subsequence converges, then $\{x_j\}_{j =
1}^\infty$ converges to the same limit, because it is a Cauchy
sequence.

        Let $E$ be a nonempty set, and consider $\ell^p(E)$, $0 < p
\le \infty$.  This is a metric space with respect to the metric
associated to the norm $\|f\|_p$ when $p \ge 1$, or the $p$-norm
$\|f\|_p$ when $0 < p \le 1$, and it is well known that this space is
complete.  For if $\{f_j\}){j = 1}^\infty$ is a Cauchy sequence in
$\ell^p(E)$, then it is easy to see that $\{f_j(x)\}_{j = 1}^\infty$
is a Cauchy sequence in ${\bf R}$ or ${\bf C}$ for each $x \in E$, as
appropriate.  This implies that $\{f_j(x)\}_{j = 1}^\infty$ converges
pointwise on $E$, since the real and complex numbers are complete.
One can check that the limit $f(x)$ is also in $\ell^p(E)$, and that
$\{f_j\}_{j = 1}^\infty$ converges to $f$ in the $\ell^p$ metric, as
desired.

\section{Infinite series}
\label{infinite series}
\setcounter{equation}{0}

        Let $V$ be a real or complex vector space equipped with a norm
or $p$-norm $\|v\|$, $0 < p \le 1$.  This determines a natural metric
on $V$, as in the previous section.  As usual, an infinite series
$\sum_{j = 1}^\infty v_j$ with terms $v_j \in V$ is said to converge
if the corresponding sequence of partial sums $\sum_{j = 1}^n v_j$
converges in $V$ as $n \to \infty$.  Let us say that $\sum_{j =
1}^\infty v_j$ converges absolutely if
\begin{equation}
        \sum_{j = 1}^\infty \|v_j\|
\end{equation}
converges when $\|v\|$ is a norm, and if
\begin{equation}
\label{sum_{j = 1}^infty ||v_j||^p}
        \sum_{j = 1}^\infty \|v_j\|^p
\end{equation}
converges when $\|v\|$ is a $p$-norm.  Note that the convergence of
(\ref{sum_{j = 1}^infty ||v_j||^p}) is more restrictive as $p$
decreases, as in Section \ref{monotonicity}.  As in the previous
section, absolute convergence of $\sum_{j = 1}^\infty v_j$ implies
that the sequence of partial sums $\sum_{j =1}^n v_j$ is a Cauchy
sequence.  In particular, absolute convergence implies convergence
when $V$ is complete.  Conversely, $V$ is complete if every absolutely
convergent series with terms in $V$ converges in $V$, by another
argument mentioned in the previous section.

\section{$c_0(E)$}
\label{c_0(E)}
\setcounter{equation}{0}

        Let $E$ be a nonempty set, and let $f(x)$ be a real or
complex-valued function on $E$.  We say that $f$ vanishes at infinity
on $E$ if for every $\epsilon > 0$, $|f(x)| \ge \epsilon$ for only
finitely many $x \in E$.  The spaces of real or complex-valued
functions on $E$ that vanish at infinity are denoted $c_0(E, {\bf
R})$, $c_0(E, {\bf C})$, respectively, and are vector spaces with
respect to pointwise addition and scalar multiplication of functions.
As usual, we may also use $c_0(E)$ to refer to both cases at the same
time.  Note that $f(x) \ne 0$ for only finitely or countably many $x
\in E$ when $f \in c_0(E)$.

        If $f$ vanishes at infinity on $E$, then $f$ is bounded, and
so $c_0(E)$ is a linear subspace of $\ell^\infty(E)$.  More precisely,
one can check that $c_0(E)$ is a closed linear subspace of
$\ell^\infty(E)$ with respect to the $\ell^\infty$ norm.  A function
$f$ on $E$ is said to have finite support if $f(x) \ne 0$ for only
finitely many $x \in E$, in which case it obviously vanishes at
infinity.  One can also check that functions with finite support
are dense in $c_0(E)$ with respect to the $\ell^\infty$ norm, so that
$c_0(E)$ is the same as the closure in $\ell^\infty(E)$ of the linear
subspace of functions with finite support.

        If a function $f$ on $E$ is $p$-summable for some $p > 0$,
then $f$ vanishes at infinity on $E$.  More precisely, the number of
$x \in E$ such that $|f(x)| \ge \epsilon$ is less than or equal to
\begin{equation}
        \epsilon^{-p} \sum_{x \in E} |f(x)|^p.
\end{equation}
Of course, a function $f$ with finite support on $E$ is $p$-summable
for every $p > 0$.  It is not difficult to show that functions with
finite support on $E$ are dense in $\ell^p(E)$ when $0 < p < \infty$.

\section{Generalized convergence, 2}
\label{generalized convergence, 2}
\setcounter{equation}{0}

        Let $E$ be a nonempty set, let $V$ be a real or complex vector
space with a norm or $p$-norm $\|v\|$, $0 < p \le 1$, and let $f(x)$
be a $V$-valued function on $E$.  We say that $\sum_{x \in E} f(x)$
converges in the generalized sense if there is a $\lambda \in V$ such
that for every $\epsilon > 0$ there is a finite set $A_\epsilon
\subseteq E$ such that
\begin{equation}
        \biggl\|\sum_{x \in B} f(x) - \lambda\biggr\| < \epsilon
\end{equation}
whenever $B \subseteq E$ is a finite set that satisfies $A_\epsilon
\subseteq B$.  It is easy to see that $\lambda$ is unique when it
exists, in which case it may be denoted $\sum_{x \in E} f(x)$.  Of
course, this is the same as the definition in Section \ref{generalized
convergence} when $V = {\bf R}$ or ${\bf C}$, and it is equivalent to
the convergence of the net of partial sums of $f(x)$ over finite
subsets of $E$ as in Section \ref{nets}.

        Similarly, we say that $\sum_{x \in E} f(x)$ satisfies the
generalized Cauchy criterion if for every $\epsilon > 0$ there is a
finite set $A_\epsilon \subseteq E$ such that
\begin{equation}
        \biggl\|\sum_{x \in B} f(x)\biggr\| < \epsilon
\end{equation}
whenever $B \subseteq E$ is a finite set with $A_\epsilon \cap B = \emptyset$.
If $\sum_{x \in E} f(x)$ converges in the generalized sense, then it is easy
to see that $\sum_{x \in E} f(x)$ satisfies the generalized Cauchy criterion.
Conversely, let us check that $\sum_{x \in E} f(x)$ converges in the
generalized sense when $\sum_{x \in E} f(x)$ satisfies the generalized
Cauchy criterion and $V$ is complete.

        If $\sum_{x \in E} f(x)$ satisfies the generalized Cauchy
criterion, then it is easy to see that $\|f(x)\|$ vanishes at infinity
on $E$, by considering sets $B \subseteq E$ with only one element in
the previous definition.  In particular, $f(x) \ne 0$ for only
finitely or countably many $x \in E$.  If $f(x) \ne 0$ for only
finitely many $x \in E$, then convergence of the sum is trivial, and
so we suppose that $f(x) \ne 0$ for countably many $x$.  Let
$\{x_j\}_{j = 1}^\infty$ be an enumeration of the set of $x \in E$
such that $f(x) \ne 0$, so that each element of this set occurs in the
sequence exactly once, and consider the infinite series $\sum_{j =
1}^\infty f(x_j)$.  Using the generalized Cauchy criterion for
$\sum_{x \in E} f(x)$, one can check that the sequence of partial sums
of $\sum_{j = 1}^\infty f(x_j)$ forms a Cauchy sequence in $V$.  If
$V$ is complete, then it follows that $\sum_{j = 1}^\infty f(x_j)$
converges in $V$.  Using the generalized Cauchy criterion for $\sum_{x
\in E} f(x)$ again, one can show that $\sum_{x \in E} f(x)$ converges
in the generalized sense, and that the sum is the same as $\sum_{j =
1}^\infty f(x_j)$.

\section{Summable functions, 2}
\label{summable functions, 2}
\setcounter{equation}{0}

        Let $E$ be a nonempty set, and let $V$ be a real or complex
vector space equipped with a norm or $p$-norm $\|v\|$ for $0 < p \le
1$.  Suppose that $f$ is a $V$-valued function on $E$ such that
$\|f(x)\|$ is summable on $E$ when $\|v\|$ is a norm on $V$, or that
$\|f(x)\|^p$ is summable on $E$ when $\|v\|$ is a $p$-norm, $0 < p \le
1$.  If $B \subseteq E$ is a finite set, then we have that
\begin{equation}
        \biggl\|\sum_{x \in B} f(x)\biggr\| \le \sum_{x \in B} \|f(x)\|
\end{equation}
in the first case, and
\begin{equation}
        \biggl\|\sum_{x \in B} f(x)\biggr\|^p \le \sum_{x \in B} \|f(x)\|^p
\end{equation}
in the second case.  In both cases, one can use these simple estimates
to check that $\sum_{x \in E} f(x)$ satisfies the generalized Cauchy
criterion.  If $V$ is complete, then it follows that $\sum_{x \in E}
f(x)$ converges in the generalized sense, as in the previous
section.

\section{A special case}
\label{special case}
\setcounter{equation}{0}

        Let $E$ be a nonempty set, and suppose that $\phi \in
\ell^p(E)$ for some $p$, $0 < p \le \infty$.  For each $x \in E$, let
$\delta_x(y)$ be the function on $E$ defined by $\delta_x(x) = 1$ and
$\delta_x(y) = 0$ when $y \ne x$.  Consider
\begin{equation}
        f(x) = \phi(x) \, \delta_x,
\end{equation}
as a function on $E$ with values in $\ell^p(E)$.  Observe that
\begin{equation}
 \sum_{x \in E} f(x)(y) = \sum_{x \in E} \phi(x) \, \delta_x(y) = \phi(y)
\end{equation}
for each $y \in E$, where these are sums over $x \in E$ of real or
complex numbers that are equal to $0$ when $x \ne y$ and hence
converge trivially.  One can also ask about the convergence of
$\sum_{x \in E} f(x)$ in the generalized sense to $\phi$, as a sum of
elements of $\ell^p(E)$.  Of course,
\begin{equation}
        \|f(x)\|_p = |\phi(x)| \, \|\delta_x\|_p = |\phi(x)|
\end{equation}
for every $x \in E$.  Thus $\|f(x)\|_p$ is $p$-summable on $E$ when $0
< p < \infty$, and bounded on $E$ when $p = \infty$.  If $0 < p \le
1$, then this is the same as the summability condition mentioned in
the previous section.  However, one can check that $\sum_{x \in E}
f(x)$ converges to $\phi$ in the generalized sense in $\ell^p(E)$ for
every positive real number $p$.  If $p = \infty$, then $\sum_{x \in E}
f(x)$ converges to $\phi$ in the generalized sense in $\ell^\infty(E)$
if and only if $\phi \in c_0(E)$.

\section{Inner product spaces}
\label{inner products}
\setcounter{equation}{0}

        An \emph{inner product} on a real or complex vector space $V$
is a real or complex-valued function $\langle v, w \rangle$, as
appropriate, defined for $v, w \in V$ and satisfying the following
three conditions.  First, $\langle v, w \rangle$ is a linear function
of $v$ for each $w \in W$.  Second,
\begin{equation}
        \langle w, v \rangle = \langle v, w \rangle
\end{equation}
for every $v, w \in V$ in the real case, and
\begin{equation}
        \langle w, v \rangle = \overline{\langle v, w \rangle}
\end{equation}
in the complex case.  In particular,
\begin{equation}
        \langle v, v \rangle = \overline{\langle v, v \rangle} \in {\bf R}
\end{equation}
for every $v \in V$ in the complex case.  Third,
\begin{equation}
        \langle v, v \rangle > 0
\end{equation}
for every $v \in V$ with $v \ne 0$.

        Put
\begin{equation}
        \|v\| = \langle v, v \rangle^{1/2}.
\end{equation}
The \emph{Cauchy--Schwarz inequality} states that
\begin{equation}
        |\langle v, w \rangle| \le \|v\| \, \|w\|
\end{equation}
for every $v, w \in V$.  Using this, one can show that
\begin{equation}
        \|v + w\| \le \|v\| + \|w\|
\end{equation}
for every $v, w \in V$, so that $\|v\|$ defines a norm on $V$.  If $V$
is complete with respect to this norm, then $V$ is said to be a
Hilbert space.

        Let $E$ be a nonempty set, and let $f, g \in \ell^2(E)$ be given.
Remember that
\begin{equation}
        a \, b \le \frac{a^2 + b^2}{2}
\end{equation}
for every $a, b \ge 0$, since $(a - b)^2 \ge 0$, so that
\begin{equation}
 \sum_{x \in E} |f(x)| \, |g(x)| \le \frac{1}{2} \sum_{x \in E} |f(x)|^2
                              + \frac{1}{2} \sum_{x \in E} |g(x)|^2 < +\infty.
\end{equation}
Thus $|f(x)| \, |g(x)|$ is summable on $E$, and it is easy to see that
\begin{equation}
        \langle f, g \rangle = \sum_{x \in E} f(x) \, g(x)
\end{equation}
defines an inner product on $\ell^2(E, {\bf R})$, and that
\begin{equation}
        \langle f, g \rangle = \sum_{x \in E} f(x) \, \overline{g(x)}
\end{equation}
defines an inner product on $\ell^2(E, {\bf C})$.  The corresponding
norm is the same as the $\ell^2$ norm discussed in Section
\ref{p-summable functions}.  These spaces are also complete, as in
Section \ref{metric spaces}, and are therefore Hilbert spaces.

        A pair of vectors $v$, $w$ in an inner product space $V$
are said to be \emph{orthogonal} if
\begin{equation}
        \langle v, w \rangle = 0.
\end{equation}
This may also be expressed by $v \perp w$.  In this case,
\begin{equation}
        \|v + w\|^2 = \langle v + w, v + w \rangle
          = \langle v, v \rangle + \langle w, w \rangle = \|v\|^2 + \|w\|^2.
\end{equation}
If $v_1, \ldots, v_n \in V$ and $v_j \perp v_l$ when $j \ne l$, then
we get that
\begin{equation}
        \biggl\|\sum_{j = 1}^n v_j\biggr\|^2 = \sum_{j = 1}^n \|v_j\|^2.
\end{equation}

\section{Inner product spaces, 2}
\label{inner products, 2}
\setcounter{equation}{0}

        Let $E$ be a nonempty set, let $(V, \langle v, w \rangle)$ be
an inner product space, and let $f$ be a $V$-valued function on $E$
such that
\begin{equation}
        f(x) \perp f(y)
\end{equation}
when $x \ne y$.  Thus
\begin{equation}
        \biggl\|\sum_{x \in B} f(x)\biggr\|^2 = \sum_{x \in B} \|f(x)\|^2
\end{equation}
for every finite set $B \subseteq E$, as in the previous section.  If
$\|f(x)\|^2$ is a summable function on $E$, then it follows that
$\sum_{x \in E} f(x)$ satisfies the generalized Cauchy criterion, and
hence converges in the generalized sense when $V$ is complete.  In
this case, one can also check that
\begin{equation}
        \biggl\|\sum_{x \in E} f(x)\biggr\|^2 = \sum_{x \in E} \|f(x)\|^2.
\end{equation}

\section{Infinite series, 2}
\label{infinite series, 2}
\setcounter{equation}{0}

        Let $V$ be a real or complex vector space equipped with a norm
or $p$-norm $\|v\|$, $0 < p \le 1$, and let $\sum_{j = 1}^\infty v_j$
be an infinite series with terms in $V$.  This can also be considered
as a sum over $E = {\bf Z}_+$, so that the notions of convergence in
the generalized sense and the generalized Cauchy criterion are
applicable.  If $\sum_{j = 1}^\infty v_j$ converges in the ordinary
sense and satisfies the generalized Cauchy criterion as a sum over
${\bf Z}_+$, then it is easy to see that $\sum_{j = 1}^\infty v_j$
converges in the generalized sense, and to the same sum.

        Suppose that $\sum_{j = 1}^\infty v_j$ does not satisfy the
generalized Cauchy criterion.  This means that there is an $\epsilon >
0$ such that for each finite set $A \subseteq {\bf Z}_+$ there is
another finite set $B \subseteq {\bf Z}_+$ such that $A \cap B =
\emptyset$ and
\begin{equation}
        \biggl\|\sum_{j \in B} v_j\biggr\| \ge \epsilon.
\end{equation}
Using this repeatedly, one can get finite subsets $A_n$, $B_n$ of
${\bf Z}_+$ such that $\{1, \ldots, n\} \subseteq A_n$, $A_n \cap B_n
= \emptyset$, $A_n \cup B_n \subseteq A_{n + 1}$, and
\begin{equation}
        \biggl\|\sum_{j \in B_n} v_j\biggr\| \ge \epsilon
\end{equation}
for each $n$.  Let $k_n$ be the number of elements of $A_n$ and $l_n$
be the number of elements of $B_n$, so that $n \le k_n < k_n + l_n \le
k_{n + 1}$ for each $n$.  Also let $\pi$ be a one-to-one mapping of
${\bf Z}_+$ onto itself such that $A_n = \{\pi(1), \ldots, \pi(k_n)\}$
and $B_n = \{\pi(k_n + 1), \ldots, \pi(k_n + l_n)\}$ for each $n$.
This is easy to arrange, because of the inclusion and disjointness
properties of the $A_n$'s and $B_n$'s.  Thus
\begin{equation}
        \biggl\|\sum_{j = k_n + 1}^{k_n + l_n} v_j\biggr\| \ge \epsilon
\end{equation}
for each $n$.  This implies that the partial sums of $\sum_{j =
1}^\infty v_{\pi(j)}$ do not form a Cauchy sequence, and in particular
that $\sum_{j = 1}^\infty v_{\pi(j)}$ does not converge in the
ordinary sense.

        If $\sum_{j = 1}^\infty v_j$ satisfies the generalized Cauchy
criterion, then it is easy to see that the partial sums of every
rearrangement $\sum_{j = 1}^\infty v_{\pi(j)}$ of $\sum_{j = 1}^\infty
v_j$ form a Cauchy sequence.  Conversely, if the partial sums of every
rearrangement of $\sum_{j = 1}^\infty v_j$ form a Cauchy sequence,
then $\sum_{j = 1}^\infty v_j$ satisfies the generalized Cauchy
criterion, by the argument in the preceding paragraph.  Similarly,
every rearrangement of $\sum_{j = 1}^\infty v_j$ converges to the same
sum when $\sum_{j = 1}^\infty v_j$ converges in the generalized sense.
Conversely, if every rearrangement of $\sum_{j = 1}^\infty v_j$
converges, then $\sum_{j = 1}^\infty v_j$ satisfies the generalized
Cauchy criterion, by the previous remarks.  Hence $\sum_{j = 1}^\infty
v_j$ converges in the generalized sense, because it converges in the
ordinary sense, as mentioned at the beginning of the section.

\section{H\"older's inequality}
\label{holder's inequality}
\setcounter{equation}{0}

        Let $E$ be a nonempty set, and suppose that $1 \le p, q \le
\infty$ are \emph{conjugate exponents} in the sense that
\begin{equation}
        \frac{1}{p} + \frac{1}{q} = 1.
\end{equation}
If $f \in \ell^p(E)$ and $g \in \ell^q(E)$, then \emph{H\"older's
inequality} states that $f \, g \in \ell^1(E)$, and that
\begin{equation}
        \|f \, g\|_1 \le \|f\|_p \, \|g\|_q.
\end{equation}
This is quite straightforward when $p = 1$, $q = \infty$ or $p =
\infty$, $q = 1$, and so we focus now on the case where $1 < p, q <
\infty$.  Note that the $p = q = 2$ case is another version of the
Cauchy--Schwarz inequality.

        If $a$, $b$ are nonnegative real numbers, then
\begin{equation}
        a \, b \le \frac{a^p}{p} + \frac{b^q}{q}.
\end{equation}
This can be seen as a consequence of the convexity of the exponential
function.  In particular,
\begin{equation}
        |f(x)| \, |g(x)| \le \frac{|f(x)|^p}{p} + \frac{|g(x)|^q}{q}
\end{equation}
for every $x \in E$.  Hence
\begin{equation}
        \sum_{x \in B} |f(x)| \, |g(x)| \le \frac{1}{p} \sum_{x \in B} |f(x)|^p
                                    + \frac{1}{q} \sum_{x \in B} |g(x)|^q
                        \le \frac{\|f\|_p^p}{p} + \frac{\|g\|_q^q}{q}
\end{equation}
for every finite set $B \subseteq E$.  This implies that $f \, g$ is
summable on $E$, with
\begin{equation}
        \|f \, g\|_1 \le \frac{\|f\|_p^p}{p} + \frac{\|g\|_q^q}{q}.
\end{equation}

        This implies H\"older's inequality when $\|f\|_p = \|g\|_q = 1$.
Otherwise, if $f, g \ne 0$, then we can apply this to
\begin{equation}
 \widetilde{f} = \frac{f}{\|f\|_p}, \quad \widetilde{g} = \frac{g}{\|g\|_q}.
\end{equation}
Thus $\widetilde{f} \in \ell^p(E)$, $\widetilde{g} \in \ell^q(E)$,
$\|\widetilde{f}\|_p = \|\widetilde{g}\|_q = 1$, and the previous
inequality implies that
\begin{equation}
        \frac{\|f \, g\|_1}{\|f\|_p \, \|g\|_q}
         = \|\widetilde{f} \, \widetilde{g}\|_1 \le 1.
\end{equation}
Of course, H\"older's inequality is trivial when either $f$ or $g$ is
identically $0$ on $E$.

\section{Bounded linear functionals}
\label{bounded linear functionals}
\setcounter{equation}{0}

        Let $V$ be a vector space over the real or complex numbers.  A
\emph{linear functional} on $V$ is simply a linear mapping from $V$
into ${\bf R}$ or ${\bf C}$, as appropriate.  Suppose now that $V$ is
also equipped with a norm $\|v\|$.  A linear functional $\lambda$ on
$V$ is said to be \emph{bounded} with respect to this norm if there is
a nonnegative real number $C$ such that
\begin{equation}
        |\lambda(v)| \le C \, \|v\|
\end{equation}
for every $v \in V$.  In this case, we put
\begin{equation}
\label{||lambda||_* = sup {|lambda(v)| : v in V, ||v|| le 1}}
        \|\lambda\|_* = \sup \{|\lambda(v)| : v \in V, \, \|v\| \le 1\},
\end{equation}
which is the same as the smallest $C \ge 0$ for which the previous
inequality holds.

        The boundedness of a linear functional $\lambda$ on $V$ implies that
\begin{equation}
        |\lambda(v) - \lambda(w)| = |\lambda(v - w)| \le C \, \|v - w\|
\end{equation}
for some $C \ge 0$ and every $v, w \in V$.  This shows that a bounded
linear functional $\lambda$ is uniformly continuous on $V$.
Conversely, if a linear functional $\lambda$ on $V$ is continuous at
$0$, then there is a $\delta > 0$ such that
\begin{equation}
        |\lambda(v)| < 1
\end{equation}
for every $v \in V$ with $\|v\| < \delta$.  This implies that
$\lambda$ is bounded, with $C = 1/\delta$.

        The space of arbitrary linear functionals on $V$ is a vector
space with respect to pointwise addition and scalar multiplication of
functions.  It is easy to see that the space $V^*$ of bounded linear
functionals on $V$ is also a vector space in this way, and that
$\|\lambda\|_*$ defines a norm on $V^*$, known as the \emph{dual
norm}.  Note that $V^*$ is automatically complete with respect to the
dual norm.  For if $\{\lambda_j\}_{j = 1}^\infty$ is a Cauchy sequence
of bounded linear functionals on $V$ with respect to the dual norm,
then $\{\lambda_j(v)\}_{j = 1}^\infty$ is a Cauchy sequence of real or
complex numbers, as appropriate, for each $v \in V$.  Hence
$\{\lambda_j(v)\}_{j = 1}^\infty$ converges in ${\bf R}$ or ${\bf C}$
for each $v \in V$, by completeness.  It is easy to see that the limit
defines a linear functional $\lambda$ on $V$, which is also bounded
because the $\lambda_j$'s have uniformly bounded dual norms.  One can
also show that $\{\lambda_j\}_{j = 1}^\infty$ converges to $\lambda$
with respect to the dual norm, using the fact that $\{\lambda_j\}_{j =
1}^\infty$ is a Cauchy sequence with respect to the dual norm.

        The definitions of bounded linear functionals and the dual
norm also make sense when $\|v\|$ is a $p$-norm on $V$.  The dual
space $V^*$ is still a vector space in this case, and the dual norm is
still a norm on $V^*$, and not just a $p$-norm.  The dual space is
also complete with respect to the dual norm, but there are some other
problems with the dual space when $\|v\|$ is not a norm, as we shall
see.

\section{H\"older's inequality, 2}
\label{holder's inequality, 2}
\setcounter{equation}{0}

        Let $E$ be a nonempty set, and let $1 \le p, q \le \infty$ be
conjugate exponents.  For each $g \in \ell^q(E)$, put
\begin{equation}
        \lambda_g(f) = \sum_{x \in E} f(x) \, g(x)
\end{equation}
when $f \in \ell^p(E)$.  This makes sense, because of H\"older's inequality,
and satisfies
\begin{equation}
        |\lambda_g(f)| \le \|f\|_p \, \|g\|_q.
\end{equation}
Thus $\lambda_g$ is a bounded linear functional on $\ell^p(E)$, with
dual norm less than or equal to $\|g\|_q$.  It is well known and not
too difficult to show that the dual norm of $\lambda_g$ on $\ell^p(E)$
is actually equal to $\|g\|_q$.  If $q = 1$ and $p = \infty$, then one
can also restrict $\lambda_g$ to $c_0(E)$.  One can also check that
the dual norm of the restriction of $\lambda_g$ to $c_0(E)$ with
respect to the $\ell^\infty$ norm is also equal to $\|g\|_1$.

        It is also well known that every bounded linear functional
$\lambda$ on $\ell^p(E)$ is of the form $\lambda_g$ for some $g \in
\ell^q(E)$ when $1 \le p < \infty$, and that every bounded linear
functional on $c_0(E)$ with respect to the $\ell^\infty$ norm is of
the form $\lambda_g$ for some $g \in \ell^1(E)$.  The basic idea is to put
\begin{equation}
\label{g(x) = lambda(delta_x)}
        g(x) = \lambda(\delta_x),
\end{equation}
where $\delta_x(x) = 1$ and $\delta_x(y) = 0$ when $y \in E$ and $y
\ne x$.  Using the boundedness of $\lambda$, one can show that $g \in
\ell^q(E)$.  By construction,
\begin{equation}
        \lambda(f) = \lambda_g(f)
\end{equation}
when $f(x) \ne 0$ for only finitely many $x \in E$.  This implies the
same relation for every $f \in \ell^p(E)$, $1 \le p < \infty$, or $f
\in c_0(E)$, as appropriate, because of the density of functions with
finite support on $E$ in these spaces.

        If $0 < p < 1$, then $\ell^p(E) \subseteq \ell^1(E)$, and
$\|f\|_1 \le \|f\|_p$ for every $f \in \ell^p(E)$.  It follows that
the restriction of a bounded linear functional on $\ell^1(E)$ to
$\ell^p(E)$ is a bounded linear functional with respect to the
$p$-norm $\|f\|_p$.  In particular, if $g \in \ell^\infty(E)$, then
the restriction of $\lambda_g$ to $\ell^p(E)$ is a bounded linear
functional with dual norm less than or equal to $\|g\|_\infty$ with
respect to $\|f\|_p$.  One can check that the dual norm of $\lambda_g$
on $\ell^p(E)$ is actually equal to $\|g\|_\infty$, because
$\lambda_g(\delta_x) = g(x)$ and $\|\delta_x\|_p = 1$ for each $x \in
E$.

        Conversely, if $\lambda$ is a bounded linear functional on
$\ell^p(E)$, $0 < p < 1$, then $\lambda = \lambda_g$ for some $g \in
\ell^\infty(E)$.  The proof is basically the same as when $p = 1$.  If
$g$ is as in (\ref{g(x) = lambda(delta_x)}), then $g$ is bounded, and
the $\ell^\infty$ norm of $g$ is less than or equal to the dual norm
of $\lambda$ on $\ell^p(E)$, because $\|\delta_x\|_p = 1$ for each $x
\in E$.  One can then use density of functions with finite support in
$\ell^p(E)$ to show that $\lambda = \lambda_g$.

\section{Hilbert spaces}
\label{hilbert spaces}
\setcounter{equation}{0}

        Let $(V, \langle v, w \rangle)$ be a real or complex inner
product space, and put
\begin{equation}
\label{lambda_w(v) = langle v, w rangle}
        \lambda_w(v) = \langle v, w \rangle
\end{equation}
for each $w \in W$.  By the Cauchy--Schwarz inequality, this is a
bounded linear functional on $V$, with $\|\lambda_w\|_* \le \|w\|$.
More precisely,
\begin{equation}
        \|\lambda\|_* = \|w\|,
\end{equation}
because $\lambda(w) = \|w\|^2$.  If $V$ is complete, then it is well
known that every bounded linear functional on $V$ is of this form.
Let us briefly review a proof of this fact.

        Let $Y \subseteq V$, $Y \ne \emptyset$, and $z \in V$ be
given, and let $\{y_j\}_{j = 1}^\infty$ be a sequence of elements of
$Y$ such that
\begin{equation}
\label{lim_{j to infty} ||y_j - z|| = inf {||y - z|| : y in Y}}
        \lim_{j \to \infty} \|y_j - z\| = \inf \{\|y - z\| : y \in Y\}.
\end{equation}
Note that
\begin{equation}
 \biggl\|\frac{u + v}{2}\biggr\|^2 + \biggl\|\frac{u - v}{2}\biggr\|^2
          = \frac{\|u\|^2}{2} + \frac{\|v\|^2}{2}
\end{equation}
for every $u, v \in V$, which is a version of the \emph{parallelogram
law}.  Applying this to $u = y_j - z$, $v = y_l - z$, we get that
\begin{equation}
\label{||y_j - y_l||^2}
 \biggl\|\frac{y_j + y_l}{2} - z\biggr\|^2 + \frac{\|y_j - y_l\|^2}{4}
        = \frac{\|y_j - z\|^2}{2} + \frac{\|y_l - z\|^2}{2}
\end{equation}
for each $j, l \ge 1$.  If $Y$ is convex, then $(y_j + y_l)/2 \in Y$
for every $j$, $l$, and hence
\begin{equation}
 \inf \{\|y - z\| : y \in Y\} \le \biggl\|\frac{y_j + y_l}{2} - z\biggr\|.
\end{equation}
Combining this with (\ref{lim_{j to infty} ||y_j - z|| = inf {||y -
z|| : y in Y}}) and (\ref{||y_j - y_l||^2}), we get that
\begin{equation}
        \lim_{j, l \to \infty} \|y_j - y_l\| = 0.
\end{equation}
Thus $\{y_j\}_{j = 1}^\infty$ is a Cauchy sequence when $Y$ is convex.
If $V$ is complete and $Y$ is also closed, then $\{y_j\}_{j =
1}^\infty$ converges to an element $y$ of $Y$ with minimal distance to
$z$.

        If $Y$ is a linear subspace of $V$, then one can show that $y
\in Y$ has minimal distance to $z \in V$ if and only if $z - y$ is
orthogonal to every element of $Y$.  One can also check that $y$ is
uniquely determined by these properties.  If $V$ is complete, $Y$ is a
closed linear subspace of $V$, and $z \in V$, then it follows from
that there is a $y \in Y$ such that $y - z$ is orthogonal to every
element of $Y$.

        Let $\lambda$ be a bounded linear functional on $V$, and let
\begin{equation}
        Y = \{v \in V : \lambda(v) = 0\}
\end{equation}
be the kernel of $\lambda$.  Thus $Y$ is a closed linear subspace of
$V$, and $Y = V$ if and only if $\lambda = 0$.  If $\lambda \ne 0$,
then there is a $w' \in V$ such that $w' \ne 0$ and $w' \perp y$ for
every $y \in Y$, by the discussion in the previous paragraphs.  In
this case, one can check that $\lambda = \lambda_w$, where $w$ is a
scalar multiple of $w'$.  This uses the observation that $Y$ has
codimension $1$ in $V$, so that every element of $V$ can be expressed
as a linear combination of $w'$ and an element of $Y$.

\section{The Hahn--Banach theorem}
\label{hahn--banach theorem}
\setcounter{equation}{0}

        Let $V$ be a real or complex vector space with a norm $\|v\|$,
and let $W$ be a linear subspace of $V$.  The \emph{Hahn--Banach
theorem} states that every bounded linear functional on $W$ can be
extended to a bounded linear functional on $V$ with the same norm.
Note that this theorem does not work for $p$-norms, $0 < p < 1$.  By
standard arguments based on uniform continuity, a bounded linear
functional on $W$ has a unique extension to a bounded linear
functional on the closure of $W$ with the same norm, and this does
work for $p$-norms on $V$.

        It follows from the Hahn--Banach theorem that for every $v \in
V$ with $v \ne 0$ there is a $\lambda \in V^*$ such that
$\|\lambda\|_* = 1$ and
\begin{equation}
\label{lambda(v) = ||v||}
        \lambda(v) = \|v\|.
\end{equation}
More precisely, (\ref{lambda(v) = ||v||}) determines a unique linear
functional on the $1$-dimensional subspace of $V$ spanned by $v$, and
the Hahn--Banach theorem implies that there is an extension of this
linear functional to $V$ with dual norm equal to $1$.  Note that this
corollary does not hold for $\ell^p(E)$ when $0 < p < 1$ and $E$ has
at least two elements.

        Let $V$ be the space of continuous real or complex-valued
functions $f$ on the unit interval $[0, 1]$.  If $0 < p < \infty$,
then put
\begin{equation}
        \|f\|_p = \Big(\int_0^1 |f(x)|^p \, dx\Big)^{1/p}.
\end{equation}
One can check that this is a norm when $p \ge 1$ and a $p$-norm when
$0 < p \le 1$, in the same way as for $\ell^p$.  The counterpart of
$\|f\|_p$ for $p = \infty$ is the supremum norm
\begin{equation}
        \|f\|_\infty = \sup \{|f(x)| : 0 \le x \le 1\}.
\end{equation}
It is well known that $V$ is complete with respect to the supremum
norm, and not with respect to $\|f\|_p$ when $0 < p < \infty$, for
which the completions of $V$ can be described in terms of Lebesgue
integrals.

        If $0 < p \le q \le \infty$, then
\begin{equation}
        \|f\|_p \le \|f\|_q
\end{equation}
for every continuous function $f$ on $[0, 1]$.  This is easy to see
when $q = \infty$, and it follows from the convexity of $r^{q/p}$ on
the nonnegative real numbers when $q < \infty$.  One can show that the
only bounded linear functional on $V$ with respect to $\|f\|_p$ is the
trivial linear functional equal to $0$ when $0 < p < 1$.  This is
because every continuous function $f$ on $[0, 1]$ can be expressed as
$\sum_{l = 1}^n f_l$ for some continuous functions $f_1, \ldots, f_n$
such that $\sum_{l = 1}^n \|f_l\|_p$ is arbitrarily small when $p <
1$.  More precisely, one can choose the $f_l$'s to be supported on
intervals of length approximately $1/n$.

\section{Weak summability}
\label{weak summability}
\setcounter{equation}{0}

        Let $E$ be a nonempty set, and let $V$ be a real or complex
vector space with a norm $\|v\|$.  Also let $f(x)$ be a $V$-valued
function on $E$ such that $\sum_{x \in E} f(x)$ converges in the
generalized sense.  If $\lambda$ is a bounded linear functional on
$V$, then $\sum_{x \in E} \lambda(f(x))$ also converges in the
generalized sense, and
\begin{equation}
\label{lambda(sum_{x in E} f(x)) = sum_{x in E} lambda(f(x))}
        \lambda\Big(\sum_{x \in E} f(x)\Big) = \sum_{x \in E} \lambda(f(x)).
\end{equation}
Of course, $\sum_{x \in E} f(x)$ automatically converges in the
generalized sense when $\|f(x)\|$ is summable on $E$, in which case
$\lambda(f(x))$ is summable on $E$ for every $\lambda \in V^*$, and
\begin{equation}
\label{sum_{x in E} |lambda(f(x))| le ||lambda||_* sum_{x in E} ||f(x)||}
 \sum_{x \in E} |\lambda(f(x))| \le \|\lambda\|_* \sum_{x \in E} \|f(x)\|.
\end{equation}
However, we have seen examples where $\sum_{x \in E} f(x)$ converges
in the generalized sense, even though $\|f(x)\|$ is not summable on
$E$.  If $\phi(x)$ is a real or complex-valued function on $E$ such
that $\sum_{x \in E} \phi(x)$ converges in the generalized sense, then
$\phi(x)$ is summable on $E$.  In particular, $\lambda(f(x))$ is a
summable function on $E$ for every $\lambda \in V^*$ when $\sum_{x \in
E} f(x)$ converges in the generalized sense.

\section{Bounded partial sums}
\label{bounded partial sums}
\setcounter{equation}{0}

        Let $V$ be a real or complex vector space with a norm or
$p$-norm $\|v\|$, $0 < p \le 1$.  Also let $X(V)$ be the space of
sequences $\{v_j\}_{j = 1}^\infty$ of elements of $V$ such that the
partial sums $\sum_{j = 1}^n v_j$ of $\sum_{j = 1}^\infty v_j$ are
uniformly bounded in $V$.  It is easy to see that $X(V)$ is a vector
space with respect to termwise addition and scalar multiplication.
Moreover,
\begin{equation}
\label{||{v_j}_{j = 1}^infty||_{X(V)}}
        \|\{v_j\}_{j = 1}^\infty\|_{X(V)}
          = \sup_{n \ge 1} \, \biggl\|\sum_{j = 1}^n v_j \biggr\|
\end{equation}
is a norm or $p$-norm on $X(V)$, as appropriate.  If $\{v_j\}_{j =
1}^\infty \in X(V)$, then the sums $\sum_{j = l}^n v_j$ are uniformly
bounded over $1 \le l \le n$, because
\begin{equation}
\label{sum_{j = l}^n v_j = sum_{j = 1}^n v_j - sum_{j = 1}^{l - 1} v_j}
        \sum_{j = l}^n v_j = \sum_{j = 1}^n v_j - \sum_{j = 1}^{l - 1} v_j.
\end{equation}
More precisely,
\begin{equation}
 \biggl\|\sum_{j = l}^n v_j\biggr\| \le \biggl\|\sum_{j = 1}^n v_j\biggr\|
                                  + \biggl\|\sum_{j = 1}^{l - 1} v_j\biggr\|
                                   \le 2 \, \|\{v_j\}_{j = 1}^\infty\|_{X(V)}
\end{equation}
when $\|v\|$ is a norm on $V$.  Similarly,
\begin{equation}
 \biggl\|\sum_{j = l}^n v_j\biggl\|^p \le \biggl\|\sum_{j = 1}^n v_j\biggr\|^p
                            + \biggl\|\sum_{j = 1}^{l - 1} v_j \biggr\|^p
                        \le 2 \, \|\{v_j\}_{j = 1}^\infty\|_{X(V)}^p
\end{equation}
when $\|v\|$ is a $p$-norm on $V$, so that
\begin{equation}
        \biggl\|\sum_{j = l}^n v_j\biggr\|
                   \le 2^{1/p} \|\{v_j\}_{j = 1}^\infty\|_{X(V)}.
\end{equation}
In particular, $\{v_j\}_{j = 1}^\infty$ is bounded, by taking $l = n$.

        An infinite series $\sum_{j = 1}^\infty v_j$ with terms in $V$
satisfies the ordinary Cauchy criterion if for every $\epsilon > 0$
there is an $L \ge 1$ such that
\begin{equation}
        \biggl\|\sum_{j = l}^n v_j \biggr\| < \epsilon
\end{equation}
when $n \ge l \ge L$.  This is equivalent to saying that the sequence
of partial sums $\sum_{j = 1}^n v_j$ is a Cauchy sequence in $V$.
Note that the partial sums are bounded in this case, so that
$\{v_j\}_{j = 1}^\infty \in X(V)$.  Put
\begin{equation}
 \quad   X_0(V) = \bigg\{\{v_j\}_{j = 1}^\infty \in X(V) :
    \sum_{j = 1}^\infty v_j \hbox{ satisfies the Cauchy criterion}\bigg\}.
\end{equation}
It is easy to see that $X_0(V)$ is a linear subspace of $X(V)$, and
that $\{v_j\}_{j = 1}^\infty$ is an element of $X_0(V)$ when $v_j = 0$
for all but finitely many $j$.  One can also check that $X_0(V)$ is
closed in $X(V)$, and in fact that $X_0(V)$ is the closure in $X(V)$
of the linear subspace of sequences $\{v_j\}_{j = 1}^\infty$ such that
$v_j = 0$ for all but finitely many $j$.  If $V$ is complete, then
$X_0(V)$ is the same as the space of sequences $\{v_j\}_{j =
1}^\infty$ such that $\sum_{j = 1}^\infty v_j$ converges in $V$.

\section{Bounded finite subsums}
\label{bounded finite subsums}
\setcounter{equation}{0}

        Let $E$ be a nonempty set, and let $V$ be a real or complex
vector space with a norm or $p$-norm $\|v\|$, $0 < p \le 1$.  Also let
$Y(E, V)$ be the space of $V$-valued functions $f(x)$ on $E$ such that
the sums $\sum_{x \in B} f(x)$ over nonempty finite subsets $B$ of $E$
are uniformly bounded in $V$.  It is easy to see that this is a vector
space with respect to pointwise addition and scalar multiplication,
and that
\begin{eqnarray}
 \|f\|_{Y(E, V)} & = & \sup \bigg\{\biggl\|\sum_{x \in B} f(x)\biggr\| :
              B \subseteq E, \, B \ne \emptyset, \hbox{ and $B$ has} \\
 & & \hskip 1.1in \hbox{only finitely many elements}\bigg\} \nonumber
\end{eqnarray}
is a norm or $p$-norm on $Y(E, V)$, as appropriate.  Note that each
$f \in Y(E, V)$ is bounded, and that
\begin{equation}
        \sup_{x \in E} \|f(x)\| \le \|f\|_{Y(E, V)}.
\end{equation}

        Let $Y_0(E, V)$ be the set of $V$-valued functions $f(x)$ on
$E$ such that $\sum_{x \in E} f(x)$ satisfies the generalized Cauchy
criterion.  It is easy to see that this is a closed linear subspace of
$Y(E, V)$.  If $f(x) = 0$ for all but finitely many $x \in E$, then $f
\in Y_0(E, V)$, and in fact $Y_0(E, V)$ is the same as the closure in
$Y(E, V)$ of the linear subspace of $V$-valued functions on $E$ with
finite support.  If $V$ is complete, then $Y_0(E, V)$ is also the same
as the collection of $V$-valued functions $f(x)$ on $E$ such that
$\sum_{x \in E} f(x)$ converges in the generalized sense.

        If $\|v\|$ is a norm on $V$ and $\|f(x)\|$ is summable on $E$,
or if $\|v\|$ is a $p$-norm on $V$ and $\|f(x)\|$ is $p$-summable on
$E$, $0 < p \le 1$, then $f \in Y(E, V)$, and
\begin{equation}
\label{||f||_{Y(E, V)}^p le sum_{x in E} ||f(x)||^p}
        \|f\|_{Y(E, V)}^p \le \sum_{x \in E} \|f(x)\|^p.
\end{equation}
Furthermore, $f \in Y_0(E, V)$ under these conditions.  Conversely, if
$V = {\bf R}$ and $f \in Y(E, {\bf R})$, then $f$ is summable on $E$, and
\begin{equation}
        \sum_{x \in E} |f(x)| \le 2 \, \|f\|_{Y(E, {\bf R})}.
\end{equation}
More precisely,
\begin{equation}
        \sum_{x \in E \atop f(x) \ge 0} f(x), \, 
   \sum_{x \in E \atop f(x) \le 0} -f(x) \le \|f\|_{Y(E, {\bf R})}.
\end{equation}
Similarly, if $V = {\bf C}$ and $f \in Y(E, {\bf C})$, then $f$ is
summable on $E$, and
\begin{equation}
        \sum_{x \in E} |f(x)| \le 4 \, \|f\|_{Y(E, {\bf C})}.
\end{equation}
In this case, the real and imaginary parts $\re f$, $\im f$ of $f$
are in $Y(E, {\bf R})$, and satisfy
\begin{equation}
\label{||re f||_{Y(E, R)}, ||im f||_{Y(E, R)} le ||f||_{Y(E, C)}}
        \|\re f\|_{Y(E, {\bf R})}, \, \|\im f\|_{Y(E, {\bf R})}
                                          \le \|f\|_{Y(E, {\bf C})}.
\end{equation}
This implies the desired estimate for the $\ell^1$ norm of $f$, which
is less than or equal to the sum of the $\ell^1$ norms of the real and
imaginary parts of $f$.

\section{Uniform boundedness}
\label{uniform boundedness}
\setcounter{equation}{0}

        Let $V$ be a real or complex vector space with a norm or
$p$-norm $\|v\|$, and take $E = {\bf Z}_+$.  Thus a $V$-valued
function on $E$ is basically the same as a sequence with terms in $V$,
and $Y({\bf Z}_+, V)$ can be identified with a linear subspace of
$X(V)$.  Also, $Y_0({\bf Z}_+, V)$ corresponds to a linear subspace of
$X_0(V)$ with respect to this identification, and the $X(V)$ norm is
less than or equal to the $Y({\bf Z}_+, V)$ norm.  By definition,
$Y({\bf Z}_+, V)$, $Y_0({\bf Z}_+, V)$, and the $Y({\bf Z}_+, V)$ norm
are invariant under one-to-one mappings of ${\bf Z}_+$ onto itself,
while $X(V)$, $X_0(V)$, and the $X(V)$ norm are not invariant under
rearrangements.

        Suppose that $\{v_j\}_{j = 1}^\infty$ is a sequence of
elements of $V$ such that $\{v_{\pi(j)}\}_{j = 1}^\infty$ is an
element of $X(V)$ for every one-to-one mapping $\pi$ from ${\bf Z}_+$
onto itself, and let us show that $\{v_j\}_{j = 1}^\infty$ corresponds
to an element of $Y({\bf Z}_+, V)$.  This would be immediate if we
also asked that the $X(V)$ norm of $\{v_{\pi(j)}\}_{j = 1}^\infty$ be
uniformly bounded, independently of $\pi$.  If $\{v_j\}_{j =
1}^\infty$ does not correspond to an element of $Y({\bf Z}_+, V)$,
then there is a sequence of finite subsets $B_1, B_2, \ldots$ of $E$
such that
\begin{equation}
 \biggl\|\sum_{j \in B_n} v_j\biggr\| \to \infty \hbox{ as } n \to \infty.
\end{equation}
One can also argue a bit more to get the $B_n$'s to be pairwise
disjoint.  This permits us to choose $\pi$ so that $B_n = \{\pi(k_n),
\ldots, \pi(l_n)\}$ for some $k_n, l_n \in {\bf Z}_+$ with $k_n \le
l_n$ and every $n$.  Hence $\{v_{\pi(j)}\}_{j = 1}^\infty \not\in
X(V)$, as desired.  Of course, the analogous statement for the
generalized Cauchy criterion was discussed in Section \ref{infinite
series, 2}.

\section{Uniform boundedness, 2}
\label{uniform boundedness, 2}
\setcounter{equation}{0}

        Let $M$ be a metric space, and let $\mathcal{A}$ be a
collection of continuous real or complex-valued functions on $M$.
Suppose that $\mathcal{A}$ is pointwise bounded on $M$, in the sense
that
\begin{equation}
        \mathcal{A}(x) = \{f(x) : f \in \mathcal{A}\}
\end{equation}
is a bounded set in ${\bf R}$ or ${\bf C}$, as appropriate, for each
$x \in M$.  Put
\begin{equation}
        A_n = \{x \in M : |f(x)| \le n \hbox{ for each } f \in \mathcal{A}\},
\end{equation}
so that $A_n$ is a closed set in $M$ for each $n$, by continuity, and
\begin{equation}
        \bigcup_{n = 1}^\infty A_n = M,
\end{equation}
by pointwise boundedness.  If $M$ is complete, then the Baire category
theorem implies that $A_n$ contains a nonempty open set in $M$ for
some $n$.

        Suppose now that $V$ is a real or complex vector space with a
norm or $p$-norm, and that $\Lambda$ is a collection of bounded linear
functionals on $V$.  If $\Lambda$ is bounded pointwise on $V$ and $V$
is complete, then $\Lambda$ is uniformly bounded on a nonempty open
set in $V$, as in the previous paragraph.  Using linearity, one can
check that the elements of $\Lambda$ have uniformly bounded dual
norms.  This is a version of the Banach--Steinhaus theorem, or uniform
boundedness principle.  Of course, $\Lambda$ is uniformly bounded on
bounded subsets of $V$ when the dual norms of the elements of
$\Lambda$ are uniformly bounded.

        Now let $W$ be a real or complex vector space with a norm
$\|w\|$, and let $K$ be a subset of $W$.  Suppose that
\begin{equation}
\label{K(lambda) = {lambda(w) : w in K}}
        K(\lambda) = \{\lambda(w) : w \in K\}
\end{equation}
is a bounded set in ${\bf R}$ or ${\bf C}$, as appropriate, for each
bounded linear functional $\lambda$ on $W$.  Each $w \in W$ determines
a bounded linear functional on $W^*$, which sends $\lambda \in W^*$ to
its value $\lambda(w)$ at $w$.  Dual spaces are automatically
complete, and so the boundedness of $K(\lambda)$ for each $\lambda \in
W^*$ implies that the linear functionals $\lambda \mapsto \lambda(w)$
corresponding to $w \in K$ have uniformly bounded dual norm on $W^*$,
as in the preceding paragraph.  It follows that $K$ is a bounded set
in $W$, by the Hahn--Banach theorem.

\section{Sums and linear functionals}
\label{sums, linear functionals}
\setcounter{equation}{0}

        Let $E$ be a nonempty set, and let $V$ be a real or complex
vector space with a norm or $p$-norm $\|v\|$.  If $f(x)$ is a
$V$-valued function on $E$ with uniformly bounded finite subsums, then
$\lambda(f(x))$ has the same property for each bounded linear
functional $\lambda$ on $V$.  Moreover,
\begin{equation}
\label{||lambda circ f||_{Y(E, R)} or ||lambda circ f||_{Y(E, C)}}
        \|\lambda \circ f\|_{Y(E, {\bf R})} \, \hbox{ or } \,
 \|\lambda \circ f\|_{Y(E, {\bf C})} \le \|\lambda\|_* \, \|f\|_{Y(E, V)},
\end{equation}
as appropriate.  This implies that $\lambda(f(x))$ is summable on $E$, with
\begin{equation}
 \sum_{x \in E} |\lambda(f(x))| \le 2 \, \|\lambda\|_* \, \|f\|_{Y(E, V)}
\end{equation}
in the real case, and
\begin{equation}
 \sum_{x \in E} |\lambda(f(x))| \le 4 \, \|\lambda\|_* \, \|f\|_{Y(E, V)}
\end{equation}
in the complex case.

        Conversely,
\begin{equation}
\label{|lambda(sum_{x in B} f(x))| = ... le sum_{x in B} |lambda(f(x))|}
        \biggl|\lambda\Big(\sum_{x \in B} f(x)\Big)\biggr| =
 \biggl|\sum_{x \in B} \lambda(f(x))\biggr| \le \sum_{x \in B} |\lambda(f(x))|
\end{equation}
for every finite set $B \subseteq E$ and $\lambda \in V^*$.  Suppose
that $\lambda(f(x))$ is summable on $E$ for each $\lambda \in V^*$, and that
\begin{equation}
\label{sum_{x in E} |lambda(f(x))| le C ||lambda||_*}
        \sum_{x \in E} |\lambda(f(x))| \le C \, \|\lambda\|_*
\end{equation}
for some $C \ge 0$ and every $\lambda \in V^*$.  If $\|v\|$ is a norm
on $V$, then the Hahn--Banach theorem implies that
\begin{equation}
\label{||sum_{x in B} f(x)|| le C}
        \biggl\|\sum_{x \in B} f(x) \biggr\| \le C
\end{equation}
for every finite set $B \subseteq E$.  Hence $f \in Y(E, V)$ and
\begin{equation}
        \|f\|_{Y(E, V)} \le C
\end{equation}
under these conditions.

        Let $K$ be the set of vectors in $V$ of the form $\sum_{x \in
B} f(x)$, where $B \subseteq E$ is a finite set.  If $\lambda(f(x))$
is summable on $E$ for some $\lambda \in V^*$, then the set
$K(\lambda)$ as in (\ref{K(lambda) = {lambda(w) : w in K}}) is
bounded.  If $\lambda(f(x))$ is summable on $E$ for every $\lambda \in
V^*$, and if $\|v\|$ is a norm on $V$, then it follows that $K$ is a
bounded set in $V$, as in the previous section.  This is the same as
saying that $f \in Y(E, V)$.

\section{Seminorms}
\label{seminorms}
\setcounter{equation}{0}

        Let $V$ be a vector space over the real or complex numbers.  A
nonnegative real-valued function $N(v)$ on $V$ is said to be a
\emph{seminorm} if
\begin{equation}
        N(t \, v) = |t| \, N(v)
\end{equation}
for every $v \in V$ and $t \in {\bf R}$ or ${\bf C}$, as appropriate, and
\begin{equation}
        N(v + w) \le N(v) + N(w)
\end{equation}
for every $v, w \in V$.  Thus a seminorm $N(v)$ is a norm exactly when
$N(v) > 0$ for every $v \in V$ with $v \ne 0$.  As another class of
examples, $N_\lambda(v) = |\lambda(v)|$ is a seminorm on $V$ when
$\lambda$ is a linear functional on $V$.  Observe that
\begin{equation}
        \{v \in V : N(v) = 0\}
\end{equation}
is a linear subspace of $V$ when $N(v)$ is a seminorm on $V$.

        Let $\mathcal{N}$ be a collection of seminorms on $V$.  Let us
say that $U \subseteq V$ is an open set with respect to $\mathcal{N}$
if for every $u \in U$ there are finitely many seminorms $N_1, \ldots,
N_l \in \mathcal{N}$ and positive real numbers $r_1, \ldots, r_l$ such that
\begin{equation}
        \{v \in V : N_j(u - v) < r_j, \, j = 1, \ldots, l\} \subseteq U.
\end{equation}
It is easy to see that this defines a topology on $V$.  Note that this
topology is Hausdorff if and only if $\mathcal{N}$ satisfies the
positivity condition that for each $v \in V$ with $v \ne 0$ there is
an $N \in \mathcal{N}$ such that $N(v) > 0$.  If $\mathcal{N}$
consists of a single norm, then this is the usual topology associated
to the norm.

        Suppose that $V$ is equipped with a norm or $p$-norm
$\|v\|_V$, and consider the collection of seminorms on $V$ of the form
$N_\lambda(v) = |\lambda(v)|$, where $\lambda \in V^*$.  The topology
on $V$ associated to this collection of seminorms is known as the
\emph{weak topology}.  If $\|v\|_V$ is a norm on $V$, then the
Hahn--Banach theorem implies that for each $v \in V$ with $v \ne 0$
there is a $\lambda \in V^*$ such that $\lambda(v) \ne 0$.  Thus
$N_\lambda(v) > 0$, and so the weak topology on $V$ is Hausdorff when
$\|v\|_V$ is a norm.  Note that open subsets of $V$ with respect to
the weak topology are open with respect to $\|v\|_V$, because the
linear functionals being used are bounded.

        Now let $W$ be a real or complex vector space with a norm or
$p$-norm $\|w\|_W$, and consider $V = W^*$.  Each $w \in W$ determines
a linear functional $\lambda \mapsto \lambda(w)$ on $W^*$, and hence a
seminorm $N_w^*(\lambda) = |\lambda(w)|$ on $W^*$.  The topology on
$W^*$ defined by this collection of seminorms is known as the
\emph{weak$^*$ topology}.  This topology is automatically Hausdorff,
but it is helpful for $\|w\|_W$ to be a norm on $W$ so that there are
plenty of bounded linear functionals on $W$.  Note that every open set
in $W^*$ with respect to the weak$^*$ topology is also open with
respect to the dual norm on $W^*$.

\section{Sums in dual spaces}
\label{sums in dual spaces}
\setcounter{equation}{0}

        Let $E$ be a nonempty set, let $W$ be a real or complex vector
space with a norm or $p$-norm $\|w\|$, and let $f$ be a function on
$E$ with values in the dual $W^*$ of $W$.  Suppose that $f(x)(w)$ is a
summable function on $E$ for every $w \in W$, where $f(x)(w)$ refers
to the value of $f(x) \in W^*$ at $w$, and that
\begin{equation}
        \sum_{x \in E} |f(x)(w)| \le C \, \|w\|
\end{equation}
for some $C \ge 0$ and every $w \in W$.  In this case, $\sum_{x \in E}
f(x)(w)$ defines a bounded linear functional on $W$ with dual norm
$\le C$.  One can also say that $\sum_{x \in E} f(x)$ converges in the
generalized sense with respect to the weak$^*$ topology on $W^*$ under
these conditions.

        This estimate also implies that
\begin{equation}
        \biggl|\sum_{x \in B} f(x)(w)\biggr| \le C \, \|w\|
\end{equation}
for every finite set $B \subseteq E$ and $w \in W$, which is to say that
\begin{equation}
        \biggl\|\sum_{x \in B} f(x)\biggr\|_* \le C
\end{equation}
for every finite set $B \subseteq E$.  Thus $f \in Y(E, W^*)$, and
\begin{equation}
        \|f\|_{Y(E, W^*)} \le C.
\end{equation}
Conversely, if $f \in Y(E, W^*)$, then $f(x)(w)$ is summable on $E$
for every $w \in W$, with $\ell^1$ norm bounded by $2 \, \|f\|_{Y(E,
V)}$ in the real case and by $4 \, \|f\|_{Y(E, V)}$ in the complex
case.  If $W$ is complete and $f(x)(w)$ is summable on $E$ for every
$w \in W$, then one can use the uniform boundedness principle to
conclude that $f \in Y(E, W^*)$.

\section{Seminorms, 2}
\label{seminorms, 2}
\setcounter{equation}{0}

        Let $V$ be a vector space over the real or complex numbers,
and let $N_1, N_2, \ldots$ be a sequence of seminorms on $V$ such that
for each $v \in V$ with $v \ne 0$ there is a positive integer $j$ for
which $N_j(v) > 0$.  Under these conditions, one can check that
\begin{equation}
        d(v, w) = \max \{\min(N_j(v - w), 1/j) : j \in {\bf Z}_+\}
\end{equation}
defines a metric on $V$, and that the topology on $V$ determined by
this metric is the same as the one associated to this sequence of
seminorms as in Section \ref{seminorms}.  Conversely, if the topology
on $V$ determined by a collection $\mathcal{N}$ of seminorms on $V$ is
metrizable, then it is Hausdorff, and there is a countable local base
for the topology at $0$.  Using the latter, one can show that there is
a subcollection of $\mathcal{N}$ with only finitely or countably many
elements that determines the same topology on $V$.

        Suppose now that $V$ is equipped with a norm or $p$-norm
$\|v\|$, and consider the weak topology on $V$.  Suppose also that for
each $v \in V$ with $v \ne 0$ there is a $\lambda \in V^*$ such that
$\lambda(v) \ne 0$, which follows from the Hahn--Banach theorem when
$\|v\|$ is a norm on $V$, and which implies that the weak topology on
$V$ is Hausdorff.  Suppose in addition that $V^*$ is separable, and
let $\lambda_1, \lambda_2, \ldots$ be a sequence of bounded linear
functionals on $V$ whose linear span is dense in $V^*$.  Let
$N_{\lambda_1}, N_{\lambda_2}, \ldots$ be the seminorms on $V$
corresponding to the $\lambda_j$'s as in Section \ref{seminorms}.
Under these conditions, one can check that the topology induced on a
bounded set in $V$ by the weak topology is the same as the topology
induced by the seminorms $N_{\lambda_1}, N_{\lambda_2}, \ldots$, and
hence is metrizable.

        Similarly, we can consider the weak$^*$ topology on the dual
of a vector space $V$ with a norm or $p$-norm.  Suppose that $V$ is
separable, so that there is a sequence of vectors $v_1, v_2, \ldots$
in $V$ whose linear span is dense in $V$.  Let $N_{v_1}^*, N_{v_2}^*,
\ldots$ be the seminorms on $V^*$ corresponding to the $v_j$'s, as in
Section \ref{seminorms}.  If $K$ is a bounded set in $V^*$ with
respect to the dual norm, then one can again check that the topology
induced on $K$ by the weak$^*$ topology is the same as the topology
induced by the seminorms $N_{v_1}^*, N_{v_2}^*, \ldots$, and is
therefore metrizable.

        Note that the unit ball
\begin{equation}
\label{B_1^* = {lambda in V^* : ||lambda||_* le 1}}
        B_1^* = \{\lambda \in V^* : \|\lambda\|_* \le 1\}
\end{equation}
in the dual $V^*$ of $V$ is closed with respect to the weak$^*$
topology.  To see this, it is convenient to describe $B_1^*$ as the
set of $\lambda \in V^*$ such that
\begin{equation}
\label{|lambda(v)| le 1}
        |\lambda(v)| \le 1
\end{equation}
for every $v \in V$ with $\|v\| \le 1$.  The Banach--Alaoglu theorem
states that $B_1^*$ is actually compact with respect to the weak$^*$
topology.  If $V$ is separable, then the topology induced on $B_1^*$
by the weak$^*$ topology on $V^*$ is metrizable, as in the previous
paragraph.  In this case, compactness of $B_1^*$ in the weak$^*$
topology is equivalent to sequential compactness.

\section{Isometric embeddings}
\label{isometric embeddings}
\setcounter{equation}{0}

        Let $(M, d(x, y))$ be a metric space.  It is easy to check that
\begin{equation}
\label{f_p(x) = d(p, x)}
        f_p(x) = d(p, x)
\end{equation}
is a continuous function on $M$ for each $p \in M$, using the triangle
inequality.  If $M$ is bounded, then $f_p$ is also a bounded function
on $M$.  Thus $p \mapsto f_p$ defines a mapping from $M$ into the
space $C_b(M)$ of bounded continuous real-valued functions on $M$.
Using the triangle inequality, one can show that this is an isometric
embedding of $M$ into $C_b(M)$ with the supremum norm.

        If $M$ is not bounded, then one can pick a basepoint $p_0 \in
M$, and put
\begin{equation}
        \widetilde{f}_p = f_p - f_{p_0}.
\end{equation}
Using the triangle inequality again, one can check that
$\widetilde{f}_p$ is a bounded function on $M$ for each $p \in M$.
Moreover, $p \mapsto \widetilde{f}_p$ is an isometric embedding of $M$
into $C_b(M)$ for the same reasons as before, since
\begin{equation}
\label{widetilde{f}_p - widetilde{f}_q = f_p - f_q}
        \widetilde{f}_p - \widetilde{f}_q = f_p - f_q
\end{equation}
for every $p, q \in M$.

        Suppose now that $V$ is a real or complex vector space with a
norm $\|v\|$, and let $B_1^*$ be the closed unit ball in the dual
space $V^*$, as in (\ref{B_1^* = {lambda in V^* : ||lambda||_* le
1}}).  Each $v \in V$ determines a bounded linear functional on $V^*$
defined by
\begin{equation}
        L_v(\lambda) = \lambda(v),
\end{equation}
which can also be considered as a bounded continuous function on
$B_1^*$ with respect to the topology induced by the weak$^*$ topology.
Thus $v \mapsto L_v$ defines a linear mapping from $V$ into the space
$C(B_1^*)$ of continuous real or complex-valued functions on $B_1^*$
with respect to the weak$^*$ topology, as appropriate.  By the
Banach--Alaoglu theorem, $B_1^*$ is a compact Hausdorff space with
respect to this topology.  Using the Hahn--Banach theorem, it is easy
to see that $v \mapsto L_v$ is also an isometry from $V$ into
$C(B_1^*)$, with respect to the supremum norm on $C(B_1^*)$.

\part{Functions, measures, and paths}

\section{Uniform boundedness, 3}
\label{uniform boundedness, 3}
\setcounter{equation}{0}

        Let $(X, \mathcal{A})$ be a measurable space, which is to say
a set $X$ with a $\sigma$-algebra $\mathcal{A}$ of measurable subsets
of $X$, and let $p$ be a nonnegative real-valued function on
$\mathcal{A}$.  Suppose that for every sequence $A_1, A_2, \ldots$ of
pairwise-disjoint measurable subsets of $X$,
\begin{equation}
 p\Big(\bigcup_{j = 1}^\infty A_j\Big) \le \sum_{j = 1}^\infty p(A_j) < \infty.
\end{equation}
This implies that $p(\emptyset) = 0$, by taking $A_j = \emptyset$ for
each $j$.

        Let $B_1, B_2, \ldots$ be a decreasing sequence of measurable
subsets of $X$, so that $B_{j + 1} \subseteq B_j$ for each $j$, and
put $B_\infty = \bigcap_{j = 1}^\infty B_j$.  Thus $A_j = B_j
\backslash B_{j + 1}$ is a sequence of pairwise-disjoint measurable
subsets of $X$ which are also disjoint from $B_\infty$, and
\begin{equation}
        B_n = \Big(\bigcup_{j = n}^\infty A_j\Big) \cup B_\infty
\end{equation}
for each $n$.  In particular, $\sum_{j = 1}^\infty p(A_j)$ converges,
which implies that $p(B_n)$ is uniformly bounded in $n$, since
\begin{equation}
\label{p(B_n) le sum_{j = n}^infty p(A_j) + p(B_infty)}
        p(B_n) \le \sum_{j = n}^\infty p(A_j) + p(B_\infty)
\end{equation}
for each $n$.  If $B_\infty = \emptyset$, then (\ref{p(B_n) le sum_{j
= n}^infty p(A_j) + p(B_infty)}) implies that $\{p(B_n)\}_{n =
1}^\infty$ converges to $0$.  If $C_1, C_2, \ldots$ is an increasing
sequence of measurable subsets of $X$, then a similar argument shows
that $p(C_n)$ is uniformly bounded in $n$, but we shall not need this
here.

        If $A \subseteq X$ is measurable, then put
\begin{eqnarray}
 p^*(A) & = & \sup \bigg\{\sum_{j = 1}^\infty p(A_j) : A_1, A_2, \ldots
                                       \hbox{ are pairwise-disjoint} \\
          & & \quad\hbox{measurable subsets of $X$ such that }
                    A = \bigcup_{j = 1}^\infty A_j \bigg\}.      \nonumber
\end{eqnarray}
We would like to show that $p^*(A) < \infty$ under these conditions.
Equivalently, one can check that
\begin{eqnarray}
 p^*(A) & = & \sup \bigg\{\sum_{j = 1}^n p(A_j) : A_1, \ldots, A_n
                                        \hbox{ are pairwise-disjoint} \\
          & & \quad\hbox{measurable subsets of $X$ such that }
                    A = \bigcup_{j = 1}^n A_j \bigg\}.           \nonumber
\end{eqnarray}
More precisely, the second definition of $p^*(A)$ is clearly less than
or equal to the first definition, because a partition of $A$ into
finitely many measurable sets can be extended to an infinite partition
using the empty set.  To show that the first definition of $p^*(A)$ is
less than or equal to the second definition, one can approximate an
infinite partition $A_1, A_2, \ldots$ of $A$ by the finite partitions
consisting of the sets $A_1, \ldots, A_n$ and $\bigcup_{j = n +
1}^\infty A_j$ for each $n$.

        If $B_1, B_2, \ldots$ is a sequence of pairwise-disjoint
measurable subsets of $X$, then
\begin{equation}
 \sum_{l = 1}^\infty p^*(B_l) \le p^*\Big(\bigcup_{l = 1}^\infty B_l\Big),
\end{equation}
because partitions of the $B_l$'s can be combined to get a partition
of $\bigcup_{l = 1}^\infty B_l$.  Similarly,
\begin{equation}
 p^*\Big(\bigcup_{l = 1}^\infty B_l\Big) \le \sum_{l = 1}^\infty p^*(B_l),
\end{equation}
because every measurable partition $\{E_j\}_{j = 1}^\infty$ of
$\bigcup_{l = 1}^\infty B_l$ can be refined to get a partition $\{E_j \cap 
B_l\}_{j, l = 1}^\infty$ which is a combination of partitions of the $B_l$'s.
Countable subadditivity implies that $p(E_j)$ is less than or equal to
the sum of $p(E_j \cap B_l)$ over $l$ for each $j$, so that the sum of
$p(E_j)$ over $j$ is less than or equal to the sum of $p(E_j \cap B_l)$
over $j$ and $l$.  The sum of $p(E_j \cap B_l)$ over $j$ is less than
or equal to $p^*(B_l)$ for each $l$, and so the sum of $p(E_j \cap B_l)$
over $j$ and $l$ is less than or equal to the sum of $p^*(B_l)$ over $l$,
as desired.  Therefore
\begin{equation}
 p^*\Big(\bigcup_{l = 1}^\infty B_l\Big) = \sum_{l = 1}^\infty p^*(B_l),
\end{equation}
which means that $p^*$ is countably additive.

        Suppose for the sake of a contradiction that $p^*(A) = \infty$
for some measurable set $A \subseteq X$.  This implies that there is a
finite sequence of pairwise-disjoint measurable subsets $A_{1, 1}, \ldots,
A_{1, n_1}$ of $X$ such that
\begin{equation}
\label{A = bigcup_{j = 1}^{n_1} A_{1, j}}
        A = \bigcup_{j = 1}^{n_1} A_{1, j}
\end{equation}
and
\begin{equation}
        \sum_{j = 1}^{n_1} p(A_{1, j}) \ge 1.
\end{equation}
We also have that $p^*(A_{1, j}) = \infty$ for some $j$, since
\begin{equation}
        p^*(A) = p^*(A_{1, 1}) + \cdots + p^*(A_{1, n_1}),
\end{equation}
and so we can relabel the indices, if necessary, to get that
\begin{equation}
        p^*(A_{1, n_1}) = \infty.
\end{equation}
This permits us to repeat the process, to get a finite sequence $A_{2, 1},
\ldots, A_{2, n_2}$ of pairwise-disjoint measurable subsets of $X$ such that
\begin{equation}
        A_{1, n_1} = \sum_{j = 1}^{n_2} A_{2, j}
\end{equation}
and
\begin{equation}
        \sum_{j = 1}^{n_2} p(A_{2, j}) \ge 2.
\end{equation}
As before, $p^*(A_{2, j}) = \infty$ for some $j$, and we can relabel
the indices if necessary to get that $p^*(A_{2, n_2}) = \infty$.
Continuing in this way, we get a finite sequence $A_{k, 1}, \ldots,
A_{k, n_k}$ of pairwise-disjoint measurable subsets of $X$ for each
positive integer $k$ such that
\begin{equation}
        \bigcup_{l = 1}^{n_k} A_{k, l} = A_{k - 1, n_{k - 1}}
\end{equation}
when $k \ge 2$,
\begin{equation}
\label{sum_{l = 1}^{n_k} p(A_{k, l}) ge k}
        \sum_{l = 1}^{n_k} p(A_{k, l}) \ge k,
\end{equation}
and $p^*(A_{k, n_k}) = \infty$.

        However,
\begin{equation}
        \sum_{k = 1}^\infty \sum_{l = 1}^{n_k - 1} p(A_{k, l}) < \infty,
\end{equation}
because the $A_{k, l}$'s are pairwise disjoint when $l < n_k$.  Hence
the sums
\begin{equation}
        \sum_{l = 1}^{n_k - 1} p(A_{k, l})
\end{equation}
are uniformly bounded in $k$, and even converge to $0$ as $k \to
\infty$.  By construction, $A_{k + 1, n_{k + 1}} \subseteq A_{k, n_k}$
for each $k$, and so $p(A_{k, n_k})$ is also uniformly bounded in $k$,
as mentioned earlier in the section.  This implies that the sums
\begin{equation}
\label{sum_{l = 1}^{n_k} p(A_{k_l}) = ...}
        \sum_{l = 1}^{n_k} p(A_{k_l}) = \sum_{l = 1}^{n_k - 1} p(A_{k, l})
                                                           + p(A_{k, n_k})
\end{equation}
are uniformly bounded in $k$ as well.  This contradicts (\ref{sum_{l =
1}^{n_k} p(A_{k, l}) ge k}), and we conclude that $p^*(A) < \infty$
for every measurable set $A \subseteq X$.

        Of course,
\begin{equation}
\label{p(A) le p^*(A)}
        p(A) \le p^*(A)
\end{equation}
for every measurable set $A \subseteq X$, and in fact $p^*$ is the
smallest countably-additive measure with this property.  More
precisely, if $\rho$ is a countably-additive measure such that $p(A)
\le \rho(A)$ for every measurable set $A \subseteq X$, then $p^*(A)
\le \rho(A)$ for each $A$.  This follows directly from the definition
of $p^*(A)$.  Observe too that the hypothesis that $\sum_{j =
1}^\infty p(A_j)$ converges when $A_1, A_2, \ldots$ is a sequence of
pairwise-disjoint measurable sets is necessary in order to have a
finite measure $\rho$ such that $p(A) \le \rho(A)$.

\section{Real and complex measures}
\label{real and complex measures}
\setcounter{equation}{0}

        Let $(X, \mathcal{A})$ be a measurable space, and let $\mu$ be a
real or complex measure on this space.  This means that $\mu$ is a real
or complex-valued function on $\mathcal{A}$ such that
\begin{equation}
 \mu\Big(\bigcup_{j = 1}^\infty A_j \Big) = \sum_{j = 1}^\infty \mu(A_j)
\end{equation}
for every sequence $A_1, A_2, \ldots$ of pairwise-disjoint measurable
subsets of $X$.  More precisely, the convergence of the series
$\sum_{j = 1}^\infty \mu(A_j)$ is part of the definition.  It follows
that the series converges absolutely, because every rearrangement of
the series is of the same type.  Note that $\mu(\emptyset) = 0$ is
also implied by the definition, by taking $A_j = \emptyset$.  If $p(A)
= |\mu(A)|$, then it is easy to see that $p(A)$ satisfies the
conditions described in the previous section.  Hence $p^*(A)$ is a
countably-additive finite measure, which is commonly denoted
$|\mu|(A)$.

        In the real case, $\mu$ is also known as a signed measure on
$X$, and it is easy to see that
\begin{equation}
        \mu^+(A) = \frac{|\mu|(A) + \mu(A)}{2},
         \quad   \mu^- = \frac{|\mu|(A) - \mu(A)}{2}
\end{equation}
are finite nonnegative measures on $X$.  Note that
\begin{equation}
        \mu(A) = \mu^+(A) - \mu^-(A)
\end{equation}
and
\begin{equation}
        |\mu|(A) = \mu^+(A) + \mu^-(A)
\end{equation}
for each measurable set $A \subseteq X$.  Similarly, if $\mu$ is a
complex measure on $X$, then $\mu$ can be expressed as a linear
combination of finite nonnegative measures on $X$, by applying this
argument to the real and imaginary parts of $\mu$.

        There are a number of simplifications that can be made in the
previous section when $p(A) = |\mu(A)|$ for a real measure $\mu$ on
$X$.  The first simplification is to replace the earlier definition of
$p^*(A)$ with
\begin{equation}
\label{p^*(A) = sup{|mu(B)| + |mu(C)| : A = B cup C, B cap C = emptyset}}
        p^*(A) = \sup\{|\mu(B)| + |\mu(C)| : B, C \in \mathcal{A}, \,
                                      A = B \cup C, \, B \cap C = \emptyset\}.
\end{equation}
The right side is clearly less than or equal to the earlier definition
of $p^*(A)$.  To show the opposite inequality, let $\{A_j\}_{j =
1}^\infty$ be any sequence of pairwise-disjoint measurable subsets of
$X$ such that $A = \bigcup_{j = 1}^\infty A_j$.  If $B$ is the union
of the $A_j$'s with $\mu(A_j) \ge 0$ and $C$ is the union of the
$A_j$'s with $\mu(A_j) < 0$, then $A = B \cup C$, $B \cap C =
\emptyset$, and
\begin{equation}
        \sum_{j = 1}^\infty |\mu(A_j)| = \mu(B) - \mu(C) = |\mu(B)| + |\mu(C)|.
\end{equation}
This implies that the earlier definition of $p^*(A)$ is less than or
equal to the right side of (\ref{p^*(A) = sup{|mu(B)| + |mu(C)| : A =
B cup C, B cap C = emptyset}}), by taking the supremum over all such
sequences $\{A_j\}_{j = 1}^\infty$.  In the same way, we also have
that
\begin{equation}
\label{p^*(A) = sup{mu(B) - mu(C) : A = B cup C, B cap C = emptyset}}
        p^*(A) = \sup\{\mu(B) - \mu(C) : B, C \in \mathcal{A}, \,
                                      A = B \cup C, \, B \cap C = \emptyset\}.
\end{equation}

        This makes it much easier to show that $p^*(A) < \infty$.  If
$p^*(A) = \infty$ for some measurable set $A \subseteq X$, then there
are disjoint measurable sets $B$, $C$ such that $A = B \cup C$ and
$\mu(B) - \mu(C)$ is as large as we want.  Of course,
\begin{equation}
\label{mu(A) = mu(B) + mu(C)}
        \mu(A) = \mu(B) + \mu(C),
\end{equation}
which implies that both $|\mu(B)|$ and $|\mu(C)|$ are as large as we
want.  Because $p^*$ is subadditive, we also have that $p^*(B) =
\infty$ or $p^*(C) = \infty$.  Put $A_1 = B$ if $p^*(B) = \infty$, and
otherwise $A_1 = C$.  Repeating the process, we get a decreasing
sequence $\{A_j\}_{j = 1}^\infty$ of measurable subsets of $X$ such
that $p^*(A_l) = \infty$ for each $l$ and $|\mu(A_l)| \to \infty$ as
$l \to \infty$.  This contradicts the fact that $|\mu(A_l)|$ is
bounded when $A_{l + 1} \subseteq A_l$ for each $l$, as in the
previous section.  One can also use the fact that $\{\mu(A)_j\}_{j =
1}^\infty$ converges under these conditions, and hence is bounded,
which is based on a similar argument.  It follows that $p^*(A) <
\infty$ when $p(A) = |\mu(A)|$ for a complex measure $\mu$, by
considering the real and imaginary parts of $\mu$.

        In the real case, we can combine (\ref{p^*(A) = sup{mu(B) -
mu(C) : A = B cup C, B cap C = emptyset}}) and (\ref{mu(A) = mu(B) +
mu(C)}) to get that
\begin{equation}
\label{mu^+(A) = sup{mu(B) : B in mathcal{A}, B subseteq A}}
        \mu^+(A) = \sup \{\mu(B) : B \in \mathcal{A}, \, B \subseteq A\}.
\end{equation}
We may restrict our attention to $B \subseteq A$ such that $\mu(B) \ge
0$ here, since $B = \emptyset$ has these properties.  Similarly,
\begin{equation}
\label{mu^-(A) = sup {-mu(C) : C in mathcal{A}, C subseteq A}}
        \mu^-(A) = \sup \{-\mu(C) : C \in \mathcal{A}, \, C \subseteq A\}.
\end{equation}

        If $\mu_1$ is a nonnegative real measure on $X$ such that
$\mu(A) \le \mu_1(A)$ for every measurable set $A \subseteq X$, then
\begin{equation}
\label{mu^+(A) le mu_1(A)}
        \mu^+(A) \le \mu_1(A)
\end{equation}
for every $A \in \mathcal{A}$.  More precisely, this uses the fact that
\begin{equation}
        \mu_1(A) = \mu_1(B) + \mu_1(A \backslash B) \ge \mu_1(B)
\end{equation}
when $B \subseteq A$, because $\mu_1(A \backslash B) \ge 0$.
Similarly, if $\mu_2$ is a nonnegative real measure on $X$ such that
$\mu(A) \ge - \mu_2(A)$ for every measurable set $A \subseteq X$, then
\begin{equation}
\label{mu^-(A) le mu_2(A)}
        \mu^-(A) \le \mu_2(A)
\end{equation}
for every $A \in \mathcal{A}$.  Of course, $\mu_1 = \mu^+$ and $\mu_2
= \mu^-$ have these properties, by construction.

        If $\mu_1$, $\mu_2$ are finite nonnegative real measures on
$X$ such that
\begin{equation}
        \mu(A) = \mu_1(A) - \mu_2(A)
\end{equation}
for every measurable set $A \subseteq X$, then
\begin{equation}
\label{-mu_2(A) le mu(A) le mu_1(A)}
        -\mu_2(A) \le \mu(A) \le \mu_1(A)
\end{equation}
for every $A \in \mathcal{A}$.  Thus $\mu_1$ and $\mu_2$ satisfy
(\ref{mu^+(A) le mu_1(A)}) and (\ref{mu^-(A) le mu_2(A)}),
respectively, as in the preceding paragraph.  As before, $\mu_1 =
\mu^+$, $\mu_2 = \mu^-$ have this property, by construction.

        Suppose that $P$, $Q$ are disjoint measurable subsets of $X$
such that $P \cup Q = X$ and
\begin{equation}
        |\mu|(X) = \mu(P) - \mu(Q).
\end{equation}
This is the same as saying that the supremum in (\ref{p^*(A) =
sup{mu(B) - mu(C) : A = B cup C, B cap C = emptyset}}) is attained
when $A = X$, with $B = P$ and $C = Q$.  If $E$ is a measurable subset
of $P$ such that $\mu(E) < 0$, then
\begin{equation}
        \mu(P) = \mu(P \backslash E) + \mu(E) < \mu(P \backslash E)
\end{equation}
and
\begin{equation}
        \mu(Q) > \mu(Q) + \mu(E) = \mu(Q \cup E),
\end{equation}
which implies that
\begin{equation}
        \mu(P) - \mu(Q) < \mu(P \backslash E) - \mu(Q \cup E),
\end{equation}
contradicting maximality.  Thus $\mu(E) \ge 0$ for every measurable
set $E \subseteq P$, and similarly $\mu(E) \le 0$ for every measurable
set $E \subseteq Q$.  Using this, one can check that
\begin{equation}
\label{mu^+(A) = mu(A cap P), mu^-(A) = mu(A cap Q)}
        \mu^+(A) = \mu(A \cap P), \quad  \mu^-(A) = \mu(A \cap Q)
\end{equation}
for every measurable set $A \subseteq X$, which is to say that the
suprema in (\ref{mu^+(A) = sup{mu(B) : B in mathcal{A}, B subseteq
A}}) and (\ref{mu^-(A) = sup {-mu(C) : C in mathcal{A}, C subseteq
A}}) are attained with $B = A \cap P$ and $C = A \cap Q$.

        The \emph{Hahn decomposition theorem} states that there are
disjoint measurable subsets $P$, $Q$ of $X$ such that $P \cup Q = X$
and (\ref{mu^+(A) = mu(A cap P), mu^-(A) = mu(A cap Q)}) holds for
every measurable set $A \subseteq X$.  One way to prove this is to
show that the supremum in (\ref{p^*(A) = sup{mu(B) - mu(C) : A = B cup
C, B cap C = emptyset}}) is attained when $A = X$, as in the next
paragraph.  Another way is to use the Radon--Nikodym theorem,
discussed in Section \ref{radon--nikodym theorem}.

        Suppose that $\{B_j\}_{j = 1}^\infty$, $\{C_j\}_{j =
1}^\infty$ are sequences of measurable subsets of $X$ such that $B_j
\cap C_j = \emptyset$ and $B_j \cup C_j = X$ for each $j$, and
\begin{equation}
\label{lim_{j to infty} (mu(B_j) - mu(C_j)) = |mu|(X)}
        \lim_{j \to \infty} (\mu(B_j) - \mu(C_j)) = |\mu|(X).
\end{equation}
Observe that
\begin{equation}
 |\mu|(X) - (\mu(B_j) - \mu(C_j)) = 2 \, (\mu^-(B_j) + \mu^+(C_j))
\end{equation}
for each $j$, because $|\mu|(X) = |\mu|(B_j) + |\mu|(C_j)$.  Hence
\begin{equation}
\label{lim_{j to infty} mu^-(B_j) = lim_{j to infty} mu^+(C_j) = 0}
 \lim_{j \to \infty} \mu^-(B_j) = \lim_{j \to \infty} \mu^+(C_j) = 0.
\end{equation}
Using this, one can show that $\{B_j\}_{j = 1}^\infty$, $\{C_j\}_{j =
1}^\infty$ are Cauchy sequences with respect to the semimetric on
$\mathcal{A}$ associated to $|\mu|$ as in Section \ref{distances
between sets}, and hence converge.  This is equivalent to saying that
the sequences of their indicator functions are Cauchy sequences in
$L^1(X, |\mu|)$, and hence converge in $L^1(X, |\mu|)$ to indicator
functions of measurable subsets of $X$.  More precisely, (\ref{lim_{j
to infty} mu^-(B_j) = lim_{j to infty} mu^+(C_j) = 0}) implies that
$\{B_j\}_{j = 1}^\infty$ converges to the empty set with respect to
$\mu^-$, and that $\{C_j\}_{j = 1}^\infty$ converges to the empty set
with respect to $\mu^+$.  This implies in turn that $\{B_j\}_{j =
1}^\infty$ converges to $X$ with respect to $\mu^+$, and that
$\{C_j\}_{j = 1}^\infty$ converges to $X$ with respect to $\mu^-$,
because $C_j = X \backslash B_j$ for each $j$.  It follows that
$\{B_j\}_{j = 1}^\infty$, $\{C_j\}_{j = 1}^\infty$ are Cauchy
sequences with respect to both $\mu^+$ and $\mu^-$, and are thus
Cauchy sequences with respect to $|\mu| = \mu^+ + \mu^-$.  The limits
of these sequences correspond to measurable subsets $P$, $Q$ of $X$
that are determined up to sets of $|\mu|$-measure $0$.  By
construction, $Q$ is the same as $X \backslash P$ up to a set of
$|\mu|$-measure $0$, and we may as well take $Q = X \backslash P$.  We
also have that $\mu^-(P) = \mu^+(Q) = 0$, $|\mu|(X) = \mu(P) -
\mu(Q)$, and so on.

\section{Vector-valued measures}
\label{vector-valued measures}
\setcounter{equation}{0}

        Let $(X, \mathcal{A})$ be a measurable space, and let $V$ be a
real or complex vector space with a norm $\|v\|$.  More precisely,
suppose that $V$ is a Banach space, which means that $V$ is complete
with respect to the metric associated to the norm.  Let $\mu$ be a
$V$-valued function on $\mathcal{A}$ such that
\begin{equation}
 \mu\Big(\bigcup_{j = 1}^\infty A_j \Big) = \sum_{j = 1}^\infty \mu(A_j)
\end{equation}
for every sequence $A_1, A_2, \ldots$ of pairwise-disjoint measurable
subsets of $X$.  Again convergence of the sum
\begin{equation}
\label{sum_{j = 1}^infty mu(A_j)}
        \sum_{j = 1}^\infty \mu(A_j)
\end{equation}
is part of the hypothesis, which implies convergence of rearrangements
of the sum.  However, in this case, absolute convergence
\begin{equation}
        \sum_{j = 1}^\infty \|\mu(A_j)\| < \infty
\end{equation}
is an additional condition.  If we have absolute convergence, then
$p(A) = \|\mu(A)\|$ satisfies the requirements of Section \ref{uniform
boundedness, 3}.  This implies that $\|\mu\|(A) = p^*(A)$ is a
countably-additive finite nonnegative measure.

        Let $\nu$ be a countably-additive finite nonnegative measure
on $(X, \mathcal{A})$, and take $V$ to be $L^q(X, \nu)$ for some $q$,
$1 \le q < \infty$.  Also let ${\bf 1}_A(x)$ be the indicator function
of $A \subseteq X$, equal to $1$ when $x \in A$ and to $0$ when $x \in
X \backslash A$.  If $\mu(A) = {\bf 1}_A$ for each measurable set $A
\subseteq X$, then $\mu$ is a $V$-valued function on $\mathcal{A}$
that satisfies the countable additivity condition described in the
previous paragraph.  If $q = 1$, then $\mu$ also satisfies the
absolute convergence condition.  This does not normally work when $q >
1$, even when $\nu$ is Lebesgue measure on the unit interval.

        Let $\mu$ be an arbitrary $V$-valued function $\mu$ on
$\mathcal{A}$ that satisfies the countable additivity condition
mentioned at the beginning of the section, not necessarily with
absolute convergence.  If $\lambda$ is a bounded linear functional on
$V$, then
\begin{equation}
        \mu_\lambda(A) = \lambda(\mu(A))
\end{equation}
defines a real or complex measure on $(X, \mathcal{A})$, as
appropriate.  In particular, $\mu_\lambda$ has finite total
variation $|\mu_\lambda|$, and
\begin{equation}
\label{|mu_lambda(A)| le |mu_lambda|(A) le |mu_lambda|(X)}
        |\mu_\lambda(A)| \le |\mu_\lambda|(A) \le |\mu_\lambda|(X)
\end{equation}
for every measurable set $A \subseteq X$.  Thus
\begin{equation}
\label{lambda(mu(A)) : A in mathcal{A}}
        \{\lambda(\mu(A)) : A \in \mathcal{A}\}
\end{equation}
is a bounded set of real or complex numbers, as appropriate, for each
$\lambda \in V^*$.  It follows that
\begin{equation}
        \{\mu(A) : A \in \mathcal{A}\}
\end{equation}
is a bounded set in $V$, as in Section \ref{uniform boundedness, 2}.

        If $\alpha$ is a real measure on $(X, \mathcal{A})$, then
\begin{equation}
\label{|alpha|(X) le 2 sup {|alpha(A)| : A in mathcal{A}}}
        |\alpha|(X) \le 2 \, \sup \{|\alpha(A)| : A \in \mathcal{A}\},
\end{equation}
because of (\ref{p^*(A) = sup{|mu(B)| + |mu(C)| : A = B cup C, B cap C
= emptyset}}).  Similarly, if $\beta$ is a complex measure on $(X,
\mathcal{A})$, then
\begin{equation}
        |\beta|(X) \le 4 \, \sup \{|\beta(A)| : A \in \mathcal{A}\},
\end{equation}
by applying (\ref{|alpha|(X) le 2 sup {|alpha(A)| : A in mathcal{A}}})
to the real and imaginary parts of $\beta$.  If $\mu$ is a
countably-additive $V$-valued function on $\mathcal{A}$ and $\lambda$
is a bounded linear functional on $V$, as in the previous paragraph, then
\begin{equation}
        |\mu_\lambda(A)| = |\lambda(\mu(A))| \le \|\lambda\|_* \, \|\mu(A)\|
\end{equation}
for every measurable set $A \subseteq X$.  Hence
\begin{equation}
        |\mu_\lambda|(X)
          \le 2 \, \|\lambda\|_* \, \sup \{\|\mu(A)\| : A \in \mathcal{A}\}
\end{equation}
in the real case, and
\begin{equation}
        |\mu_\lambda|(X)
          \le 4 \, \|\lambda\|_* \, \sup \{\|\mu(A)\| : A \in \mathcal{A}\}
\end{equation}
in the complex case.

        If $\mu$ is a countably-additive $V$-valued function on
$\mathcal{A}$ and $B_1, B_2, \ldots$ is an increasing sequence of
measurable subsets of $X$, then
\begin{equation}
\label{lim_{j to infty} mu(B_j) = mu(bigcup_{j = 1}^infty B_j)}
        \lim_{j \to \infty} \mu(B_j) = \mu\Big(\bigcup_{j = 1}^\infty B_j\Big).
\end{equation}
This follows from countable additivity by taking $A_1 = B_1$ and $A_j
= B_j \backslash B_{j - 1}$ when $j \ge 2$, as usual.  Conversely,
this continuity condition implies countable additivity when $\mu$ is
finitely additive, by taking $B_n = \bigcup_{j = 1}^n A_j$.
Similarly, if $C_1, C_2, \ldots$ is a decreasing sequence of
measurable subsets of $X$, then
\begin{equation}
        \lim_{l \to \infty} \mu(C_l) = \mu\Big(\bigcap_{l = 1}^\infty C_l\Big).
\end{equation}
This is equivalent to (\ref{lim_{j to infty} mu(B_j) = mu(bigcup_{j =
1}^infty B_j)}) when $\mu$ is finitely additive, with $B_j = X
\backslash C_j$.

        Let us use these continuity conditions to give another proof
of the fact that $\mu$ is bounded, like the one for real measures in
the previous section.  Put
\begin{equation}
 \widehat{\mu}(A) = \sup \{\|\mu(B)\| : B \in \mathcal{A}, \ B \subseteq A\}
\end{equation}
for each measurable set $A \subseteq X$, which may be $+\infty$ a priori.
Observe that
\begin{equation}
        \widehat{\mu}(A \cup A') \le \widehat{\mu}(A) + \widehat{\mu}(A')
\end{equation}
for any measurable sets $A, A' \subseteq X$.  This is because any
measurable subset $B$ of $A \cup A'$ can be expressed as the union of
$B \cap A \subseteq A$ and $B \backslash A \subseteq A'$, which are
automatically disjoint.  Thus $\mu(B)$ is the sum of $\mu(B \cap A)$
and $\mu(B \backslash A)$, so that $\|\mu(B)\|$ is less than or equal
to the sum of $\|\mu(B \cap A)\|$ and $\|\mu(B \backslash A)\|$, which
is less than or equal to the sum of $\widehat{\mu}(A)$ and
$\widehat{\mu}(A')$, as desired.

        Suppose for the sake of a contradiction that $\widehat{\mu}(A)
= +\infty$ for some measurable set $A \subseteq X$.  Hence there are
measurable sets $B \subseteq A$ such that $\|\mu(B)\|$ is as large as
we want.  Because $\mu(A)$ is equal to the sum of $\mu(B)$ and $\mu(A
\backslash B)$, it follows that $\|\mu(B)\|$ and $\|\mu(A \backslash
B)\|$ can both be as large as we want at the same time.  Using the
finite subadditivity of $\widehat{\mu}$ discussed in the previous
paragraph, we get that $\widehat{\mu}(B) = +\infty$ or
$\widehat{\mu}(A \backslash B) = +\infty$.  By taking $C_1 = B$ or $A
\backslash B$, as appropriate, we get a measurable subset of $A$ such
that $\widehat{C_1} = +\infty$ and $\|\mu(C_1)\|$ is as large as we
like.  Repeating the process, we get a decreasing sequence of
measurable sets $C_1, C_2, \ldots$ such that $\widehat{\mu}(C_l) =
+\infty$ for each $l \ge 1$ and $\|\mu(C_l)\| \to \infty$ as $l \to
\infty$.  This contradicts the fact that $\{\mu(C_l)\}_{l = 1}^\infty$
converges in $V$ to $\mu\big(\bigcap_{l = 1}^\infty C_l\big)$, by the
continuity condition that follows from countable additivity.

        Let $E$ be a nonempty set, and let $f(x)$ be a $V$-valued
function on $E$ such that $\sum_{x \in E} f(x)$ converges in the
generalized sense.  In particular, $\sum_{x \in E} f(x)$ satisfies the
generalized Cauchy condition, and so for each $\epsilon > 0$ there is
a finite set $B_\epsilon \subseteq E$ such that
\begin{equation}
        \biggl\|\sum_{x \in C} f(x)\biggr\| < \epsilon
\end{equation}
for every nonempty finite set $C \subseteq X \backslash B_\epsilon$.
It follows that $\sum_{x \in A} f(x)$ satsfies the generalized Cauchy
condition for every nonempty set $A \subseteq X$, since we can use $A
\cap B_\epsilon$ in place of $B_\epsilon$ for the sum over $A$.  Hence
$\sum_{x \in A} f(x)$ converges in the generalized sense for every
nonempty set $A \subseteq E$, because $V$ is complete.  Put
\begin{equation}
\label{mu(A) = sum_{x in A} f(x)}
        \mu(A) = \sum_{x \in A} f(x)
\end{equation}
for each $A \subseteq E$, which is interpreted as being $0$ when $A =
\emptyset$.  It is easy to see that this is a finitely-additive
$V$-valued measure on the algebra of all subsets of $E$.  Note that
\begin{equation}
        \|\mu(C)\| \le \epsilon
\end{equation}
for every $C \subseteq X \backslash B_\epsilon$, since we can reduce
to the previous case by approximating $C$ by finite sets.  Using this,
one can check that $\mu$ is countably-additive.  If $\|f(x)\|$ is a
summable function on $E$, then $\mu$ satisfies the additional absolute
convergence condition mentioned at the beginning of the section.

\section{The Radon--Nikodym theorem}
\label{radon--nikodym theorem}
\setcounter{equation}{0}

        Let $(X, \mathcal{A})$ be a measurable space, and let $\mu$,
$\nu$ be a finite nonnegative measures on $(X, \mathcal{A})$ such that
\begin{equation}
\label{mu(A) le C nu(A)}
        \mu(A) \le C \, \nu(A)
\end{equation}
for some $C \ge 0$ and every measurable set $A \subseteq X$.
A special case of the Radon--Nikodym theorem states that there is a
bounded nonnegative measurable function $h$ on $X$ such that
\begin{equation}
\label{mu(A) = int_A h d nu}
        \mu(A) = \int_A h \, d\nu
\end{equation}
for every measurable set $A \subseteq X$.  Von Neumann's trick for
showing this is to observe first that
\begin{equation}
\label{lambda(f) = int_X f d mu}
        \lambda(f) = \int_X f \, d\mu
\end{equation}
is a bounded linear functional on $L^2(\nu)$.  More precisely,
\begin{equation}
        |\lambda(f)| \le \int_X |f| \, d\mu \le C \, \int_X |f| \, d\nu
           \le C \, \nu(X)^{1/2} \, \Big(\int_X |f|^2 \, d\nu\Big)^{1/2},
\end{equation}
using our hypothesis on $\mu$ and $\nu$ in the second step, and the
Cauchy--Schwarz inequality in the third step.  Because $L^2(X, \nu)$
is a Hilbert space, the Riesz representation theorem implies that
there is an $h \in L^2(X, \nu)$ such that
\begin{equation}
        \lambda(f) = \int_X f \, h \, d\nu
\end{equation}
for every $f \in L^2(X, \nu)$.  Hence
\begin{equation}
        \mu(A) = \lambda({\bf 1}_A) = \int_A h \, d\nu
\end{equation}
for every measurable set $A \subseteq X$.  It follows that
\begin{equation}
        h(x) \le C
\end{equation}
almost everywhere on $X$ with respect to $\nu$ under these conditions.

        Instead of (\ref{mu(A) le C nu(A)}), suppose now that $\mu(A)
= 0$ for every measurable set $A \subseteq X$ such that $\nu(A) = 0$,
In this case, $\mu$ is said to be \emph{absolutely continuous} with
respect to $\nu$, denoted $\mu \ll \nu$.  The Radon--Nikodym theorem
states that there is then a nonnegative measurable function $h$ on $X$
such that (\ref{mu(A) = int_A h d nu}) holds for every measurable set
$A \subseteq X$.  More precisely, $h$ is also integrable with respect
to $\nu$, because
\begin{equation}
        \int_X h \, d\nu = \mu(X) < \infty.
\end{equation}
To see this, we apply the previous version to $\mu$ and $\nu_1 = \mu + \nu$,
since
\begin{equation}
\label{mu(A) le mu(A) + nu(A) = nu_1(A)}
        \mu(A) \le \mu(A) + \nu(A) = \nu_1(A)
\end{equation}
for every measurable set $A \subseteq X$ trivially.  This leads to a
real-valued measurable function $h_1$ on $X$ such that $0 \le h_1 \le 1$ and
\begin{equation}
\label{mu(A) = int_A h_1 d nu_1}
        \mu(A) = \int_A h_1 \, d\nu_1
\end{equation}
for every measurable set $A$.  If
\begin{equation}
\label{B = {x in X : h_1(x) = 1}}
        B = \{x \in X : h_1(x) = 1\},
\end{equation}
then $B$ is measurable, and
\begin{equation}
\label{mu(B) = nu_1(B) = mu(B) + nu(B)}
        \mu(B) = \nu_1(B) = \mu(B) + \nu(B),
\end{equation}
which implies that $\nu(B) = 0$, and hence $\mu(B) = 0$.  Thus $h_1 <
1$ $\nu$-almost everywhere, and one may as well take $h_1$ so that $0
\le h_1 < 1$ everywhere on $X$.  If $A \subseteq X$ is measurable,
then
\begin{equation}
\label{mu(A) = int_A h_1 d mu + int_A h_1 d nu}
        \mu(A) = \int_A h_1 \, d\mu + \int_A h_1 \, d\nu
\end{equation}
implies that
\begin{equation}
\label{int_A (1 - h_1) d mu = int_A h_1 d nu}
        \int_A (1 - h_1) \, d\mu = \int_A h_1 \, d\nu,
\end{equation}
and one can show that (\ref{mu(A) = int_A h d nu}) holds with
$h = h_1 / (1 - h_1)$.  More precisely,
\begin{equation}
\label{int_X g (1 - h_1) d mu = int_X g h_1 d nu}
        \int_X g \, (1 - h_1) \, d\mu = \int_X g \, h_1 \, d\nu
\end{equation}
for every bounded measurable function $g$ on $X$, because of
(\ref{int_A (1 - h_1) d mu = int_A h_1 d nu}).  If $h_1 \le 1 -
\delta$ on $A$ for some $\delta > 0$, then one can take $g = 1/(1 -
h_1)$ on $A$, $g = 0$ on $X \backslash A$, to get (\ref{mu(A) = int_A
h d nu}).  One can then use countable additivity to get (\ref{mu(A) =
int_A h d nu}) for arbitrary measurable sets $A$.

        If $\mu$ is a real or complex measure on $(X, \mathcal{A})$,
and not necessarily positive, then $\mu$ is still said to be
absolutely continuous with respect to $\nu$ when $\mu(A) = 0$ for
every measurable set $A \subseteq X$ such that $\nu(A) = 0$.  This is
equivalent to the condition that the total variation measure $|\mu|$
be absolutely continuous with respect to $\nu$, which implies that
$\mu$ can be expressed as a linear combination of finite nonnegative
measures on $X$ that are absolutely continuous with respect to $\nu$.
It follows from the previous case that there is a real or
complex-valued integrable function $h$ on $X$ with respect to $\nu$
for which (\ref{mu(A) = int_A h d nu}) holds.  One can also allow
$\nu$ to be $\sigma$-finite, by decomposing the domain into a
countable union of pairwise-disjoint measurable sets of finite
$\nu$-measure.  It is better to do this first when $\mu$ is
nonnegative, to get the integrability of the density $h$, and then
deal with real or complex measures $\mu$.

        Note that $h$ is determined $\nu$-almost everywhere by $\mu$.
More precisely, if $h$ is a real or complex-valued integrable function
on $X$ with respect to $\nu$ such that
\begin{equation}
        \int_A h \, d\nu = 0
\end{equation}
for every measurable set $A \subseteq X$, then $h(x) = 0$ for almost
every $x \in X$ with respect to $\nu$.  In the real case, one can
simply take $A$ to be the set where $h(x) > 0$ or $h(x) < 0$.  The
complex case follows from the real case, by considering the real and
imaginary parts of $h$ separately.  If $h'$, $h''$ are integrable
functions on $X$ with respect to $\nu$ such that
\begin{equation}
        \int_A h' \, d\nu = \int_A h'' \, d\nu
\end{equation}
for every measurable set $A \subseteq X$, then it follows that $h = h'
- h''$ is equal to $0$ almost everywhere on $X$ with respect to $\nu$.

        Of course, any real or complex measure $\mu$ on $X$ is
absolutely continuous with respect to the corresponding total
variation measure $|\mu|$.  The Radon--Nikodym theorem implies that
there is an integrable function $h$ on $X$ with respect to $|\mu|$
such that
\begin{equation}
        \mu(A) = \int_A h \, d|\mu|
\end{equation}
for every measurable set $A \subseteq X$.  Clearly
\begin{equation}
        |\mu(A)| \le \int_A |h| \, d|\mu|
\end{equation}
for every measurable set $A \subseteq X$, which implies that
\begin{equation}
        |\mu|(A) \le \int_A |h| \, d|\mu|,
\end{equation}
since the right side is a nonnegative measure on $X$.  It follows that
$|h(x)| \ge 1$ for almost every $x \in X$ with respect to $|\mu|$, and
we would like to check that $|h(x)| = 1$ almost everywhere on $X$.

        If $\mu$ is real and $A_1 = \{x \in X : h(x) > 1\}$ has
positive $|\mu|$-measure, then
\begin{equation}
\label{mu(A_1) = int_{A_1} h d|mu| > |mu|(A_1) ge mu(A_1)}
        \mu(A_1) = \int_{A_1} h \, d|\mu| > |\mu|(A_1) \ge \mu(A_1),
\end{equation}
a contradiction.  Thus $|\mu|(A_1) = 0$, and $|\mu|(\{x \in X : h(x) <
- 1\}) = 0$ for similar reasons.  In the complex case, put $A_\alpha =
\{x \in X : \re (\alpha \, h(x)) > 1\}$ for each $\alpha \in {\bf C}$
with $|\alpha| = 1$.  If $|\mu|(A_\alpha) > 0$ for some $\alpha$, then
\begin{equation}
        |\mu(A_\alpha)| \ge \re (\alpha \, \mu(A_\alpha))
                           = \int_{A_\alpha} \re (\alpha \, h) \, d|\mu|
                              > |\mu|(A_\alpha) \ge |\mu(A_\alpha)|,
\end{equation}
which is a contradiction again.  This shows that $|\mu|(A_\alpha) = 0$
for every complex number $\alpha$ with $|\alpha| = 1$.  Let
$\{\alpha_j\}_{j = 1}^\infty$ be a sequence of complex numbers with
$|\alpha_j| = 1$ for each $j$ which is dense in the unit circle in
${\bf C}$, such as an enumeration of the points on the circle that
correspond to angles that are rational multiples of $2 \, \pi$.  If $x
\in X$ and $|h(x)| > 1$, then $x \in A_{\alpha_j}$ when $\alpha_j$ is
sufficiently close to $\overline{h(x)}/|h(x)|$.  Equivalently,
\begin{equation}
        \{x \in X : |h(x)| > 1\} = \bigcup_{j = 1}^\infty A_{\alpha_j},
\end{equation}
and so $|\mu|(\{x \in X : |h(x)| > 1\}) = 0$, as desired.  In
particular, $h(x) = \pm 1$ almost everywhere on $X$ with respect to
$|\mu|$ in the real case, which implies the Hahn decomposition, as in
Section \ref{real and complex measures}.

\section{The Lebesgue decomposition}
\label{lebesgue decomposition}
\setcounter{equation}{0}

        Let $(X, \mathcal{A})$ be a measurable space, and let $\mu$
and $\nu$ be positive finite measures on $X$.  If $\nu_1 = \mu + \nu$,
then $\mu \le \nu_1$, and there is a real-valued measurable function
$h_1$ on $X$ that satisfies $0 \le h_1 \le 1$ and (\ref{mu(A) = int_A
h_1 d nu_1}), as before.  Let $B$ be as in (\ref{B = {x in X : h_1(x)
= 1}}), so that $B$ is measurable and satisfies (\ref{mu(B) = nu_1(B)
= mu(B) + nu(B)}), which implies that $\nu(B) = 0$.  However, without
the additional hypothesis of absolute continuity of $\mu$ with respect
to $\nu$, we do not necessarily have that $\mu(B) = 0$.  Instead, let
$\mu'$, $\mu''$ be the measures defined by
\begin{equation}
        \mu'(A) = \mu(A \cap B), \quad \mu''(A) = \mu(A \cap (X \backslash B)).
\end{equation}
By construction, $\mu'$ and $\nu$ are mutually singular, in the sense
that $\nu(B) = 0$ and $\mu'(X \backslash B) = 0$.  We still have
(\ref{mu(A) = int_A h_1 d mu + int_A h_1 d nu}), (\ref{int_A (1 - h_1)
d mu = int_A h_1 d nu}), and (\ref{int_X g (1 - h_1) d mu = int_X g
h_1 d nu}), which imply that
\begin{equation}
\label{mu''(A) = int_{A cap (X backslash B)} frac{h_1}{1 - h_1} d nu}
 \mu''(A) = \int_{A \cap (X \backslash B)} \frac{h_1}{1 - h_1} \, d\nu
\end{equation}
for every measurable set $A \subseteq X$.  In particular, $\mu''$ is
absolutely continuous with respect to $\nu$.  Of course, $\mu = \mu' +
\mu''$, which is known as the \emph{Lebesgue decomposition} of $\mu$.
If $\mu$ is a real or complex measure on $X$, then an analogous
decomposition can be obtained by applying this argument to $|\mu|$ in
place of $\mu$.

\section{The Riesz representation theorem}
\label{riesz representation theorem}
\setcounter{equation}{0}

       Let $(X, \mathcal{A}, \mu)$ be a measure space, and let $1 \le
p, q \le \infty$ be conjugate exponents, so that $1/p + 1/q = 1$.
If $f \in L^p(X)$ and $g \in L^q(X)$, then the integral version of
H\"older's inequality implies that $f \, g \in L^1(X)$, and that
\begin{equation}
        \|f \, g\|_1 \le \|f\|_p \, \|g\|_q.
\end{equation}
The proof is basically the same as for sums, as in Section
\ref{holder's inequality}.  It follows that
\begin{equation}
\label{lambda_g(f) = int_X f g d mu}
        \lambda_g(f) = \int_X f \, g \, d\mu
\end{equation}
defines a bounded linear functional on $L^p(X)$ when $g \in L^q(X)$,
with dual norm less than or equal to $\|g\|_q$.  If $p = \infty$, then
it is easy to see that the dual norm of $\lambda_g$ is equal to
$\|g\|_1$, by choosing $f \in L^\infty(X)$ such that $\|f\|_\infty =
1$ and $f \, g = |g|$.  Similarly, if $1 < p < \infty$, then the dual
norm of $\lambda_g$ on $L^p(X)$ is equal to $\|g\|_q$, because there
is an $f \in L^p(X)$ such that $f \, g = |f|^p = |g|^q$.  The dual
norm of $\lambda_g$ on $L^1(X)$ is also equal to $\|g\|_\infty$, under
an additional hypothesis.  More precisely, we should ask that for each
measurable set $A \subseteq X$ with $\mu(A) > 0$ there is a measurable
set $B \subseteq A$ such that $0 < \mu(B) < \infty$.  This condition
holds when $\mu$ is $\sigma$-finite on $X$, and for counting measure
on any set $X$.  If $0 \le t < \|g\|_\infty$, then we can apply this
to $A_t = \{x \in X : |g(x)| \ge t\}$ to get a measurable set $B_t
\subseteq A_t$ with $0 < \mu(B_t) < \infty$.  Put $f_t(x) = g(x) /
|g(x)|$ for every $x \in B_t$ when $g$ is real-valued, $f_t(x) =
\overline{g(x)} / |g(x)|$ for every $x \in B_t$ when $g$ is
complex-valued, and $f_t(x) = 0$ for every $x \in X \backslash B_t$ in
both cases.  It is easy to see that $f_t \in L^1(X)$, $\|f_t\|_1 =
\mu(B)$, and $\lambda_g(f_t) \ge t \, \mu(B)$, which implies that the
dual norm of $\lambda_g$ on $L^1(X)$ is greater than or equal to $t$.
It follows that the dual norm of $\lambda_g$ on $L^1(X)$ is greater
than or equal to $\|g\|_\infty$, since this holds for every
nonnegative real number $t$ such that $t < \|g\|_\infty$.  Hence the
dual norm of $\lambda_g$ on $L^1(X)$ is equal to $\|g\|_\infty$, since
we already know that it is less than or equal to $\|g\|_\infty$.

        Conversely, every bounded linear functional on $L^p(X)$ can be
realized in this way when $1 < p < \infty$, and also when $p = 1$ and
$X$ has $\sigma$-finite $\mu$-measure.  To see this, let us begin with
the case where $\mu(X) < \infty$.  Let $\lambda$ be a bounded linear
functional on $L^p(X)$, $1 \le p < \infty$, and put
\begin{equation}
\label{nu(A) = lambda({bf 1}_A)}
        \nu(A) = \lambda({\bf 1}_A)
\end{equation}
for every measurable set $A \subseteq X$.  Here ${\bf 1}_A$ denotes
the indicator function on $X$ associated to $A$, equal to $1$ on $A$
and to $0$ on $X \backslash A$.  If $A_1, A_2, \ldots$ is a sequence
of pairwise-disjoint measurable subsets of $X$, then $\sum_{j =
1}^\infty {\bf 1}_{A_j}$ converges in $L^p(X)$ to the indicator function
associated to $\bigcup_{j = 1}^\infty A_j$ when $p < \infty$, and hence
\begin{equation}
        \nu\Big(\bigcup_{j = 1}^\infty A_j\Big) = \sum_{j = 1}^\infty \nu(A_j).
\end{equation}
Thus $\nu$ is a real or complex measure on $X$, as appropriate.  This
measure is also absolutely continuous with respect to $\mu$, since
${\bf 1}_A = 0$ in $L^p(X)$ when $\mu(A) = 0$.  The Radon--Nikodym
theorem implies that there is a $g \in L^1(X)$ such that
\begin{equation}
\label{nu(A) = int_A g d mu}
        \nu(A) = \int_A g \, d\mu
\end{equation}
for every measurable set $A \subseteq X$.  By linearity, it follows
that
\begin{equation}
\label{lambda(f) = int_X f g d mu}
        \lambda(f) = \int_X f \, g \, d\mu
\end{equation}
for every measurable simple function $f$ on $X$.  This also holds when
$f$ is a bounded measurable function on $X$, by approximating $f$ by
simple functions.  If $p = 1$, then one can use this to show that $g
\in L^\infty(X)$, with $L^\infty$ norm less than or equal to the dual
norm of $\lambda$ on $L^1(X)$, in the same way as in the previous
paragraph.  If $p > 1$, then one can first show that the $L^q$ norm of
the restriction of $g$ to any set on which it is bounded is less than
or equal to the dual norm of $\lambda$ on $L^p(X)$, by the same type
of argument as in the previous paragraph.  This implies that $g \in
L^q(X)$, with $L^q$ norm less than or equal to the dual norm of
$\lambda$ on $L^p(X)$.  In both cases, one can then use the
boundedness of $\lambda$ on $L^p(X)$ and the fact that that $g \in
L^q(X)$ to show that (\ref{lambda(f) = int_X f g d mu}) holds for
every $f \in L^p(X)$, because simple functions are dense in $L^p(X)$.

        Suppose now that $X$ has $\sigma$-finite $\mu$-measure, so
that there is a sequence of measurable subsets $E_1, E_2, \ldots$ of
$X$ such that $\mu(E_l) < \infty$ for each $l \ge 1$ and $\bigcup_{l =
1}^\infty E_l = X$.  We may also suppose that $E_k \cap E_l =
\emptyset$ when $k \ne l$, by replacing $E_l$ with $E_l \backslash
(E_1 \cup \cdots E_{l - 1})$ when $l > 1$.  If $\lambda$ is a bounded
linear functional on $L^p(X)$, then the restriction of $\lambda$ to $f
\in L^p(X)$ such that $f = 0$ on $X \backslash E_l$ defines a bounded
linear functional on $L^p(E_l)$ for each $l$.  By the previous
argument, for each positive integer $l$, there is a $g_l \in L^q(E_l)$
such that
\begin{equation}
        \lambda(f) = \int_{E_l} f \, g_l \, d\mu
\end{equation}
for every $f \in L^p(E_l)$.  Let $g$ be the function on $X$ defined by
$g = g_l$ on $E_l$ for each $l$.  Thus the restriction of $g$ to
$\bigcup_{j = 1}^n E_l$ is in $L^q$ for each $n$, and $\lambda(f)$ is
equal to the integral of $f$ times $g$ when $f \in L^p(X)$ and $f = 0$
on $X \backslash \Big(\bigcup_{l = 1}^n E_l\Big)$.  In particular, the
$L^q$ norm of the restriction of $g$ to $\bigcup_{l = 1}^n E_l$ is
less than or equal to the dual norm of the restriction of $\lambda$ to
$L^p\Big(\bigcup_{l = 1}^n E_l\Big)$ for each $n$, which is bounded by
the dual norm of $\lambda$ on $L^p(X)$.  This implies that $g \in
L^q(X)$, with $L^q$ norm less than or equal to the dual norm of
$\lambda$ on $L^p(X)$.  Every $f \in L^p(X)$ can be approximated in
the $L^p$ norm by functions that are equal to $0$ on $X \backslash
\Big(\bigcup_{l = 1}^n E_l\Big)$ for some $n$, because $q < \infty$,
and so $\lambda(f)$ is given by the integral of $f$ times $g$ for
every $f \in L^p(X)$.

        If $1 < p < \infty$, then we can drop the hypothesis that $X$
be $\sigma$-finite.  To see this, let a bounded linear functional
$\lambda$ on $L^p(X)$ be given.  We may as well suppose that $\lambda
\ne 0$, since otherwise there is nothing to do.  In particular,
$L^p(X) \ne \{0\}$, which is to say that there are measurable subsets
of $X$ with positive finite measure.  If $Y \subseteq X$ is measurable
and $\sigma$-finite, then there is a $g_Y \in L^q(Y)$ such that
\begin{equation}
        \lambda(f) = \int_Y f \, g_Y \, d\mu
\end{equation}
for every $f \in L^p(X)$ with $f = 0$ on $X \backslash Y$, by the
previous argument.  Moreover, the $L^q$ norm of $g_Y$ is equal to the
dual norm of the restriction of $\lambda$ to $L^p(Y)$, which is less
than or equal to the dual norm of $\lambda$ on $L^p(X)$.  Let $f_1,
f_2, \ldots$ be a sequence of elements of $L^p(X)$ such that
$\|f_j\|_p = 1$ for each $j$ and $\{|\lambda(f_j)|\}_{j = 1}^\infty$
converges to the dual norm of $\lambda$ on $L^p(X)$.  Observe that
\begin{equation}
        Y_0 = \bigcup_{j = 1}^\infty \{x \in X : f_j(x) \ne 0\}
\end{equation}
is a measurable set with $\sigma$-finite measure, because the set
where $f_j \ne 0$ has this property for each $j$.  Hence there is a
$g_{Y_0} \in L^q(Y_0)$ with the properties mentioned earlier.  By
construction, the dual norm of $\lambda$ on $L^p(X)$ is equal to the
dual norm of the restriction of $\lambda$ to $L^p(Y_0)$, which is
equal to the $L^q$ norm of $g_{Y_0}$.  If $Y \subseteq X$ is
measurable and $\sigma$-finite, and if $Y_0 \subseteq Y$, then $g_Y =
g_{Y_0}$ almost everywhere on $Y_0$, by uniqueness of the
representation.  However, the $L^q$ norm of $g_Y$ is less than or
equal to the dual of norm of $\lambda$ on $L^p(X)$, which is equal to
the $L^q$ norm of $g_{Y_0}$.  This implies that $g_Y = 0$ almost
everywhere on $Y \backslash Y_0$, since $q < \infty$.  Let $g$ be the
function on $X$ equal to $g_{Y_0}$ on $Y_0$ and to $0$ on $X
\backslash Y_0$.  If $f \in L^p(X)$, then the previous argument can be
applied to
\begin{equation}
        Y = Y_0 \cup \{x \in X : f(x) \ne 0\},
\end{equation}
to get that $\lambda(f)$ is equal to the integral of $f$ times $g$, as
desired.

\section{Lengths of paths}
\label{lengths of paths}
\setcounter{equation}{0}

        Let $(M, d(x, y))$ be a metric space, and let $f$ be a
function on a closed interval $[a, b]$ in the real line with values in
$M$.  If $\mathcal{P} = \{t_j\}_{j = 0}^n$ is a partition of $[a, b]$,
in the sense that
\begin{equation}
        a = t_0 < t_1 < \cdots < t_n = b,
\end{equation}
then we put
\begin{equation}
        \Lambda_a^b(\mathcal{P}) = \sum_{j = 1}^n d(f(t_j), f(t_{j - 1})).
\end{equation}
Note that
\begin{equation}
        d(f(a), f(b)) \le \Lambda_a^b(\mathcal{P}),
\end{equation}
because of the triangle inequality. Similarly,
\begin{equation}
        \Lambda_a^b(\mathcal{P}) \le \Lambda_a^b(\mathcal{P}')
\end{equation}
when $\mathcal{P}'$ is another partition of $[a, b]$ that is a
refinement of $\mathcal{P}$, which means that $\mathcal{P}'$ includes
the points in $\mathcal{P}$.  The length $\Lambda_a^b$ of the path
$f(t)$, $a \le t \le b$, is defined to be the supremum of
$\Lambda_a^b(\mathcal{P})$ over all partitions $\mathcal{P}$ of $[a,
b]$, which may be infinite.

        Suppose that $a \le r \le b$, and that $\mathcal{P}_1$,
$\mathcal{P}_2$ are partitions of $[a, r]$, $[r, b]$, respectively.
We can combine $\mathcal{P}_1$, $\mathcal{P}_2$ to get a partition
$\mathcal{P}$ of $[a, b]$ that satisfies
\begin{equation}
\label{Lambda_a^r(mathcal{P}_1) + Lambda_r^b(mathcal{P}_2) = ...}
        \Lambda_a^r(\mathcal{P}_1) + \Lambda_r^b(\mathcal{P}_2)
         = \Lambda_a^b(\mathcal{P}).
\end{equation}
Thus
\begin{equation}
 \Lambda_a^r(\mathcal{P}_1) + \Lambda_r^b(\mathcal{P}_2) \le \Lambda_a^b,
\end{equation}
which implies that
\begin{equation}
\label{Lambda_a^r + Lambda_r^b le Lambda_a^b}
        \Lambda_a^r + \Lambda_r^b \le \Lambda_a^b,
\end{equation}
by taking the supremum over all partitions $\mathcal{P}_1$,
$\mathcal{P}_2$ of $[a, r]$, $[r, b]$.  In the other direction, if
$\mathcal{P}$ is any partition of $[a, b]$, then $\mathcal{P}$ may or
may not include $r$, but we can add $r$ to $\mathcal{P}$ if necessary
to get a refinement $\mathcal{P}'$ of $\mathcal{P}$ that does contain
$r$.  This permits $\mathcal{P}'$ to be expressed as the combination
of partitions $\mathcal{P}_1$, $\mathcal{P}_2$ of $[a, r]$, $[r, b]$,
respectively, so that
\begin{equation}
        \Lambda_a^b(\mathcal{P}) \le \Lambda_a^b(\mathcal{P}')
         = \Lambda_a^r(\mathcal{P}_1) + \Lambda_r^b(\mathcal{P}_2).
\end{equation}
Hence
\begin{equation}
        \Lambda_a^b(\mathcal{P}) \le \Lambda_a^r + \Lambda_r^b
\end{equation}
for every partition $\mathcal{P}$ of $[a, b]$, and therefore
\begin{equation}
\label{Lambda_a^b le Lambda_a^r + Lambda_r^b}
        \Lambda_a^b \le \Lambda_a^r + \Lambda_r^b.
\end{equation}
Combining this with (\ref{Lambda_a^r + Lambda_r^b le Lambda_a^b}), we get that
\begin{equation}
\label{Lambda_a^b = Lambda_a^r + Lambda_r^b}
        \Lambda_a^b = \Lambda_a^r + \Lambda_r^b.
\end{equation}
In particular,
\begin{equation}
        \Lambda_a^r \le \Lambda_a^b
\end{equation}
when $a \le r \le b$, which can be seen more directly by extending any
partition of $[a, r]$ to a partition of $[a, b]$.

        The diameter of a nonempty set $E \subseteq M$ is defined by
\begin{equation}
        \diam E = \sup \{d(x, y) : x, y \in E\},
\end{equation}
which is finite exactly when $E$ is bounded.  If $a \le r \le t \le b$
and $\mathcal{P}$ is a partition of $[a, b]$ consisting of these points, then
\begin{equation}
        d(f(r), f(t)) \le \Lambda_a^b(\mathcal{P}) \le \Lambda_a^b.
\end{equation}
It follows that
\begin{equation}
        \diam f([a, b]) \le \Lambda_a^b.
\end{equation}
Note that $\Lambda_a^b = 0$ if and only if $f$ is constant.

        Consider the special case where $M = {\bf R}$ and $f : [a, b]
\to {\bf R}$ is monotone increasing. If $\mathcal{P} = \{t_j\}_{j =
0}^n$ is any partition of $[a, b]$, then
\begin{equation}
\label{Lambda_a^b(mathcal{P}) = ...  = f(b) - f(a)}
        \Lambda_a^b(\mathcal{P}) = \sum_{j = 1}^n (f(t_j) - f(t_{j - 1}))
                                 = f(b) - f(a).
\end{equation}
This implies that
\begin{equation}
        \Lambda_a^b = f(b) - f(a).
\end{equation}

\section{Lipschitz mappings}
\label{lipschitz mappings}
\setcounter{equation}{0}

        Let $(M_1, d_1(x, y))$ and $(M_2, d_2(u, v))$ be metric
spaces.  A mapping $f : M_1 \to M_2$ is said to be \emph{Lipschitz} if
there is a constant $k \ge 0$ such that
\begin{equation}
\label{d_2(f(x), f(y)) le k d_1(x, y)}
        d_2(f(x), f(y)) \le k \, d_1(x, y)
\end{equation}
for every $x, y \in M_1$.  Thus Lipschitz mappings are automatically
uniformly continuous, and $f$ is Lipschitz with $k = 0$ if and only if
$f$ is constant.

        If $M_2$ is the real line with the standard metric, then
$f : M_1 \to {\bf R}$ is Lipschitz with constant $k$ if and only if
\begin{equation}
\label{f(x) le f(y) + k d_1(x, y)}
        f(x) \le f(y) + k \, d_1(x, y)
\end{equation}
for every $x, y \in M_1$.  More precisely, (\ref{d_2(f(x), f(y)) le k
d_1(x, y)}) implies (\ref{f(x) le f(y) + k d_1(x, y)}) directly, and
to get the converse, one can apply the latter both to $x$, $y$ and
with the roles of $x$, $y$ exchanged.  In particular,
\begin{equation}
\label{f_p(x) = d_1(x, p)}
        f_p(x) = d_1(x, p)
\end{equation}
is Lipschitz with constant $1$ on $M_1$ for every $p \in M_1$, by the
triangle inequality.  For example, $f(x) = |x|$ is Lipschitz with
constant $1$ on the real line.

        Suppose now that $f$ is a Lipschitz mapping with constant $k$
from a closed interval $[a, b]$ in the real line with the standard
metric into a metric space $(M_2, d_2(u, v))$.  If $\mathcal{P} =
\{t_j\}_{j = 0}^n$ is a partition of $[a, b]$, then
\begin{equation}
\label{Lambda_a^b(mathcal{P}) le k (b - a)}
 \Lambda_a^b(\mathcal{P}) = \sum_{j = 1}^n d_2(f(t_j), f(t_{j - 1}))
                    \le \sum_{j = 1}^n k \, (t_j - t_{j - 1}) = k \, (b - a).
\end{equation}
Thus $f$ has length $\Lambda_a^b \le k \, (b - a)$.

        If $M_1$, $M_2$, and $M_3$ are metric spaces, and $f_1 : M_1
\to M_2$, $f_2 : M_2 \to M_3$ are Lipschitz mappings with constants
$k_1$, $k_2$, respectively, then their composition $f_2 \circ f_1$ is
a Lipschitz mapping from $M_1$ into $M_2$ with constant $k_1 \, k_2$.
Similarly, if $f_1 : [a, b] \to M_2$ has length $\Lambda_a^b$ and $f_2
: M_2 \to M_3$ is Lipschitz with constant $k_2$, then $f_2 \circ f_1 :
[a, b] \to M_3$ has length $\le k_2 \, \Lambda_a^b$.

\section{Bounded variation}
\label{bounded variation}
\setcounter{equation}{0}

        A real-valued function $f$ on a closed interval $[a, b]$ in the real
line is said to have \emph{bounded variation} if it has finite length as a
mapping into ${\bf R}$ with the standard metric.  In this case, the length
of $f$ is also known as its \emph{total variation}.  We can also consider the
positive and negative variations of $f$ separately, as follows.

        For each real number $x$, put $x_+ = x$ when $x \ge 0$, $x_+ =
0$ when $x \le 0$, $x_- = -x$ when $x \le 0$, and $x_- = 0$ when $x \ge 0$.
Thus
\begin{equation}
\label{x_+ + x_- = |x|, x_+ - x_- = x}
        x_+ + x_- = |x|, \quad x_+ - x_- = x
\end{equation}
and
\begin{equation}
\label{(x + y)_+ le x_+ + y_+, (x + y)_- le x_- + y_-}
        (x + y)_+ \le x_+ + y_+, \quad (x + y)_- \le x_- + y_-
\end{equation}
for every $x, y \in {\bf R}$.  If $\mathcal{P} = \{t_j\}_{j = 0}^n$ is
a partition of $[a, b]$, then put
\begin{equation}
\label{P_a^b(mathcal{P}) = sum_{j = 1}^n (f(t_j) - f(t_{j - 1}))_+}
        P_a^b(\mathcal{P}) = \sum_{j = 1}^n (f(t_j) - f(t_{j - 1}))_+
\end{equation}
and
\begin{equation}
\label{N_a^b(mathcal{P}) = sum_{j = 1}^n (f(t_j) - f(t_{j - 1}))_-}
        N_a^b(\mathcal{P}) = \sum_{j = 1}^n (f(t_j) - f(t_{j - 1}))_-.
\end{equation}
Note that
\begin{equation}
\label{P_a^b(mathcal{P}) + N_a^b(mathcal{P}) = Lambda_a^b(mathcal{P})}
        P_a^b(\mathcal{P}) + N_a^b(\mathcal{P}) = \Lambda_a^b(\mathcal{P})
\end{equation}
and
\begin{equation}
\label{P_a^b(mathcal{P}) - N_a^b(mathcal{P}) = f(b) - f(a)}
        P_a^b(\mathcal{P}) - N_a^b(\mathcal{P}) = f(b) - f(a),
\end{equation}
by (\ref{x_+ + x_- = |x|, x_+ - x_- = x}).  If $\mathcal{P}'$ is
another partition of $[a, b]$ which is a refinement of $\mathcal{P}$,
then it is easy to see that
\begin{equation}
\label{P_a^b({P}) le P_a^b({P}'), N_a^b({P}) le N_a^b({P}')}
        P_a^b(\mathcal{P}) \le P_a^b(\mathcal{P}'), \quad 
         N_a^b(\mathcal{P}) \le N_a^b(\mathcal{P}'),
\end{equation}
using (\ref{(x + y)_+ le x_+ + y_+, (x + y)_- le x_- + y_-}).

        Let $P_a^b$, $N_a^b$ be the suprema of $P_a^b(\mathcal{P})$,
$N_a^b(\mathcal{P})$ over all partitions $\mathcal{P}$ of $[a, b]$,
respectively.  Clearly
\begin{equation}
        \Lambda_a^b \le P_a^b + N_a^b,
\end{equation}
by (\ref{P_a^b(mathcal{P}) + N_a^b(mathcal{P}) = Lambda_a^b(mathcal{P})}).
To get the opposite inequality
\begin{equation}
\label{P_a^b + N_a^b le Lambda_a^b}
        P_a^b + N_a^b \le \Lambda_a^b,
\end{equation}
one should be a bit more careful, because the partitions $\mathcal{P}$
of $[a, b]$ for which $P_a^b(\mathcal{P})$ approaches $P_a^b$ may not
be the same as the partitions for which $N_a^b(\mathcal{P})$
approaches $N_a^b$.  However, using common refinements of such
partitions, one can get partitions $\mathcal{P}$ such that
$P_a^b(\mathcal{P})$, $N_a^b(\mathcal{P})$ approach $P_a^b$, $N_a^b$
at the same time.  This implies (\ref{P_a^b + N_a^b le Lambda_a^b}),
from which it follows that
\begin{equation}
        P_a^b + N_a^b = \Lambda_a^b.
\end{equation}
Observe also that
\begin{equation}
        P_a^r + P_r^b = P_a^b, \quad N_a^r + N_r^b = N_a^b
\end{equation}
for each $r$, $a \le r \le b$.  This uses the same arguments as for
$\Lambda_a^b$, in Section \ref{lengths of paths}.

        Suppose now that $f$ has bounded variation, so that
$\Lambda_a^b < \infty$, and hence $P_a^b, N_a^b < \infty$.  Using
(\ref{P_a^b(mathcal{P}) - N_a^b(mathcal{P}) = f(b) - f(a)}), one can
check that
\begin{equation}
\label{P_a^b - N_a^b = f(b) - f(a)}
        P_a^b - N_a^b = f(b) - f(a).
\end{equation}
More precisely, one should be careful to use partitions $\mathcal{P}$
of $[a, b]$ such that $P_a^b$, $N_a^b$ are simultaneously approximated
by $P_a^b(\mathcal{P})$, $N_a^b(\mathcal{P})$, respectively, as in the
previous paragraph.  Similarly,
\begin{equation}
\label{P_a^r - N_a^r = f(r) - f(a)}
        P_a^r - N_a^r = f(r) - f(a)
\end{equation}
for each $r \in [a, b]$, since the restriction of $f$ to $[a, r]$ also
has bounded variation.  Of course, $P_a^r$ and $N_a^r$ are monotone
increasing on $[a, b]$.

\section{Functions and measures}
\label{functions, measures}
\setcounter{equation}{0}

        Let $\alpha(x)$ be a monotone increasing real-valued function
on the real line.  As usual, the one-sided limits $\alpha(x+) =
\lim_{y \to x+} \alpha(y)$, $\alpha(x-) = \lim_{z \to x-} \alpha(z)$
exist for every $x \in {\bf R}$, and are given by
\begin{eqnarray}
        \alpha(x+) & = & \sup \{\alpha(y) : y \in {\bf R}, \, y < x\}, \\
        \alpha(x-) & = & \inf \{\alpha(z) : z \in {\bf R}, \, x < z\}.
\end{eqnarray}
Thus
\begin{equation}
\label{alpha(x-) le alpha(x) le alpha(x+)}
        \alpha(x-) \le \alpha(x) \le \alpha(x+)
\end{equation}
for every $x \in {\bf R}$, and $\alpha(x+) = \alpha(x-)$ exactly when
$\alpha$ is continuous at $x$.  Moreover,
\begin{equation}
        \alpha(x+) \le \alpha(y-)
\end{equation}
for every $x, y \in {\bf R}$ with $x < y$.  Remember that the set of
$x \in {\bf R}$ at which $\alpha$ is not continuous has only finitely
or countably many elements.

        It is well known that there is a unique positive Borel measure
$\mu_\alpha$ on ${\bf R}$ that satisfies
\begin{equation}
        \mu_\alpha((a, b)) = \alpha(b-) - \alpha(a+), \quad  
         \mu_\alpha([a, b]) = \alpha(b+) - \alpha(a-)
\end{equation}
for every $a, b \in R$ with $a < b$.  The expression for closed
intervals also makes sense when $a = b$, in which case it reduces to
\begin{equation}
        \alpha(\{a\}) = \alpha(a+) - \alpha(a-).
\end{equation}
Of course, this is equal to $0$ when $\alpha$ is continuous at $a$.
Alternatively, if $f$ is a continuous real-valued function on the real
line with compact support, then one can define the Riemann--Stieltjes
integral
\begin{equation}
\label{int_{-infty}^infty f(x) d alpha(x)}
        \int_{-\infty}^\infty f(x) \, d\alpha(x).
\end{equation}
This is a nonnegative linear functional on the space of continuous
functions with compact support on ${\bf R}$, and the Riesz
representation theorem leads to a positive Borel measure that is the
same as $\mu_\alpha$.  As another approach, if $\alpha$ is a strictly
increasing continuous function on ${\bf R}$, then one can get
$\mu_\alpha$ from Lebesgue measure using a change of variables.  If
$\alpha$ is monotone increasing and continuous, but perhaps not
strictly increasing, then
\begin{equation}
        \beta(x) = \alpha(x) + x
\end{equation}
is continuous and strictly increasing, the previous argument can be
used to get $\mu_\beta$, and one can get $\mu_\alpha$ by subtracting
Lebesgue measure from $\mu_\beta$.  If $\alpha$ is not continuous,
then one can account for the discontinuities directly with sums of
multiples of Dirac masses.

        Let us say that a real-valued function $\alpha$ on ${\bf R}$
has bounded variation if it has bounded variation on every closed
interval $[a, b]$, and if the total variation $\Lambda_a^b$ of
$\alpha$ on $[a, b]$ is uniformly bounded.  This implies that $\alpha$
is bounded on ${\bf R}$, since
\begin{equation}
\label{|alpha(a) - alpha(b)| le Lambda_a^b}
        |\alpha(a) - \alpha(b)| \le \Lambda_a^b
\end{equation}
for every $a, b \in {\bf R}$ with $a \le b$.  It is easy to see that
bounded monotone functions on ${\bf R}$ have bounded variation.
Conversely, one can check that a function with bounded variation on
${\bf R}$ can be expressed as a difference of monotone increasing
functions that are bounded.  Complex-valued functions of bounded
variation on ${\bf R}$ can be defined analogously, and represented as
linear combinations of bounded monotone real-valued functions.

        If $\alpha$ is a real or complex-valued function of bounded
variation on ${\bf R}$, then there is a real or complex measure Borel
measure $\mu_\alpha$ on ${\bf R}$ associated to $\alpha$ as before.
More precisely, if $\alpha$ is given as a linear combination of
bounded monotone increasing real-valued functions, then $\mu_\alpha$
is the same as the corresponding linear combination of positive finite
measures.  In this case, the Riemann-Stieltjes integral
(\ref{int_{-infty}^infty f(x) d alpha(x)}) defines a bounded linear
functional on the space of continuous functions on ${\bf R}$ with
compact support with respect to the supremum norm, which leads to a
real or complex Borel measure on ${\bf R}$, as appropriate.

\section{Continuity conditions}
\label{continuity conditions}
\setcounter{equation}{0}

        Let $(M, d(x, y))$ be a complete metric space, and let $f :
[a, b] \to M$ be a path of finite length $\Lambda_a^b$.  If
$\{t_j\}_{j = 1}^\infty$ is a monotone sequence of elements of $[a,
b]$, then it is easy to see that
\begin{equation}
        \sum_{j = 1}^n d(f(t_j), f(t_{j + 1})) \le \Lambda_a^b
\end{equation}
for every positive integer $n$.  This implies that $\sum_{j =
1}^\infty d(f(t_j), f(t_{j + 1}))$ converges, and hence that
$\{f(t_j)\}_{j = 1}^\infty$ converges in $M$, as in Section
\ref{metric spaces}.  Using this, one can check that $f(r+) = \lim_{t
\to r+} f(t)$ exists for every $r \in [a, b)$, and similarly that
$f(r-) = \lim_{t \to r-} f(t)$ exists for every $r \in (a, b]$.  More
precisely, this also uses the observation that two strictly increasing
or two strictly decreasing sequences with the same limit can be
combined into a single monotone sequence, and hence that the
corresponding sequences of values of $f$ have the same limit in $M$.

        Alternatively, let $\Lambda_u^v$ be the length of the
restriction of $f$ to $[u, v]$ when $a \le u \le v \le b$.  Of course,
$\Lambda_a^r$ is monotone increasing in $r$, and hence
\begin{equation}
\label{lim_{t to r-} Lambda_a^t = sup_{a le t < r} Lambda_a^t}
        \lim_{t \to r-} \Lambda_a^t = \sup_{a \le t < r} \Lambda_a^t
\end{equation}
when $a < r \le b$.  Let $\epsilon > 0$ be given, and choose $u \in
[a, r)$ so that
\begin{equation}
        \Lambda_a^u > \sup_{a \le t < r} \Lambda_a^t - \epsilon.
\end{equation}
Because $\Lambda_a^t = \Lambda_a^u + \Lambda_u^t$ when $u \le t < r$,
we get that
\begin{equation}
        \sup_{u \le t < r} \Lambda_u^t < \epsilon.
\end{equation}
One can also use this to deal with $f(r-)$, and similarly for $f(r+)$
when $a \le r < b$.

        If $a \le r \le t \le b$, then
\begin{equation}
        d(f(r), f(t)) \le \Lambda_r^t,
\end{equation}
as usual.  It follows that $f$ is continuous on the right at $r \in
[a, b)$ when
\begin{equation}
        \lim_{t \to r+} \Lambda_r^t = 0,
\end{equation}
and that $f$ is continuous from the left at $r \in (a, b]$ when
\begin{equation}
        \lim_{t \to r-} \Lambda_t^r = 0.
\end{equation}
Equivalently, continuity of $\Lambda_a^r$ from the right or the left
implies continuity of $f(r)$ from the right or the left at the same
point, respectively.  In particular, $f(r)$ is continuous at every
point where $\Lambda_a^r$ is continuous, which includes all but at
most finitely or countably many elements of $[a, b]$, because
$\Lambda_a^r$ is monotone increasing in $r$.

        Conversely, $\Lambda_a^r$ is continuous from the right or left
at any point where $f$ is continuous from the right or left.  To see
this, let $r \in (a, b]$ and $\epsilon > 0$ be given, and let
$\mathcal{P} = \{t_j\}_{j = 0}^n$ be a partition of $[a, r]$ such that
\begin{equation}
\label{Lambda_a^r(mathcal{P}) > Lambda_a^r - epsilon}
        \Lambda_a^r(\mathcal{P}) > \Lambda_a^r - \epsilon.
\end{equation}
If $t_{n - 1} < t < t_n = r$, then let $\mathcal{P}_t$ be the
partition of $[a, r]$ obtained by adding $t$ between $t_{n - 1}$ and
$t_n = r$ in $\mathcal{P}$.  Thus $\mathcal{P}_t$ is a refinement of
$\mathcal{P}$, so that
\begin{equation}
        \Lambda_a^r(\mathcal{P}_t) \ge \Lambda_a^r(\mathcal{P}).
\end{equation}
We can also consider $\mathcal{P}_t$ as the combination of a partition
of $[0, t]$ with a single step from $t$ to $r$, which implies that
\begin{equation}
        \Lambda_a^r(\mathcal{P}_t) \le \Lambda_a^t + d(f(t), f(r)).
\end{equation}
Hence
\begin{equation}
        \Lambda_a^t + d(f(t), f(r)) > \Lambda_a^r - \epsilon
\end{equation}
when $t_{n - 1} < t < r$.  This shows that $\Lambda_a^r$ is continuous
from the left at $r$ when $f(r)$ is continuous from the left at $r$,
using also the fact that $\Lambda_a^t \le \Lambda_a^r$ when $a \le t
\le r$.  The argument for continuity on the right is very similar.

\section{Maximal functions}
\label{maximal functions}
\setcounter{equation}{0}

        Let $\mu$ be a positive finite Borel measure on the real line.
The \emph{Hardy--Littlewood maximal function} associated to $\mu$ is defined by
\begin{equation}
\label{mu^*(x) = sup_{x in I} frac{mu(I)}{|I|}}
        \mu^*(x) = \sup_{x \in I} \frac{\mu(I)}{|I|},
\end{equation}
where the supremum is taken over all open intervals $(a, b)$ that
contain $x$, and $|I| = b - a$ is the length of $I$.  Put
\begin{equation}
        E_t = \{x \in {\bf R} : \mu^*(x) > t\}
\end{equation}
for each $t > 0$.  Thus $x \in E_t$ if and only if there is an open
interval $I$ such that $x \in I$ and
\begin{equation}
        \mu(I) > t \, |I|.
\end{equation}
In this case, $I \subseteq E_t$, and it follows that $E_t$ is an open
set in ${\bf R}$.

        Suppose that $K \subseteq E_t$ is compact.  This implies that
there are finitely many open intervals $I_1, \ldots, I_n$ in ${\bf R}$
such that
\begin{equation}
        K \subseteq \bigcup_{j = 1}^n I_j
\end{equation}
and
\begin{equation}
        \mu(I_j) > t \, |I_j|
\end{equation}
for each $j$.  A basic property of the real line is that for any three
intervals with a point in common, one of the intervals is contained in
the union of the other two.  This permits us to reduce the collection
of intervals $I_1, \ldots, I_n$ in such a way that no element of ${\bf
R}$ is contained in more than two of these intervals.

        It follows that
\begin{equation}
        \sum_{j = 1}^n |I_j| < t^{-1} \sum_{j = 1}^n \mu(I_j)
         \le 2 \, t^{-1} \, \mu\Big(\bigcup_{j = 1}^n I_j\Big).
\end{equation}
More precisely, if ${\bf 1}_A$ is the indicator function on ${\bf R}$
associated to $A \subseteq {\bf R}$, then
\begin{equation}
        \sum_{j = 1}^n \mu(I_j)
         = \int_{\bf R} \Big(\sum_{j = 1}^n {\bf 1}_{I_j}\Big) \, d\mu
         \le  \int_{\bf R} 2 \, {\bf 1}_{\bigcup_{j = 1}^n I_j} \, d\mu
          = 2 \, \mu\Big(\bigcup_{j = 1}^n I_j\Big).
\end{equation}
If $|K|$ denotes the Lebesgue measure of $K$, then we get that
\begin{equation}
        |K| \le 2 \, t^{-1} \, \mu({\bf R}).
\end{equation}
Hence
\begin{equation}
        |E_t| \le 2 \, t^{-1} \, \mu({\bf R}),
\end{equation}
because $K$ is an arbitrary compact subset of $E_t$.

        If $f$ is an integrable function on ${\bf R}$, then we put
\begin{equation}
\label{f^*(x)}
        f^*(x) = \sup_{x \in I} \frac{1}{|I|} \int_I |f(y)| \, dy.
\end{equation}
This is the same as the maximal function $\mu^*(x)$ associated to the measure
\begin{equation}
        \mu(A) = \int_A |f(y)| \, dy.
\end{equation}
Thus the estimate in the previous paragraph can be re-expressed in this case as
\begin{equation}
 |\{x \in {\bf R} : f^*(x) > t\}| \le 2 \, t^{-1} \, \int_{\bf R} |f(y)| \, dy
\end{equation}
for each $t > 0$.

\section{Lebesgue's theorem}
\label{lebesgue's theorem}
\setcounter{equation}{0}

        Let $f$ be a locally integrable function on the real line.
A famous theorem of Lebesgue implies that
\begin{equation}
 \lim_{r \to 0} \frac{1}{2 r} \int_{x - r}^{x + r} |f(y) - f(x)| \, dy = 0
\end{equation}
for almost every $x \in {\bf R}$.  We may as well suppose that $f$ is
integrable on ${\bf R}$, since the problem is local.

        Put
\begin{eqnarray}
 L(f)(x) & = & \limsup_{r \to 0} \frac{1}{2 r} \int_{x - r}^{x + r}
                                                 |f(y) - f(x)| \, dy \\
         & = & \lim_{\epsilon \to 0} \, \sup_{0 < r < \epsilon} \frac{1}{2 r}
                         \int_{x - r}^{x + r} |f(y) - f(x)| \, dy.  \nonumber
\end{eqnarray}
Observe that
\begin{equation}
        L(f_1 + f_2)(x) \le L(f_1)(x) + L(f_2)(x),
\end{equation}
and that
\begin{equation}
        L(g)(x) = 0
\end{equation}
when $g$ is continuous at $x$.  It follows that
\begin{equation}
        L(f) = L(f - g)
\end{equation}
for every continuous function $g$.

        We also have that
\begin{equation}
        L(f)(x) \le f^*(x) + |f(x)|,
\end{equation}
where $f^*(x)$ is as in (\ref{f^*(x)}).  This implies that
\begin{equation}
        L(f)(x) \le (f - g)^*(x) + |f(x) - g(x)|
\end{equation}
for every continuous function $g$ on ${\bf R}$.  Hence
\begin{eqnarray}
\lefteqn{\{x \in {\bf R} : L(f)(x) > t\}} \\
         & \subseteq & \{x \in {\bf R} : (f - g)^*(x) > t/2\}
          \cup \{x \in {\bf R} : |f(x) - g(x)| > t/2\} \nonumber
\end{eqnarray}
for every $t > 0$.

        As in the previous section,
\begin{equation}
        \quad |\{x \in {\bf R} : (f - g)^*(x) > t/2\}|
          \le 2 \, (t / 2)^{-1} \, \|f - g\|_1 = 4 \, t^{-1} \, \|f - g\|_1.
\end{equation}
Similarly,
\begin{eqnarray}
\lefteqn{|\{x \in {\bf R} : |f(x) - g(x)| > t/2\}|} \\
        & \le & (t/2)^{-1} \int_{\bf R} |f(y) - g(y)| \, dy 
         = 2 \, t^{-1} \, \|f - g\|_1 \nonumber
\end{eqnarray}
for every $t > 0$.  Of course, we can choose $g$ so that $\|f - g\|_1$ is
arbitrarily small, because continuous functions are dense in $L^1({\bf R})$.
Using this, one can show that $L(f)(x) = 0$ almost everywhere, as desired.

\section{Singular measures}
\label{singular measures}
\setcounter{equation}{0}

        Let $\mu$ be a positive finite Borel measure on the real line
which is singular with respect to Lebesgue measure.  This means that
there is a Borel set $B \subseteq {\bf R}$ whose Lebesgue measure
$|B|$ is $0$ while $\mu({\bf R} \backslash B) = 0$.  Let us check that
\begin{equation}
        \lim_{r \to 0} \frac{\mu((x - r, x + r))}{2 r} = 0
\end{equation}
for almost every $x \in {\bf R}$ with respect to Lebesgue measure.  If
$B$ happens to be a closed set in ${\bf R}$, then this holds trivially
for every $x \in {\bf R} \backslash B$.  The idea is to use the
maximal function to make an approximation by this type of situation.

        Consider
\begin{eqnarray}
 L(\mu)(x) & = & \limsup_{r \to 0} \frac{\mu((x - r, x+ r))}{2 \, r} \\
           & = & \lim_{\epsilon \to 0} \, \sup_{0 < r < \epsilon}
                            \frac{\mu((x - r, x + r))}{2 r}, \nonumber
\end{eqnarray}
in analogy with the previous section.  Thus 
\begin{equation}
        L(\mu)(x) \le \mu^*(x)
\end{equation}
and
\begin{equation}
        L(\mu_1 + \mu_2)(x) \le L(\mu_1)(x) + L(\mu_2)(x)
\end{equation}
for every pair of positive Borel measures $\mu_1$, $\mu_2$ on ${\bf R}$.

        Let $U$ be an open set in ${\bf R}$ such that $B \subseteq U$,
and let $K$ be a compact set in ${\bf R}$ such that $K \subseteq U$.
Also let $\mu_1$, $\mu_2$ be the Borel measures on ${\bf R}$ defined by
\begin{equation}
 \mu_1(A) = \mu(A \cap K), \quad \mu_2(A) = \mu(A \cap ({\bf R} \backslash K)).
\end{equation}
Thus $L(\mu_1)(x) = 0$ when $x \in {\bf R} \backslash K$, which
implies that
\begin{equation}
        L(\mu)(x) \le L(\mu_2)(x) \le \mu_2^*(x)
\end{equation}
for every $x \in {\bf R} \backslash K$, and hence for every $x \in
{\bf R} \backslash U$.  The main point now is to choose $K \subseteq
U$ so that
\begin{equation}
\label{mu(U backslash K) = mu({bf R} backslash K) = mu_2({bf R})}
        \mu(U \backslash K) = \mu({\bf R} \backslash K) = \mu_2({\bf R})
\end{equation}
is arbitrarily small.  This is easy to do, using the fact that open
subsets of the real line are $\sigma$-compact.  This implies that
$L(\mu)(x) = 0$ for Lebesgue almost every $x \in {\bf R} \backslash
U$, by the maximal function estimates in Section \ref{maximal
functions}.  More precisely,
\begin{equation}
\label{{x in {bf R} - U : L(mu)(x) > t} subseteq {x in {bf R} : mu_2^*(x) > t}}
        \{x \in {\bf R} \backslash U : L(\mu)(x) > t\}
         \subseteq \{x \in {\bf R} : \mu_2^*(x) > t\}
\end{equation}
for every $t > 0$, and the Lebesgue measure of the set on the right
can be made arbitrarily small, by choosing $K$ so that (\ref{mu(U
backslash K) = mu({bf R} backslash K) = mu_2({bf R})}) is small.  This
implies that $L(\mu)(x) \le t$ almost everywhere on ${\bf R}
\backslash U$ with respect to Lebesgue measure for each $t > 0$, and
hence that $L(\mu)(x) = 0$ almost everywhere on ${\bf R} \backslash
U$, by taking $t = 1/n$, where $n$ is a positive integer.  It follows
that $L(\mu)(x) = 0$ for Lebesgue almost every $x \in {\bf R}$, as
desired, since we can also choose $U$ to have arbitrarily small
Lebesgue measure, because $|B| = 0$.

\section{Differentiability almost everywhere}
\label{differentiability almost everywhere}
\setcounter{equation}{0}

        Let $\alpha$ be a bounded real-valued monotone increasing
function on the real line, and let $\mu_\alpha$ be the corresponding
positive Borel measure on ${\bf R}$, as in Section \ref{functions,
measures}.  Using the Lebesgue decomposition and Radon--Nikodym
theorem, we get an integrable function $f$ with respect to Lebesgue
measure and a Borel measure $\nu$ that is singular with respect to
Lebesgue measure such that
\begin{equation}
        \mu_\alpha(A) = \int_A f(y) \, dy + \nu(A).
\end{equation}
We would like to show that $\alpha(x)$ is differentiable almost
everywhere on ${\bf R}$ with respect to Lebesgue measure, and more
precisely that $\alpha'(x) = f(x)$ almost everywhere.

        Thus we would like to show that
\begin{equation}
        \lim_{h \to 0} \frac{\alpha(x + h) - \alpha(x)}{h} = f(x)
\end{equation}
for almost every $x \in {\bf R}$.  As a first approximation, we have that
\begin{equation}
        \lim_{h \to 0} \frac{1}{h} \int_x^{x + h} f(y) \, dy = f(x)
\end{equation}
for almost every $x \in {\bf R}$, by Lebesgue's theorem.  More
precisely, this integral is supposed to be oriented, as in calculus,
so that the integral from $x$ to $x + h$ is $-1$ times the integral
from $x + h$ to $x$.  This means that we are looking at the average of
$f$ over the interval $[x, x + h]$ when $h > 0$, and over $[x + h, x]$
when $h < 0$.

        It remains to show that
\begin{equation}
 \frac{\alpha(x + h) - \alpha(x)}{h} - \frac{1}{h} \int_x^{x + h} f(y) \, dy
\end{equation}
converges to $0$ as $h \to 0$ for almost every $x \in {\bf R}$.  If
$\alpha$ is continuous at $x$ and $x + h$, then this difference is
equal to $\nu([x, x + h])/ h$ when $h > 0$, and similarly when $h <
0$.  In any case, this difference is nonnegative, bounded by $\nu([x,
x + h])/ h$ when $h > 0$, and similarly for $h < 0$.  Hence the
difference converges to $0$ almost everywhere, as in the previous
section.

        Of course, it is not important that $\alpha$ be bounded or
defined on the whole line, since the problem is local.  If $\alpha$ is
a real or complex-valued function of bounded variation on ${\bf R}$,
then $\alpha$ can be expressed as a linear combination of monotone
functions, and is therefore differentiable almost everywhere too.

\section{Maximal functions, 2}
\label{maximal functions, 2}
\setcounter{equation}{0}

        The maximal function of a positive Borel measure $\mu$ on
${\bf R}$ can also be given by
\begin{equation}
        \mu^*(x) = \sup_{x \in I} \frac{\mu(I)}{|I|},
\end{equation}
where now the supremum is taken over all closed intervals $I = [a, b]$
that contain $x$ and have positive length $|I| = b - a$.  The previous
definition is clearly less than or equal to this one, since every open
interval $(a, b)$ is contained in a closed interval $[a, b]$ with the
same length, and
\begin{equation}
        \mu((a, b)) \le \mu([a, b]).
\end{equation}
In the other direction, one can approximate closed intervals by open
intervals that contain them.

        Let $\alpha$ be a bounded monotone increasing real-valued
function on the real line.  If $\mu_\alpha$ is the corresponding
measure, as in Section \ref{functions, measures}, then its maximal
function can be expressed directly in terms of $\alpha$, by
\begin{equation}
        \mu_\alpha^*(x) = \sup_{a \le x \le b \atop a < b}
                                \frac{\alpha(b) - \alpha(a)}{b - a}.
\end{equation}
More precisely, the supremum is taken over $a, b \in {\bf R}$ with $a
\le x \le b$ and $a < b$, and this expression for the maximal function
is trapped between the previous two, by (\ref{alpha(x-) le alpha(x) le
alpha(x+)}).  If $E_t = \{x \in {\bf R} : \mu_\alpha^*(x) > t\}$, then
the main estimate from Section \ref{maximal functions} can be
reformulated as
\begin{equation}
\label{|E_t| le 2 t^{-1} (sup_{x in R} alpha(x) - inf_{x in R} alpha(x))}
        |E_t| \le 2 \, t^{-1} \, (\sup_{x \in {\bf R}} \alpha(x)
                                     - \inf_{x \in {\bf R}} \alpha(x)).
\end{equation}

        Now let $(M, d(x, y))$ be a metric space, and let $f : [a, b]
\to M$ be a path of finite length.  Let $\alpha(r)$ be the length
$\Lambda_a^r$ of the restriction of $f$ to $[a, r]$ when $a \le r \le
b$, and put $\alpha(r) = 0$ when $r < a$, $\alpha(r) = \Lambda_a^b$
when $r > b$.  Thus $\alpha$ is a bounded monotone increasing function
on ${\bf R}$, and
\begin{equation}
\label{d(f(r), f(r')) le Lambda_r^{r'} = alpha(r') - alpha(r)}
        d(f(r), f(r')) \le \Lambda_r^{r'} = \alpha(r') - \alpha(r)
\end{equation}
when $a \le r \le r' \le b$.  If $[r, r']$ contains an element of
${\bf R} \backslash E_t$, where $t > 0$ and $E_t$ is as in the
previous paragraph, then
\begin{equation}
        d(f(r), f(r')) \le \alpha(r') - \alpha(r) \le t \, (r' - r).
\end{equation}
In particular, the restriction of $f$ to $[a, b] \backslash E_t$ is
Lipschitz with constant $t$.  Note that $[a, b] \backslash E_t$ is a
closed set, because $E_t$ is open.  Also, (\ref{|E_t| le 2 t^{-1}
(sup_{x in R} alpha(x) - inf_{x in R} alpha(x))}) reduces to
\begin{equation}
        |E_t| \le 2 \, t^{-1} \, \Lambda_a^b.
\end{equation}

        If our metric space is a real or complex vector space with a
norm, then we can extend the restriction of $f$ to $[a, b] \backslash
E_t$ to a $t$-Lipschitz function $f_t$ on $[a, b]$.  Remember that
$E_t$ can be expressed as the union of finitely or countably many
pairwise-disjoint open intervals, since $E_t$ is an open set in ${\bf
R}$.  If $I$ is one of these open intervals and $I \subseteq [a, b]$,
then $f_t$ is defined on $I$ as the affine function that agrees with
$f$ on the endpoints.  If $a$ or $b$ is an element of $E_t$, and $I$
is an open interval in $E_t$ that contains $a$ or $b$ and whose other
endpoint is in $[a, b]$, then we can take $f_t$ to be the constant on
$I \cap [a, b]$ that agrees with $f$ at the other endpoint of $I$.  Of
course, if $[a, b] \subseteq E_t$, then there is nothing to do.

\section{Vector-valued functions}
\label{vector-valued functions}
\setcounter{equation}{0}

        Let $V$ be a real or complex vector space with a norm.  As
usual, a function $F : [a, b] \to V$ is said to be differentiable at
$x \in (a, b)$ if
\begin{equation}
        \lim_{h \to 0} \frac{F(x + h) - F(x)}{h}
\end{equation}
exists in $V$.  One can also consider one-sided limits at the endpoints.

        For example, let $V$ be $L^1([0, 1])$, with respect to
Lebesgue measure.  Let $F(x)$ be the indicator function of $[0, x]$ as
an element of $L^1([0, 1])$ for each $x \in [0, 1]$.  It is easy to
see that
\begin{equation}
        \|F(x) - F(y)\|_1 = |x - y|
\end{equation}
for every $x, y \in [0, 1]$, so that $F$ is actually an isometric
embedding of $[0, 1]$ in $L^1([0, 1])$.  However, one can also check
that $F$ is not differentiable at any point in $[0, 1]$.  The
derivative of $F$ at $x \in [0, 1]$ is basically a Dirac mass at $x$,
in a weak sense that we shall discuss later.

        Now let $V$ be $L^\infty({\bf R})$.  If $f$ is a bounded real
or complex-valued Lipschitz function on ${\bf R}$, then let $F : {\bf
R} \to L^\infty({\bf R})$ be the mapping that sends $x \in {\bf R}$ to
the translate $f_x(\cdot) = f(\cdot - x)$ of $f$ by $x$.  It is easy
to see that this is a Lipschitz mapping from the real line into
$L^\infty({\bf R})$, because $f$ is a Lipschitz function on ${\bf R}$.
If $F$ is differentiable at any point in ${\bf R}$ as a mapping into
$L^\infty({\bf R})$, then the difference quotient for $f$ would
converge uniformly on ${\bf R}$.  This would imply that $f$ is
continuously differentiable on ${\bf R}$, with uniformly continuous
derivative.  Conversely, if $f$ is continuously differentiable on
${\bf R}$, with uniformly continuous derivative, then the difference
quotient for $f$ does converge uniformly to the derivative of $f$, and
$F$ is differentiable at every point in ${\bf R}$.  More precisely,
the derivative of $F$ at $x \in {\bf R}$ corresponds to $-1$ times the
derivative of $f$ translated by $x$ in this case.  If $F$ is not
bounded, then one can take $F(x) = f_x - f$, and get similar
conclusions.

        Let $V$ be any vector space with a norm $\|v\|$ again, and
suppose that $F$, $G$ are $V$-valued functions on an interval $[a, b]$
with finite length.  One can check that $F - G$ also has finite length
on $[a, b]$, which is less than or equal to the sum of the lengths of
$F$ and $G$.  It follows that $\|F - G\|$ has finite length as a
real-valued function on $[a, b]$, which is to say that it has bounded
variation.  In particular, $\|F - G\|$ is differentiable almost
everywhere as a real-valued function on $[a, b]$.  If $x \in [a, b]$
is a limit point of the set where $F = G$, and hence a limit point of
the set where $\|F - G\| = 0$, and if $\|F - G\|$ is differentiable at
$x$, then the derivative of $\|F - G\|$ at $x$ is equal to $0$.  This
implies that the derivative of $F - G$ exists at $x$ and is equal to
$0$, under these conditions.  In particular, this can be applied to
Lipschitz approximations $G$ of $F$ as in the previous section.

\section{Uniform boundedness, 4}
\label{uniform boundedness, 4}
\setcounter{equation}{0}

        Let $W$ be a real or complex vector space with a norm $\|w\|$,
and let $\{\lambda_j\}_{j = 1}^\infty$ be a sequence of bounded linear
functionals on $W$.  Suppose that the dual norms of the $\lambda_j$'s
are uniformly bounded, so that
\begin{equation}
\label{||lambda_j||_* le L}
        \|\lambda_j\|_* \le L
\end{equation}
for some $L \ge 0$ and each $j$.  Under these conditions, one can
check that the set of $w \in W$ such that $\{\lambda_j(w)\}_{j =
1}^\infty$ is a Cauchy sequence in ${\bf R}$ or ${\bf C}$, as
appropriate, is closed.  Because of the completeness of the real and
complex numbers, this is the same as saying that the set of $w \in W$
such that $\{\lambda_j(w)\}_{j = 1}^\infty$ converges in ${\bf R}$ or
${\bf C}$ is closed.  It is easy to see that this is also a linear
subspace of $W$.

        In particular, $\{\lambda_j(w)\}_{j = 1}^\infty$ converges for
every $w \in W$ if it converges for a set of $w$'s whose linear span
is dense in $W$.  In this case,
\begin{equation}
        \lambda(w) = \lim_{j \to \infty} \lambda_j(w)
\end{equation}
defines a linear functional on $W$.  More precisely, $\lambda$ is a
bounded linear functional on $W$, with
\begin{equation}
\label{||lambda||_* le L}
        \|\lambda\|_* \le L,
\end{equation}
because of (\ref{||lambda_j||_* le L}).

        Conversely, if $\{\lambda_j(w)\}_{j = 1}^\infty$ converges for
every $w \in W$, then $\{\lambda_j(w)\}_{j = 1}^\infty$ is bounded for
every $w \in W$.  The Banach--Steinhaus theorem implies that the
$\lambda_j$'s have uniformly bounded dual norms when $W$ is complete,
as in Section \ref{uniform boundedness, 2}.

        Suppose now that $E$ is a set of real numbers, and that for
each $t \in E$ we have a bounded linear functional $\lambda_t$ on $W$.
Suppose also that $0$ is a limit point of $E$ in ${\bf R}$, and that
the $\lambda_t$'s have uniformly bounded dual norms.  If
\begin{equation}
        \lim_{t \to 0 \atop t \in E} \lambda_t(w)
\end{equation}
exists in ${\bf R}$ or ${\bf C}$, as appropriate, for a set of $w \in
W$ whose linear span is dense in $W$, then this limit exists for every
$w \in W$, and determines a bounded linear functional on $W$.  This is
a variant of the earlier discussion for sequences.  One can also apply
the previous remarks to sequences of elements of $E$ that converge to
$0$.

\section{Weak$^*$ derivatives}
\label{weak* derivatives}
\setcounter{equation}{0}

        Let $W$ be a real or complex vector space with a norm $\|w\|$,
and let $F(x)$ be a function on a closed interval $[a, b]$ in the real
line with values in the dual $W^*$ of $W$.  Thus $F(x)(w)$ is a real
or complex-valued function of $x$ on $[a, b]$ for each $w \in W$, as
appropriate.  If $F(x)$ has finite length as a mapping from $[a, b]$
into $W^*$, then $F(x)(w)$ has bounded variation as a real or
complex-valued function of $x$ on $[a, b]$ for every $w \in W$.  This
implies that for each $w \in W$ there is a set $Z(w) \subseteq [a, b]$
of Lebesgue measure $0$ such that $F(x)(w)$ is differentiable for
every $x \in [a, b] \backslash Z(w)$.

        Suppose that $W$ is separable, so that there is a collection
$\{w_l\}_l$ of finitely or countably many elements of $W$ whose linear
span is dense in $W$.  Thus $Z = \bigcup_l Z(w_l)$ also has Lebesgue
measure $0$.  If $x \in [a, b] \backslash Z$, then $F(x)(w_l)$ is
differentiable at $x$ for each $l$.

        We also know that
\begin{equation}
 \sup_{a \le y \le b \atop y \ne x} \frac{\|F(x) - F(y)\|_*}{|x - y|} < \infty
\end{equation}
for almost every $x \in [a, b]$.  This follows from the finiteness
almost everywhere of the maximal function associated to the function
$\Lambda_a^r$ that measures the length of $F$ on $[a, r]$, as in
Section \ref{maximal functions, 2}.  If $x$ has this property and $x
\not\in Z$, then one can check that the derivative
\begin{equation}
        \lim_{h \to 0} \frac{F(x + h)(w) - F(x)(w)}{h}
\end{equation}
of $F(x)(w)$ at $x$ exists for every $w \in W$, using the remarks in
the previous section.  Hence the derivative
\begin{equation}
\label{lim_{h to 0} frac{F(x + h) - F(x)}{h}, weak*}
        \lim_{h \to 0} \frac{F(x + h) - F(x)}{h}
\end{equation}
exists for almost every $x \in [a, b]$ in the weak$^*$ topology under
these conditions.

        Let $W$ be the space of continuous real or complex-valued
functions on $[0, 1]$ with the supremum norm, so that $W^*$ can be
identified with the space of real or complex Borel measures on $[0,
1]$, as appropriate.  Also let $F(x)$ be the function on $[0, 1]$ with
values in $W^*$ that assigns to $x \in [0, 1]$ the measure on $[0, 1]$
that is Lebesgue measure on $[0, x]$.  This is basically the same as
the function on $[0, 1]$ with values in $L^1([0, 1])$ discussed in
Section \ref{vector-valued functions}, by identifying integrable
functions on $[0, 1]$ with absolutely continuous measures with respect
to Lebesgue measure.  Now that we consider $F$ to take values in
$W^*$, it is easy to see that the derivative of $F$ exists with
respect to the weak$^*$ topology on $W^*$ at every $x \in [0, 1]$, and
corresponds to a Dirac mass at $x$.

\section{Lipschitz functions}
\label{lipschitz functions}
\setcounter{equation}{0}

        Let $f$ be a real or complex-valued Lipschitz function on the
real line.  Thus $f$ is differentiable almost everywhere, since it has
bounded variation on any bounded interval.  In particular,
\begin{equation}
\label{lim_{j to infty} frac{f(x + h_j) - f(x)}{h_j} = f'(x)}
        \lim_{j \to \infty} \frac{f(x + h_j) - f(x)}{h_j} = f'(x)
\end{equation}
almost everywhere for every sequence $\{h_j\}_{j = 1}^\infty$ of
nonzero real numbers that converges to $0$.  This implies that
\begin{equation}
\lim_{j \to \infty} \int_{\bf R} \frac{f(x + h_j) - f(x)}{h_j} \, \phi(x) \, dx
                                    = \int_{\bf R} f'(x) \, \phi(x) \, dx
\end{equation}
for every integrable function $\phi$ on ${\bf R}$, by the dominated
convergence theorem.  More precisely, this also uses the fact that the
difference quotients are uniformly bounded, because $f$ is Lipschitz.
Hence
\begin{equation}
 \lim_{h \to 0} \int_{\bf R} \frac{f(x + h) - f(x)}{h} \, \phi(x) \, dx
                                   = \int_{\bf R} f'(x) \, \phi(x) \, dx.
\end{equation}
This is the same as saying that
\begin{equation}
        \lim_{h \to 0} \frac{f(x + h) - f(x)}{h} = f'(x)
\end{equation}
in the weak$^*$ topology on $L^\infty({\bf R})$, as the dual of $L^1({\bf R})$.

        Alternatively, we can start with the identity
\begin{equation}
\label{int_R (f(x + h) - f(x))/h phi(x) dx = ...}
 \int_{\bf R} \frac{f(x + h) - f(x)}{h} \, \phi(x) \, dx =
        - \int_{\bf R} f(x) \, \frac{\phi(x) - \phi(x - h)}{h} \, dx,
\end{equation}
which uses the change of variables $x \mapsto x - h$.  This implies that
\begin{equation}
 \lim_{h \to 0} \int_{\bf R} \frac{f(x + h) - f(x)}{h} \, \phi(x) \, dx
                                      = - \int_{\bf R} f(x) \, \phi'(x) \, dx
\end{equation}
when $\phi$ is a continuously-differentiable function with compact
support on ${\bf R}$, for instance.  Thus
\begin{equation}
\label{lambda_h(phi) = int_{bf R} frac{f(x + h) - f(x)}{h} phi(x) dx}
 \lambda_h(\phi) = \int_{\bf R} \frac{f(x + h) - f(x)}{h} \, \phi(x) \, dx
\end{equation}
defines a bounded family of linear functionals on $L^1({\bf R})$ that
converges as $h \to 0$ on a dense linear subspace of $L^1({\bf R})$,
and hence converges on all of $L^1({\bf R})$, as in Section
\ref{uniform boundedness, 4}.  The limit is a bounded linear
functional on $L^1({\bf R})$ that can be expressed by integration with
an element of $L^\infty({\bf R})$, that corresponds to the derivative
of $f$.

        If $f$ is a bounded Lipschitz function on ${\bf R}$, then we
can take $F : {\bf R} \to L^\infty({\bf R})$ to be the function that
sends are real number to the corresponding translate of $f$, as in
Section \ref{vector-valued functions}.  Otherwise, we can take a
difference between $f$ and its translate to get an element of
$L^\infty({\bf R})$, as before.  This defines a Lipschitz mapping from
${\bf R}$ into $L^\infty({\bf R})$, with a weak$^*$ derivative at
every point.

\section{Averages}
\label{averages}
\setcounter{equation}{0}

        Let $f$ be a locally integrable function on the real line, and put
\begin{equation}
        A_h(f)(x) = \frac{1}{h} \int_x^{x + h} f(y) \, dy
\end{equation}
for every $h, x \in {\bf R}$ with $h \ne 0$.  As before, the integral
in this expression is considered to be oriented, as in ordinary
calculus, so that 
\begin{equation}
        A_h(f)(x) = \frac{1}{|h|} \int_{x - |h|}^x f(y) \, dy
\end{equation}
when $h < 0$.  In particular,
\begin{equation}
\label{|A_h(f)(x)| le A_h(|f|)(x)}
        |A_h(f)(x)| \le A_h(|f|)(x).
\end{equation}

        If $f \in L^p({\bf R})$$, 1 \le p \le \infty$, then $A_h(f)
\in L^p({\bf R})$ for every $h \ne 0$, and
\begin{equation}
        \|A_h(f)\|_p \le \|f\|_p.
\end{equation}
This is very easy to see when $p = \infty$.  If $p = 1$, then one can
integrate (\ref{|A_h(f)(x)| le A_h(|f|)(x)}) in $x$, and the use
Fubini's theorem.  If $1 < p < \infty$, then
\begin{equation}
\label{|A_h(f)(x)|^p le A_h(|f|^p)(x)}
        |A_h(f)(x)|^p \le A_h(|f|^p)(x),
\end{equation}
by the convexity of $r^p$ on the nonnegative real numbers, as in
Jensen's inequality.  One can then integrate in $x$ and apply Fubini's
theorem, as when $p = 1$.

        If $f$ is continuous at $x$, then
\begin{equation}
        \lim_{h \to 0} A_h(f)(x) = f(x).
\end{equation}
If $f$ is uniformly continuous, then this holds with uniform convergence.
If $f$ is a continuous function on ${\bf R}$, then $f$ is uniformly continuous
on bounded sets, and we get uniform convergence on bounded sets.

        If $f \in L^p({\bf R})$, $1 \le p < \infty$, then
\begin{equation}
\label{lim_{h to 0} ||A_h(f) - f||_p = 0}
        \lim_{h \to 0} \|A_h(f) - f\|_p = 0.
\end{equation}
To see this, observe first that this holds for every continuous
function $f$ with compact support on the real line.  More precisely,
$f$ is uniformly continuous in this case, so that $A_h(f)$ converges
to $f$ uniformly as $h \to 0$, as in the previous paragraph.  Also,
the support of $A_h(f)$ is contained in a single compact set when $|h|
\le 1$, say, and hence uniform convergence implies convergence in the
$L^p({\bf R})$ norm.  Any $f \in L^p({\bf R})$ can be approximated in
the $L^p$ norm by a continuous function with compact support when $p <
\infty$, and one can get (\ref{lim_{h to 0} ||A_h(f) - f||_p = 0})
using this approximation and the uniform bounds for $A_h$ on $L^p({\bf
R})$.

\section{$L^p$ derivatives}
\label{L^p derivatives}
\setcounter{equation}{0}

        If $f$, $g$ are locally integrable functions on the real
line, then we say that $f' = g$ in the sense of distributions if
\begin{equation}
 \int_{\bf R} f(x) \, \phi'(x) \, dx = - \int_{\bf R} g(x) \, \phi(x) \, dx
\end{equation}
for every continuously-differentiable function $\phi$ with compact
support on ${\bf R}$.  If $f$ is continuously differentiable on ${\bf
R}$, then the ordinary derivative of $f$ has this property, by
integration by parts.  Similarly, if
\begin{equation}
\label{lim_{h to 0} frac{f(x + h) - f(x)}{h} = g(x)}
        \lim_{h \to 0} \frac{f(x + h) - f(x)}{h} = g(x)
\end{equation}
with respect to the $L^1$ norm on any bounded interval in the real
line, then $f' = g$ in the sense of distributions.  This follows from
(\ref{int_R (f(x + h) - f(x))/h phi(x) dx = ...}), by taking the limit
as $h \to 0$.

        Suppose that
\begin{equation}
\label{f(x + h) - f(x) in L^p({bf R})}
        f(x + h) - f(x) \in L^p({\bf R})
\end{equation}
for some $p$, $1 \le p < \infty$, and every $h \in {\bf R}$, which
holds in particular when $f \in L^p({\bf R})$.  We say that $f$ is
differentiable in the $L^p$ sense, with derivative equal to $g$, if $g
\in L^p({\bf R})$, and one has convergence in (\ref{lim_{h to 0}
frac{f(x + h) - f(x)}{h} = g(x)}) in the $L^p$ norm.  This implies
that $f' = g$ in the sense of distributions, as in the previous
paragraph.  If $g \in L^p({\bf R})$ and
\begin{equation}
\label{f(x) = int_a^x g(y) dy}
        f(x) = \int_a^x g(y) \, dy
\end{equation}
for some $a \in {\bf R}$, then
\begin{equation}
\label{frac{f(x + h) - f(x)}{h} = A_h(g)(x)}
        \frac{f(x + h) - f(x)}{h} = A_h(g)(x)
\end{equation}
converges to $g$ as $h \to 0$ in the $L^p$ norm, as in the previous
section, and so the derivative of $f$ is equal to $g$ in the $L^p$
sense.  If $g$ is locally integrable, then $A_h(g) \to g$ as $h \to 0$
in the $L^1$ norm on every bounded interval, and we still have that
$f' = g$ in the sense of distributions.

        Note that $f' = 0$ in the sense of distributions when
\begin{equation}
        \int_{\bf R} f(x) \, \phi'(x) \, dx = 0
\end{equation}
for every continuously-differentiable function $\phi$ with compact support.
If $\psi$ is a continuous function with compact support on ${\bf R}$ such that
\begin{equation}
        \int_{\bf R} \psi(y) \, dy = 0,
\end{equation}
then
\begin{equation}
        \phi(x) = \int_{-\infty}^x \psi(y) \, dy
\end{equation}
is continuously differentiable and has compact support, and $\phi' = \psi$.
Thus $f' = 0$ in the sense of distributions if and only if
\begin{equation}
        \int_{\bf R} f(x) \, \psi(x) \, dx = 0
\end{equation}
for every continuous function $\psi$ with compact support and integral $0$.
One can show that this happens if and only if $f$ is constant almost 
everywhere.

        If $f' = g$ in the sense of distributions, then it follows
that that the difference between $f$ and (\ref{f(x) = int_a^x g(y)
dy}) is constant almost everywhere, since they have the same
derivative.  In particular, (\ref{frac{f(x + h) - f(x)}{h} =
A_h(g)(x)}) holds for each $h \ne 0$ and almost every $x$.  If $g \in
L^p({\bf R})$, then we get that the derivative of $f$ is equal to $g$
in the $L^p$ sense, as before.

\section{$L^p$ Lipschitz conditions}
\label{L^p lipschitz conditions}
\setcounter{equation}{0}

        Let $f$ be a locally integrable function on the real line that
satisfies (\ref{f(x + h) - f(x) in L^p({bf R})}) for some $p$, $1 \le
p < \infty$, and every $h \in {\bf R}$, such as an $L^p$ function.
Suppose that
\begin{equation}
\label{(int_{bf R} |f(x + h) - f(x)|^p dx)^{1/p} le C |h|}
         \Big(\int_{\bf R} |f(x + h) - f(x)|^p \, dx\Big)^{1/p} \le C \, |h|
\end{equation}
for some $C \ge 0$ and every $h \in {\bf R}$, which is the same as
saying that
\begin{equation}
        \frac{f(x + h) - f(x)}{h}
\end{equation}
is uniformly bounded in $L^p({\bf R})$.  Note that this happens when
$f' = g \in L^p({\bf R})$ in the sense of distributions, since the
difference quotient is equal to $A_h(g)$.

        Suppose also that $1 < p < \infty$, and let $q$ be the
conjugate exponent to $p$, $1/p + 1/q = 1$.  If $\lambda_h$ is as in
(\ref{lambda_h(phi) = int_{bf R} frac{f(x + h) - f(x)}{h} phi(x) dx})
for $h \ne 0$, then $\lambda_h$ is a uniformly bounded family of
linear functionals on $L^q({\bf R})$, by H\"older's inequality.
As in Section \ref{lipschitz functions},
\begin{equation}
\label{lim_{h to 0} lambda_h(phi) = - int_{bf R} f(x) phi'(x) dx}
        \lim_{h \to 0} \lambda_h(\phi) = - \int_{\bf R} f(x) \, \phi'(x) \, dx
\end{equation}
for every continuously-differentiable function $\phi$ with compact
support on ${\bf R}$.  Because these functions are dense in $L^q({\bf
R})$, it follows that
\begin{equation}
\label{lim_{h to 0} lambda_h(phi)}
        \lim_{h \to 0} \lambda_h(\phi)
\end{equation}
exists for every $\phi \in L^q({\bf R})$, as in Section \ref{uniform
boundedness, 4}.  The limit determines a bounded linear functional on
$L^q({\bf R})$, and so there is a function $g \in L^p({\bf R})$ such that
\begin{equation}
\label{lim_{h to 0} lambda_h(phi) = int_{bf R} g(x) phi(x) dx}
        \lim_{h \to 0} \lambda_h(\phi) = \int_{\bf R} g(x) \, \phi(x) \, dx
\end{equation}
for every $\phi \in L^q({\bf R})$.  In particular, this holds when
$\phi$ is a continuously-differentiable function with compact support
on ${\bf R}$, for which we have (\ref{lim_{h to 0} lambda_h(phi) = -
int_{bf R} f(x) phi'(x) dx}).  This shows that $f' = g$ in the sense
of distributions.

        If $p = 1$, then it is better to think of $\lambda_h$ as a
uniformly bounded family of linear functionals on the space $C_0({\bf
R})$ of continuous functions on the real line that vanish at infinity,
equipped with the supremum norm.  We still have (\ref{lim_{h to 0}
lambda_h(phi) = - int_{bf R} f(x) phi'(x) dx}) for every
continuously-differentiable function $\phi$ with compact support on
${\bf R}$, and hence that (\ref{lim_{h to 0} lambda_h(phi)}) exists
for every $\phi \in C_0({\bf R})$, as in Section \ref{uniform boundedness, 4}.
The limit determines a bounded linear functional on $C_0({\bf R})$,
and so there is a real or complex Borel measure $\mu$ on ${\bf R}$ such
that
\begin{equation}
        \lim_{h \to 0} \lambda_h(\phi) = \int_{\bf R} \phi \, d\mu
\end{equation}
for every $\phi \in C_0({\bf R})$.  Combining this with (\ref{lim_{h
to 0} lambda_h(phi) = - int_{bf R} f(x) phi'(x) dx}), we get that
\begin{equation}
        \int_{\bf R} f(x) \, \phi'(x) \, dx = - \int_{\bf R} \phi \, d\mu
\end{equation}
for every continuously-differentiable function $\phi$ with compact
support on ${\bf R}$.  This can be expressed by saying that $f' = \mu$
in the sense of distributions.

        If $\alpha$ is a function of bounded variation on ${\bf R}$,
and if $\mu_\alpha$ is the corresponding real or complex Borel measure
as in Section \ref{functions, measures}, then $\alpha' = \mu$ in the
sense of distributions.  This is basically another version of
integration by parts.  One can also show that every real or complex
Borel measure on the real line is of this form.  If $f$ is a locally
integrable function on ${\bf R}$ such that $f' = \mu$ in the sense of
distributions for some real or complex Borel measure $\mu$, then it
follows that $f$ is equal almost everywhere to a function of bounded
variation.  Conversely, one can check that such functions satisfy the
integrated Lipschitz condition (\ref{(int_{bf R} |f(x + h) - f(x)|^p
dx)^{1/p} le C |h|}) with $p = 1$.

\section{Dyadic intervals}
\label{dyadic intervals}
\setcounter{equation}{0}

        In this section, it will be convenient to use $[0, 1)$ as the
unit interval, consisting of $x \in {\bf R}$ with $0 \le x < 1$.  By a
\emph{dyadic subinterval} of the unit interval we mean an interval of
the form $[j \, 2^{-l}, (j + 1) \, 2^{-1})$, where $j$, $l$ are
nonnegative integers and $j < 2^l$.  Thus the unit interval is the
disjoint union of these dyadic intervals at level $l$.  If $I$, $I'$
are dyadic intervals of arbitrary lengths, then either $I \subseteq
I'$, $I' \subseteq I$, or $I \cap I' = \emptyset$.  More precisely, if
$|I| \le |I'|$, where $|I|$ denotes the length of $I$, then either $I
\subseteq I'$ or $I \cap I' = \emptyset$.

        Let $\mu$ be a positive Borel measure on $[0, 1)$.  The
\emph{dyadic maximal function} associated to $\mu$ is defined by
\begin{equation}
        \mu^*_\delta(x) = \sup_{x \in I} \frac{\mu(I)}{|I|},
\end{equation}
where now the supremum is taken over all dyadic intervals that contain
a given point $x \in [0, 1)$.  Similarly, if $f$ is an integrable
function on $[0, 1)$, then we put
\begin{equation}
        f^*_\delta(x) = \sup_{x \in I} \frac{1}{|I|} \int_I |f(y)| \, dy,
\end{equation}
where again the supremum is taken over all dyadic intervals $I$ such
that $x \in I$.  This is the same as $\mu^*_\delta(x)$, where $\mu$ is
the Borel measure on $[0, 1)$ defined by
\begin{equation}
\label{mu(A) = int_A |f(y)| dy}
        \mu(A) = \int_A |f(y)| \, dy,
\end{equation}
as in Section \ref{maximal functions}.

        Consider
\begin{equation}
        E_{\delta, t} = \{x \in [0, 1) : \mu^*_\delta(x) > t\}
\end{equation}
for each $t > 0$.  Thus $x \in E_{\delta, t}$ if and only if there is
a dyadic interval $I$ such that $x \in I$ and
\begin{equation}
\label{mu(I) > t |I|}
        \mu(I) > t \, |I|,
\end{equation}
in which case $I \subseteq E_{\delta, t}$.  Let $I(x)$ be the maximal
dyadic interval that contains $x$ and satisfies (\ref{mu(I) > t |I|})
for each $x \in E_{\delta, t}$.  If $x, y \in E_{\delta, t}$, then
either $I(x) = I(y)$ or $I(x) \cap I(y) = \emptyset$, by maximality
and the nesting properties of dyadic intervals mentioned before.

        Let $\mathcal{M}_t$ be the collection of dyadic intervals of
the form $I(x)$ for some $x$ in $E_{\delta, t}$.  Note that the
elements of $\mathcal{M}_t$ are pairwise disjoint, and
\begin{equation}
        \bigcup_{I \in \mathcal{M}_t} I = E_{\delta, t}.
\end{equation}
Hence
\begin{equation}
 |E_{\delta, t}| = \sum_{I \in \mathcal{M}_t} |I|
 < t^{-1} \, \sum_{I \in \mathcal{M}_t} \mu(I) = t^{-1} \, \mu(E_{\delta, t}).
\end{equation}
This is almost the same as the estimate in Section \ref{maximal
functions}, but without the additional factor of $2$.  Although we
have focused on dyadic subintervals of the unit interval for
simplicity, there is an analogous discussion for arbitrary dyadic
intervals in the real line, and the corresponding maximal functions.

\section{Dyadic averages}
\label{dyadic averages}
\setcounter{equation}{0}

        Let $f$ be an integrable function on $[0, 1)$, and put
\begin{equation}
        A_l(f)(x) = 2^l \int_{j \, 2^{-l}}^{(j + 1) \, 2^{-l}} f(y) \, dy
\end{equation}
when $j \, 2^{-l} \le x < (j + 1) \, 2^{-l}$.  Thus $A_l(f)(x)$ is the
average of $f$ over the dyadic interval of length $2^{-l}$ that
contains $x$.  In particular, $A_l(f)$ is constant on dyadic intervals
of length $2^{-l}$, by construction.  Also,
\begin{eqnarray}
 \int_0^1 A_l(f)(x) \, dx & = & \sum_{j = 0}^{2^l - 1}
                    \int_{j \, 2^{-l}}^{(j + 1) \, 2^{-l}} A_l(f)(x) \, dx  \\
 & = & \sum_{j = 0}^{2^l - 1} \int_{j \, 2^{-l}}^{(j + 1) \, 2^{-l}} f(x) \, dx
                              =  \int_0^1 f(x) \, dx.               \nonumber
\end{eqnarray}

        If $f \in L^p([0, 1))$, $1 \le p \le \infty$, then
\begin{equation}
\label{||A_l(f)||_p le ||f||_p}
        \|A_l(f)\|_p \le \|f\|_p.
\end{equation}
This is immediate when $p = \infty$.  Note that
\begin{equation}
        |A_l(f)(x)| \le A_l(|f|)(x)
\end{equation}
for every $x \in [0, 1)$ and $l \ge 0$, and that
\begin{equation}
        |A_l(f)(x)|^p \le A_l(|f|^p)(x)
\end{equation}
when $f \in L^p([0, 1))$, $1 < p < \infty$, by Jensen's inequality.
To estimate $\|A_l(f)\|_p$, one can integrate these inequalities
using the identity in the previous paragraph.

        As in Section \ref{averages},
\begin{equation}
\label{lim_{l to infty} A_l(f)(x) = f(x)}
        \lim_{l \to \infty} A_l(f)(x) = f(x)
\end{equation}
when $f$ is continuous at $x$, and with uniform convergence when $f$
is uniformly continuous on $[0, 1)$.  If $f$ is a continuous function
on $[0, 1]$, then $f$ is uniformly continuous, by compactness.  If $f
\in L^p([0, 1))$, $1 \le p < \infty$, then
\begin{equation}
        \lim_{l \to \infty} \|A_l(f) - f\|_p = 0.
\end{equation}
This follows from uniform convergence when $f$ is a continuous
function on $[0, 1]$, and otherwise one can approximate by continuous
functions using the uniform bound (\ref{||A_l(f)||_p le ||f||_p}).  Of
course, $L^p([0, 1))$ is the same as $L^p([0, 1])$, and so continuous
functions on $[0, 1]$ are still dense in this space when $p < \infty$.

        If $f \in L^1([0, 1))$, then Lebesgue's theorem implies that
(\ref{lim_{l to infty} A_l(f)(x) = f(x)}) holds almost everywhere on
$[0, 1)$.  More precisely,
\begin{equation}
\label{lim_{l to infty} 2^l int_{I_l(x)} |f(y) - f(x)| dy = 0}
        \lim_{l \to \infty} 2^l \int_{I_l(x)} |f(y) - f(x)| \, dy = 0
\end{equation}
for almost every $x \in [0, 1)$, where $I_l(x)$ denotes the dyadic
interval of length $2^{-l}$ that contains $x$.  This follows from
Lebesgue's theorem, as in Section \ref{lebesgue's theorem}, and one
can also establish it a bit more directly.  Specifically, one can use
the estimate for the dyadic maximal function in the previous section,
instead of the estimate for the Hardy--Littlewood maximal function in
Section \ref{maximal functions}.

\section{Rademacher functions}
\label{rademacher functions}
\setcounter{equation}{0}

        Let $r_1(x), r_2(x), \ldots$ be the functions defined on $[0, 1)$ by
\begin{eqnarray}
 r_l(x) & = & 1 \quad\enspace\hbox{when } j \, 2^{-l} \le x < (j + 1) \, 2^{-l}
                                    \hbox{ and $j$ is even} \\
        & = & -1 \quad\hbox{when } j \, 2^{-l} \le x < (j + 1) \, 2^{-l}
                                    \hbox{ and $j$ is odd}. \nonumber
\end{eqnarray}
Thus $r_l(x)$ is constant on each dyadic interval of length $2^{-l}$,
\begin{equation}
\label{int_I r_l(x) dx = 0}
        \int_I r_l(x) \, dx = 0
\end{equation}
for each dyadic interval $I$ of length $2^{-l + 1}$, and $|r_l(x)| =
1$ for every $x \in [0, 1)$ and positive integer $l$.  These are known
as the \emph{Rademacher functions} on the unit interval.

        Let $X$ be the set of sequences $x = \{x_k\}_{k = 1}^\infty$
with $x_k = 1$ or $-1$ for each $k$.  Equivalently, $X$ is the
Cartesian product of a sequence of copies of $\{1, -1\}$.  This is a
compact Hausdorff topological space with respect to the product
topology, which is homeomorphic to the usual middle-thirds Cantor set.
There is a natural continuous mapping from $X$ onto the closed unit
interval $[0, 1]$, defined by
\begin{equation}
        \beta(x) = \sum_{k = 1}^\infty \Big(\frac{x_k + 1}{2}\Big) \, 2^{-k}.
\end{equation}
Each element $x$ of $X$ corresponds to an infinite binary sequence
$\{(x_k + 1)/2\}_{k = 1}^\infty$, and $\beta$ sends $x$ to the real
number with that binary expansion.  Every real number in $[0, 1]$ has
a binary expansion, and the binary expansion is unique for all but a
countable set of real numbers.  Dydadic rational numbers of the form
$j \, 2^{-l}$, $0 < j < 2^l$, have two binary expansions, which agree
up to a point where one has a $1$ followed by all $0$'s, and the other
has a $0$ followed by all $1$'s.

        There is a natural Borel probablility measure on $X$, which is
the product measure associated to $1$, $-1$ having probability $1/2$
in each coordinate.  This probability measure corresponds exactly to
Lebesgue measure on $[0, 1]$ under the mapping $\beta$.  That $\beta$
fails to be one-to-one on a countable set does not really matter here,
since countable sets have measure $0$.  Thus $[0, 1)$ and $X$ are
basically the same as probability spaces.  The Rademacher functions
$r_l$ on $[0, 1)$ correspond to the coordinate functions $x \mapsto
x_l$ on $X$, which are independent identically distributed random
variables.

        In particular,
\begin{equation}
        \int_0^1 r_{l_1}(x) \, r_{l_2}(x) \cdots r_{l_n}(x) \, dx = 0
\end{equation}
when $1 \le l_1 < l_2 < \cdots < l_n$.  Because of independence, the
integral of the product should be the same as the product of the
individual integrals, each of which is $0$, by (\ref{int_I r_l(x) dx =
0}).  One can see this more directly by observing that the integral
over each dyadic interval of length $2^{-l_n + 1}$ is $0$, because the
integral of $r_{l_n}$ over such an interval is $0$, as in (\ref{int_I
r_l(x) dx = 0}), while the other functions in the integral are
constant over these intervals.

\section{$L^p$ estimates}
\label{L^p estimates}
\setcounter{equation}{0}

        The Rademacher functions are orthonormal in $L^2([0, 1))$, since 
\begin{equation}
        \|r_l\|_2 = \Big(\int_0^1 |r_l(x)|^2 \, dx\Big)^{1/2} = 1
\end{equation}
for each $l$, and
\begin{equation}
        \langle r_k, r_l \rangle = \int_0^1 r_k(x) \, r_l(x) \, dx = 0
\end{equation}
when $k \ne l$.  This implies that
\begin{equation}
\label{||sum_{l = 1}^n a_l r_l||_2 = (sum_{l = 1}^n a_l^2)^{1/2}}
        \biggl\|\sum_{l = 1}^n a_l \, r_l\biggr\|_2
                      = \Big(\sum_{l = 1}^n a_l^2\Big)^{1/2}
\end{equation}
for every $a_1, \ldots, a_n \in {\bf R}$.  Let us check that
\begin{equation}
\label{||sum_{l = 1}^n a_l r_l||_infty = sum_{l = 1}^n |a_l|}
 \biggl\|\sum_{l = 1}^n a_l \, r_l\biggr\|_\infty = \sum_{l = 1}^n |a_l|.
\end{equation}
The left side of (\ref{||sum_{l = 1}^n a_l r_l||_infty = sum_{l = 1}^n
  |a_l|}) is clearly less than or equal to the right side, by the
triangle inequality.  To get the opposite inequality, one can choose a
dyadic interval of length $2^{-n}$ on which $a_l \, r_l = |a_l|$ for
$l = 1, \ldots, n$.

        Before proceeding, it will be helpful to remember two basic
facts about $L^p$ norms.  The first is that
\begin{equation}
        \|f\|_p = \Big(\int_0^1 |f(x)|^p \, dx\Big)^{1/p}
\end{equation}
is monotone increasing in $p$, by Jensen's inequality.  The second
fact is that the $L^p$ norm is logarithmically convex in $1/p$,
which means that
\begin{equation}
\label{||f||_r le ||f||_p^t ||f||_q^{1 - t}}
        \|f\|_r \le \|f\|_p^t \, \|f\|_q^{1 - t}
\end{equation}
when $p, q, r > 0$, $0 < t < 1$, and
\begin{equation}
        \frac{1}{r} = \frac{t}{p} + \frac{1 - t}{q}.
\end{equation}
This can be derived from H\"older's inequality.  It is a little
simpler to start with the $r = 1$ case, and then get (\ref{||f||_r le
||f||_p^t ||f||_q^{1 - t}}) by applying the $r = 1$ case to $|f|^r$.

        If $2 < p < \infty$, then there is a constant $C(p) > 0$ such that
\begin{equation}
\label{||sum_{l = 1}^n a_l r_l||_p le C(p) (sum_{l = 1}^n a_l^2)^{1/2}}
        \biggl\|\sum_{l = 1}^n a_l \, r_l\biggr\|_p
                \le C(p) \, \Big(\sum_{l = 1}^n a_l^2\Big)^{1/2}
\end{equation}
for every $a_1, \ldots, a_n \in {\bf R}$.  Of course, it is very
important here that $C(p)$ does not depend on $n$.  To prove
(\ref{||sum_{l = 1}^n a_l r_l||_p le C(p) (sum_{l = 1}^n
a_l^2)^{1/2}}), it suffices to restrict our attention to $p = 2^k$ for
some positive integer $k \ge 2$, because of the monotonicity of the
$L^p$ norm.  One can get better constants for the intermediate
exponents using (\ref{||sum_{l = 1}^n a_l r_l||_2 = (sum_{l = 1}^n
a_l^2)^{1/2}}) and (\ref{||f||_r le ||f||_p^t ||f||_q^{1 - t}}).  If
$p = 2^k$, then one can expand
\begin{equation}
        \biggl\|\sum_{l = 1}^n a_l \, r_l\biggr\|_{2^k}^{2^k}
          = \int_0^1 \Big(\sum_{l = 1}^n a_l \, r_l\Big)^{2^k} \, dx
\end{equation}
into a $2^k$-fold sum, where each term has the product of $2^k$
coefficients $a_l$ times the integral of the product of $2^k$
Rademacher functions $r_l$.  As in the previous section, most of these
integrals are equal to $0$.  The only way that the integral is not
equal to $0$ is to have $r_l$ occur an even number of times for each
$l$.  In this case, the integral is equal to $1$, and the coefficients
are products of $2^{k - 1}$ factors of $r_l^2$, $1 \le l \le n$.
This permits one to estimate the $2^k$-fold sum by a constant multiple of
\begin{equation}
        \Big(\sum_{l = 1}^n a_l^2\Big)^{2^{k - 1}},
\end{equation}
as desired.  The $k = 2$ case is already a nice exercise.

        If $0 < p < 2$, then there is a constant $C(p) > 0$ such that
\begin{equation}
\label{(sum_{l = 1}^n a_l^2)^{1/2} le C(p) ||sum_{l = 1}^n a_l r_l||_p}
        \Big(\sum_{l = 1}^n a_l^2\Big)^{1/2}
               \le C(p) \, \biggl\|\sum_{l = 1}^n a_l \, r_l\biggr\|_p
\end{equation}
for every $a_1, \ldots, a_n \in {\bf R}$.  Again, it is very important
that $C(p)$ not depend on $n$.  This time, we can apply (\ref{||f||_r
le ||f||_p^t ||f||_q^{1 - t}}) to $f = \sum_{l = 1}^n a_l \, r_l$, $r
= 2$, and $q = 4$ to get that
\begin{equation}
        \Big(\sum_{l = 1}^n a_l^2\Big)^{1/2}
              \le \biggl\|\sum_{l = 1}^n a_l \, r_l\biggr\|_p^t
                \, \biggl\|\sum_{l = 1}^n a_l \, r_l\biggr\|_4^{1 - t}
\end{equation}
for some $t$, $0 < t < 1$.  Using the previous estimate with $p = 4$,
we get that
\begin{equation}
        \Big(\sum_{l = 1}^n a_l^2\Big)^{1/2}
          \le C(4) \, \Big(\sum_{l = 1}^n a_l^2\Big)^{(1 - t)/2}
                    \, \biggl\|\sum_{l = 1}^n a_l \, r_l\biggr\|_p^t.
\end{equation}
This implies (\ref{(sum_{l = 1}^n a_l^2)^{1/2} le C(p) ||sum_{l = 1}^n
a_l r_l||_p}), by dividing both sides by $\Big(\sum_{l = 1}^n
a_l^2\Big)^{(1 - t)/2}$, at least when $a_l \ne 0$ for some $l$.

\section{Rademacher sums}
\label{rademacher sums}
\setcounter{equation}{0}

        Let $a_1, a_2, \ldots$ be a sequence of real numbers such that
$\sum_{l = 1}^\infty a_l^2$ converges, and consider
\begin{equation}
        f(x) = \sum_{l = 1}^\infty a_l \, r_l(x).
\end{equation}
This series converges in $L^2([0, 1))$, by the orthonormality of the
Rademacher functions.  Moreover, the series converges in $L^p([0, 1))$
for every $p < \infty$, by the estimates in the previous section.
Using these estimates, one can also check that this series converges
in $L^p([0, 1))$ in the generalized sense for every $p < \infty$, as
in Section \ref{generalized convergence, 2}.

        Observe that
\begin{equation}
\label{A_n(f)(x) = sum_{l = 1}^n a_l r_l(x)}
        A_n(f)(x) = \sum_{l = 1}^n a_l \, r_l(x)
\end{equation}
for every $n$, where $A_n$ is the dyadic averaging operator in Section
\ref{dyadic averages}.  This follows from the fact that $A_n(r_l) = 0$
when $l > n$.  By Lebesgue's theorem,
\begin{equation}
\label{lim_{n to infty} A_n(f)(x) = f(x)}
        \lim_{n \to \infty} A_n(f)(x) = f(x)
\end{equation}
almost everywhere on $[0, 1)$, which implies that the series defining
$f$ converges almost everywhere.  However, if $\sum_{l = 1}^\infty a_l
\, r_l(x)$ converges in the generalized sense as a sum of real numbers
for any $x \in [0, 1)$, then
\begin{equation}
        \sum_{l = 1}^\infty |a_l \, r_l(x)| = \sum_{l = 1}^\infty |a_l|
\end{equation}
converges, as in Section \ref{generalized convergence}.  Similarly, if
$f \in L^\infty([0, 1))$, then $A_n(f)$ is uniformly bounded, and
hence $\sum_{l = 1}^\infty |a_l|$ converges, by (\ref{||sum_{l = 1}^n
a_l r_l||_infty = sum_{l = 1}^n |a_l|}).

        Let $\pi$ be a one-to-one mapping from the set ${\bf Z}_+$ of
positive integers onto itself, and let $X$ be the space of all
sequences $\{x_k\}_{k = 1}^\infty$ with $x_k = \pm 1$ for each $k$, as
in Section \ref{rademacher functions}.  Thus $\pi$ determines a
measure-preserving homeomorphism from $X$ onto itself, which sends
$\{x_k\}_{k = 1}^\infty$ to $\{x_{\pi(k)}\}_{k = 1}^\infty$.  Using
this transformation, one can check that $\sum_{l = 1}^\infty
a_{\pi(l)} \, r_{\pi(l)}(x)$ also converges almost everywhere.  More
precisely, this rearrangement of the series corresponds to the
composition of $\sum_{l = 1}^\infty a_{\pi(l)} \, r_l(x)$ with the
automorphism on $X$ just mentioned.  This new series is of the same
type as the previous one, and so converges almost everywhere for the
same reasons as before.

\section{Lacunary series}
\label{lacunary series}
\setcounter{equation}{0}

        Let ${\bf T}$ be the unit circle in the complex plane,
consisting of $z \in {\bf C}$ with $|z| = 1$.  It is well known that
\begin{equation}
\label{int_{bf T} z^j |dz| = 0}
        \int_{\bf T} z^j \, |dz| = 0
\end{equation}
for every nonzero integer $j$, where $|dz|$ denotes the element of arc
length.  If $j = 0$, then $z^j$ is interpreted as being equal to $1$,
and the integral is equal to $2 \pi$, the circumference of the circle.
The usual integral inner product for complex-valued functions in
$L^2({\bf T})$ is defined by
\begin{equation}
        \langle f, g \rangle
          = \frac{1}{2 \pi} \int_{\bf T} f(z) \, \overline{g(z)} \, |dz|,
\end{equation}
and the corresponding norm is given by
\begin{equation}
        \|f\|_2 = \Big(\frac{1}{2 \pi}\int_{\bf T} |f(z)|^2\Big)^{1/2}.
\end{equation}
The functions $z^j$, $j \in {\bf Z}$, are orthonormal with respect to
this inner product, because of (\ref{int_{bf T} z^j |dz| = 0}) and the
fact that the integral is equal to $2 \pi$ when $j = 0$.  It is well
known that the linear span of these functions is dense in $L^2({\bf
T})$, and more precisely that their linear span is dense in the space
of continuous functions on ${\bf T}$ with respect to the supremum
norm.  This implies that $z^j$, $j \in {\bf Z}$, is an orthonormal
basis for $L^2({\bf T})$.

        Let $n_1 < n_2 < \cdots$ be a strictly increasing sequence of
positive integers, and let $a_1, a_2, \ldots$ be a sequence of complex
numbers such that $\sum_{j = 1}^\infty |a_j|^2$ converges.  Thus 
\begin{equation}
\label{f(z) = sum_{j = 1}^infty a_j z^{n_j}}
        f(z) = \sum_{j = 1}^\infty a_j \, z^{n_j}
\end{equation}
converges in $L^2({\bf T})$, since the $z^{n_j}$'s are orthonormal in
$L^2({\bf T})$.  We say that (\ref{f(z) = sum_{j = 1}^infty a_j z^{n_j}})
is a \emph{lacunary} or \emph{gap series} if there is a $q > 1$ such that
\begin{equation}
\label{n_{j + 1} ge q n_j}
        n_{j + 1} \ge q \, n_j
\end{equation}
for each $j$.  In this case, (\ref{f(z) = sum_{j = 1}^infty a_j
z^{n_j}}) actually converges in $L^p({\bf T})$ for each $p < \infty$.
One can also show that the series converges in the generalized sense
in $L^p({\bf T})$, as in Section \ref{generalized convergence, 2}.

        To see this, it suffices to show that for each $p \in (2, \infty)$
there is a constant $C'(p) > 0$ such that
\begin{equation}
        \biggl\|\sum_{j = 1}^L a_j \, z^{n_j}\biggr\|_p
           \le C'(p) \, \Big(\sum_{j = 1}^L |a_j|^2\Big)^{1/2}
\end{equation}
for every $a_1, \ldots a_L \in {\bf C}$ and $L \ge 1$.  It is also
enough to do this when $p = 2^k$ for some integer $k \ge 2$.  In this
case, the $p$th power of the $L^p$ norm can be expanded into a
$2^k$-fold sum, as before.  More precisely,
\begin{equation}
 \biggl|\sum_{j = 1}^L a_j \, z^{n_j}\biggr|^{2^k} =
        \Big(\sum_{j = 1}^L a_j \, z^{n_j}\Big)^{2^{k - 1}} \,
 \Big(\sum_{j = 1}^L \overline{a_j} \, \overline{z}^{n_j}\Big)^{2^{k - 1}},
\end{equation}
since $|a|^2 = a \, \overline{a}$ for every $a \in {\bf C}$.  Thus
each term in the $2^k$-fold sum has $2^{k - 1}$ $a_j$'s and
$z^{n_j}$'s, and $2^{k - 1}$ $\overline{a_j}$'s and
$\overline{z}^{n_j}$'s.

        Each term is also integrated over ${\bf T}$, and so includes
an expression of the form
\begin{equation}
 \int_{\bf T} \Big(\prod_{l = 1}^{2^{k - 1}} z^{n_{j_l}}\Big) \,
   \Big(\prod_{l' = 1}^{2^{k - 1}} \overline{z}^{n_{j'_{l'}}}\Big) \, |dz|,
\end{equation}
where the $j_l$'s and $j'_{l'}$'s are integers between $1$ and $L$.
Because of (\ref{int_{bf T} z^j |dz| = 0}), this integral is equal to
$0$ unless
\begin{equation}
 \sum_{l = 1}^{2^{k - 1}} n_{j_l} - \sum_{l' = 1}^{2^{k - 1}} n_{j'_{l'}} = 0.
\end{equation}
If $q$ is large enough, depending on $k$, then the only way that this
can happen is if the largest of the $n_{j_l}$'s is equal to the
largest of the $n_{j'_{l'}}$'s.  One can then repeat the argument to
get that the $n_{j_l}$'s and $n_{j'_{l'}}$'s are permutations of each
other.  This permits the $2^k$-fold sum to be estimated in terms of
$\Big(\sum_{j = 1}^L |a_j|^2\Big)^{2^{k - 1}}$, as in Section \ref{L^p
estimates}.  If $q$ is not sufficiently large for this argument, then
one can express (\ref{f(z) = sum_{j = 1}^infty a_j z^{n_j}}) as a sum
of finitely many lacunary series with larger gaps.  More precisely,
(\ref{f(z) = sum_{j = 1}^infty a_j z^{n_j}}) can be expressed as the
sum of $r$ lacunary series with gaps of size $q^r$ for each positive
integer $r$, by taking every $r$th term in the series.

\section{Walsh functions}
\label{walsh functions}
\setcounter{equation}{0}

        If $I = \{l_1, \ldots, l_n\}$ is a finite set of positive integers,
then the corresponding \emph{Walsh function} $w_I$ on $[0, 1)$ is defined by
\begin{equation}
\label{w_I(x) = r_{l_1}(x) r_{l_2}(x) r_{l_n}(x)}
        w_I(x) = r_{l_1}(x) \, r_{l_2}(x) \cdots r_{l_n}(x),
\end{equation}
where the $r_l$'s are Rademacher functions.  If $I = \emptyset$, then
we take $w_I$ to be the constant function $1$.  Thus
\begin{equation}
        |w_I(x)| = 1
\end{equation}
for every $x \in [0, 1)$ and finite set $I$ of positive integers, and
\begin{equation}
        \int_0^1 w_I(x) \, dx = 0
\end{equation}
when $I \ne \emptyset$, as in Section \ref{rademacher functions}.
This implies that
\begin{equation}
        \int_0^1 w_I(x) \, w_{I'}(x) \, dx = 0
\end{equation}
when $I \ne I'$, so that the Walsh functions are orthonormal in $L^2([0, 1))$.

        The Walsh functions actually form an orthonormal basis for
$L^2([0, 1))$.  To see this, it suffices to show that the linear span
of the Walsh functions is dense in $L^2([0, 1))$.  Note that $w_I(x)$
is constant on dyadic intervals of length $2^{-n}$ when $I \subseteq
\{1, \ldots, n\}$, because of the corresponding property of the
Rademacher functions.  One can check that the linear span of the Walsh
functions $w_I$ with $I \subseteq \{1, \ldots, n\}$ is exactly the
same as the space of functions on $[0, 1)$ that are constant on dyadic
intervals of length $2^{-n}$.  Both spaces have dimension $2^n$, for
instance, since there are $2^n$ subsets of $\{1, \ldots, n\}$, and
$2^n$ dyadic intervals of length $2^{-n}$.  It follows that the linear
span of all Walsh functions is the space of dyadic step functions on
$[0, 1)$, which are the functions that are constant on dyadic
intervals of length $2^{-n}$ for some $n$.  Hence the Walsh functions
form an orthonormal basis of $L^2([0, 1))$, because the dyadic step
functions are dense in $L^2([0, 1))$.

        There is another description of the Walsh functions in terms
of harmonic analysis.  Let $X$ be the space of sequences $\{x_k\}_{k =
1}^\infty$ with $x_k = \pm 1$ for each $k$, as in Section
\ref{rademacher functions}.  It is easy to see that $X$ is a
commutative group with respect to coordinatewise multiplication.  More
precisely, $X$ is a topological group with respect to the product
topology, because the group operations are continuous with respect to
this topology.  Note that the probability measure on $X$ described
before is invariant under translations defined by this group
structure, and hence corresponds to Haar measure on $X$.  The
Rademacher functions may be identified with the coordinate functions
on $X$, and so the Walsh functions may be identified with products of
coordinate functions on $X$.  One can check that these are continuous
homomorphisms from $X$ into the multiplicative group of nonzero
complex numbers, and that every such homomorphism arises in this way.

\section{Independent random variables}
\label{independent random variables}
\setcounter{equation}{0}

        Let $(X_1, \mu_1), \ldots, (X_n, \mu_n)$ be probability
spaces, and let $X = X_1 \times \cdots \times X_n$ be their product,
with the product measure $\mu = \mu_1 \times \cdots \times \mu_n$.
Also let $f_1, \ldots, f_n$ be real or complex-valued functions on
$X_1, \ldots, X_n$, respectively, which can be identified with
functions on $X$ that are constant in the other variables.
Suppose that $f_j \in L^2(X_j, \mu_j)$,
\begin{equation}
\label{int_{X_j} f_j d mu_j = 0}
        \int_{X_j} f_j \, d\mu_j = 0,
\end{equation}
and
\begin{equation}
\label{||f_j||_{L^2(X_j, mu_j)} = (int_{X_j} |f_j|^2 d mu_j)^{1/2} = 1}
 \|f_j\|_{L^2(X_j, \mu_j)} = \Big(\int_{X_j} |f_j|^2 \, d\mu_j\Big)^{1/2} = 1
\end{equation}
for each $j$.  It may be that the $(X_j, \mu_j)$'s are copies of the
same space, for instance, and that the $f_j$'s are copies of the same
function on this space.  As functions on $X$, it is easy to see that
$f_1, \ldots, f_n$ are orthonormal in $L^2(X, \mu)$.  This is because
\begin{equation}
        \int_X f_j \, f_l \, d\mu = \Big(\int_{X_j} f_j \, d\mu_j\Big) \,
                                     \Big(\int_{X_l} f_l \, d\mu_l\Big) = 0
\end{equation}
when $j \ne l$ in the real case, and
\begin{equation}
 \int_X f_j \, \overline{f_l} \, d\mu = \Big(\int_{X_j} f_j \, d\mu_j\Big) \,
                            \Big(\int_{X_l} \overline{f_l} \, d\mu_l\Big) = 0
\end{equation}
in the complex case.  Hence
\begin{equation}
\label{||sum_{j = 1}^n a_j f_j||_{L^2(X, mu)} = (sum_{j = 1}^n |a_j|^2)^{1/2}}
        \biggl\|\sum_{j = 1}^n a_j \, f_j\biggr\|_{L^2(X, \mu)}
                             = \Big(\sum_{j = 1}^n |a_j|^2\Big)^{1/2}
\end{equation}
for any real or complex numbers $a_1, \ldots, a_n$, as appropriate.

        Let $k$ be a positive integer, and put $p = 2^k$.  Suppose in
addition that $f_j \in L^p(X_j, \mu_j)$ for each $j$, and that
\begin{equation}
        \|f_j\|_{L^p(X_j, \mu_j)}
             = \Big(\int_{X_j} |f_j|^p \, d\mu_j\Big)^{1/p} \le L_p
\end{equation}
for some $L_p \ge 0$ and $j = 1, \ldots, n$.  In this case, one can show that
\begin{equation}
        \biggl\|\sum_{j = 1}^n a_j \, f_j\biggr\|_{L^p(X, \mu)}
          \le C(p, L_p) \, \Big(\sum_{j = 1}^n |a_j|^2\Big)^{1/2}
\end{equation}
for some constant $C(p, L_p) \ge 0$ and all $a_1, \ldots, a_n \in {\bf
R}$ or ${\bf C}$, as appropriate.  As usual, it is very important that
$C(p, L_p)$ does not depend on $n$ here.  To see this, one can expand
\begin{equation}
        \biggl\|\sum_{j = 1}^n a_j \, f_j\biggr\|_{L^p(X, \mu)}^p
          = \int_X \biggl|\sum_{j = 1}^n a_j \, f_j\biggr|^p \, d\mu
\end{equation}
into a $2^k$-fold sum, where each term is a product of $2^k$ $a_j$'s
and perhaps their complex conjugates times the integral of a product
of $2^k$ $f_j$'s and perhaps their complex conjugates, as in Sections
\ref{L^p estimates} and \ref{lacunary series}.  The integrals can be
estimated individually using H\"older's inequality and the hypothesis
that the $f_j$'s have bounded $L^p$ norms.  The main point is that the
integral is equal to $0$ whenever an $f_j$ occurs exactly once for
some $j$, because the integral over $X$ of a product of $f_j$'s and
perhaps their complex conjugates is equal to the product of the integrals
over the $X_j$'s of the corresponding $f_j$'s for $j = 1, \ldots, n$.
In the remaining terms, there is a product of $2^k$ $a_j$'s and perhaps
their complex conjugates, in which each $a_j$ either does not occur or occurs
more than once.  This permits one to estimate the sum by a constant multiple of
\begin{equation}
        \Big(\sum_{j = 1}^n |a_j|^2\Big)^{2^{k - 1}},
\end{equation}
as before.  This is a bit more complicated than in the context of
Rademacher functions, where the integrals are equal to $0$ when any $f_j$
occurs an odd number of times.  However, one can use the monotonicity
of $\ell^p$ norms as in Section \ref{monotonicity} to deal with this.

        These estimates for $p = 2^k$ imply analogous estimates for $2
\le p \le 2^k$, as in Section \ref{L^p estimates}.  In particular,
there are analogous estimates for every $p \in (2, \infty)$ when the
$f_j$'s have bounded $L^p$ norms for each $p \in (2, \infty)$.  Using
the upper bound for $k = 2$, one also gets that
\begin{equation}
        \Big(\sum_{j = 1}^n |a_j|^2\Big)^{1/2}
 \le C(p, L_4) \, \bigg\|\sum_{j = 1}^n a_j \, f_j\biggr\|_{L^p(X, \mu)}
\end{equation}
for $0 < p < 2$, as in Section \ref{L^p estimates}.  Here $C(p, L_4)$
is a positive constant that does not depend on $n$, but does depend on
$p$ and the upper bound $L_4$ for the $L^4$ norms of the $f_j$'s.

        Suppose now that $(X_1, \mu_1), (X_2, \mu_2), \ldots$ is an
infinite sequence of probability spaces, $X = \prod_{j = 1}^\infty
X_j$ is their product, and $\mu$ is the corresponding product measure
on $X$.  Let $f_1, f_2, \ldots$ be real or complex-valued functions on
$X_1, X_2, \ldots$, respectively, which can be identified with
functions on $X$ that are constant in the other variables.  As before,
suppose also that $f_j \in L^2(X_j, \mu_j)$ satisfies (\ref{int_{X_j}
f_j d mu_j = 0}) and (\ref{||f_j||_{L^2(X_j, mu_j)} = (int_{X_j}
|f_j|^2 d mu_j)^{1/2} = 1}) for each $j$, so that the $f_j$'s are
orthonormal in $L^2(X, \mu)$.  If $a_1, a_2, \ldots$ is a sequence of
real or complex numbers such that $\sum_{j = 1}^\infty |a_j|^2$
converges, then $\sum_{j = 1}^\infty a_j \, f_j$ converges in $L^2(X,
\mu)$.  If $k \in {\bf Z}_+$, $p = 2^k$, and $f_j \in L^p(X_j, \mu_j)$
for each $j$, with uniformly bounded $L^p$ norm, then it follows from
the previous estimates that $\sum_{j = 1}^\infty a_j \, f_j$ converges
in $L^p(X, \mu)$.  More precisely, $\sum_{j = 1}^\infty a_j \, f_j$
converges in $L^p(X, \mu)$ in the generalized sense, as in Section
\ref{generalized convergence, 2}.  In particular, if $f_j \in L^p(X_j,
\mu_j)$ for every $j \ge 1$ and $p \in (2, \infty)$, with
$\|f_j\|_{L^p(X_j, \mu_j)}$ uniformly bounded in $j$ for each $p > 2$,
then $\sum_{j = 1}^\infty a_j \, f_j$ converges in $L^p(X, \mu)$ in
the generalized sense for each $p \in (2, \infty)$.  If the $(X_j,
\mu_j)$'s are copies of the same space, and the $f_j$ are copies of
the same function on this space, then of course the $f_j$'s have the
same $L^p$ norm for each $j$.

\section{Linear functions on ${\bf R}^n$}
\label{linear functions on R^n}
\setcounter{equation}{0}

        Let $\mu$ be a Borel probability measure on ${\bf R}^n$ that
is not the Dirac mass at $0$, so that
\begin{equation}
        \mu({\bf R}^n \backslash \{0\}) > 0.
\end{equation}
Remember that a linear transformation $T$ from ${\bf R}^n$ onto itself
is said to be an \emph{orthogonal transformation} if $T$ preserves the
standard inner product on ${\bf R}^n$, and hence the standard
Euclidean norm on ${\bf R}^n$.  Suppose that $\mu$ is invariant under
orthogonal transformations, in the sense that
\begin{equation}
        \mu(T(E)) = \mu(E)
\end{equation}
for every Borel set $E \subseteq {\bf R}^n$ and every orthogonal
transformation $T$ on ${\bf R}^n$.  For example, $\mu$ might be
surface measure on the unit sphere normalized to have total measure
$1$, or $\mu$ could be absolutely continuous with respect to Lebesgue
measure, with a radial density.  Also let $p$ be a positive real
number, and suppose that
\begin{equation}
        \int_{{\bf R}^n} |x|^p \, d\mu(x) < \infty.
\end{equation}
Note that this integral is positive, by hypothesis.  If $\mu$ is
normalized surface measure on the unit sphere, then this condition
halds for every $p > 0$.  If $\mu$ is given by a radial density times
Lebesgue measure, then this condition depends on the integrability
properties of the density.

        Consider
\begin{equation}
        \lambda_v(x) = \sum_{j = 1}^n x_j \, v_j
\end{equation}
for each $v \in {\bf R}^n$.  This is a linear function on ${\bf R}^n$,
and every real-valued linear function on ${\bf R}^n$ is of this form.
By hypothesis, $\lambda_v \in L^p({\bf R}^n, \mu)$ for each $v \in
{\bf R}^n$.  Because of invariance under orthogonal transformations,
\begin{equation}
        \|\lambda_v\|_{L^p({\bf R}^n, \mu)}
 = \Big(\int_{{\bf R}^n} |\lambda_v(x)|^p \, d\mu(x)\Big)^{1/p}
           = C(p, \mu) \, |v|,
\end{equation}
where
\begin{equation}
        C(p, \mu) = \Big(\int_{{\bf R}^n} |x_1|^p \, d\mu(x)\Big)^{1/p}
\end{equation}
and
\begin{equation}
        |v| = \Big(\sum_{j = 1}^n v_j^2\Big)^{1/2}
\end{equation}
is the standard norm on ${\bf R}^n$.  Note that $0 < C(p, \mu) < \infty$.

        Remember that
\begin{equation}
\label{int_{-infty}^infty exp (- t^2) dt = sqrt{pi}}
        \int_{-\infty}^\infty \exp (- t^2) \, dt = \sqrt{\pi}.
\end{equation}
To see this, one can begin with
\begin{eqnarray}
        \Big(\int_{-\infty}^\infty \exp (-t^2) \, dt\Big)^2
           & = & \Big(\int_{-\infty}^\infty \exp (-t^2) \, dt\Big)
                  \Big(\int_{-\infty}^\infty \exp (-u^2) \, du\Big) \\
           & = & \int_{{\bf R}^2} \exp (-t^2 - u^2) \, dt du. \nonumber
\end{eqnarray}
Using polar coordinates, we get that
\begin{equation}
        \Big(\int_{-\infty}^\infty \exp (-t^2) \, dt\Big)^2
          = 2 \pi \int_0^\infty r \, \exp (-r^2) \, dr.
\end{equation}
The derivative of $\exp (-r^2)$ is $- 2 \, r \, \exp (-r^2)$, and so
\begin{equation}
        \int_0^\infty 2 \, r \, \exp (-r^2) \, dr = 1.
\end{equation}
This implies (\ref{int_{-infty}^infty exp (- t^2) dt = sqrt{pi}}), as desired.

        Let $\mu_n$ be the measure on ${\bf R}^n$ given by $\pi^{-n/2}
\, \exp (-|x|^2)$ times Lebesgue measure.  Thus $\mu_n({\bf R}^n) =
1$, by the previous computations, and $\mu_n$ is clearly invariant
under orthogonal transformations.  Also, $|x|^p \in L^p({\bf R}^n,
\mu_n)$ for every $p > 0$.  Moreover, $\mu_n$ is the same as the
product of n copies of $\mu_1$ on $n$ copies of ${\bf R}$, as in the
previous section.

\section{Countability conditions}
\label{countability conditions}
\setcounter{equation}{0}

        Remember that a collection $\beta$ of open subsets of a
topological space $X$ is said to be a \emph{base} for the topology of
$X$ if for every open set $U$ in $X$ and every point $p \in U$ there
is an open set $V \in \beta$ such that $p \in V$ and $V \subseteq U$.
In this case,
\begin{equation}
        U = \bigcup \{V : V \in \beta, \, V \subseteq U\}
\end{equation}
for every open set $U$ in $X$.  Conversely, $\beta$ is a base for the
topology of $X$ if every open set in $X$ can be expressed as a union
of elements of $\beta$.  It is especially nice to have a base $\beta$
for the topology of $X$ with only finitely or countably many elements.
This implies that there is a dense set in $X$ with only finitely or
countably many elements, by picking an element in each nonempty open
set in the base.  Conversely, if the topology on $X$ is determined by
a metric, and if there is a dense set in $X$ with only finitely or
countably many elements, then there is a base for the topology of $X$
with only finitely or countably many elements.  More precisely, the
collection of open balls in $X$ with centers contained in a dense
subset of $X$ and radii of the form $1/n$, $n \in {\bf Z}_+$, is a
base for the topology of $X$.

        Suppose that $\beta$ is a base for the topology of $X$ with
only finitely or countably many elements, and let $\{U_i\}_{i \in I}$
be a collection of open subsets of $X$.  For each $i \in I$, let
$\beta_i$ be the set of $V \in \beta$ such that $V \subseteq U_i$.
Thus
\begin{equation}
        U_i = \bigcup \{V : V \in \beta_i\}
\end{equation}
for each $i \in I$, because $\beta$ is a base for the topology of $X$.
If $\beta' = \bigcup_{i \in I} \beta_i$, then it follows that
\begin{equation}
        \bigcup_{i \in I} U_i = \bigcup \{V : V \in \beta'\}.
\end{equation}
For each $V \in \beta'$, let $i(V)$ be an element of $I$ such that $V
\subseteq U_{i(V)}$.  Also let $I'$ be the set of $i(V)$, $V \in
\beta'$.  Note that $I'$ has only finitely or countably many elements,
because $\beta' \subseteq \beta$ has only finitely or countably many
elements.  In addition,
\begin{equation}
        \bigcup_{i \in I'} U_i \subseteq \bigcup_{i \in I} U_i
   = \bigcup \{V : V \in \beta'\} \subseteq \bigcup U_{i(V)} : V \in \beta'\}
                                              = \bigcup_{i \in I'} U_i,
\end{equation}
which implies that $\bigcup_{i \in I'} U_i = \bigcup_{i \in I} U_i$.

        A set $E \subseteq X$ is said to be \emph{$\sigma$-compact} if
there is a sequence $K_1, K_2, \ldots$ of compact subsets of $X$ such
that $E = \bigcup_{n = 1}^\infty K_n$.  Suppose that $X$ is a locally
compact Hausdorff space, and that $U$ is an open set in $X$.  For each
$p \in U$, let $U(p)$ be an open set in $X$ such that $p \in U(p)$,
$\overline{U(p)}$ is compact, and $\overline{U(p)} \subseteq U$.  If
there is a base for the topology of $X$ with only finitely or
countably many elements, then it follows that there is a set $A
\subseteq U$ with only finitely or countably many elements such that
$U = \bigcup_{p \in A} U(p)$.  Hence $U = \bigcup_{p \in A}
\overline{U(p)}$, so that $U$ is $\sigma$-compact.

        Suppose that $X$ is a locally compact Hausdorff space in which
every open set is $\sigma$-compact.  As in Theorem 2.18 in \cite{r2},
every positive Borel measure $\mu$ on $X$ such that $\mu(K) < \infty$
when $K \subseteq X$ is compact automatically satisfies strong
regularity properties.  It is easy to see that the real line has this
property, for instance, as well as ${\bf R}^n$ for every positive
integer $n$.  If $X$ is a locally compact Hausdorff space, and there
is a base for the topology of $X$ with only finitely or countably many
elements, then $X$ has this property, by the remarks in the previous
paragraph.

\section{Separation conditions}
\label{separation conditions}
\setcounter{equation}{0}

        Remember that a topological space $X$ satisfies the
\emph{first separation condition} if for every pair of distinct
elements $p$, $q$ of $X$ there is an open set $U \subseteq X$ such
that $p \in U$ and $q \not\in U$.  This is equivalent to asking that
every set $A \subseteq X$ with exactly one element be closed, which
implies that finite subsets of $X$ are closed.  Similarly, $X$
satisfies the \emph{second separation condition} if for every pair
$p$, $q$ of distinct elements of $X$ there are disjoint open subsets
$U$, $V$ of $X$ such that $p \in U$, $q \in V$.  In this case, $X$ is
said to be a \emph{Hausdorff} topological space, and $X$ clearly
satisfies the first separation condition.  If $X$ satisfies the first
separation condition and for every point $p \in X$ and closed set $B
\subseteq X$ with $p \not\in B$ there are disjoint open subsets $U$,
$V$ of $X$ such that $p \in U$ and $B \subseteq V$, then $E$ satisfies
the \emph{third separation condition}, and is also said to be
\emph{regular}.  Note that regular topological spaces are Hausdorff,
since one can take $B = \{q\}$ when $q \in X$ and $q \ne p$.  If $X$
satisfies the first separation condition and for every pair $A$, $B$
of disjoint closed subsets of $X$ there are disjoint open sets $U$,
$V$ such that $A \subseteq U$, $B \subseteq V$, then $X$ satisfies the
\emph{fourth separation condition}, and is also said to be
\emph{normal}.  As before, normal spaces are automatically Hausdorff
and regular.  It is well known that metric spaces are normal.

        Equivalently, $X$ is Hausdorff if for every pair of distinct
elements $p$, $q$ of $X$ there is an open set $U \subseteq X$ such
that $p \in U$ and $q$ is not in the closure $\overline{U}$ of $U$.
Similarly, $X$ satisfies the third separation condition if and only if
it satisfies the first separation condition and for every point $p \in
X$ and open set $W \subseteq X$ with $p \in W$ there is an open set $U
\subseteq X$ such that $p \in U$ and $\overline{U} \subseteq W$.  This
formulation of regularity makes it clear that it is a local property.
In the same way, $X$ is normal if and only if for every closed set $A
\subseteq X$ and open set $W \subseteq X$ with $A \subseteq W$ there
is an open set $U \subseteq X$ such that $A \subseteq U$ and
$\overline{U} \subseteq W$.

        If $X$ is Hausdorff, then compact subsets of $X$ are closed,
and one can show that $X$ satisfies the analogues of regularity and
normality for compact sets instead of closed sets.  This implies that
compact Hausdorff spaces are normal, because closed sets of compact
spaces are compact.  If $X$ is regular, then one can show that $X$
satisfies the analogue of normality in which at least one of the
closed sets is compact.  One can also show that locally compact
Hausdorff spaces are regular.

        It is easy to see that the Cartesian product of a family of
topological spaces that satisfy the first or second separation
condition has the same property with respect to the product topology.
This is because a pair of distinct elements of the product are
different in at least one coordinate, and the appropriate separation
condition can then be applied in the corresponding space.  One can
also check that a product of regular spaces is regular.  This uses the
local characterization of regularity mentioned before.

\section{Metrizability}
\label{metrizability}
\setcounter{equation}{0}

        Let $(X, d(x, y))$ be a metric space, and put
\begin{equation}
        B(p, r) = \{x \in X : d(p, x) < r\}
\end{equation}
for each $p \in X$ and $r > 0$.  This is the open ball in $X$ with
center $p$ and radius $r$, which is well known to be an open set in
$X$, by the triangle inequality.  If $A \subseteq X$ and $r > 0$, then
\begin{equation}
\label{A_r}
        A_r = \bigcup_{p \in A} B(p, r)
            = \{x \in X : d(x, p) < r \hbox{ for some } p \in A\}
\end{equation}
is an open set in $X$ that contains $A$.  It is easy to check that
\begin{equation}
        \overline{A} = \bigcap_{r > 0} A_r = \bigcap_{n = 1}^\infty A_{1/n},
\end{equation}
where $\overline{A}$ denotes the closure of $A$ in $X$.  In
particular, every closed set in $X$ can be expressed as the
intersection of a sequence of open sets.  This implies that every open
set in $X$ can be expressed as the union of a sequence of closed sets.
If $X$ is compact, then every closed set in $X$ is compact, and hence
every open set in $X$ is $\sigma$-compact.  If $X$ is
$\sigma$-compact, then every closed set in $X$ is $\sigma$-compact,
and it follows that every open set in $X$ is $\sigma$-compact as well.

        Now let $(X_1, d_1), (X_2, d_2), \ldots$ be a sequence of
metric spaces, and let $X = \prod_{j = 1}^\infty X_j$ be their
Cartesian product, with the product topology.  One can check that
\begin{equation}
\label{d(x, y) = max_{j ge 1} (min(d_j(x_j, y_j), 1/j))}
        d(x, y) = \max_{j \ge 1} (\min(d_j(x_j, y_j), 1/j))
\end{equation}
defines a metric on $X$ for which the corresponding topology is the product
topology, where $x = \{x_j\}_{j = 1}^\infty$, $y = \{y_j\}_{j = 1}^\infty$.
In particular, $X$ may be considered as a compact metric space when $X_j$
is compact for each $j$.

        Uhrysohn's famous metrization theorem implies that there is a
metric on a topological space $X$ that determines the same topology
when $X$ is regular and there is a countable base for the topology of
$X$.  If $X$ is compact, and the topology on $X$ is determined by a
metric, then it is easy to show that there is a dense set in $X$ with
only finitely or countably many elements, which implies that there is
a base for the topology of $X$ with only finitely or countably many
elements.  This also works when $X$ is $\sigma$-compact.  Thus a base
for the topology of $X$ with only finitely or countably many elements
is necessary for metrizability of a compact or $\sigma$-compact
topological space.

\section{Partitions of unity}
\label{partitions of unity}
\setcounter{equation}{0}

        Let $X$ be a compact Hausdorff topological space.  Suppose
that for each $p \in X$, we have an open set $U(p)$ in $X$ such that
$p \in U(p)$.  By Uhryson's lemma, there is a nonnegative continuous
real-valued function $\phi_p(x)$ on $X$ such that $\phi(p) > 0$ and
the support of $\phi_p$ is contained in $U(p)$.  If
\begin{equation}
\label{U_1(p) = {x in X : phi_p(x) > 0}}
        U_1(p) = \{x \in X : \phi_p(x) > 0\},
\end{equation}
then $U_1(p)$ is an open set in $X$ such that $p \in U_1(p)$ and
$U_1(p) \subseteq U(p)$.  By compactness, there are finitely many
elements $p_1, \ldots, p_n$ of $X$ such that
\begin{equation}
        X = \bigcup_{j = 1}^n U_1(p_j).
\end{equation}
This implies that $\sum_{j = 1}^n \phi_{p_j}(x) > 0$ for every $x \in
X$.  Hence
\begin{equation}
        \psi_j(x) = \frac{\phi_{p_j}(x)}{\sum_{l = 1}^n \phi_{p_l}(x)}
\end{equation}
defines a nonnegative continuous real-valued function on $X$.  Also,
\begin{equation}
        \sum_{j = 1}^n \psi_j(x) = 1
\end{equation}
for every $x \in X$, and $\psi_j(x) > 0$ if and only if $\phi_{p_j}(x) > 0$.

        As an application, let $V$ be a real or complex vector space
equipped with a norm $\|v\|$, and let $f$ be a continuous mapping from
$X$ into $V$.  Let $\epsilon > 0$ be given, and let $U(p)$ be an open
set in $X$ such that $p \in U(p)$ and
\begin{equation}
\label{||f(x) - f(p)|| < epsilon}
        \|f(x) - f(p)\| < \epsilon
\end{equation}
for every $x \in U(p)$.  Put
\begin{equation}
\label{g(x) = sum_{j = 1}^n psi_j(x) f(p_j)}
        g(x) = \sum_{j = 1}^n \psi_j(x) \, f(p_j),
\end{equation}
where $p_1, \ldots, p_n$ and $\psi_1, \ldots, \psi_n$ are as in the
previous paragraph.  Thus 
\begin{equation}
\label{||f(x) - g(x)|| le sum_{j = 1}^n psi_j(x) ||f(x) - f(p_j)|| < epsilon}
 \|f(x) - g(x)\| \le \sum_{j = 1}^n \psi_j(x) \, \|f(x) - f(p_j)\| < \epsilon
\end{equation}
for every $x \in X$, using (\ref{||f(x) - f(p)|| < epsilon}) and the
fact that $x \in U(p_j)$ when $\psi_j(x) > 0$.  The same argument
works when the topology on $V$ is determined by a collection
$\mathcal{N}$ of seminorms, and $\|v\|$ is replaced by the maximum of
finitely many seminorms in $\mathcal{N}$.

\section{Product spaces}
\label{product spaces}
\setcounter{equation}{0}

        Let $X$, $Y$ be compact Hausdorff topological spaces, and let
$X \times Y$ be their Cartesian product, equipped with the product
topology.  Thus $X \times Y$ is also a compact Hausdorff space.  Also
let $f(x, y)$ be a continuous real or complex-valued function on $X
\times Y$, and let $\epsilon > 0$ be given.  For each $x \in X$ and $y
\in Y$, there are open sets $U(x, y) \subseteq X$, $V(x, y) \subseteq
Y$ such that $x \in U(x, y)$, $y \in V(x, y)$, and
\begin{equation}
\label{|f(u, v) - f(w, z)| < epsilon}
        |f(u, v) - f(w, z)| < \epsilon
\end{equation}
for every $u, w \in U(x, y)$ and $v, z \in V(x, y)$, by the continuity
of $f$ at $(x, y)$ and the definition of the product topology.  If we
fix $x \in X$ for a moment, and apply this to each $y \in Y$, then the
open sets $V(x, y)$, $y \in Y$, form an open covering of $Y$.  By compactness
of $Y$, there are finitely many elements $y_1, \ldots, y_n$ of $Y$ such that
\begin{equation}
        Y = \bigcup_{j = 1}^n V(x, y_j).
\end{equation}
Put $U(x) = \bigcap_{j = 1}^n U(x, y_j)$, so that $U(x)$ is an open
set in $X$ that contains $x$.  Moreover,
\begin{equation}
        |f(u, y) - f(w, y)| < \epsilon
\end{equation}
for every $u, w \in U(x)$ and $y \in Y$, by applying (\ref{|f(u, v) -
f(w, z)| < epsilon}) to $v = z = y$, which is contained in $V(y_j)$
for some $j$.  Similarly, one can use compactness of $X$ to show that
for every $y \in Y$ there is an open set $V(y) \subseteq Y$ such that
$y \in V(y)$ and
\begin{equation}
        |f(x, v) - f(x, z)| < \epsilon
\end{equation}
for every $v, z \in V(y)$ and $x \in X$.

        Let $\mu$, $\nu$ be regular Borel probability measures on $X$,
$Y$, respectively.  By the Riesz representation theorem, this is
equivalent to having positive linear functionals on the spaces of
continuous functions on $X$, $Y$ that take the value $1$ on the
constant functions identically equal to $1$ on these spaces.  If $f(x,
y)$ is a continuous function on $X \times Y$, then it follows from
the uniform continuity properties in the previous paragraph that
\begin{equation}
        \int_X f(x, y) \, d\mu(x), \quad \int_Y f(x, y) \, d\nu(y)
\end{equation}
are continuous functions on $Y$, $X$, respectively.  Thus
\begin{equation}
        \int_Y \Big(\int_X f(x, y) \, d\mu(x)\Big) d\nu(y), \quad
         \int_X \Big(\int_Y f(x, y) \, d\nu(y)\Big) d\mu(x)
\end{equation}
define nonnegative linear functionals on the space of continuous
functions on $X \times Y$ that take the value $1$ on the constant
function $1$.  One can also show that these two linear functionals are
the same, because they are the same when $f$ is a linear combination
of products of continuous functions on $X$ and $Y$, and because these
functions are dense in the space of all continuous functions on $X
\times Y$ with respect to the supremum norm.  The latter statement can
be verified using partitions of unity on $X$ and uniform continuity
over $Y$, for instance, as in the preceding section and paragraph.
The Riesz representation theorem implies that there is a unique
regular Borel probability measure $\mu \times \nu$ on $X \times Y$
such that this linear functional on the space of continuous functions
on $X \times Y$ is given by
\begin{equation}
        \int_{X \times Y} f(x, y) \, d(\mu \times \nu)(x, y).
\end{equation}
There are analogous arguments for nonnegative Borel measures with
suitable regularity properties on locally compact Hausdorff spaces,
which correspond to nonnegative linear functionals on continuous
functions with compact support on these spaces.  If the measures are
finite, then one can simply compactify the spaces using one-point
compactifications.

        Let $\beta_X$, $\beta_Y$ be bases for the topologies of $X$,
$Y$, respectively.  It is easy to see that
\begin{equation}
 \beta_{X \times Y} = \{U \times V : U \in \beta_X, \, V \in \beta_Y\}
\end{equation}
is a base for the topology of $X \times Y$.  In particular, $\beta_{X
\times Y}$ has only finitely or countably many elements when
$\beta_X$, $\beta_Y$ have only finitely or countably many elements.
In this case, it follows that every open set in $X \times Y$ is the
union of finitely or countably many products of open subsets of $X$
and $Y$.  Otherwise, one can check that an open set in $X \times Y$
that is also $\sigma$-compact is the union of finitely or countably
many products of open subsets of $X$ and $Y$.

\section{Product spaces, 2}
\label{product spaces, 2}
\setcounter{equation}{0}

        Let $I$ be a nonempty set, and suppose that for each $i \in I$
we have a topological space $X_i$.  In practice, we shall be
interested in sets $I$ with only finitely or countably many elements.
Let $X = \prod_{i \in I} X_i$ be the corresponding Cartesian product,
equipped with the product topology.

        Suppose that $\beta_i$ is a base for the topology of $X_i$ for
each $i \in I$, and let $\beta$ be the collection of subsets of $X$ of
the form $\prod_{i \in I} U_i$, where $U_i \in \beta_i$ for each $i
\in I$, and $U_i = X_i$ for all but finitely many $i$.  It is easy to
check that $\beta$ is a base for the product topology on $X$.  If $I$
has only finitely or countably many elements, and each $\beta_i$ has
only finitely or countably many elements, then $\beta$ has only
finitely or countably many elements too.  This follows from the fact
that the Cartesian product of finitely many countable sets is
countable when $I$ has only finitely many elements.  If $I$ is a
countably infinite set, then one can use the same argument for finite
subsets of $I$, and apply this to an increasing sequence of finite
subsets of $I$ whose union is all of $I$.

        If $X_i$ is Hausdorff for each $i \in I$, then $X$ is
Hausdorff.  If $X_i$ is compact for each $i \in I$, then $X$ is
compact, by Tychonoff's theorem.  Of course, this is much more
elementary when $I$ has only finitely many elements.  If $I$ has only
finitely or countably many elements and each $X_i$ is metrizable, then
$X$ is metrizable, and compactness can be handled in a simpler way
using sequential compactness.  This approach can also be applied
directly when $I$ has only finitely or countably many elements and
there is a base for the topology of $X_i$ with only finitely or
countably many elements for each $i \in I$, so that there is also a
base for the topology of $X$ with only finitely or countably many elements.

        Let $f$ be a continuous real or complex-valued function on $X$.
For each $\epsilon > 0$ and $x \in X$, there is an open set $U(x)$ in $X$
such that $x \in U(x)$ and
\begin{equation}
\label{|f(y) - f(z)| < epsilon}
        |f(y) - f(z)| < \epsilon
\end{equation}
for every $y, z \in U(x)$.  More precisely, we can take $U(x)$ to be a
basic open set in the product topology, so that there is a finite set
$I(x) \subseteq I$ such that $U(x) = \prod_{i \in I} U_i(x)$ for some
open sets $U_i(x) \subseteq X$, where $U_i(x) = X_i$ for every $i \in
I \backslash I(x)$.  In particular, if $y \in U(x)$, $z \in X$, and
$y_i = z_i$ for each $i \in I(x)$, then it follows that $z \in U(x)$,
and hence (\ref{|f(y) - f(z)| < epsilon}) holds.

        If $X_i$ is compact for each $i \in I$, so that $X$ is compact,
then there are finitely many elements $x(1), \ldots, x(n)$ of $X$ such that
\begin{equation}
        X = \bigcup_{j = 1}^n U(x(n)).
\end{equation}
Put $I_\epsilon = \bigcup_{j = 1}^n I(x(j))$, so that $I_\epsilon
\subseteq I$ has only finitely many elements.  If $y, z \in X$ satisfy
$y_i = z_i$ for every $i \in I_\epsilon$, then it is easy to see that
(\ref{|f(y) - f(z)| < epsilon}) holds.  This is because $y \in
U(x(j))$ for some $j = 1, \ldots, n$, and so $z \in U(x(j))$ too.
Thus continuous functions on $X$ may be approximated uniformly by
functions of finitely many variables under these conditions.

        Suppose that $\mu_i$ is a regular Borel probability measure on
$X_i$ for each $i$.  If $A \subseteq I$ is a nonempty set with only
finitely many elements, then let $L_A(f)$ be the function on $X$ which
is constant in $x_i$ for each $i \in A$ obtained by integrating $f$ in
$x_i$ with respect to $\mu_i$ for each $i \in A$.  If $A \cap
I_\epsilon = \emptyset$, then
\begin{equation}
        |L_A(f)(y) - f(y)| < \epsilon
\end{equation}
for every $y \in X$, since (\ref{|f(y) - f(z)| < epsilon}) holds for
every $z \in X$ such that $y_i = z_i$ when $i \in I \backslash A$.  If
$A, B \subseteq I$ are finite sets such that $I_\epsilon \subseteq A,
B$, then
\begin{equation}
        |L_A(f)(y) - L_B(f)(y)| < 2 \, \epsilon
\end{equation}
for every $y \in X$.  This uses the previous estimate applied to $A
\backslash B$ and $B \backslash A$, to estimate the difference between
each of $L_A(f)$, $L_B(f)$ and $L_{A \cap B}(f)$.

        Let $\mathcal{A}$ be the collection of all finite subsets of
$I$, ordered by inclusion.  This is a directed system, because for
every $A, B \in \mathcal{A}$ we have that $A \cup B \in \mathcal{A}$
and $A, B \subseteq A \cup B$.  If $f$ is a continuous function on
$X$, then one can think of $\{L_A(f)\}_{A \in \mathcal{A}}$ as a net
of functions on $X$ indexed by $\mathcal{A}$.  One can show that this
net converges uniformly to a constant on $X$ for every continuous
function on $X$.  This uses the fact that the net satisfies a uniform
Cauchy condition on $X$, as in the previous paragraph.

        In the limit, we get a positive linear functional on the space
of continuous functions on $X$ which takes the value $1$ on the
constant function $1$.  The Riesz representation theorem implies that
this linear functional can be expressed in terms of a unique regular
Borel probability measure on $X$, which corresponds to the product of
the $\mu_i$'s.  As usual, the situation is especially nice when $I$ is
countably infinite, and each $X_i$ has a base $\beta_i$ for its
topology with only finitely or countably many elements.  This leads to
a base $\beta$ for the topology of $X$ consisting of only finitely or
countably many basic open sets in $X$, as before, which implies in
particular that every open set in $X$ is the union of finitely or
countably many basic open sets.  Otherwise, every open set in $X$ that
is also $\sigma$-compact is the union of finitely or countably many
basic open sets, as in the previous section.

\part{Conditional expectation and martingales}

\section{$\sigma$-Subalgebras}
\label{sigma-subalgebras}
\setcounter{equation}{0}

        Let $(X, \mathcal{A}, \mu)$ be a probability space, and let
$\mathcal{B}$ be a $\sigma$-subalgebra of $\mathcal{A}$.  Thus $(X,
\mathcal{B}, \mu)$ is also a probability space, where the measure
$\mu$ is restricted to $\mathcal{B}$.  If a real or complex-valued
function $f$ on $X$ is measurable with respect to $\mathcal{B}$, then
it is automatically measurable with respect to $\mathcal{A}$ as well.
If $f$ is measurable with respect to $\mathcal{B}$ and integrable with
respect to $\mu$, then $f$ is also integrable as a function which is
measurable with respect to $\mathcal{A}$, and the integral
\begin{equation}
        \int_X f \, d\mu
\end{equation}
is the same with respect to both $\mathcal{A}$ and $\mathcal{B}$.

        For example, $\mathcal{B}$ might consist of only the empty set
$\emptyset$ and $X$ itself, in which case the only functions on $X$
that are measurable with respect to $\mathcal{B}$ are constant
functions.  As another example, one might take $X$ to be the closed
unit interval $[0, 1]$, $\mathcal{A}$ to be the $\sigma$-algebra of
Lebesgue measurable subsets of $[0, 1]$, $\mu$ to be Lebesgue measure
on $[0, 1]$, and $\mathcal{B}$ to be the $\sigma$-algebra of Borel
subsets of $[0, 1]$.  It is well known that for each Lebesure
measurable set $A \subseteq [0, 1]$ there are Borel sets $B_1, B_2
\subseteq [0, 1]$ such that $B_1 \subseteq A \subseteq B_2$ and
$\mu(B_2 \backslash B_1) = 0$.  More precisely, one can take $B_1$ to
be a countable union of compact sets, and $B_2$ to be a countable
intersection of relatively open sets in $[0, 1]$.

        Let $(X_1, \mathcal{A}_1, \mu_1)$, $(X_2, \mathcal{A}_2,
\mu_2)$ be probability spaces, and let $X = X_1 \times X_2$ be their
Cartesian product, with the corresponding product measure $\mu_1
\times \mu_2$ and $\sigma$-algebra $\mathcal{A}$.  Let $\mathcal{B}_1$
be the collection of subsets of $X$ of the form $E \times X_2$ with $E
\in \mathcal{A}_1$, and let $\mathcal{B}_2$ be the collection of
subsets of $X$ of the form $X_1 \times E$ with $E \in \mathcal{A}_2$.
It is easy to see that $\mathcal{B}_1$, $\mathcal{B}_2$ are
$\sigma$-subalgebras of $\mathcal{A}$, and that a function $f(x_1,
x_2)$ on $X$ is measurable with respect to $\mathcal{B}_1$ or
$\mathcal{B}_2$ if and only if it is measurable with respect to
$\mathcal{A}$ and constant in $x_2$ or $x_1$, respectively.  Thus
measurable functions on $X$ with respect to $\mathcal{B}_1$,
$\mathcal{B}_2$ may be identified with functions on $X_1$, $X_2$ that
are measurable with respect to $\mathcal{A}_1$, $\mathcal{A}_2$,
respectively.

        As a variant of this, suppose that $X_1$, $X_2$ are
topological spaces, and let $X = X_1 \times X_2$ be equipped with the
product topology.  If $A_1 \subseteq X_1$, $A_2 \subseteq X_2$ are
Borel sets, then $A_1 \times X_2$, $X_1 \times A_2$ are Borel sets in
$X$, by standard reasoning.  In particular,
\begin{equation}
        A_1 \times A_2 = (A_1 \times X_2) \cap (X_1 \times A_2)
\end{equation}
is a Borel set in $X$.  At any rate, the collections of subsets of $X$
of the form $A_1 \times X_2$, $X_1 \times A_2$, where $A_1$, $A_2$ are
Borel subsets of $X_1$, $X_2$, respectively, are $\sigma$-subalgebras
of the Borel sets in $X$.  As in Section \ref{product spaces}, if
there are bases for the topologies of $X_1$, $X_2$ with only finitely
or countably many elements, then every open set in $X$ is the union of
finitely or countably many products of open subsets of $X_1$ and
$X_2$.  This implies that every open set in $X$ is in the
$\sigma$-algebra generated by products of Borel sets in $X_1$, $X_2$,
and hence that every Borel set in $X$ is in this $\sigma$-algebra.  It
follows that the $\sigma$-algebra of subsets of $X$ generated by
products of Borel sets in $X_1$, $X_2$ is the same as the
$\sigma$-algebra of Borel sets in $X$ under these conditions.

\section{$L^p$ Spaces}
\label{L^p spaces}
\setcounter{equation}{0}

        Let $(X, \mathcal{A}, \mu)$ be a probability space, and let
$\mathcal{B}$ be a $\sigma$-subalgebra of $\mathcal{A}$.  If $f$, $g$
are measurable functions on $X$ with respect to $\mathcal{A}$, then
\begin{equation}
\label{x in X : f(x) = g(x)}
        \{x \in X : f(x) = g(x)\}
\end{equation}
is a measurable set in $X$ with respect to $\mathcal{A}$.  If $f$, $g$
are measurable with respect to $\mathcal{B}$, then (\ref{x in X : f(x)
= g(x)}) is measurable with respect to $\mathcal{B}$.  Of course, $f$
and $g$ are said to be equal almost everywhere with respect to $\mu$ when
\begin{equation}
        \mu(\{x \in X : f(x) \ne g(x)\}) = 0.
\end{equation}

        Let $L^p(X, \mathcal{A})$, $L^p(X, \mathcal{B})$ be the $L^p$
spaces of measurable functions on $X$ with respect to $\mathcal{A}$,
$\mathcal{B}$, for $0 < p \le \infty$.  These spaces also involve the
measure $\mu$, but we omit this from the notation when it is
unambiguous.  Because measurable functions on $X$ with respect to
$\mathcal{B}$ are also measurable with respect to $\mathcal{A}$, we
get an isometric linear embedding of $L^p(X, \mathcal{B})$ into
$L^p(X, \mathcal{A})$ for each $p$, $0 < p \le \infty$.

        Note that $L^p(X, \mathcal{B})$ corresponds to a closed linear
subspace of $L^p(X, \mathcal{A})$ for each $p$, $0 < p \le \infty$.
One way to see this is to use the completeness of $L^p(X,
\mathcal{B})$ and the fact that the embedding into $L^p(X,
\mathcal{A})$ is isometric.  Basically the same argument can be given
more explicitly as follows.  Suppose that $\{f_j\}_{j = 1}^\infty$ is
a sequence of elements of $L^p(X, \mathcal{B})$ that converges in the
$L^p$ norm to $f \in L^p(X, \mathcal{A})$.  By passing to a
subsequence, we may suppose that $\{f_j\}_{j = 1}^\infty$ converges
pointwise almost everywhere to $f$.  It is well known that the set of
$x \in X$ such that $\{f_j(x)\}_{j = 1}^\infty$ converges in ${\bf R}$
or ${\bf C}$, as appropriate, is measurable with respect to
$\mathcal{B}$, because each $f_j$ is measurable with respect to
$\mathcal{B}$.  The complement of this set has measure $0$ by
hypothesis, and we may suppose that $\{f_j(x)\}_{j = 1}^\infty$
converges in ${\bf R}$ or ${\bf C}$ for every $x \in X$, by setting
$f_j(x) = 0$ on the set where the sequence does not converge
initially.  The limit is automatically measurable with respect to
$\mathcal{B}$, and equal to $f$ almost everywhere.  This shows that
$f$ is in the image of $L^p(X, \mathcal{B})$ in $L^p(X, \mathcal{A})$,
as desired.

\section{Conditional expectation}
\label{conditional expectation}
\setcounter{equation}{0}

        Let $(X, \mathcal{A}, \mu)$ be a probability space, and let
$\mathcal{B}$ be a $\sigma$-subalgebra of $\mathcal{A}$.  If $f \in
L^1(X, \mathcal{A})$, then
\begin{equation}
\label{mu_f(A) = int_A f d mu}
        \mu_f(A) = \int_A f \, d\mu
\end{equation}
defines a real or complex measure on $\mathcal{A}$, as appropriate.
By construction, $\mu_f$ is absolutely continuous with respect to
$\mu$.  Hence the restriction of $\mu_f$ to $\mathcal{B}$ is
absolutely continuous with respect to the restriction of $\mu$ to
$\mathcal{B}$.  The Radon--Nikodym theorem implies that there is a
measurable function $f_\mathcal{B}$ on $X$ with respect to
$\mathcal{B}$ which is integrable with respect to $\mu$ and satisfies
\begin{equation}
        \mu_f(B) = \int_B f_\mathcal{B} \, d\mu
\end{equation}
for every $B \in \mathcal{B}$.  If $f'_\mathcal{B}$ is another
measurable function on $X$ with respect to $\mathcal{B}$ which is
integrable with respect to $\mu$ and satisfies
\begin{equation}
        \mu_f(B) = \int_B f'_\mathcal{B} \, d\mu
\end{equation}
for every $B \in \mathcal{B}$, then it is easy to see that
$f'_\mathcal{B} = f_\mathcal{B}$ almost everywhere with respect to
$\mu$.  Thus $f_\mathcal{B}$ is uniquely determined as an element of
$L^1(X, \mathcal{B})$.  This function $f_\mathcal{B}$ is known as the
\emph{conditional expectation} of $f$ with respect to $\mathcal{B}$,
and may be denoted $E(f \mid \mathcal{B})$.

        For example, if $\mathcal{B} = \{\emptyset, X\}$, so that only
constant functions are measurable with respect to $\mathcal{B}$, then
$E(f \mid \mathcal{B})$ reduces to the ordinary expectation
\begin{equation}
\label{E(f) = int_X f d mu}
        E(f) = \int_X f \, d\mu.
\end{equation}
If $\mathcal{A} = \mathcal{B}$, then $f_\mathcal{B} = f$.  For any
$\mathcal{A}$, $\mathcal{B}$, we can take $f_\mathcal{B} = f$ when $f$
is measurable with respect to $\mathcal{B}$.

        Let $(X_1, \mathcal{A}_1, \mu_1)$, $(X_2, \mathcal{A}_2,
\mu_2)$ be probability spaces, and let $X = X_1 \times X_2$ with the
product measure $\mu = \mu_1 \times \mu_2$ and corresponding
$\sigma$-algebra $\mathcal{A}$.  Also let $\mathcal{B}_1$,
$\mathcal{B}_2$ be the $\sigma$-subalgebras of $\mathcal{A}$ defined
in Section \ref{sigma-subalgebras}.  If $f(x_1, x_2) \in L^1(X,
\mathcal{A})$, then
\begin{eqnarray}
        f_1(x_1) & = & \int_{X_2} f(x_1, x_2) \, d\mu_2(x_2), \\
        f_2(x_2) & = & \int_{X_1} f(x_1, x_2) \, d\mu_1(x_1)
\end{eqnarray}
are defined almost everywhere on $X_1$, $X_2$, respectively, and
determine integrable functions on these spaces, as in Fubini's
theorem.  In this case,
\begin{equation}
        f_{\mathcal{B}_1}(x_1, x_2) = f_1(x_1), \quad 
           f_{\mathcal{B}_2}(x_1, x_2) = f_2(x_2)
\end{equation}
are measurable functions on $X$ with respect to $\mathcal{B}_1$,
$\mathcal{B}_2$, respectively, and satisfy the requirements of the
conditional expectation, again by Fubini's theorem.

\section{Product spaces, 3}
\label{product spaces, 3}
\setcounter{equation}{0}

        Let $X_1$, $X_2$ be compact Hausdorff topological spaces, and
let $X = X_1 \times X_2$ be their Cartesian product, with the product
topology.  Also let $\mu_1$, $\mu_2$ be regular Borel probability
measures on $X_1$, $X_2$, respectively, which may be given by positive
linear functionals on the spaces of continuous functions on $X_1$,
$X_2$ that take the value $1$ on the constant functions equal to $1$
on these spaces, by the Riesz representation theorem.  If $f(x_1,
x_2)$ is a continuous function on $X$, then
\begin{equation}
        f_1(x_1) = \int_{X_2} f(x_1, x_2) \, d\mu_2(x_2), \quad
         f_2(x_2) = \int_{X_1} f(x_1, x_2) \, d\mu_1(x_1)
\end{equation}
are continuous functions on $X_1$, $X_2$, respectively, by the uniform
continuity properties of $f(x_1, x_2)$ in each variable separately
discussed in Section \ref{product spaces}.  In addition,
\begin{equation}
        \int_{X_1} f_1(x_1) \, d\mu_1(x_1) = \int_{X_2} f_2(x_2) \, d\mu_2(x_2)
\end{equation}
defines a positive linear functional on the space of continuous
functions on $X$ that takes the value $1$ on the constant $1$, and
hence determines a regular Borel probability measure $\mu$ on $X$ by
the Riesz representation theorem, as in Section \ref{product spaces}
again.

        In this context, one can think of $\mu_f$ as the regular Borel
measure on $X$ determined by
\begin{equation}
        \phi \mapsto \int_X \phi \, f \, d\mu,
\end{equation}
as a bounded linear functional on the space of continuous functions on
$X$.  If $\psi$ is a continuous function on $X_1$, which can also be
considered as a continuous function on $X$ that is constant in $x_2$,
then this linear functional applied to $\phi(x_1, x_2) = \psi(x_1)$
reduces to
\begin{equation}
        \int_{X_1} \psi \, f_1 \, d\mu_1 = \int_X \psi \, f_1 \, d\mu.
\end{equation}
Of course, there is an analogous statement for continuous functions on
$X_2$.  In this way, conditional expectation can be expressed more
directly in terms of linear functionals on continuous functions.

\section{Measurable partitions}
\label{partitions}
\setcounter{equation}{0}

        Let $(X, \mathcal{A}, \mu)$ be a probability space, and let
$\mathcal{P}$ be a partition of $X$ consisting of finitely or
countably many measurable subsets of $X$.  Thus the elements of
$\mathcal{P}$ are pairwise-disjoint measurable subsets of $X$ whose
union is all of $X$.  Let $\mathcal{B} = \mathcal{B}(\mathcal{P})$ be
the collection of subsets of $X$ that can be expressed as unions of
elements of $\mathcal{P}$, including the empty set.  It is easy to see
that $\mathcal{B}$ is a $\sigma$-subalgebra of $\mathcal{A}$, and that
a function $f$ on $X$ is measurable with respect to $\mathcal{P}$ if
and only if $f$ is constant on each of the elements of $\mathcal{P}$.

        If $f \in L^1(X, \mathcal{A})$, then one can check that
\begin{equation}
        f_\mathcal{B}(x) = \frac{1}{\mu(A)} \int_A f \, d\mu
\end{equation}
when $x \in A \in \mathcal{P}$ and $\mu(A) > 0$.  Let us ask that
$\mu(A) > 0$ for every $A \in \mathcal{P}$, for the sake of
simplicity.  Thus $f_\mathcal{B}(x)$ is defined for every $x \in X$ by
this expression, and is constant on elements of $\mathcal{P}$, and
hence is measurable with respect to $\mathcal{B}$.

        If $\nu$ is a real or complex measure on $\mathcal{A}$, then
the restriction of $\nu$ to a $\sigma$-subalgebra $\mathcal{B}$ of
$\mathcal{A}$ may be absolutely continuous with respect to the
restriction of $\mu$ to $\mathcal{B}$, even if $\nu$ is not absolutely
continuous with respect to $\mu$ on $\mathcal{A}$.  In this case, the
Radon--Nikodym theorem implies that there is a unique $f_\mathcal{B}
\in L^1(X, \mathcal{B})$ such that
\begin{equation}
        \nu(B) = \int_B f_\mathcal{B} \, d\mu
\end{equation}
for every $B \in \mathcal{B}$, as before.  If $\mathcal{B} =
\mathcal{B}(\mathcal{P})$ and $\mu(A) > 0$ for every $A \in
\mathcal{P}$, then any measure on $\mathcal{B}$ is absolutely
continuous with respect to the restriction of $\mu$ to $\mathcal{B}$.
As in the previous situation,
\begin{equation}
        f_\mathcal{B}(x) = \frac{\nu(A)}{\mu(A)}
\end{equation}
for every $x \in A \in \mathcal{P}$.

\section{Basic properties}
\label{basic properties}
\setcounter{equation}{0}

        Let $(X, \mathcal{B}, \mu)$ be a measure space, and let
$g_0$ be a real-valued integrable function on $X$.  If
\begin{equation}
        \int_B g_0 \, d\mu \ge 0
\end{equation}
for every $B \in \mathcal{B}$, then $g_0 \ge 0$ almost everywhere on $X$.
To see this, put
\begin{equation}
        B_0 = \{x \in X : g_0(x) < 0\}.
\end{equation}
If $\mu(B_0) > 0$, then it follows that
\begin{equation}
        \int_{B_0} g_0 \, d\mu < 0,
\end{equation}
a contradiction.

        Suppose now that $g$ is a real or complex-valued integrable
function on $X$, and that $h$ is a nonnegative real-valued integrable
function on $X$ such that
\begin{equation}
       \biggl|\int_B g \, d\mu\biggr| \le \int_B h \, d\mu
\end{equation}
for every $B \in \mathcal{B}$.  We would like to check that $|g| \le
h$ almost everywhere on $X$ under these conditions.  If $g$ is
real-valued, then we can apply the previous argument to $h \pm g$, to
get that $h \pm g \ge 0$ almost everywhere on $X$.  If $g$ is
complex-valued, then the same argument shows that $\re \alpha \, g \le
h$ almost everywhere on $X$ for every $\alpha \in {\bf C}$ with
$|\alpha| = 1$.  This implies that $|g| \le h$ almost everywhere, by
using a countable dense set of $\alpha$'s in the unit circle.

        Now let $(X, \mathcal{A}, \mu)$ be a probability space, and
let $\mathcal{B}$ be a $\sigma$-subalgebra of $\mathcal{A}$.  If $f
\in L^1(X, \mathcal{A})$ is real-valued and nonnegative, then
\begin{equation}
\label{int_B f_mathcal{B} d mu = int_B f d mu ge 0}
        \int_B f_\mathcal{B} \, d\mu = \int_B f \, d\mu \ge 0
\end{equation}
for every $B \in \mathcal{B}$.  This implies that $f_\mathcal{B} \ge
0$ almost everywhere on $X$, by the argument at the beginning of the
section.  Of course, it is important here that $f_\mathcal{B}$ is also
measurable with respect to $\mathcal{B}$.  Similarly, if $f > 0$
almost everywhere on $X$, then
\begin{equation}
        \int_B f_\mathcal{B} \, d\mu = \int_B f \, d\mu > 0
\end{equation}
for every $B \in \mathcal{B}$ with $\mu(B) > 0$, and one can use this
to show that $f_\mathcal{B} > 0$ almost everywhere on $X$.

        If $f$ is any integrable function on $X$ that is measurable
with respect to $\mathcal{A}$, then we can apply the preceding
observation to $|f|$ to get that
\begin{equation}
        |f|_\mathcal{B} = E(|f| \mid \mathcal{B}) \ge 0
\end{equation}
almost everywhere on $X$.  Moreover,
\begin{equation}
 \biggl|\int_B f_\mathcal{B} \, d\mu\biggr| = \biggl|\int_B f \, d\mu\biggr|
                      \le \int_B |f| \, d\mu = \int_B |f|_\mathcal{B} \, d\mu
\end{equation}
for every $B \in \mathcal{B}$, which implies that
\begin{equation}
\label{|f_mathcal{B}| le |f|_mathcal{B}}
        |f_\mathcal{B}| \le |f|_\mathcal{B}
\end{equation}
almost everywhere on $X$, by the earlier remarks.  As before, it is
important here that both $f_\mathcal{B}$ and $|f|_\mathcal{B}$ are
measurable with respect to $\mathcal{B}$, to apply the arguments at
the beginning of the section.  In particular,
\begin{equation}
\label{int_X |f_mathcal{B}| d mu le ... = int_X |f| d mu}
        \int_X |f_\mathcal{B}| \, d\mu \le \int_X |f|_\mathcal{B} \, d\mu
                                                      = \int_X |f| \, d\mu,
\end{equation}
using the fact that $X \in \mathcal{B}$ in the last step.

        Alternatively, let $\nu$ be a real or complex measure on
$\mathcal{A}$, and let $|\nu|$ be the corresponding total variation
measure on $\mathcal{A}$.  Also let $\nu_\mathcal{B}$ be the
restriction of $\nu$ to $\mathcal{B}$, and let $|\nu_\mathcal{B}|$ be
its total variation, as a measure on $\mathcal{B}$.  It is easy to
see that
\begin{equation}
        |\nu_\mathcal{B}|(B) \le |\nu|(B)
\end{equation}
for every $B \in \mathcal{B}$, so that $|\nu_\mathcal{B}|$ is less
than or equal to the restriction of $|\nu|$ to $\mathcal{B}$.  If $f
\in L^1(X, \mathcal{A})$ and $\mu_f$ is as in (\ref{mu_f(A) = int_A f
d mu}), then one can show that $|\mu_f| = \mu_{|f|}$.  This gives
another way to look at (\ref{|f_mathcal{B}| le |f|_mathcal{B}}), since
the restriction of $\mu_f$ to $\mathcal{B}$ is given by integrating
$f_\mathcal{B}$.

        Note that $f \mapsto f_\mathcal{B}$ defines a linear mapping
from $L^1(X, \mathcal{A})$ into $L^1(X, \mathcal{B})$, because of the
uniqueness of the conditional expectation.  More precisely, this
mapping sends $L^1(X, \mathcal{A})$ onto $L^1(X, \mathcal{B})$,
because $f_\mathcal{B} = f$ when $f$ is measurable with respect to
$\mathcal{B}$.  If $f$, $f'$ are real-valued integrable functions on
$X$ that are measurable with respect to $\mathcal{A}$ and satisfy $f
\le f'$ almost everywhere on $X$, then
\begin{equation}
        f_\mathcal{B} \le f'_\mathcal{B}
\end{equation}
almost everywhere on $X$.  This follows from the linearity of the
conditional expectation and the fact that $f' - f \ge 0$ almost
everywhere, so that $(f' - f)_\mathcal{B} \ge 0$ almost everywhere on
$X$.  If $f$ is a real or complex-valued integrable function on $X$
and $f'$ is a nonnegative real-valued integrable function on $X$ such
that $|f| \le f'$ almost everywhere, then we get that
\begin{equation}
        |f_\mathcal{B}| \le |f|_\mathcal{B} \le f'_\mathcal{B}
\end{equation}
almost everywhere on $X$.  In particular, this holds when $f'$ is a
constant, in which case $f'_\mathcal{B}$ is the same constant.  This
implies that $f_\mathcal{B} \in L^\infty(X, \mathcal{B})$ when $f \in
L^\infty(X, \mathcal{A})$, with
\begin{equation}
\label{||f_mathcal{B}||_infty le ||f||_infty}
        \|f_\mathcal{B}\|_\infty \le \|f\|_\infty.
\end{equation}

        Let $f$ be a real-valued integrable function on $X$ that is
measurable with respect to $\mathcal{A}$ and takes values in an
interval $I \subseteq {\bf R}$ almost everywhere.  This interval may
be open, closed, or half-open and half-closed, and it may also be
unbounded, such as a half-line or the whole real line.  One can check
that $f_\mathcal{B}$ takes values in $I$ almost everywhere as well, by
comparing $f$ with constant functions.  If $\phi : I \to {\bf R}$ is
convex and $\phi \circ f$ is integrable on $X$, then Jensen's
inequality implies that
\begin{equation}
\label{phi(frac{1}{mu(A)} int_A f d mu) le frac{1}{mu(A)} int_A phi circ f dmu}
        \phi\Big(\frac{1}{\mu(A)} \int_A f \, d\mu\Big)
           \le \frac{1}{\mu(A)} \int_A \phi \circ f \, d\mu
\end{equation}
for every $A \in \mathcal{A}$ with $\mu(A) > 0$.  Hence
\begin{equation}
\label{phi(frac{1}{mu(B)} int_B f_mathcal{B} d mu) le ...}
        \phi\Big(\frac{1}{\mu(B)} \int_B f_\mathcal{B} \, d\mu\Big)
            \le \frac{1}{\mu(B)} \int_B (\phi \circ f)_\mathcal{B} \, d\mu
\end{equation}
for every $B \in \mathcal{B}$ with $\mu(B) > 0$, because these
averages can be reduced to those in (\ref{phi(frac{1}{mu(A)} int_A f d
mu) le frac{1}{mu(A)} int_A phi circ f dmu}).  Using this, one can
check that
\begin{equation}
\label{phi(f_mathcal{B}) le (phi circ f)_mathcal{B}}
        \phi(f_\mathcal{B}) \le (\phi \circ f)_\mathcal{B}
\end{equation}
almost everywhere on $X$.  More precisely, one can apply the previous
inequality for averages to sets $B \in \mathcal{B}$ where
$f_\mathcal{B}$, $(\phi \circ f)_\mathcal{B}$ are approximately
constant.

        Of course, $\phi(t) = |t|^p$ is a convex function on the real
line when $1 \le p < \infty$.  If $f \in L^p(X, \mathcal{A})$ is
real-valued, then we get that
\begin{equation}
\label{|f_mathcal{B}|^p le (|f|^p)_mathcal{B} = E(|f|^p mid mathcal{B})}
        |f_\mathcal{B}|^p \le (|f|^p)_\mathcal{B} = E(|f|^p \mid \mathcal{B})
\end{equation}
almost everywhere on $X$, as in the previous paragraph.  If $f$ is
complex-valued, then one can apply this to $|f|$, to get that
\begin{equation}
        |f_\mathcal{B}|^p \le (|f|_\mathcal{B})^p \le (|f|^p)_\mathcal{B},
\end{equation}
using (\ref{|f_mathcal{B}| le |f|_mathcal{B}}) in the first step.
It follows that
\begin{equation}
 \int_X |f_\mathcal{B}|^p \, d\mu \le \int_X (|f|^p)_\mathcal{B} \, d\mu
                                    = \int_X |f|^p \, d\mu,
\end{equation}
because $X \in \mathcal{B}$, and that $f_\mathcal{B} \in L^p(X,
\mathcal{B})$ in particular.  Equivalently,
\begin{equation}
\label{||f_mathcal{B}||_p le ||f||_p}
        \|f_\mathcal{B}\|_p \le \|f\|_p,
\end{equation}
which also holds when $p = \infty$, as in (\ref{||f_mathcal{B}||_infty
le ||f||_infty}).

        Remember that ${\bf 1}_E(x)$ denotes the indicator function of
a set $E \subseteq X$, equal to $1$ when $x \in E$ and to $0$ when $x
\in X \backslash E$.  If $f \in L^1(X, \mathcal{A})$ and $A, E \in
\mathcal{A}$, then of course
\begin{equation}
        \int_A f \, {\bf 1}_E \, d\mu = \int_{A \cap E} f \, d\mu.
\end{equation}
If $B, E \in \mathcal{B}$, then $B \cap E \in \mathcal{B}$, and
\begin{eqnarray}
\int_B (f \, {\bf 1}_E)_\mathcal{B} \, d\mu & = & \int_B f \, {\bf 1}_E \, d\mu
                                              = \int_{B \cap E} f \, d\mu \\
           & = & \int_{B \cap E} f_\mathcal{B} \, d\mu
                     = \int_B f_\mathcal{B} \, {\bf 1}_E \, d\mu. \nonumber
\end{eqnarray}
This implies that
\begin{equation}
        (f \, {\bf 1}_E)_\mathcal{B} = f_\mathcal{B} \, {\bf 1}_E,
\end{equation}
since $f_\mathcal{B} \, {\bf 1}_E$ is measurable with respect to $\mathcal{B}$.
Similarly, if $g \in L^\infty(X, \mathcal{B})$, then
\begin{equation}
\label{(f g)_mathcal{B} = f_mathcal{B} g}
        (f \, g)_\mathcal{B} = f_\mathcal{B} \, g.
\end{equation}
This follows from the previous statement by approximating $g$ by
simple functions that are measurable with respect to $\mathcal{B}$.
If $f \in L^p(X, \mathcal{A})$, $1 \le p \le \infty$, then (\ref{(f
g)_mathcal{B} = f_mathcal{B} g}) also works for $g \in L^q(X,
\mathcal{B})$, where $1/p + 1/q = 1$, by the same argument.

        Note that $f_\mathcal{B} = 0$ almost everywhere on $X$ if and
only if
\begin{equation}
\label{int_B f d mu = 0}
        \int_B f \, d\mu = 0
\end{equation}
for every $B \in \mathcal{B}$.  If $f \in L^p(X, \mathcal{A})$, $1 \le
p \le \infty$, then this implies that
\begin{equation}
\label{int_X f g d mu = 0}
        \int_X f \, g \, d\mu = 0
\end{equation}
for every $g \in L^q(X, \mathcal{B})$, where $1/p + 1/q = 1$ again.
This uses the fact that simple functions are dense in $L^q(X,
\mathcal{B})$.  If $p = 2$, then the collection of $f \in L^2(X,
\mathcal{A})$ such that $f_\mathcal{B} = 0$ is the same as the
orthogonal complement of $L^2(X, \mathcal{B})$ as a linear subspace of
$L^2(X, \mathcal{A})$, and $f \mapsto f_\mathcal{B}$ is the same as the
orthogonal projection of $L^2(X, \mathcal{A})$ onto $L^2(X, \mathcal{B})$.

        Suppose now that $\mathcal{B}_1$, $\mathcal{B}_2$ are
$\sigma$-subalgebras of $\mathcal{A}$, with $\mathcal{B}_1 \subseteq
\mathcal{B}_2$.  If $f$ is an integrable function on $X$ with respect
to $\mathcal{A}$, then
\begin{equation}
        (f_{\mathcal{B}_2})_{\mathcal{B}_1} = f_{\mathcal{B}_1}.
\end{equation}
To see this, let $B \in \mathcal{B}_1$ be given, and observe that
\begin{equation}
        \int_B (f_{\mathcal{B}_2})_{\mathcal{B}_1} \, d\mu
              = \int_B f_{\mathcal{B}_2} \, d\mu = \int_B f \, d\mu
                                           = \int_B f_{\mathcal{B}_1} \, d\mu,
\end{equation}
because $B \in \mathcal{B}_2$ as well.  This corresponds to the fact
that restricting a measure $\nu$ on $\mathcal{A}$ to $\mathcal{B}_1$
is the same as restricting $\nu$ to $\mathcal{B}_2$, and then to
$\mathcal{B}_2$.

\section{Distances between measurable sets}
\label{distances between sets}
\setcounter{equation}{0}

        Remember that the symmetric difference $A \bigtriangleup B$ of
two sets $A$, $B$ is defined by
\begin{equation}
        A \bigtriangleup B = (A \backslash B) \cup (B \backslash A).
\end{equation}
If $C$ is another set, then it is easy to see that
\begin{equation}
\label{A bigtriangleup C subseteq (A bigtriangleup B) cup (B bigtriangleup C)}
 A \bigtriangleup C \subseteq (A \bigtriangleup B) \cup (B \bigtriangleup C).
\end{equation}
Let $(X, \mathcal{A}, \mu)$ be a probability space, and define $d(A,
B)$ for $A, B \in \mathcal{A}$ by
\begin{equation}
        d(A, B) = \mu(A \bigtriangleup B).
\end{equation}
Thus $d(A, A) = 0$, $d(A, B) = d(B, A) \ge 0$, and
\begin{equation}
\label{d(A, C) le d(A, B) + d(B, C)}
        d(A, C) \le d(A, B) + d(B, C)
\end{equation}
for every $A, B, C \in \mathcal{A}$, by (\ref{A bigtriangleup C
subseteq (A bigtriangleup B) cup (B bigtriangleup C)}).  This shows
that $d(A, B)$ is a semimetric on $\mathcal{A}$, which means that it
satisfies all of the requirements of a metric, except that $d(A, B) =
0$ may not imply that $A = B$.  In this case, $d(A, B) = 0$ when $A$
and $B$ are the same up to sets of measure $0$.  Equivalently, $d(A,
B)$ is equal to the distance between the indicator functions ${\bf
1}_A$, ${\bf 1}_B$ in $L^1$.

        Observe that $(X \backslash A) \bigtriangleup (X \backslash B)
= A \bigtriangleup B$ for every $A, B \subseteq X$, and hence
\begin{equation}
\label{d(X backslash A, X backslash B) = d(A, B)}
        d(X \backslash A, X \backslash B) = d(A, B)
\end{equation}
when $A, B \in \mathcal{A}$.  Moreover,
\begin{eqnarray}
            & & (A_1 \cup A_2) \bigtriangleup (B_1 \cup B_2) \\
          & = & ((A_1 \cup A_2) \backslash (B_1 \cup B_2)) \cup
                 ((B_1 \cup B_2) \backslash (A_1 \cup A_2)) \nonumber \\
& = & (A_1 \backslash (B_1 \cup B_2)) \cup (A_2 \backslash (B_1 \cup B_2)) \cup
(B_1 \backslash (A_1 \cup A_2)) \cup (B_2 \backslash (A_1 \cup A_2))\nonumber\\
 & \subseteq & (A_1 \backslash B_1) \cup (A_2 \backslash B_2) \cup
                (B_1 \backslash A_1) \cup (B_2 \backslash A_2) \nonumber \\
 & = & (A_1 \bigtriangleup B_1) \cup (A_2 \bigtriangleup B_2) \nonumber
\end{eqnarray}
for every $A_1, A_2, B_1, B_2 \subseteq X$.  Therefore
\begin{equation}
        d(A_1 \cup A_2, B_1 \cup B_2) \le d(A_1, B_1) + d(A_2, B_2)
\end{equation}
when $A_1, A_2, B_1, B_2 \in \mathcal{A}$.  This implies that
\begin{equation}
        d(A_1 \cap A_2, B_1 \cap B_2) \le d(A_1, B_1) + d(A_2, B_2)
\end{equation}
for every $A_1, A_2, B_1, B_2 \in \mathcal{A}$, because
\begin{equation}
 X \backslash (A_1 \cap A_2) = (X \backslash A_1) \cup (X \backslash A_2),
\end{equation}
and similarly for $X \backslash (B_1 \cap B_2)$.  This also uses
(\ref{d(X backslash A, X backslash B) = d(A, B)}) applied to $A_1 \cap
A_2$, $B_1 \cap B_2$ instead of $A$, $B$, and then to $A_1$, $B_1$ and
$A_2$, $B_2$.

        If $A_1 \subseteq A_2 \subseteq \cdots$ is an increasing
sequence of measurable subsets of $X$, then $\{A_j\}_{j = 1}^\infty$
converges to their union $\bigcup_{j = 1}^\infty A_j$ with respect to
$d(A, B)$, in the sense that
\begin{equation}
        \lim_{n \to \infty} d\Big(A_n, \bigcup_{j = 1}^\infty A_j\Big) = 0.
\end{equation}
To see this, note that $A_n \subseteq \bigcup_{j = 1}^\infty A_j$ for
each $n$, so that
\begin{equation}
        A_n \bigtriangleup \Big(\bigcup_{j = 1}^\infty A_j\Big)
         = \Big(\bigcup_{j = 1}^\infty A_j\Big) \backslash A_n
         = \bigcup_{j = n}^\infty (A_{j + 1} \backslash A_j).
\end{equation}
Hence
\begin{equation}
        d\Big(A_n, \bigcup_{j = 1}^\infty A_j\Big)
         = \sum_{j = n}^\infty \mu(A_{j + 1} \backslash A_j).
\end{equation}
Of course, the sets $A_{j + 1} \backslash A_j$ are pairwise disjoint,
and so $\sum_{j = 1}^\infty \mu(A_{j + 1} \backslash A_j)$ converges,
by countable additivity.  This implies that
\begin{equation}
 \lim_{n \to \infty} \sum_{j = n}^\infty \mu(A_{j + 1} \backslash A_j) = 0,
\end{equation}
as desired.  Similarly, if $B_1 \supseteq B_2 \supseteq \cdots$ is a
decreasing sequence of measurable sets, then $\{B_j\}_{j = 1}^\infty$
converges to $\bigcap_{j = 1}^\infty B_j$ with respect to $d(A, B)$,
in the sense that
\begin{equation}
        \lim_{n \to \infty} d\Big(B_n, \bigcap_{j = 1}^\infty B_j\Big) = 0.
\end{equation}
This follows from the previous case applied to $A_j = X \backslash B_j$.

        Let $\{A_j\}_{j = 1}^\infty$ be a sequence of subsets of $X$, and put
\begin{equation}
 B_k = \bigcup_{j = k}^\infty A_j, \quad C_l = \bigcap_{j = l}^\infty A_j
\end{equation}
for each $k, l \ge 1$.  Thus 
\begin{equation}
\label{B_{k + 1} subseteq B_k, C_l subseteq C_{l + 1}, and C_k subseteq B_k}
        B_{k + 1} \subseteq B_k, \, C_l \subseteq C_{l + 1},
                                \hbox{ and } C_k \subseteq B_k
\end{equation}
for each $k$, $l$.  The upper and lower limits of $\{A_j\}_{j =
1}^\infty$ are the subsets of $X$ defined by
\begin{equation}
\label{limsup_{j to infty} A_j, liminf_{j to infty} A_j}
        \limsup_{j \to \infty} A_j = \bigcap_{k = 1}^\infty B_k, \quad
          \liminf_{j \to \infty} A_j = \bigcup_{l = 1}^\infty C_l.
\end{equation}
In particular,
\begin{equation}
        \liminf_{j \to \infty} A_j \subseteq \limsup_{j \to \infty} A_j.
\end{equation}

        Suppose that $A_j \in \mathcal{A}$ for each $j$, so that $B_k,
C_l \in \mathcal{A}$ for every $k$, $l$, and hence
\begin{equation}
 \limsup_{j \to \infty} A_j, \, \liminf_{j \to \infty} A_j \in \mathcal{A}.
\end{equation}
Because of monotonicity,
\begin{equation}
 \lim_{k \to \infty} \mu(B_k) = \mu\Big(\limsup_{j \to \infty} A_j\Big), \quad
  \lim_{l \to \infty} \mu(C_l) = \mu\Big(\liminf_{j \to \infty} A_j\Big).
\end{equation}
It follows that
\begin{equation}
\label{mu(limsup_{j to infty} A_j) = mu(liminf_{j to infty} A_j)}
        \mu\Big(\limsup_{j \to \infty} A_j\Big) =
           \mu\Big(\liminf_{j \to \infty} A_j\Big)
\end{equation}
if and only if
\begin{equation}
\label{lim_{n to infty} mu(B_n backslash C_n) = 0}
        \lim_{n \to \infty} \mu(B_n \backslash C_n) = 0.
\end{equation}
If this condition holds and $A \in \mathcal{A}$ satisfies
\begin{equation}
 \liminf_{j \to \infty} A_j \subseteq A \subseteq \limsup_{j \to \infty} A_j,
\end{equation}
then it is easy to see that
\begin{equation}
        \lim_{n \to \infty} d(A_n, A) = 0.
\end{equation}
More precisely,
\begin{equation}
        A_n \bigtriangleup A = (A_n \backslash A) \cup (A \backslash A_n)
 \subseteq (B_n \backslash A) \cup (A \backslash C_n) = B_n \backslash C_n,
\end{equation}
and so
\begin{equation}
        d(A_n, A) \le \mu(B_n \backslash C_n) \to 0 \hbox{ as } n \to \infty.
\end{equation}

        Let us check that (\ref{lim_{n to infty} mu(B_n backslash C_n)
= 0}) holds when $\sum_{j = 1}^\infty d(A_{j + 1}, A_j)$ converges.
The main point is that
\begin{equation}
\label{B_n backslash A_n, A_n backslash C_n}
B_n \backslash A_n \subseteq \bigcup_{j = n}^\infty (A_{j + 1} \backslash A_j),
 \quad  
A_n \backslash C_n \subseteq \bigcup_{j = n}^\infty (A_j \backslash A_{j + 1})
\end{equation}
for each $n$.  More precisely, if $x \in B_n \backslash A_n$, then $x
\in A_{j + 1}$ for some $j \ge n + 1$, and $x \not\in A_n$.  If $j$ is
the smallest integer such that $j \ge n$ and $x \in A_{j + 1}$, then
$x \not\in A_j$, and so $x \in A_{j + 1} \backslash A_j$, as desired.
Similarly, if $y \in A_n \backslash C_n$, then $y \not\in A_{j + 1}$
for some $j \ge n$.  If $j$ is the smallest integer such that $j \ge
n$ and $y \not\in A_{j + 1}$, then $y \in A_j$, and so $y \in A_j
\backslash A_{j + 1}$.  This proves (\ref{B_n backslash A_n, A_n
backslash C_n}).

        It follows that
\begin{equation}
 \mu(B_n \backslash A_n) \le \sum_{j = n}^\infty \mu(A_{j + 1} \backslash A_j),
 \quad 
 \mu(A_n \backslash C_n) \le \sum_{j = n}^\infty \mu(A_j \backslash A_{j + 1})
\end{equation}
for each $n$.  Hence
\begin{equation}
 \mu(B_n \backslash C_n) = \mu(B_n \backslash A_n) + \mu(A_n \backslash C_n)
                          \le \sum_{j = n}^\infty d(A_{j + 1}, A_j),
\end{equation}
using the fact that $C_n \subseteq A_n \subseteq B_n$ in the first
step.  If $\sum_{j = 1}^\infty d(A_{j + 1}, A_j)$ converges, then the
right side tends to $0$ as $n \to \infty$, and so (\ref{lim_{n to
infty} mu(B_n backslash C_n) = 0}) holds.  This implies that there is
an $A \in \mathcal{A}$ such that $\lim_{n \to \infty} d(A_n, A) = 0$,
by the earlier remarks.  If instead $\{A_j\}_{j = 1}^\infty$ satisfies
the Cauchy condition
\begin{equation}
        \lim_{j, l \to \infty} d(A_j, A_l) = 0,
\end{equation}
then there is a subsequence $\{A_{j_n}\}_{n = 1}^\infty$ of
$\{A_j\}_{j = 1}^\infty$ such that $\sum_{n = 1}^\infty d(A_{j_{n +
1}}, A_{j_n})$ converges.  This implies that there is an $A \in
\mathcal{A}$ such that $\lim_{n \to \infty} d(A_{j_n}, A) = 0$, as
before.  Using the Cauchy condition, one can check that $\lim_{j \to
\infty} d(A_j, A) = 0$.

        If $\mathcal{E} \subseteq \mathcal{A}$, then let
$\overline{\mathcal{E}}$ be the collection of $A \in \mathcal{A}$ such
that for each $\epsilon > 0$ there is an $E \in \mathcal{A}$ that
satisfies $d(A, E) < \epsilon$.  This is basically the same as the
closure of a set in a metric space, except that $d(A, B)$ is only a
semimetric.  In particular, note that $\overline{\mathcal{E}}$
automatically contains every $A \in \mathcal{A}$ for which there is an
$E \in \mathcal{E}$ such that $d(A, E) = 0$.  As in the context of
metric spaces, one can check that
\begin{equation}
        \overline{\overline{\mathcal{E}}} = \mathcal{E}.
\end{equation}
If $\mathcal{E}$ is a subalgebra of $\mathcal{A}$, then it is easy to
see that $\overline{\mathcal{E}}$ is also a subalgebra of
$\mathcal{A}$, using the properties of the distance related to unions,
intersections, and complements discussed earlier in this section.

        Let us check that $\overline{\mathcal{E}}$ is actually a
$\sigma$-algebra when $\mathcal{E}$ is an algebra.  It suffices to
show that $\bigcup_{j = 1}^\infty A_j \in \overline{\mathcal{E}}$ for
every sequence $A_1, A_2, \ldots$ of elements of
$\overline{\mathcal{E}}$.  Of course, $\bigcup_{j = 1}^n A_j \in
\overline{\mathcal{E}}$ for each $n$, because $\overline{\mathcal{E}}$
is an algebra.  We also know that $\bigcup_{j = 1}^n A_j$ converges to
$\bigcup_{j = 1}^\infty A_j$ as $n \to \infty$ with respect to $d(A,
B)$, because of monotonicity.  It follows that $\bigcup_{j = 1}^\infty
A_j \in \overline{\mathcal{E}}$, by combining these two facts.

        If $A \in \overline{\mathcal{E}}$, then there is a sequence
$\{A_j\}_{j = 1}^\infty$ of elements of $\mathcal{E}$ such that
$\sum_{j = 1}^\infty d(A_j, A)$ converges.  This implies that $\sum_{j
= 1}^\infty d(A_{j + 1}, A_j)$ converges, by the triangle inequality.
Thus $\{A_j\}_{j = 1}^\infty$ converges to $\limsup_{j \to \infty}
A_j$, $\liminf_{j \to \infty} A_j$ with respect to $d(A, B)$, by the
earlier discussion, and $A$ differs from these limits by sets of
measure $0$.  In particular, $A \in \mathcal{E}$ when $\mathcal{E}$ is
a $\sigma$-subalgebra of $\mathcal{A}$ that contains all elements of
$\mathcal{A}$ with measure $0$.  It follows that
$\overline{\mathcal{E}} = \mathcal{E}$ when $\mathcal{E}$ is a
$\sigma$-subalgebra of $\mathcal{A}$ that contains the sets of measure
$0$.

\section{Sequences of $\sigma$-subalgebras}
\label{sequences of sigma-subalgebras}
\setcounter{equation}{0}

        Let $(X, \mathcal{A}, \mu)$ be a probability space, and let
$\mathcal{B}_1 \subseteq \mathcal{B}_2 \subseteq \cdots$ be an
increasing sequence of $\sigma$-subalgebras of $\mathcal{A}$.  Thus
$\mathcal{E} = \bigcup_{j = 1}^\infty \mathcal{B}_j$ is a subalgebra
of $\mathcal{A}$, but not necessarily a $\sigma$-subalgebra.  If
$\mathcal{C} = \overline{\mathcal{E}}$ is the closure of $\mathcal{E}$
with respect to the semimetric $d(A, B)$, then $\mathcal{C}$ is the
smallest $\sigma$-subalgebra of $\mathcal{A}$ that contains
$\mathcal{E}$ and the sets of measure $0$, as in the previous section.

        Put
\begin{equation}
        f_n = f_{\mathcal{B}_n} = E(f \mid \mathcal{B}_n)
\end{equation}
for each $f \in L^1(X, \mathcal{A})$ and $n \ge 1$, and
\begin{equation}
        f_\infty = f_\mathcal{C} = E(f \mid \mathcal{C}).
\end{equation}
Note that
\begin{equation}
        f_n = E(f_\infty \mid \mathcal{B}_n)
\end{equation}
for each $n$, since $\mathcal{B}_n \subseteq \mathcal{C}$.  If $f \in
L^p(X, \mathcal{A})$ for some $p$, $1 \le p \le \infty$, then $f_n \in
L^p(X, \mathcal{B}_n)$ for each $n$, $f_\infty \in L^p(X, \mathcal{C})$, and
\begin{equation}
\label{||f_n||_p le ||f_infty||_p le ||f||_p}
        \|f_n\|_p \le \|f_\infty\|_p \le \|f\|_p.
\end{equation}
If $f$ happens to be measurable with respect to $\mathcal{B}_l$ for
some $l \ge 1$, then
\begin{equation}
\label{f_n = f_infty = f}
        f_n = f_\infty = f
\end{equation}
for every $n \ge l$.

        If $1 \le p < \infty$, then
\begin{equation}
        \bigcup_{l = 1}^\infty L^p(X, \mathcal{B}_l)
\end{equation}
is dense in $L^p(X, \mathcal{C})$.  To see this, one can first
approximate elements of $L^p(X, \mathcal{C})$ by simple functions that
are measurable with respect to $\mathcal{C}$.  The latter can then be
approximated by simple functions that are measurable with respect to
$\mathcal{B}_l$ for some $l$, using the definition of $\mathcal{C}$.
This implies that
\begin{equation}
        \lim_{n \to \infty} f_n = f_\infty
\end{equation}
in the $L^p$ norm when $f \in L^p(X, \mathcal{A})$, $1 \le p <
\infty$.  More precisely, one may as well take $f = f_\infty$, so that
$f$ is already measurable with respect to $\mathcal{C}$.  If $f$ is
measurable with respect to $\mathcal{B}_l$ for some $l$, then one can
apply (\ref{f_n = f_infty = f}).  Otherwise, one can approximate $f$
by $g \in L^p(X, \mathcal{B}_l)$ for some $l$, by previous remarks
about density in $L^p(X, \mathcal{C})$.  The main point is that $f_n$
is also approximated by $g$ when $n \ge l$, uniformly in $n$, because
of (\ref{||f_n||_p le ||f_infty||_p le ||f||_p}).

        Suppose that $X_1, X_2, \ldots$ is a sequence of compact
Hausdorff spaces, and that $X = \prod_{j = 1}^\infty X_j$ is their
Cartesian product, with the product topology.  Let $\mu_j$ be a
regular Borel probability measure on $X_j$ for each $j$, and let $\mu$
be the corresponding product measure on $X$.  Also let $\mathcal{B}_n$
be the collection of subsets of $X$ of the form $B \times \prod_{j = n
+ 1}^\infty X_j$, where $B$ is a Borel set in $\prod_{j = 1}^n X_j$.
If $f$ is a continuous real or complex-valued function on $X$, then
$f_n$ is the function of $x_1, \ldots, x_n$ obtained by integrating
$f$ in the variables $x_j$ for $j \ge n + 1$.  In this case,
$\{f_n\}_{n = 1}^\infty$ converges to $f$ uniformly on $X$, because of
the uniform continuity properties discussed in Section \ref{product
spaces, 2}.

\section{Martingales}
\label{martingales}
\setcounter{equation}{0}

        Let $(X, \mathcal{A}, \mu)$ be a probability space, and let
$\mathcal{B}_1 \subseteq \mathcal{B}_2 \subseteq \cdots$ be an
increasing sequence of $\sigma$-subalgebras of $\mathcal{A}$, also
known as a \emph{filtration}.  A sequence $\{f_j\}_{j = 1}^\infty$ of
functions on $X$ is said to be a \emph{martingale} with respect to
this filtration if $f_j \in L^1(X, \mathcal{B}_j)$ for each $j$, and
\begin{equation}
        f_j = E(f_l \mid \mathcal{B}_j)
\end{equation}
when $1 \le j \le l$.  In particular, this implies that
\begin{equation}
        \|f_j\|_1 \le \|f_l\|_1
\end{equation}
for each $j \le l$.  If $f_j \in L^p(X, \mathcal{B}_j)$ for some $p$,
$1 \le p \le \infty$, and every $j$, then
\begin{equation}
        \|f_j\|_p \le \|f_l\|_p
\end{equation}
for each $j \le l$.  If $f \in L^1(X, \mathcal{A})$ and $f_j = E(f
\mid \mathcal{B}_j)$ for each $j$, then $\{f_j\}_{j = 1}^\infty$ is a
martingale.

        Let $(X_1, \mathcal{A}_1, \mu_1), (X_2, \mathcal{A}_2, \mu_2),
\ldots$ be a sequence of probability spaces, and let $X = \prod_{j =
1}^\infty X_j$ be their Cartesian product, with the product measure
$\mu$ on the corresponding $\sigma$-algebra $\mathcal{A}$.  Also let
$\mathcal{B}_n$ be the collection of subsets of $X$ of the form $B
\times \prod_{j = n + 1}^\infty X_j$, where $B$ is a measurable subset
of $\prod_{j = 1}^n X_j$.  This defines an increasing sequence of
$\sigma$-subalgebras of $\mathcal{A}$.  Let $a_j$ be an integrable
function on $X_j$ such that
\begin{equation}
        \int_{X_j} a_j \, d\mu_j = 0
\end{equation}
for each $j$, which can also be considered as an integrable function
on $X$ that does not depend on $x_l$ when $j \ne l$.  In this case,
\begin{equation}
        f_n = \sum_{j = 1}^n a_j
\end{equation}
defines a martingale with respect to this filtration.

        Let $(X, \mathcal{A}, \mu)$ be any probability space again,
with an increasing sequence $\mathcal{B}_j$ of $\sigma$-algebras of
$\mathcal{A}$.  Also let $\{f_j\}_{j = 1}^\infty$ be a martingale with
respect to this filtration, with $f_j \in L^2(X, \mathcal{B}_j)$ for
each $j$.  Thus
\begin{equation}
        \int_B f_j \, d\mu = \int_B f_{j + 1} \, d\mu
\end{equation}
for each $B \in \mathcal{B}_j$, which implies that
\begin{equation}
        \int_X b \, f_j \, d\mu = \int_X b \, f_{j + 1} \, d\mu
\end{equation}
for every $b \in L^2(X, \mathcal{B}_j)$.  Equivalently,
\begin{equation}
        \int_X b \, (f_j - f_{j + 1}) \, d\mu = 0
\end{equation}
for every $b \in L^2(X, \mathcal{B}_j)$.  It follows that the
functions $f_1$ and $f_{j + 1} - f_j$, $j \ge 1$, are all orthogonal
to each other in $L^2(X, \mathcal{A})$.

\section{$L^p$ Boundedness}
\label{L^p boundedness}
\setcounter{equation}{0}

        Let $(X, \mathcal{A}, \mu)$ be a probability space, and let
$\mathcal{B}_1 \subseteq \mathcal{B}_2 \subseteq \cdots$ be an
increasing sequence of $\sigma$-subalgebras of $\mathcal{A}$.  As
before, put $\mathcal{E} = \bigcup_{j = 1}^\infty \mathcal{B}_j$, and
let $\mathcal{C} = \overline{\mathcal{E}}$ be the closure of
$\mathcal{E}$ with respect to the semimetric $d(A, B)$.  Let $1 < p
\le \infty$ be given, and let $\{f_j\}_{j = 1}^\infty$ be a martingale
on $X$ with respect to the $\mathcal{B}_j$'s such that $f_j \in L^p(X,
\mathcal{B}_j)$ for each $j$, and the $L^p$ norms $\|f_j\|_p$ are
uniformly bounded.

        If $B \in \mathcal{B}_l$ for some $l$, then
\begin{equation}
        \int_B f_l \, d\mu = \int_B f_n \, d\mu
\end{equation}
when $n \ge l$.  This implies that
\begin{equation}
        \int_X f_l \, g \, d\mu = \int_X f_n \, g \, d\mu
\end{equation}
when $g \in L^q(X, \mathcal{B}_l)$, where $1/p + 1/q = 1$.  In
particular,
\begin{equation}
\label{lim_{n to infty} int_X f_n g d mu}
        \lim_{n \to \infty} \int_X f_n \, g \, d\mu
\end{equation}
exists for every $g \in L^q(X, \mathcal{B}_l)$, $l \ge 1$.  Note that
$\bigcup_{l = 1}^\infty L^q(X, \mathcal{B}_l)$ is dense in $L^q(X,
\mathcal{C})$, as in Section \ref{sequences of sigma-subalgebras},
because $1 \le q < \infty$.  It follows that the limit (\ref{lim_{n to
infty} int_X f_n g d mu}) exists for every $g \in L^q(X,
\mathcal{C})$, using also the uniform boundedness of the $L^p$ norms
of the $f_j$'s, as in Section \ref{uniform boundedness, 4}.

        More precisely, (\ref{lim_{n to infty} int_X f_n g d mu})
defines a bounded linear functional on $L^q(X, \mathcal{C})$ under
these conditions.  The Riesz representation theorem implies that there
is an $f \in L^p(X, \mathcal{C})$ such that
\begin{equation}
\label{lim_{n to infty} int_X f_n g d mu = int_X f g d mu}
        \lim_{n \to \infty} \int_X f_n \, g \, d\mu = \int_X f \, g \, d\mu
\end{equation}
for every $g \in L^q(X, \mathcal{C})$ under these conditions.  If $g
\in L^q(X, \mathcal{B}_l)$ for some $l$, then we get that
\begin{equation}
        \int_X f_l \, g \, d\mu = \int_X f \, g \, d\mu.
\end{equation}
In particular,
\begin{equation}
        \int_B f_l \, d\mu = \int_X f \, d\mu
\end{equation}
for each $B \in \mathcal{B}_l$, which implies that
\begin{equation}
        f_l = E(f \mid \mathcal{B}_l)
\end{equation}
for each $l$.

        If $1 < p < \infty$, then it follows that $\{f_l\}_{l =
1}^\infty$ converges to $f$ in the $L^p$ norm, as in Section
\ref{sequences of sigma-subalgebras}.  If $p = 2$, then
\begin{equation}
 \|f_n\|_2^2 = \|f_1\|_2^2 + \sum_{j = 1}^{n - 1} \|f_{j + 1} - f_j\|_2^2
\end{equation}
for each $n$, because of orthogonality, as in the previous section.
The boundedness of the $L^2$ norms $\|f_n\|_2$ is equivalent to the
convergence of the series
\begin{equation}
        \sum_{j = 1}^\infty \|f_{j + 1} - f_j\|_2^2,
\end{equation}
which implies the convergence of the series $\sum_{j = 1}^\infty (f_{j
+ 1} - f_j)$ in $L^2(X, \mathcal{C})$.  This gives a more direct proof
of the convergence of $\{f_j\}_{j = 1}^\infty$ in $L^2(X,
\mathcal{C})$ in this case.  Of course, if $\{f_j\}_{j = 1}^\infty$ is
a martingale such that $f_j \in L^p(X, \mathcal{B}_l)$ converges to $f
\in L^p(X, \mathcal{A})$ in the $L^p$ norm for any $p$, $1 \le p \le
\infty$, then $f \in L^p(X, \mathcal{C})$ and $f_l = E(f \mid
\mathcal{B}_l)$ for each $l$, for basically the same reasons as before.

\section{Uniform integrability}
\label{uniform integrability}
\setcounter{equation}{0}

        Let $(X, \mathcal{A}, \mu)$ be a probability space, and let
$\mathcal{B}_1 \subseteq \mathcal{B}_2 \subseteq \cdots$ be an
increasing sequence of $\sigma$-subalgebras of $\mathcal{A}$.  Also
let $\{f_j\}_{j = 1}^\infty$ be a martingale with respect to this
filtration with bounded $L^1$ norms, so that there is a $C \ge 0$ such
that
\begin{equation}
\label{||f_n||_1 le C}
        \|f_n\|_1 \le C
\end{equation}
for each $n$.  Note that this holds automatically when $f_j \ge 0$ for
each $j$, because
\begin{equation}
\label{||f_j||_1 = int_X f_j d mu = int_X f_1 d mu}
        \|f_j\|_1 = \int_X f_j \, d\mu = \int_X f_1 \, d\mu
\end{equation}
for each $j \ge 1$ in this case.

        Suppose that the $f_j$'s are \emph{uniformly integrable} as well,
in the sense that for each $\epsilon > 0$ there is a $\delta > 0$ such that
\begin{equation}
\label{int_A |f_n| d mu < epsilon}
        \int_A |f_n| \, d\mu < \epsilon
\end{equation}
for every $A \in \mathcal{A}$ with $\mu(A) < \delta$ and every $n \ge 1$.
It is well known that this condition holds automatically for a single
integrable function, by approximating that function by bounded
functions in the $L^1$ norm, for instance.  Similarly, any finite collection
of integrable functions has this property.  Using this, it is easy to
check that a sequence of integrable functions that converges in the
$L^1$ norm is uniformly integrable.  If there is a $p > 1$ such that
$f_n \in L^p$ for each $n$ and $\|f_n\|_p$ is uniformly bounded, then
$\{f_n\}_{n = 1}^\infty$ is uniformly integrable, because of H\"older's
inequality.

        If $\{f_n\}_{n = 1}^\infty$ satisfies (\ref{||f_n||_1 le C}), then
\begin{equation}
        \mu(\{x \in X : |f_n(x)| > t\}) \le t^{-1} \, C
\end{equation}
for each $t > 0$, by Tchebychev's inequality.  If $\{f_j\}_{j =
1}^\infty$ is uniformly integrable too, then it follows that
\begin{equation}
        \int_{\{x \in X : |f_n(x)| > t\}} |f_n(x)| \, d\mu(x) \to 0
                                      \quad\hbox{as}\quad t \to \infty,
\end{equation}
uniformly in $n$.  Conversely, the latter condition implies that
$\{f_n\}_{n = 1}^\infty$ has bounded $L^1$ norms and is uniformly integrable.

        As usual, put $\mathcal{E} = \bigcup_{j = 1}^\infty
\mathcal{B}_j$, and let $\mathcal{C} = \overline{\mathcal{E}}$ be the
closure of $\mathcal{E}$ with respect to the semimetric $d(A, B)$.
Note that
\begin{equation}
        \int_A f_n \, d\mu = \int_A f_l \, d\mu
\end{equation}
for every $A \in \mathcal{B}_l$ and $n \ge l$.  We would like to show that
\begin{equation}
\label{{int_A f_n d mu}_{n = 1}^infty}
        \bigg\{\int_A f_n \, d\mu\bigg\}_{n = 1}^\infty
\end{equation}
is a Cauchy sequence in ${\bf R}$ or ${\bf C}$, as appropriate, for
every $A \in \mathcal{C}$, and hence converges.  This is obvious when
$A \in \mathcal{E}$, and one can deal with $A \in \mathcal{C}$ by
approximation, using uniform integrability.  The main point is that
\begin{equation}
        \int_A f_n \, d\mu,  \, n \in {\bf Z}_+,
\end{equation}
is an equicontinuous family of functions of $A \in \mathcal{A}$ with
respect to the semimetric $d(A, B)$, since
\begin{equation}
        \biggl|\int_A f_n \, d\mu - \int_B f_n \, d\mu\biggr|
         \le \int_{A \bigtriangleup B} |f_n| \, d\mu
\end{equation}
for every $A, B \in \mathcal{A}$.

        Put
\begin{equation}
        \nu(A) = \lim_{n \to \infty} \int_A f_n \, d\mu
\end{equation}
for each $A \in \mathcal{C}$.  Uniform integrability implies that for
each $\epsilon > 0$ there is a $\delta > 0$ such that
\begin{equation}
        |\nu(A)| \le \epsilon
\end{equation}
for every $A \in \mathcal{C}$ such that $\mu(A) < \delta$.  This
follows by taking the limit as $n \to \infty$ in the definition of
uniform integrability of $\{f_n\}_{n = 1}^\infty$, using the same
$\delta$ as before.

        Clearly $\nu(A)$ is finitely additive on $\mathcal{C}$, and
countable additivity follows from this continuity condition.  For if
$A_1, A_2, \ldots$ is a sequence of pairwise-disjoint subsets of $X$
in $\mathcal{C}$, then countable additivity of $\mu$ implies that
\begin{equation}
        \lim_{k \to \infty} \mu\Big(\bigcup_{j = k + 1}^\infty A_j\Big) = 0,
\end{equation}
and hence
\begin{equation}
        \lim_{k \to \infty} \nu\Big(\bigcup_{j = k + 1}^\infty A_j\Big) = 0
\end{equation}
too, by the continuity condition.  Because of finite additivity,
we also have that
\begin{equation}
        \nu\Big(\bigcup_{j = 1}^\infty A_j\Big) = \sum_{j = 1}^k \nu(A_j) +
                                   \nu\Big(\bigcup_{j = k + 1}^\infty A_j\Big)
\end{equation}
for each $k \ge 1$.  It follows that $\sum_{j = 1}^\infty \nu(A_j)$
converges to $\nu\Big(\bigcup_{j = 1}^\infty A_j\Big)$, as desired.

        Thus $\nu$ is a countably-additive real or complex measure on
$\mathcal{C}$, as appropriate.  Moreover, $\nu$ is absolutely
continuous with respect to the restriction of $\mu$ to $\mathcal{C}$.
The Radon--Nikodym theorem implies that there is an $f \in L^1(X,
\mathcal{C})$ such that
\begin{equation}
\label{nu(A) = int_A f d mu}
        \nu(A) = \int_A f \, d\mu
\end{equation}
for every $A \in \mathcal{C}$.  In particular,
\begin{equation}
        \int_A f_l \, d\mu = \int_A f \, d\mu
\end{equation}
when $A \in \mathcal{B}_l$, which implies that
\begin{equation}
        f_l = E(f \mid \mathcal{B}_l)
\end{equation}
for each $l$.  Conversely, this implies that $\{f_l\}_{l = 1}^\infty$
converges to $f$ in the $L^1$ norm, as in Section \ref{sequences of
sigma-subalgebras}, which implies that $\{f_l\}_{l = 1}^\infty$ is
uniformly integrable.

\section{Maximal functions, 3}
\label{maximal functions, 3}
\setcounter{equation}{0}

        Let $(X, \mathcal{A}, \mu)$ be a probability space, let
$\mathcal{B}_1 \subseteq \mathcal{B}_2 \subseteq \cdots$ be an
increasing sequence of $\sigma$-subalgebras of $\mathcal{A}$, and let
$\{f_j\}_{j = 1}^\infty$ be a martingale on $X$ with respect to this
filtration.  Consider the maximal functions
\begin{equation}
        f_n^*(x) = \max_{1 \le j \le n} |f_j(x)|
\end{equation}
and
\begin{equation}
        f^*(x) = \sup_{j \ge 1} |f_j(x)|.
\end{equation}
Note that $f_n^*$ is measurable with respect to $\mathcal{B}_n$, and that
\begin{equation}
        f^*(x) = \lim_{n \to \infty} f_n^*(x)
\end{equation}
is measurable with respect to the smallest $\sigma$-algebra
$\mathcal{B}_\infty$ that contains $\mathcal{E} = \bigcup_{j =
1}^\infty \mathcal{B}_j$.  If $X = [0, 1)$, $\mu$ is Lebesgue measure,
and $\mathcal{B}_j$ consists of unions of dyadic intervals of length
$2^{-j}$, then this is a variant of the dyadic maximal function, as in
Section \ref{dyadic intervals}.

        Put
\begin{equation}
        E(t) = \{x \in X : f^*(x) > t\}
\end{equation}
for each $t > 0$, as well as
\begin{equation}
        E_1(t) = \{x \in X : |f_1(x)| > t\}
\end{equation}
and
\begin{equation}
        E_l(t) = \{x \in X : |f_l(x)| > t, \, f_{l - 1}^*(x) \le t\}
\end{equation}
when $l \ge 2$.  Thus $E_l(t) \in \mathcal{B}_l$ for each $l$, $t$,
$E_l(t) \cap E_n(t) = \emptyset$ when $l < n$, and
\begin{equation}
        E(t) = \bigcup_{l = 1}^\infty E_l(t).
\end{equation}
Similarly,
\begin{equation}
        \bigcup_{l = 1}^n E_l(t) = \{x \in X : f_n^*(x) > t\}
\end{equation}
for each $n \ge 1$.  If $l \le n$, then
\begin{equation}
        t \, \mu(E_l(t)) \le \int_{E_l(t)} |f_l| \, d\mu
                          \le \int_{E_l(t)} |f_n| \, d\mu,
\end{equation}
because $f_l = E(f_n \mid \mathcal{B}_l)$ and hence $|f_l| \le E(|f_n|
\mid \mathcal{B}_l)$, as in Section \ref{basic properties}.  This
implies that
\begin{eqnarray}
 t \, \mu\Big(\bigcup_{l = 1}^n E_l(t)\Big) = \sum_{l = 1}^n t \, \mu(E_l(t))
          & \le & \sum_{l = 1}^n \int_{E_l(t)} |f_n| \, d\mu \\
           & = & \int_{\bigcup_{l = 1}^n E_l(t)} |f_n| \, d\mu. \nonumber
\end{eqnarray}

        Suppose now that the $f_n$'s have bounded $L^1$ norms, so that
\begin{equation}
        \|f_n\|_1 \le C
\end{equation}
for some $C \ge 0$ and every $n \ge 1$.  The previous estimate implies that
\begin{equation}
        t \, \mu\Big(\bigcup_{l = 1}^n E_l(t)\Big) \le C
\end{equation}
for each $n$.  Hence
\begin{equation}
\label{t mu(E(t)) = t mu(bigcup_{l = 1}^infty E_l(t)) le C}
 t \, \mu(E(t)) = t \, \mu\Big(\bigcup_{l = 1}^\infty E_l(t)\Big) \le C.
\end{equation}
This is basically the same as the estimates in Sections \ref{maximal
functions} and \ref{dyadic intervals}, except that the measure $\mu$
here corresponds to Lebesgue measure before, and the martingale
$\{f_j\}_{j = 1}^\infty$ corresponds to the measure $\mu$ or function
$f$ before.  The martingale may be generated by a function or measure
on $X$, through conditional expectation.

        We also have that
\begin{eqnarray}
        \sum_{l = 1}^n \int_{E_l(t)} |f_l| \, d\mu
         & \le & \sum_{l = 1}^n \int_{E_l(t)} |f_n| \, d\mu \\
           & = & \int_{\bigcup_{l = 1}^n E_l(t)} |f_n| \, d\mu \le C \nonumber
\end{eqnarray}
for each $n$, since $|f_l| \le E(|f_n| \mid \mathcal{B}_l)$ when $l \le n$.
Hence
\begin{equation}
        \sum_{l = 1}^\infty \int_{E_l(t)} |f_l| \, d\mu \le C.
\end{equation}
This shows that the function $h$ defined on $X$ by $h = f_l$ on $E_l(t)$,
$h = 0$ on $X \backslash E(t)$, is integrable, with $\|h\|_1 \le C$.

\section{Convergence almost everywhere}
\label{convergence almost everywhere}
\setcounter{equation}{0}

        Let $(X, \mathcal{A}, \mu)$ be a probability space, let
$\mathcal{B}_1 \subseteq \mathcal{B}_2 \subseteq \cdots$ be an
increasing sequence of $\sigma$-subalgebras of $\mathcal{A}$, and let
$\{f_j\}_{j = 1}^\infty$ be a martingale on $X$ with respect to this
filtration.  Observe that
\begin{equation}
\label{{f_n - f_l}_{n = l}^infty}
        \{f_n - f_l\}_{n = l}^\infty
\end{equation}
is a martingale with respect to the filtration $\mathcal{B}_l
\subseteq \mathcal{B}_{l + 1} \subseteq \cdots$ for each $l \ge 1$.
Put
\begin{equation}
        A_l(t) = \bigg\{x \in X : \sup_{n \ge l} |f_n(x) - f_l(x)| > t\bigg\}
\end{equation}
for every $l \ge 1$ and $t > 0$.

        If $\|f_n\|_1$ is bounded, then
\begin{equation}
        t \, \mu(A_l(t)) \le \sup_{n \ge l} \|f_n - f_l\|_1
\end{equation}
for every $t > 0$, as in the previous section.  This implies that
\begin{equation}
        t \, \mu\Big(\bigcap_{l = 1}^\infty A_l(t)\Big)
         \le \inf_{l \ge 1} \Big(\sup_{n \ge l} \|f_n - f_l\|_1\Big),
\end{equation}
for each $t > 0$, and hence
\begin{equation}
        \mu\Big(\bigcap_{l = 1}^\infty A_l(t)\Big) = 0
\end{equation}
for every $t > 0$ when $\{f_n\}_{n = 1}^\infty$ is a Cauchy sequence
in $L^1(X, \mathcal{A})$.  Thus
\begin{equation}
 \mu\Big(\bigcup_{k = 1}^\infty \bigcap_{l = 1}^\infty A_l(1/k)\Big) = 0.
\end{equation}
Of course,
\begin{equation}
 X \backslash \Big(\bigcup_{k = 1}^\infty \bigcap_{l = 1}^\infty A_l(1/k)\Big)
  = \bigcap_{k = 1}^\infty \bigcup_{l = 1}^\infty (X \backslash A_l(1/k)).
\end{equation}

        If $x$ is in this set, then it is easy to see that
$\{f_n(x)\}_{n = 1}^\infty$ is a Cauchy sequence in ${\bf R}$ or ${\bf
C}$, as appropriate.  It follows that $\{f_n\}_{n = 1}^\infty$
converges pointwise almost everywhere on $X$ when it converges in the
$L^1$ norm.  As in Section \ref{sequences of sigma-subalgebras}, this
happens when there is an $f \in L^1(X, \mathcal{A})$ such that $f_n =
E(f \mid \mathcal{B}_n)$ for each $n$.  In particular, this happens
when $\{f_n\}_{n = 1}^\infty$ is uniformly integrable, as in Section
\ref{uniform integrability}.  This includes the case where there is a
$p > 1$ such that $f_n \in L^p(X, \mathcal{B}_n)$ for each $n$ and
$\|f_n\|_p$ is bounded, as in Section \ref{L^p boundedness}.

        Suppose that we simply know that $\|f_n\|_1$ is uniformly
bounded in $n$.  Let $t > 0$ be given, and put $g_1 = f_1$, and
\begin{eqnarray}
        g_n(x) & = & f_n(x) \quad\hbox{when }
                x \in X \backslash \Big(\bigcup_{l = 1}^{n - 1} E_l(t)\Big) \\
 & = & f_l(x) \quad\hbox{ when } x \in E_l(t), \, 1 \le l \le n - 1 \nonumber
\end{eqnarray}
for $n \ge 2$, where $E_l(t)$ is as in the previous section.  Note
that $g_n$ is measurable with respect to $\mathcal{B}_n$ for each $n$,
because $E_l(t) \in \mathcal{B}_l \subseteq \mathcal{B}_n$ when $l \le
n$, as in the previous section, and $f_l$ is measure with respect to
$\mathcal{B}_l$ and hence $\mathcal{B}_n$ when $l \le n$.  Moreover,
\begin{eqnarray}
 \int_X |g_n| \, d\mu & = &
 \int_{X \backslash \Big(\bigcup_{l = 1}^n E_l(t)\Big)} |f_n| \, d\mu
            + \sum_{l = 1}^{n - 1} \int_{E_l(t)} |f_l| \, d\mu \\
 & \le & \int_{X \backslash \Big(\bigcup_{l = 1}^{n - 1} E_l(t)\Big)} |f_n| \,
            d\mu + \sum_{l = 1}^{n - 1} \int_{E_l(t)} |f_n| \, d\mu \nonumber
\end{eqnarray}
when $n \ge 2$, using the fact that $|f_l| \le E(|f_n| \mid
\mathcal{B}_l)$ in the second step.  This implies that
\begin{equation}
        \int_X |g_n| \, d\mu \le \int_X |f_n| \, d\mu,
\end{equation}
which obviously holds when $n = 1$ as well.

        Let us check that $\{g_n\}_{n = 1}^\infty$ is a martingale on
$X$ with respect to the $\mathcal{B}_n$'s.  It suffices to show that
\begin{equation}
\label{int_A g_n d mu = int_A g_{n + 1} d mu}
        \int_A g_n \, d\mu = \int_A g_{n + 1} \, d\mu
\end{equation}
for each $A \in \mathcal{B}_n$ and $n \ge 1$, so that $g_n = E(g_{n +
1} \mid \mathcal{B}_n)$.  If $A \subseteq X \backslash \Big(\bigcup_{l
= 1}^n E_l(t)\Big)$, then $g_n = f_n$ and $g_{n + 1} = f_{n + 1}$ on
$A$, and so
\begin{equation}
        \int_A g_n \, d\mu = \int_A f_n \, d\mu = \int_A f_{n + 1} \, d\mu
                                                = \int_A g_{n + 1} \, d\mu.
\end{equation}
This uses the facts that $f_n = E(f_{n + 1} \mid \mathcal{B}_n)$ and
$A \in \mathcal{B}_n$ in the middle step.  If $A \subseteq E_n(t)$,
then $g_n = f_n$ on $A$ because $A \subseteq X \backslash
\Big(\bigcup_{l = 1}^{n - 1} E_l(t)\Big)$, and $g_{n + 1} = f_n$ on
$A$ by definition of $g_{n + 1}$.  Hence
\begin{equation}
        \int_A g_n \, d\mu = \int_A f_n \, d\mu = \int_A g_{n + 1} \, d\mu.
\end{equation}
Similarly, if $A \subseteq E_l(t)$ for some $l = 1, \ldots, n - 1$,
then $g_n = g_{n + 1} = f_l$ on $A$, and so
\begin{equation}
        \int_A g_n \, d\mu = \int_A f_l \, d\mu = \int_A g_{n + 1} \, d\mu.
\end{equation}
Every $A \in \mathcal{B}_n$ can be expressed as the disjoint union of
its intersections with $X \backslash \Big(\bigcup_{l = 1}^n
E_l(t)\Big)$ and $E_l(t)$, $1 \le l \le n$, each of which is in
$\mathcal{B}_n$.  Thus (\ref{int_A g_n d mu = int_A g_{n + 1} d mu})
follows by combining the previous cases.

        Now let us check that $\{g_n\}_{n = 1}^\infty$ is uniformly
integrable.  Let $h$ be the function on $X$ defined by $h = f_n$ on
$E_n(t)$ and $h = 0$ on $X \backslash E(t)$, as in the previous
section.  Observe that $g_n = h$ on $\bigcup_{l = 1}^n E_l(t)$, while
$g_n = f_n$ on $X \backslash \Big(\bigcup_{l = 1}^n E_l(t)\Big)$.
Moreover,
\begin{equation}
\label{|g_n| = |f_n| le t}
        |g_n| = |f_n| \le t
\end{equation}
on $X \backslash \Big(\bigcup_{l = 1}^n E_l(t)\Big)$, by definition of
$E_l(t)$.  This implies that
\begin{equation}
\label{|g_n| le max (|h|, t)}
        |g_n| \le \max (|h|, t)
\end{equation}
on $X$ for each $n$, so that the uniform integrability of $\{g_n\}_{n
= 1}^\infty$ follows from the integrability of $h$.

        Thus $\{g_n\}_{n = 1}^\infty$ converges pointwise almost
everywhere on $X$, as mentioned earlier in the section.  By
construction, $g_n = f_n$ on $X \backslash E(t)$ for each $n$, and so
$\{f_n\}_{n = 1}^\infty$ converges pointwise almost everywhere on $X
\backslash E(t)$ for each $t > 0$.  It follows that $\{f_n\}_{n =
1}^\infty$ converges pointwise almost everywhere on
\begin{equation}
        \bigcup_{k = 1}^\infty (X \backslash E(k))
         = X \backslash \Big(\bigcap_{k = 1}^\infty E(k)\Big).
\end{equation}
Of course,
\begin{equation}
        \mu\Big(\bigcap_{k = 1}^\infty E(k)\Big) \le \inf_{k \ge 1} \mu(E(k)),
\end{equation}
and $\mu(E(t)) \le t^{-1} \sup_{n \ge 1} \|f_n\|_1 \to 0$ as $t \to
\infty$, by (\ref{t mu(E(t)) = t mu(bigcup_{l = 1}^infty E_l(t)) le
C}).  Hence
\begin{equation}
\label{mu(bigcap_{k = 1}^infty E(k)) = 0}
        \mu\Big(\bigcap_{k = 1}^\infty E(k)\Big) = 0,
\end{equation}
which implies that $\{f_n\}_{n = 1}^\infty$ converges pointwise almost
everywhere on $X$.

\section{Other measures}
\label{other measures}
\setcounter{equation}{0}

        Let $(X, \mathcal{A}, \mu)$ be a probability space, and let
$\mathcal{B}_1 \subseteq \mathcal{B}_2 \subseteq \cdots$ be an
increasing sequence of $\sigma$-subalgebras of $\mathcal{A}$.  Also
let $\nu$ be a real or complex measure on a $\sigma$-algebra
$\mathcal{B} \subseteq \mathcal{A}$ that contains each
$\mathcal{B}_j$.  Suppose that the restriction of $\nu$ to
$\mathcal{B}_j$ is absolutely continuous with respect to the
restriction of $\mu$ to $\mathcal{B}_j$ for each $j$.  In particular,
this happens when each $\mathcal{B}_j$ is associated to a partition of
$X$ by finitely or countably many sets of positive $\mu$-measure, as
in Section \ref{partitions}.  Under these conditions, the
Radon--Nikodym theorem implies that there is an $f_j \in L^1(X,
\mathcal{B}_j)$ for each $j \ge 1$ such that
\begin{equation}
        \int_B f_j \, d\mu = \nu(B)
\end{equation}
for every $B \in \mathcal{B}_j$.

        By construction, $\{f_j\}_{j = 1}^\infty$ is a martingale on
$X$ with respect to the $\mathcal{B}_j$'s.  Moreover,
\begin{equation}
\label{int_X |f_j| d mu le |nu|(X)}
        \int_X |f_j| \, d\mu \le |\nu|(X)
\end{equation}
for each $j$, where $|\nu|$ denotes the total variation measure
associated to $\nu$.  As in Section \ref{basic properties},
(\ref{int_X |f_j| d mu le |nu|(X)}) basically corresponds to the
statement that the total variation of the restriction of $\nu$ to
$\mathcal{B}_j$ is less than or equal to the restriction of $|\nu|$ to
$\mathcal{B}_j$.  If $\nu$ is absolutely continuous with respect to
the restriction of $\mu$ to $\mathcal{B}$, so that there is an $f \in
L^1(X, \mathcal{B})$ such that
\begin{equation}
\label{nu(B) = int_B f d mu}
        \nu(B) = \int_B f \, d\mu
\end{equation}
for every $B \in \mathcal{B}$, then $f_j = E(f \mid \mathcal{B}_j)$
for each $j$.

        Put
\begin{equation}
        d'(A, B) = \mu(A \bigtriangleup B) + |\nu|(A \bigtriangleup B)
\end{equation}
for every $A, B \in \mathcal{B}$.  This defines a semimetric on
$\mathcal{B}$, as in Section \ref{distances between sets}, and the
closure $\mathcal{C}'$ of $\mathcal{E} = \bigcup_{j = 1}^\infty
\mathcal{B}_j$ with respect to $d'(A, B)$ is a $\sigma$-subalgebra of
$\mathcal{B}$ that contains $\mathcal{E}$.  More precisely,
$\mathcal{C}'$ is the smallest $\sigma$-subalgebra of $\mathcal{B}$
that contains $\mathcal{E}$ and the sets $A \in \mathcal{B}$ such that
$\mu(A) = |\nu|(A) = 0$.  In particular, $\mathcal{C}'$ contains the
smallest $\sigma$-algebra $\mathcal{B}_\infty$ that contains
$\mathcal{E}$, and $\mathcal{C}'$ is contained in the closure
$\mathcal{C}$ of $\mathcal{E}$ with respect to $d(A, B) = \mu(A
\bigtriangleup B)$.

        Suppose that $\{f_j\}_{j = 1}^\infty$ converges to a function
$f \in L^1(X, \mathcal{B})$ in the $L^1$ norm.  If $A \in
\mathcal{B}_l$ for some $l$, so that
\begin{equation}
        \int_A f_j \, d\mu = \int_A f_l \, d\mu = \nu(A)
\end{equation}
when $j \ge l$, then
\begin{equation}
        \int_A f \, d\mu = \lim_{j \to \infty} \int_A f_j \, d\mu = \nu(A).
\end{equation}
Thus
\begin{equation}
        \int_A f \, d\mu = \nu(A)
\end{equation}
for every $A \in \mathcal{E}$, and hence for every $A \in
\mathcal{C}'$, because both sides of the equation are continuous with
respect to $d'(A, B)$.  This uses the analogue of uniform
integrability for the single integrable function $f$.  It follows that
the restriction of $\nu$ to $\mathcal{B}_\infty \subseteq
\mathcal{C}'$ is absolutely continuous with respect to the restriction
of $\mu$ to $\mathcal{B}_\infty$ under these conditions.

\section{Finitely-additive measures}
\label{finite additivity}
\setcounter{equation}{0}

        Let $(X, \mathcal{A}, \mu)$ be a probability space, let
$\mathcal{B}_1 \subseteq \mathcal{B}_2 \subseteq \cdots$ be an
increasing sequence of $\sigma$-subalgebras of $\mathcal{A}$, and let
$\{f_j\}_{j = 1}^\infty$ be a martingale on $X$ with respect to this
filtration.  If we put
\begin{equation}
\label{nu(A) = int_A f_j d mu}
        \nu(A) = \int_A f_j \, d\mu
\end{equation}
when $A \in \mathcal{B}_j$, then $\nu$ is well-defined on $\mathcal{E}
= \bigcup_{j = 1}^\infty \mathcal{B}_j$, because
\begin{equation}
        \int_A f_j \, d\mu = \int_A f_l \, d\mu
\end{equation}
when $A \in \mathcal{B}_l$ and $j \ge l$.  It is easy to see that
$\nu$ is finitely additive on $\mathcal{E}$.

        Suppose that the $f_j$'s have bounded $L^1$ norms, so that
there is a $C \ge 0$ with the property that $\|f_j\|_1 \le C$ for
every $j \ge 1$.  Let $A_1, \ldots, A_n$ be finitely many
pairwise-disjoint subsets of $X$ that are contained in $\mathcal{E}$.
Thus $A_1, \ldots, A_n \in \mathcal{B}_l$ for some $l$, and hence
\begin{eqnarray}
\sum_{k = 1}^n |\nu(A_k)| = \sum_{k = 1}^n \biggl|\int_{A_k} f_l \, d\mu\biggr|
         & \le & \sum_{k = 1}^n \int_{A_k} |f_l| \, d\mu \\
          & = & \int_{\bigcup_{k = 1}^n A_k} |f_l| \, d\mu \le C. \nonumber
\end{eqnarray}
Conversely, if
\begin{equation}
\label{sum_{k = 1}^n |nu(A_k)| le C}
        \sum_{k = 1}^n |\nu(A_k)| \le C
\end{equation}
for every collection of finitely many pairwise disjoint elements $A_1,
\ldots, A_n$ of $\mathcal{B}_j$, then $\|f_j\|_1 \le C$.  If $\nu$ has
an extension to a countably-additive real or complex measure on a
$\sigma$-algebra that contains $\mathcal{E}$, then (\ref{sum_{k = 1}^n
|nu(A_k)| le C}) holds for each $j$, with $C$ equal to the total variation
of the extension of $\nu$ on $X$.

        For example, let $X$ be $[0, 1)$ equipped with Lebesgue
measure, and let $\mathcal{B}_j$ be the collection of subsets of $[0,
1)$ that are unions of dyadic intervals of length $2^{-j}$.  In this
case, $\mathcal{E}$ is the algebra of subsets of $[0, 1)$ that can be
expressed as the union of finitely many dyadic intervals.  Put
\begin{eqnarray}
        f_j(x) & = & 0 \quad\hbox{ when } 0 \le x < 1 - 2^{-j}, \\
               & = & 2^j \quad\hbox{when } 1 - 2^{-j} \le x < 1. \nonumber
\end{eqnarray}
Thus
\begin{equation}
        \int_I f_j(x) \, dx = 0
\end{equation}
when $I = [l \, 2^{-j}, (l + 1) \, 2^{-j})$, $0 \le l \le 2^j - 2$, and
\begin{equation}
        \int_I f_j(x) \, dx = 1
\end{equation}
when $I = [1 - 2^{-j}, 1)$, which corresponds to $l = 2^j - 1$.  It is
easy to see that $\{f_j\}_{j = 1}^\infty$ is a martingale on $[0, 1)$
with respect to this filtration.  The finitely-additive measure $\nu$
on $\mathcal{E}$ is characterized by $\nu(I) = 1$ when $I$ is a dyadic
interval with $1$ as an endpoint, and $\nu(I) = 0$ for every other
dyadic interval $I$.  Note that $\|f_j\|_1 = 1$ for each $j$, and that
(\ref{sum_{k = 1}^n |nu(A_k)| le C}) holds with $C = 1$, as it should.
If $I_j = [1 - 2^{-j}, 1)$, then $I_{j + 1} \subseteq I_j$ and
$\nu(I_j) = 1$ for each $j \ge 1$, but $\bigcap_{j = 1}^\infty I_j =
\emptyset$.  Basically, this martingale corresponds to a Dirac mass at
the point $1$.  Since $1$ is not included as an element of $X = [0,
1)$, there is no countably-additive measure on $X$ from which the
martingale is obtained.

        Let $(X, \mathcal{A}, \mu)$ be any probability space again,
with an increasing sequence $\mathcal{B}_1 \subseteq \mathcal{B}_2
\subseteq \cdots$ of $\sigma$-subalgebras of $\mathcal{A}$, and let
$\{f_j\}_{j = 1}^\infty$ be a martingale on $X$ with respect to this
filtration with bounded $L^1$ norms.  As in Section \ref{convergence
almost everywhere}, $\{f_j\}_{j = 1}^\infty$ converges pointwise
almost everywhere on $X$.  The limit determines an element $g$ of
$L^1(X, \mathcal{C})$, where $\mathcal{C}$ is the closure of
$\mathcal{E}$ with respect to the usual semimetric $d(A, B) = \mu(A
\bigtriangleup B)$ on $\mathcal{A}$, as in Section \ref{distances
between sets}.  Equivalently, $\mathcal{C}$ is the smallest
$\sigma$-subalgebra of $\mathcal{A}$ that contains $\mathcal{E}$ and
every $A \in \mathcal{A}$ with $\mu(A) = 0$.  If $g_j = E(g \mid
\mathcal{B}_j)$ for each $j$, then $\{g_j\}_{j = 1}^\infty$ is a
martingale on $X$ with respect to this filtration that converges to
$g$ in the $L^1$ norm, as in Section \ref{sequences of
sigma-subalgebras}.  Hence $\{g_j\}_{j = 1}^\infty$ also converges to
$g$ pointwise almost everywhere on $X$, as in Section \ref{convergence
almost everywhere}.  If $h_j = f_j - g_j$, then $\{h_j\}_{j =
1}^\infty$ is also a martingale on $X$ with respect to this
filtration, and with bounded $L^1$ norms.  By construction,
$\{h_j\}_{j = 1}^\infty$ converges to $0$ pointwise almost everywhere
on $X$.  One can think of $\{g_j\}_{j = 1}^\infty$ as the ``regular
part'' of the martingale $\{f_j\}_{j = 1}^\infty$, and of $\{h_j\}_{j
= 1}^\infty$ as the ``singular part'' of $\{f_j\}_{j = 1}^\infty$.

\section{Maximal functions, 4}
\label{maximal functions, 4}
\setcounter{equation}{0}

        Let $(X, \mathcal{A}, \mu)$ be a probability space, and let
$\mathcal{B}_1 \subseteq \mathcal{B}_2 \subseteq \cdots$ be an
increasing sequence of $\sigma$-subalgebras of $\mathcal{A}$.
If $f \in L^1(X, \mathcal{A})$, then $f_j = E(f \mid \mathcal{B}_j)$
defines a martingale  on $X$ with respect to this filtration,
and we get the corresponding maximal function
\begin{equation}
        f^*(x) = \sup_{j \ge 1} |f_j(x)|,
\end{equation}
as before.  Note that $f \mapsto f^*$ is sublinear, in the sense that
\begin{equation}
        (a \, f)^* = |a| \, f^*
\end{equation}
and
\begin{equation}
        (f + g)^* \le f^* + g^*
\end{equation}
for every $f, g \in L^1(X, \mathcal{A})$ and $a \in {\bf R}$ or ${\bf C}$.

        If $f \in L^\infty(X, \mathcal{A})$, then $f_j \in L^\infty(X,
\mathcal{A})$ for each $j$, and
\begin{equation}
\label{||f_j||_infty le ||f||_infty}
        \|f_j\|_\infty \le \|f\|_\infty,
\end{equation}
as in Section \ref{basic properties}.  This implies that $f^* \in
L^\infty(X, \mathcal{A})$, and that
\begin{equation}
        \|f^*\|_\infty \le \|f\|_\infty.
\end{equation}

        If $f \in L^1(X, \mathcal{A})$, then
\begin{equation}
        \|f_j\|_1 \le \|f\|
\end{equation}
for each $j$, as in Section \ref{basic properties}.  Put
\begin{equation}
\label{E(t) = {x in X : f^*(x) > t}}
        E(t) = \{x \in X : f^*(x) > t\}
\end{equation}
for each $t > 0$, so that
\begin{equation}
\label{mu(E(t)) le t^{-1} ||f||_1}
        \mu(E(t)) \le t^{-1} \, \|f\|_1,
\end{equation}
as in Section \ref{maximal functions, 3}.

        Let $g$ be the function defined on $X$ by
\begin{eqnarray}
        g(x) & = & f(x) \quad\hbox{when } |f(x)| \le t/2 \\
             & = & 0    \qquad\ \hbox{when } |f(x)| > t/2. \nonumber
\end{eqnarray}
Thus $g \in L^\infty(X, \mathcal{A})$, and hence $g^* \in L^\infty(X,
\mathcal{A})$, with
\begin{equation}
        \|g^*\|_\infty \le \|g\|_\infty \le \frac{t}{2}.
\end{equation}
This implies that
\begin{equation}
        f^*(x) \le (f - g)^*(x) + g^*(x) \le (f - g)^*(x) + \frac{t}{2}
\end{equation}
for almost every $x \in X$, so that
\begin{equation}
        (f - g)^*(x) > t/2
\end{equation}
for almost every $x \in E(t)$.

        It follows that
\begin{equation}
 \mu(E(t)) \le \mu(\{x \in X : (f - g)^*(x) > t/2\}) \le t^{-1} \, \|f - g\|_1.
\end{equation}
Using the definition of $g$, we get that
\begin{equation}
\label{mu(E(t)) le t^{-1} int_{{x in X : |f(x)| > t/2}} |f(x)| d mu(x)}
 \mu(E(t)) \le t^{-1} \int_{\{x \in X : |f(x)| > t/2\}} |f(x)| \, d\mu(x).
\end{equation}

        If $h$ is a nonnegative measurable function on $X$, then
\begin{equation}
        A(h) = \{(x, r) \in X \times {\bf R} : 0 < r < h(x)\}
\end{equation}
is a measurable subset of $X \times {\bf R}$.  This is easy to see
when $h$ is a measurable simple function, and otherwise $h$ can be
approximated by an increasing sequence of measurable simple functions.
Integrating $p \, r^{p - 1}$ over $A(h)$ with respect to the product
of $\mu$ on $X$ and Lebesgue measure on ${\bf R}$, we get that
\begin{equation}
\label{int_X h^p d mu = int_0^infty p r^{p - 1} mu({x in X : h(x) > r}) dr}
        \int_X h^p \, d\mu =
          \int_0^\infty p \, r^{p - 1} \, \mu(\{x \in X : h(x) > r\}) \, dr.
\end{equation}
More precisely, the left side of (\ref{int_X h^p d mu = int_0^infty p
r^{p - 1} mu({x in X : h(x) > r}) dr}) obtained by integrating $p \,
r^{p - 1}$ over $A(h)$ in $r$ and then $x$, while the right side is
obtained by integrating in $x$ and then $r$.

        In particular, if $1 < p < \infty$, then
\begin{eqnarray}
\int_X (f^*)^p \, d\mu & = & \int_0^\infty p \, t^{p - 1} \, \mu(E(t)) \, dt \\
 & \le & \int_0^\infty p \, t^{p - 2} \int_{\{x \in X : |f(x)| > t/2\}}
                                             |f(x)| \, d\mu(x) \, dt, \nonumber
\end{eqnarray}
by (\ref{mu(E(t)) le t^{-1} int_{{x in X : |f(x)| > t/2}} |f(x)| d mu(x)}).
Interchanging the order of integration, we get that
\begin{eqnarray}
 \int_X (f^*)^p \, d\mu & \le & \int_X \int_0^{2 \, |f(x)|}
                              |f(x)| \, p \, t^{p - 2} \, dt \, d\mu(x) \\
 & = & \frac{p \, 2^{p - 1}}{p - 1} \int_X |f(x)|^p \, d\mu(x). \nonumber
\end{eqnarray}
This shows that $f^* \in L^p(X, \mathcal{A})$ when $f \in L^p(X,
\mathcal{A})$ and $p > 1$.

        By constrast, if $f \in L^p(X, \mathcal{A})$, then
\begin{equation}
        (f^*(x))^p \le (|f|^p)^*(x),
\end{equation}
by (\ref{|f_mathcal{B}|^p le (|f|^p)_mathcal{B} = E(|f|^p mid mathcal{B})}).
As before,
\begin{equation}
 \mu(\{x \in X : (|f|^p)^*(x) > t\}) \le t^{-1} \int_X |f(x)|^p \, d\mu(x)
\end{equation}
for every $t > 0$.  This implies that
\begin{equation}
 \mu(\{x \in X : (f^*(x))^p > t\}) \le t^{-1} \int_X |f(x)|^p \, d\mu(x),
\end{equation}
or equivalently
\begin{equation}
        \mu(\{x \in X : f^*(x) > t\}) \le t^{-p} \int_X |f(x)|^p \, d\mu(x)
\end{equation}
for every $t > 0$.  This is not strong enough to imply that $f^* \in
L^p$, by integrating over $t$ as in the previous paragraph.  However,
it does have the advantage of working uniformly over $p \ge 1$.

        Note that we get the same estimates for the dyadic maximal
function, as in Section \ref{dyadic intervals}, which corresponds to
$X = [0, 1)$ with Lebesgue measure, and where $\mathcal{B}_j$ consists
of unions of dyadic intervals of length $2^{-j}$.  There are also
similar estimates for the Hardy--Littlewood maximal function on the
real line, as in Section \ref{maximal functions}, but with an extra
factor of $2$ in (\ref{mu(E(t)) le t^{-1} ||f||_1}), and in the later
steps.

\section{Decreasing sequences of $\sigma$-algebras}
\label{decreasing sequences of sigma-algebras}
\setcounter{equation}{0}

        Let $(X, \mathcal{A}, \mu)$ be a probability space, and
suppose that $\mathcal{A}_1 \supseteq \mathcal{A}_2 \supseteq \cdots$
is a decreasing sequence of $\sigma$-subalgebras of $\mathcal{A}$.  
As a basic scenario, it may be that $X = \prod_{j = 1}^\infty X_j$ is
the Cartesian product of a sequence of probability spaces $X_1, X_2,
\ldots$, and that $\mathcal{A}_n$ consists of subsets of $X$ of the
form $\prod_{j = 1}^n X_j \times A$, where $A$ is a measurable subset
of $\prod_{j = n + 1}^\infty X_j$.  In this case, conditional expectation
with respect to $\mathcal{A}_n$ corresponds to integrating a function on $X$
in $x_1, \ldots, x_n$.  Basically, conditional expectation with respect to
smaller $\sigma$-algebras corresponds to averaging functions over larger
sets.

        Note that $\mathcal{A}_\infty = \bigcap_{j = 1}^\infty
\mathcal{A}_j$ is automatically a $\sigma$-subalgebra of
$\mathcal{A}$.  If $A_j \in \mathcal{A}_j$ satisfies $A_j \subseteq
A_{j + 1}$ for each $j$, then $\bigcup_{j = 1}^\infty A_j \in
\mathcal{A}_\infty$, because
\begin{equation}
 \bigcup_{j = 1}^\infty A_j = \bigcup_{j = n}^\infty A_j \in \mathcal{A}_n
\end{equation}
for each $n$.  Similarly, if $B_j \in \mathcal{A}_j$ satisfies $B_{j +
1} \subseteq B_j$ for each $j$, then
\begin{equation}
 \bigcap_{j = 1}^\infty B_j = \bigcap_{j = n}^\infty B_j \in \mathcal{A}_n
\end{equation}
for each $n$, and so $\bigcap_{j = 1}^\infty B_j \in
\mathcal{A}_\infty$.  If $E_j \in \mathcal{A}_j$ for each $j$, then it
follows that
\begin{equation}
        \limsup_{j \to \infty} E_j = \bigcap_{l = 1}^\infty
                                   \Big(\bigcup_{j = l}^\infty E_j\Big), \quad
        \liminf_{j \to \infty} E_j = \bigcup_{l = 1}^\infty
                                   \Big(\bigcap_{j = l}^\infty E_j\Big)
\end{equation}
are also elements of $\mathcal{A}_\infty$, by taking $A_l = \bigcup_{j
= l}^\infty E_j$ and $B_l = \bigcap_{j = l}^\infty E_j$.

        If $f$ is a measurable function on $X$ with respect to
$\mathcal{A}$, and if $f_j$ is a measurable function on $X$ with
respect to $\mathcal{A}_j$ such that $f = f_j$ almost everywhere for
each $j$, then there is a measurable function $f_\infty$ on $X$ with
respect to $\mathcal{A}_\infty$ such that $f = f_\infty$ almost everywhere.
To see this, put
\begin{equation}
        E_j = \{x \in X : f_j(x) = f_{j + 1}(x)\},
\end{equation}
so that $E_j \in \mathcal{A}_j$ for each $j$.  Thus $B_l = \bigcap_{j
= l}^\infty E_j \in \mathcal{A}_l$, and $f_j(x) = f_l(x)$ for every $x
\in B_l$ and $j \ge l$.  By hypothesis, $\mu(X \backslash E_j) = 0$
for each $j$, and so $\mu(X \backslash B_l) = 0$ for each $l$, since
$X \backslash B_l = \bigcup_{j = l}^\infty (X \backslash E_j)$.  We
also have that $\bigcup_{l = 1}^\infty B_l \in \mathcal{A}_\infty$, as
in the previous paragraph.  Put
\begin{eqnarray}
        f_\infty(x) & = & 0 \qquad\ \hbox{ when }
                  x \in X \backslash \Big(\bigcup_{l = 1}^\infty B_l\Big) \\
 & = & f_l(x) \quad\hbox{when } x \in B_l \hbox{ for some } l \ge 1. \nonumber
\end{eqnarray}
This is well defined, because $f_j(x) = f_l(x)$ when $x \in B_l$ and
$j \ge l$.  Moreover, $f_\infty$ is measurable with respect to
$\mathcal{A}_l$ for every $l$, because $f_l$ is measurable with
respect to $\mathcal{A}_l$.  This implies that $f_\infty$ is
measurable with respect to $\mathcal{A}_\infty$.  It is easy to see
that $f = f_\infty$ almost everywhere, since $f = f_l$ almost
everywhere.

        Let $\{f_j\}_{j = 1}^\infty$ be a sequence of real-valued
functions on $X$ such that $f_j$ is measurable with respect to
$\mathcal{A}_j$ for each $j$.  Thus
\begin{equation}
        \sup_{j \ge l} f_j(x), \quad \inf_{j \ge l} f_j(x)
\end{equation}
are measurable with respect to $\mathcal{A}_l$ for each $l$.  This
implies that
\begin{equation}
        \limsup_{j \to \infty} f_j(x), \quad \liminf_{j \to \infty} f_j(x)
\end{equation}
are measurable with respect to $\mathcal{A}_l$ for each $l$, and hence
are measurable with respect to $\mathcal{A}_\infty$.  In particular,
the set of $x \in X$ on which $\{f_j(x)\}_{j = 1}^\infty$ converges is
measurable with respect to $\mathcal{A}_\infty$, and the limit defines
a measurable function with respect to $\mathcal{A}_\infty$ on this
set.  The analogous statement for complex-valued functions follows by
considering the real and imaginary parts separately.

        Let $f \in L^1(X, \mathcal{A})$ be given, and put $f_j = E(f
\mid \mathcal{A}_j)$ for each $j \ge 1$, and $f_0 = f$.  Thus
\begin{equation}
        f = \sum_{j = 1}^n (f_{j - 1} - f_j) + f_n
\end{equation}
for each $n \ge 1$.  If $f \in L^2(X, \mathcal{A})$, then the
functions $f_{j - 1} - f_j$, $1 \le j \le n$, and $f_n$ are pairwise
orthogonal in $L^2(X, \mathcal{A})$, as in Section \ref{martingales}.
This implies that
\begin{equation}
\label{||f||_2^2 = sum_{j = 1}^n ||f_{j - 1} - f_j||_2^2 + ||f_n||_2^2}
        \|f\|_2^2 = \sum_{j = 1}^n \|f_{j - 1} - f_j\|_2^2 + \|f_n\|_2^2
\end{equation}
for each $n$, and hence that $\sum_{j = 1}^\infty \|f_{j - 1} -
f_j\|_2^2$ converges.  Therefore
\begin{equation}
\label{sum_{j = 1}^infty (f_{j - 1} - f_j)}
        \sum_{j = 1}^\infty (f_{j - 1} - f_j)
\end{equation}
converges in $L^2(X, \mathcal{A})$, by orthogonality, which implies
that $\{f_n\}_{n = 1}^\infty$ converges in $L^2(X, \mathcal{A})$.  Of
course, $\{f_n\}_{n = l}^\infty$ converges in $L^2(X, \mathcal{A}_l)$
for each $l$, and the limits correspond to the same element of $L^2(X,
\mathcal{A})$ for each $l$.  Thus the limit may be represented by an
element $f_\infty$ of $L^2(X, \mathcal{A}_\infty)$, by the earlier
remarks.  In particular, $f_\infty = E(f_\infty \mid \mathcal{A}_\infty)$,
which implies that
\begin{equation}
        f_\infty = E(f \mid \mathcal{A}_\infty).
\end{equation}
This uses the fact that $E(f \mid \mathcal{A}_\infty) = E(f_j \mid
\mathcal{A}_\infty)$ for each $j$, since $f_j = E(f \mid
\mathcal{A}_j)$ and $\mathcal{A}_\infty \subseteq \mathcal{A}_j$, and
the convergence of $\{f_j\}_{j = 1}^\infty$ to $f_\infty$ in $L^2(X,
\mathcal{A})$.

        If $f \in L^p(X, \mathcal{A})$, $1 \le p < 2$, then $L^2(X,
\mathcal{A})$ is a dense linear subspace of $L^p(X, \mathcal{A})$, and
one can use this to show that $\{f_j\}_{j = 1}^\infty$ converges to
$E(f \mid \mathcal{A}_\infty)$ in the $L^p$ norm.  This also uses the
fact that the conditional expectation operators have operator norm $1$
on $L^p$ for each $p$.  If $f \in L^\infty(X, \mathcal{A})$, then $f_j
\in L^\infty(X, \mathcal{A}_j)$ with $\|f_j\|_\infty \le \|f\|_\infty$
for each $j$.  This together with convergence in $L^2(X, \mathcal{A})$
implies convergence in $L^p(X, \mathcal{A})$ for every $p < \infty$.
If $f \in L^p(X, \mathcal{A})$, $2 < p < \infty$, then one can show
again that $\{f_j\}_{j = 1}^\infty$ converges to $E(f \mid
\mathcal{A}_\infty)$ in the $L^p$ norm, since this holds on the
dense linear subspace $L^\infty(X, \mathcal{A})$ of $L^p(X, \mathcal{A})$,
and because the expectation operators are uniformly bounded on $L^p$.

        There are also maximal function estimates in this context.  To
see this, one can begin by observing that
\begin{equation}
        f_n^*(x) = \max_{1 \le j \le n} |f_j(x)|
\end{equation}
is basically the same as before, because one can simply rearrange the
indices to get an increasing sequence of $n$ $\sigma$-algebras.  Hence
the estimates for $f_n^*$ are the same as before, and the
corresponding estimates for
\begin{equation}
        f^*(x) = \sup_{j \ge 1} |f_j(x)|
\end{equation}
can be obtained by passing to the limit as $n \to \infty$.
Convergence almost everywhere then follows from convergence in the
$L^1$ norm, as in Section \ref{convergence almost everywhere}.

\section{Doubly-infinite sequences}
\label{doubly-infinite sequences}
\setcounter{equation}{0}

        A probability space $(X, \mathcal{A}, \mu)$ may also have a
doubly-infinite sequence
\begin{equation}
        \cdots \subseteq \mathcal{B}_{-1} \subseteq \mathcal{B}_0
                           \subseteq \mathcal{B}_1 \subseteq \cdots
\end{equation}
of $\sigma$-subalgebras of $\mathcal{A}$.  In particular, this occurs
very naturally in the context of doubly-infinite products.  Let $(X_j,
\mathcal{A}_j, \mu_j)$, $j \in {\bf Z}$ be a family of probability
spaces indexed by the integers, and let $X = \prod_{j =
-\infty}^\infty X_j$ be their Cartesian product, equipped with the
product measure $\mu$.  Thus $X$ consists of the doubly-infinite
sequences $x = \{x_j\}_{j = -\infty}^\infty$ such that $x_j \in X_j$
for each $j$.  If $\mathcal{B}_n$ is the collection of subsets of $X$
of the form $A \times \prod_{j = n + 1}^\infty X_j$, where $A$ is a
measurable subset of $\prod_{j = -\infty}^n X_j$, then $\mathcal{B}_n$
is a $\sigma$-subalgebra of the $\sigma$-algebra of measurable subsets
of $X$, and $\mathcal{B}_n \subseteq \mathcal{B}_{n + 1}$ for each $n$.

        Suppose that $(X_j, \mathcal{A}_j, \mu_j)$ is a copy of the
same probability space for each $j$.  In this case, we can define the
shift mapping $T : X \to X$ by $T(x) = y$, where $x = \{x_j\}_{j =
-\infty}^\infty, y = \{y_j\}_{j = -\infty}^\infty \in X$ satisfy
\begin{equation}
        y_j = x_{j - 1}
\end{equation}
for each $j$.  If $A \subseteq X$ is measurable, then $T(A)$ is also
measurable, and
\begin{equation}
        \mu(T(A)) = \mu(A).
\end{equation}
Similarly, $T$ maps $\mathcal{B}_n$ onto $\mathcal{B}_{n + 1}$ for each $n$.

        If the $X_j$'s are compact Hausdorff topological spaces, then
$X$ is too, with respect to the product topology.  If the $X_j$'s are
all copies of the same topological space, then $T$ is a homeomorphism.
If the $X_j$'s are all metrizable, then $X$ is as well, as in Section
\ref{metrizability}.  However, this does not mean that there is a
metric $d(x, y)$ on $X$ that determines the product topology and which
is invariant under $T$ in the sense that
\begin{equation}
        d(T(x), T(y)) = d(x, y)
\end{equation}
for every $x, y \in X$.  If $x, y \in X$ satisfy $x_j = x_l$ for every
$j, l \in {\bf Z}$ and $x_j = y_j$ for all but exactly one $j \in {\bf
Z}$, then $T(x) = x$ and $\lim_{n \to \infty} T^n(y) = x$, which would
not be possible if there were an invariant metric.

\section{Submartingales}
\label{submartingales}
\setcounter{equation}{0}

        Let $(X, \mathcal{A}, \mu)$ be a probability space, and let
$\mathcal{B}_1 \subseteq \mathcal{B}_2 \subseteq \cdots$ be an
increasing sequence of $\sigma$-subalgebras of $\mathcal{A}$.  Also
let $\{f_j\}_{j = 1}^\infty$ be a sequence of real-valued functions on
$X$ such that $f_j \in L^1(X, \mathcal{B}_j)$ for each $j$.  We say
that $\{f_j\}_{j = 1}^\infty$ is a \emph{submartingale} on $X$ with
respect to this filtration if
\begin{equation}
        f_j \le E(f_{j + 1} \mid \mathcal{B}_j)
\end{equation}
almost everywhere on $X$ with respect to $\mu$ for each $j$.
Similarly, $\{f_j\}_{j = 1}^\infty$ is a \emph{supermartingale} if
\begin{equation}
        f_j \ge E(f_{j + 1} \mid \mathcal{B}_j)
\end{equation}
almost everywhere on $X$ for each $j$.  Thus $\{f_j\}_{j = 1}^\infty$
is a martingale if and only if it is both a submartingale and a
supermartingale, and $\{f_j\}_{j = 1}^\infty$ is a supermartingale if
and only if $\{-f_j\}_{j = 1}^\infty$ is a submartingale.

        If $\{g_j\}_{j = 1}^\infty$ is a real or complex martingale on
$X$ with respect the $\mathcal{B}_j$'s, then $\{|g_j|\}_{j =
1}^\infty$ is a submartingale on $X$.  If in addition $g_j \in L^p(X,
\mathcal{B}_j)$ for some $p$, $1 < p < \infty$, and each $j$, then
$\{|g_j|^p\}_{j = 1}^\infty$ is a submartingale as well.  More
generally, if $\phi$ is a convex function on an interval $I$ in the
real line, which may be unbounded, and if $g_j$ takes valued in $I$
and $\phi \circ g_j \in L^1(X, \mathcal{B}_j)$ for each $j$, then
$\{\phi \circ g_j\}_{j = 1}^\infty$ is a submartingale.  These
statements use the remarks in Section \ref{basic properties}.  The
latter also works when $\{g_j\}_{j = 1}^\infty$ is a submartingale and
$\phi$ is both convex and monotone increasing on $I$.

        If $\{f_j\}_{j = 1}^\infty$ is a submartingale on $X$ and $a$
is a nonnegative real number, then $\{a \, f_j\}_{j = 1}^\infty$ is a
submartingale.  If $\{f_j\}_{j = 1}^\infty$, $\{g_j\}_{j = 1}^\infty$
are submartingales, then their sum $\{f_j + g_j\}_{j = 1}^\infty$ is a
martingale too.  Their maximum $\{\max(f_j, g_j)\}_{j = 1}^\infty$ is
a submartingale as well, because
\begin{equation}
        f_j \le E(f_{j + 1} \mid \mathcal{B}_j)
                  \le E(\max(f_{j + 1}, g_{j + 1}) \mid \mathcal{B}_j)
\end{equation}
and
\begin{equation}
        g_j \le E(g_{j + 1} \mid \mathcal{B}_j)
                  \le E(\max(f_{j + 1}, g_{j + 1}) \mid \mathcal{B}_j) 
\end{equation}
imply that
\begin{equation}
        \max(f_j, g_j) \le E(\max(f_{j + 1}, g_{j + 1}) \mid \mathcal{B}_j).
\end{equation}
Of course, $\{f_j + g_j\}_{j = 1}^\infty$ is a martingale when
$\{f_j\}_{j = 1}^\infty$, $\{g_j\}_{j = 1}^\infty$ are martingales,
but $\{\max(f_j, g_j)\}_{j = 1}^\infty$ is not normally a martingale
in this case.

        Let $\{f_j\}_{j = 1}^\infty$ be a sequence of real-valued
functions on $X$ with $f_j \in L^1(X, \mathcal{B}_j)$ for each $j$, as
before.  Thus $\{f_j\}_{j = 1}^\infty$ is determined by the initial
function $f_1$ and the sequence of differences $f_{j + 1} - f_j$.  The
condition that $\{f_j\}_{j = 1}^\infty$ be a martingale can be
expressed by
\begin{equation}
        E(f_{j + 1} - f_j \mid \mathcal{B}_j) = 0
\end{equation}
for each $j$, while the condition that $\{f_j\}_{j = 1}^\infty$ be a
submartinagle is expressed by
\begin{equation}
        E(f_{j + 1} - f_j \mid \mathcal{B}_j) \ge 0.
\end{equation}

        Suppose that $\{f_j\}_{j = 1}^\infty$ is a submartingale, and
put
\begin{equation}
\label{a_j = E(f_{j + 1} - f_j mid mathcal{B}_j) ge 0}
        a_j = E(f_{j + 1} - f_j \mid \mathcal{B}_j) \ge 0
\end{equation}
for each $j$.  Also put $A_l = \sum_{j = 1}^{l - 1} a_j$ when $l \ge
2$, and $A_1 = 0$.  Note that $A_l \in L^1(X, \mathcal{B}_{l - 1})$
when $l \ge 2$, and $A_l(x)$ is monotone increasing in $l$ for each $x
\in X$.  By construction, $\{f_l - A_l\}_{l = 1}^\infty$ is a
martingale, because
\begin{equation}
        (f_{l + 1} - A_{l + 1}) - (f_l - A_l) = f_{l + 1} - f_l - a_l
\end{equation}
and
\begin{eqnarray}
        E(f_{l + 1} - f_l - a_l \mid \mathcal{B}_l)
 & = & E(f_{l + 1} - f_l \mid \mathcal{B}_l) - E(a_l \mid \mathcal{B}_l) \\
 & = & E(f_{l + 1} - f_l \mid \mathcal{B}_l) - a_l = 0. \nonumber
\end{eqnarray}

        Conversely, if $\{\phi_j\}_{j = 1}^\infty$ is any sequence of
real-valued functions on $X$ such that $\phi_j \in L^1(X,
\mathcal{B}_j)$ and $\phi_j \le \phi_{j + 1}$ for each $j$, then
$\{\phi_j\}_{j = 1}^\infty$ is a submartingale on $X$.  If
$\{\psi_j\}_{j = 1}^\infty$ is a martingale on $X$, then $\{\phi_j +
\psi_j\}_{j = 1}^\infty$ is also a submartingale.  Every submartingale
on $X$ can be represented in this way, by the remarks in the previous
paragraph.

        Suppose that $f_j = \phi_j + \psi_j$ is a submartingale on
$X$, where $\{\psi_j\}_{j = 1}^\infty$ is a martingale, and $\phi_j
\le \phi_{j + 1}$ for each $j$.  If the integrals
\begin{equation}
\label{int_X f_j d mu}
        \int_X f_j \, d\mu
\end{equation}
have an upper bound in ${\bf R}$, then the integrals
\begin{equation}
\label{int_X phi_j d mu}
        \int_X \phi_j \, d\mu
\end{equation}
also have an upper bound in ${\bf R}$, because $\int_X \psi_j \, d\mu$
is constant in $j$, by hypothesis.  This implies that $\{\phi_j\}_{j =
1}^\infty$ converges pointwise almost everywhere on $X$ and in the
$L^1$ norm, by the monotone convergence theorem.  In particular, the
$\phi_j$'s have bounded $L^1$ norms.  If the $f_j$'s have bounded
$L^1$ norms, then it follows that the $\psi_j$'s have bounded $L^1$
norms too.  This implies that $\{\psi_j\}_{j = 1}^\infty$ converges
pointwise almost everywhere on $X$, as in Section \ref{convergence
almost everywhere}, and hence that $\{f_j\}_{j = 1}^\infty$ converges
pointwise almost everywhere on $X$ as well.  Similarly, $\{\psi_j\}_{j
= 1}^\infty$ converges in the $L^1$ norm when $\{f_j\}_{j = 1}^\infty$
converges in the $L^1$ norm.  Conversely, $\{f_j\}_{j = 1}^\infty$
converges in the $L^1$ norm when $\{\psi_j\}_{j = 1}^\infty$ converges
in the $L^1$ norm and the integrals (\ref{int_X f_j d mu}) have an
upper bound in ${\bf R}$.  If $\{f_j\}_{j = 1}^\infty$ is uniformly
integrable, then $\{\psi_j\}_{j = 1}^\infty$ is uniformly integrable,
because $\{\phi_j\}_{j = 1}^\infty$ converges in $L^1$ and hence is
uniformly integrable.  This implies that $\{\psi_j\}_{j = 1}^\infty$
converges in $L^1$ too, as in Section \ref{uniform integrability}, so
that $\{f_j\}_{j = 1}^\infty$ converges in $L^1$ as well, as in the
case of martingales.

        Let $\{f_j\}_{j = 1}^\infty$ be a submartingale on $X$, and
observe that
\begin{equation}
        \int_X f_j \, d\mu \le \int_X E(f_{j + 1} \mid \mathcal{B}_j) \, d\mu
                            = \int_X f_{j + 1} \, d\mu
\end{equation}
for each $j$.  If $j \ge l$, then
\begin{equation}
        E(f_j \mid \mathcal{B}_l)
         \le E(E(f_{j + 1} \mid \mathcal{B}_j) \mid \mathcal{B}_l)
          = E(f_{j + 1} \mid \mathcal{B}_l),
\end{equation}
and
\begin{equation}
        \int_X (E(f_j \mid \mathcal{B}_l) - f_l) \, d\mu
         = \int_X f_j \, d\mu - \int_X f_l \, d\mu.
\end{equation}
Suppose that $\int_X f_j \, d\mu$ has an upper bound in ${\bf R}$, and
hence converges in ${\bf R}$, by monotonicity.  The monotone
convergence theorem implies that $\{E(f_j \mid \mathcal{B}_l)\}_{j =
l}^\infty$ converges in $L^1(X, \mathcal{B}_l)$ for each $l$.  It is
easy to check that the limit $g_l$ satisfies
\begin{equation}
        g_l = E(g_{l + 1} \mid \mathcal{B}_l)
\end{equation}
for each $l$, because
\begin{equation}
        E(E(f_j \mid \mathcal{B}_{l + 1}) \mid \mathcal{B}_l)
         = E(f_j \mid \mathcal{B}_l)
\end{equation}
for each $j$, $l$.  Thus $\{g_l\}_{l = 1}^\infty$ is a martingale, and
\begin{equation}
        f_l \le E(f_j \mid \mathcal{B}_l) \le g_l
\end{equation}
when $j \ge l$, by construction.  Moreover,
\begin{equation}
        \int_X g_l \, d\mu
         = \lim_{j \to \infty} \int_X E(f_j \mid \mathcal{B}_l) \, d\mu
         = \lim_{j \to \infty} \int_X f_j \, d\mu
\end{equation}
for each $l$, which implies that
\begin{equation}
        \lim_{l \to \infty} \int_X (g_l - f_l) \, d\mu = 0,
\end{equation}
since $\int_X g_l \, d\mu$ is constant in $l$.  

        Conversely, if $\{g'_j\}_{j = 1}^\infty$ is a martingale on
$X$ such that $f_j \le g'_j$ for each $j$, then
\begin{equation}
        \int_X f_j \, d\mu \le \int_X g'_j \, d\mu
\end{equation}
has an upper bound in ${\bf R}$, because $\int_X g'_j \, d\mu$ is
constant in $j$.  In addition,
\begin{equation}
\label{E(f_j mid mathcal{B}_l) le E(g'_j mid mathcal{B}_l) = g'_l}
        E(f_j \mid \mathcal{B}_l) \le E(g'_j \mid \mathcal{B}_l) = g'_l
\end{equation}
when $j \ge l$, which implies that $g_l \le g'_l$ for each $l$, where
$g_l$ is as in the preceding paragraph.

        Let $\{f_j\}_{j = 1}^\infty$ be a submartingale on $X$ again, and put
\begin{equation}
\label{f_n^*(x) = max(f_1(x), ldots, f_n(x))}
        f_n^*(x) = \max(f_1(x), \ldots, f_n(x)).
\end{equation}
This is a bit different from the situation for martingales discussed
in Section \ref{maximal functions, 3}, since we do not take the
absolute values of the functions.  However, if $\{g_j\}_{j =
1}^\infty$ is a martingale, then $f_j = |g_j|$ is a submartingale, and
\begin{equation}
\label{f_n^*(x) = max(|g_1(x)|, ldots, |g_n(x)|)}
        f_n^*(x) = \max(|g_1(x)|, \ldots, |g_n(x)|)
\end{equation}
is the same as before.  Note that $f_n^*$ is measurable with respect
to $\mathcal{B}_n$, as before.

        Put
\begin{equation}
        A_n(t) = \{x \in X : f_n^*(x) > t\}
\end{equation}
for each $n \ge 1$ and $t \in {\bf R}$, and $A_0(t) = \emptyset$.
Thus $A_n(t) \in \mathcal{B}_n$ for each $n \ge 1$, and $A_n(t)
\subseteq A_{n + 1}(t)$.  Observe that
\begin{equation}
        A_l(t) \backslash A_{l - 1}(t)
               = \{x \in X : f_{l - 1}^*(x) \le t, \, f_l(x) > t\}
\end{equation}
when $l \ge 2$, and that
\begin{equation}
        A_1(t) \backslash A_0(t) = A_1(t) = \{x \in X : f_1(x) > t\}.
\end{equation}
In particular, $f_l > t$ on $A_l(t) \backslash A_{l - 1}(t)$, and so
\begin{equation}
        t \, \mu(A_l(t) \backslash A_{l - 1}(t))
              \le \int_{A_l(t) \backslash A_{l - 1}(t)} f_l \, d\mu.
\end{equation}
This implies that
\begin{eqnarray}
        t \, \mu(A_l(t) \backslash A_{l - 1}(t)) & \le &
 \int_{A_l(t) \backslash A_{l - 1}(t)} E(f_n \mid \mathcal{B}_l) \, d\mu \\
           & = & \int_{A_l(t) \backslash A_{l - 1}(t)} f_n \, d\mu \nonumber
\end{eqnarray}
when $1 \le l \le n$, because $f_l \le E(f_n \mid \mathcal{B}_l)$,
since $\{f_j\}_{j = 1}^\infty$ is a submartingale, and $A_l(t)
\backslash A_{l - 1}(t) \in \mathcal{B}_l$.  Of course, the sets
$A_l(t) \backslash A_{l - 1}(t)$, $1 \le l \le n$, are pairwise
disjoint, and their union is $A_n(t)$.  Hence
\begin{eqnarray}
 t \, \mu(A_n(t)) = \sum_{l = 1}^n t \, \mu(A_l(t) \backslash A_{l - 1}(t))
 & \le & \sum_{l = 1}^n \int_{A_l(t) \backslash A_{l - 1}(t)} f_n \, d\mu \\
                         & = & \int_{A_n(t)} f_n \, d\mu \nonumber
\end{eqnarray}
for each $n \ge 1$ and $t \in {\bf R}$.

\section{Another variant}
\label{another variant}
\setcounter{equation}{0}

        Let $(X, \mathcal{A}, \mu)$ be a probability space, let
$\mathcal{B}_1 \subseteq \mathcal{B}_2 \subseteq \cdots$ be an
increasing sequence of $\sigma$-subalgebras of $\mathcal{A}$, and let
$\{f_j\}_{j = 1}^\infty$ be a sequence of functions on $X$ such that
$f_j \in L^1(X, \mathcal{B}_j)$ for each $j$.  As in the previous
section, put
\begin{equation}
        a_j = E(f_{j + 1} - f_j \mid \mathcal{B}_j)
\end{equation}
for each $j$, $A_l = \sum_{j = 1}^{l - 1} a_j$ when $l \ge 2$, and
$A_1 = 0$.  Thus $A_l \in L^1(X, \mathcal{B}_{l - 1})$ when $l \ge 2$,
and $\{f_l - A_l\}_{l = 1}^\infty$ is a martingale, as before.  If
$f_j \in L^p(X, \mathcal{B}_j)$ for some $p \ge 1$ and each $j$, and
\begin{equation}
        \sum_{j = 1}^\infty \|f_{j + 1} - f_j\|_p
\end{equation}
converges, then $\{f_j\}_{j = 1}^\infty$ converges in the $L^p$ norm
and pointwise almost everywhere on $X$.  Suppose instead that $f_j \in
L^p(X, \mathcal{B}_j)$ for each $j$, $\|f_j\|_p$ is bounded, and
\begin{equation}
\label{sum_{j = 1}^infty ||a_j||_p}
        \sum_{j = 1}^\infty \|a_j\|_p
\end{equation}
converges.  This implies that $\{A_l\}_{l = 1}^\infty$ converges in
the $L^p$ norm and pointwise almost everywhere on $X$, and that $\|f_l
- A_l\|_p$ is bounded.  Because $\{f_l - A_l\}_{l = 1}^\infty$ is a
martingale, it follows that $\{f_l - A_l\}_{l = 1}^\infty$ converges
pointwise almost everywhere on $X$, and in the $L^p$ norm when $1 < p
< \infty$.

\section{Averaging functions}
\label{averaging functions}
\setcounter{equation}{0}

        Let $(X_1, \mathcal{A}_1, \mu_1), (X_2, \mathcal{A}_2, \mu_2),
\ldots$ be a sequence of probability spaces, and let $X = \prod_{j =
1}^\infty X_j$ be their Cartesian product, with the product measure
$\mu$.  Also let $\mathcal{B}_n$ be the $\sigma$-algebra of measurable
subsets of $X$ of the form $A \times \prod_{j = n + 1}^\infty X_j$,
where $A$ is a measurable subset of $\prod_{j = 1}^n X_j$.  Suppose
that $\phi_j \in L^2(X_j, \mathcal{A}_j)$ satisfies
\begin{equation}
        \int_{X_j} \phi_j \, d\mu_j = 0
\end{equation}
and
\begin{equation}
        \Big(\int_{X_j} |\phi_j|^2 \, d\mu_j\Big)^{1/2} \le C
\end{equation}
for some $C \ge 0$ and each $j$, and consider
\begin{equation}
        f_n(x) = \frac{1}{n} \sum_{j = 1}^n \phi_j(x_j),
\end{equation}
$x = \{x_j\}_{j = 1}^\infty \in X$.  Thus $f_n \in L^2(X,
\mathcal{B}_n)$ for each $n$, and
\begin{equation}
\label{||f_n||_2^2 = frac{1}{n^2} sum_{j = 1}^n ||phi_j||_2^2 le frac{C^2}{n}}
 \|f_n\|_2^2 = \frac{1}{n^2} \sum_{j = 1}^n \|\phi_j\|_2^2 \le \frac{C^2}{n},
\end{equation}
because of orthogonality.  In particular, $f_n \to 0$ in $L^2(X)$ as
$n \to \infty$.

        Observe that
\begin{eqnarray}
 f_{n + 1}(x) - f_n(x) & = & \frac{1}{n + 1} \sum_{j = 1}^{n + 1} \phi_j(x_j)
                               - \frac{1}{n} \sum_{j = 1}^n \phi_j(x_j) \\
 & = & \frac{\phi_{n + 1}(x_{n + 1})}{n + 1}
              - \frac{1}{n (n + 1)} \sum_{j = 1}^n \phi_j(x_j). \nonumber
\end{eqnarray}
If $a_n = E(f_{n + 1} - f_n \mid \mathcal{B}_n)$, as in the previous
section, then
\begin{equation}
         a_n(x) = -\frac{1}{n (n + 1)} \sum_{j = 1}^n \phi_j(x_j).
\end{equation}
This is because $\phi_j(x_j)$ is measurable with respect to
$\mathcal{B}_n$ when $j \le n$, while the conditional expectation of
$\phi_{n + 1}(x_{n + 1})$ with respect to $\mathcal{B}_n$ is equal to
$0$.  Thus $a_n = -(1/(n + 1)) \, f_n$,
\begin{equation}
        \|a_n\|_2 \le \frac{C}{\sqrt{n} (n + 1)},
\end{equation}
and so $\sum_{n = 1}^\infty \|a_n\|_2$ converges.  It follows that
$\{f_n\}_{n = 1}^\infty$ converges pointwise almost everywhere on $X$,
as in the previous section.

\section{Shift mappings}
\label{shift mappings}
\setcounter{equation}{0}

        Let $(X_0, \mathcal{A}_0, \mu_0)$, be a probability space, and
let $X$ be the space of doubly-infinite sequences $x = \{x_j\}_{j =
-\infty}^\infty$ with $x_j \in X_0$ for each $j$.  Thus $X$ is the
Cartesian product of a family of copies of $X_0$ indexed by the
integers, which is also a probability space with respect to the
product measure $\mu$.  Let $T$ be the shift mapping on $X$ defined in
Section \ref{doubly-infinite sequences}, which preserves the measure
$\mu$.  Also let $f$ be an integrable function on $X$, and consider
\begin{equation}
\label{frac{f(x) + f(T(x)) + f(T^2(x)) + cdots + f(T^n(x))}{n + 1}}
        \frac{f(x) + f(T(x)) + f(T^2(x)) + \cdots + f(T^n(x))}{n + 1}.
\end{equation}
If $f$ is constant, then (\ref{frac{f(x) + f(T(x)) + f(T^2(x)) + cdots
+ f(T^n(x))}{n + 1}}) is the same constant for each $n$.  Suppose
instead that the integral of $f$ is equal to $0$.  If $f$ is
square-integrable and depends only on one variable, then
(\ref{frac{f(x) + f(T(x)) + f(T^2(x)) + cdots + f(T^n(x))}{n + 1}})
converges to $0$ as $n \to \infty$ in the $L^2$ norm and pointwise
almoste everywhere on $X$, as in the previous section.  These are
consequences of well-known ergodic theorems as well.  One can also
deal with other $L^p$ spaces, but let us focus here on $p = 2$ for
simplicity.  If $f$ depends on only finitely many variables, then one
can get the same conclusions from analogous arguments.  More
precisely, one can begin with averages like (\ref{frac{f(x) + f(T(x))
+ f(T^2(x)) + cdots + f(T^n(x))}{n + 1}}), but using powers of $T^r$
for sufficiently large $r$ in place of powers of $T$.  An average like
(\ref{frac{f(x) + f(T(x)) + f(T^2(x)) + cdots + f(T^n(x))}{n + 1}})
with arbitrary powers of $T$ can then be estimated in terms of $r$
smaller averages involving $T^{j r + l}$, $l = 0, \ldots, r - 1$.
After that, an arbitrary function $f$ can be approximated by functions
depending on only finitely many variables.  There are also maximal
function estimates for the averages (\ref{frac{f(x) + f(T(x)) +
f(T^2(x)) + cdots + f(T^n(x))}{n + 1}}) like those that have been
discussed in other contexts.

\section{Families of $\sigma$-subalgebras}
\label{families of sigma-subalgebras}
\setcounter{equation}{0}

        Let $(X, \mathcal{A}, \mu)$ be a probability space, and let
$(\mathcal{I}, \prec)$ be a directed system.  Thus $\mathcal{I}$ is a
set, $\prec$ is a partial ordering on $\mathcal{I}$, and for each $a,
b \in \mathcal{I}$ there is a $c \in \mathcal{I}$ such that $a, b
\prec c$.  Suppose that for each $a \in \mathcal{I}$ we have a
$\sigma$-subalgebra $\mathcal{B}_a$ of $\mathcal{A}$, and that
\begin{equation}
\label{mathcal{B}_a subseteq mathcal{B}_b}
        \mathcal{B}_a \subseteq \mathcal{B}_b
\end{equation}
when $a, b \in \mathcal{I}$ and $a \prec b$.  If $\mathcal{I}$ is the
set ${\bf Z}_+$ of positive integers with the usual ordering, then
this is the same as an increasing sequence of $\sigma$-subalgebras of
$\mathcal{A}$, as in Section \ref{sequences of sigma-subalgebras}.

        Alternatively, let $I$ be a nonempty set, and let $(X_i,
\mathcal{A}_i, \mu_i)$ be a probability space for each $i$.  Consider
the Cartesian product $X = \prod_{i \in I} X_i$ of the $X_i$'s, with
the product measure $\mu$.  If $\mathcal{I}$ is the collection of
nonempty finite subsets of $I$, then $\mathcal{I}$ is partially
ordered by inclusion, and a directed system.  More precisely, if $a, b
\in \mathcal{I}$, then $a \cup b \in \mathcal{I}$, and $a, b \subseteq
a \cup b$.  Let $\mathcal{B}_a$ be the collection of subsets of $X$
that correspond to the Cartesian product of a measurable set $A
\subseteq \prod_{i \in a} X_i$ and $\prod_{i \in I \backslash a} X_i$
for each $a \in \mathcal{I}$.  It is easy to see that $\mathcal{B}_a$
is a $\sigma$-subalgebra of the $\sigma$-algebra of measurable subsets
of $X$, and that (\ref{mathcal{B}_a subseteq mathcal{B}_b}) holds.  If
the $X_i$'s are compact Hausdorff topological spaces, so that $X$ is
also a compact Hausdorff space with respect to the product topology,
then one may wish to use Borel sets.

        In this product situation, suppose that $\phi_i \in L^1(X_i,
\mathcal{A}_i)$ satisfies
\begin{equation}
\label{int_{X_i} phi_i d mu_i = 0}
        \int_{X_i} \phi_i \, d\mu_i = 0
\end{equation}
for each $i \in I$.  Put
\begin{equation}
        \Phi_a(x) = \sum_{i \in a} \phi_i(x_i)
\end{equation}
for each $a \in \mathcal{I}$, where $x = \{x_i\}_{i \in I} \in X$.
Thus $\Phi_a \in L^1(X, \mathcal{B}_a)$, and
\begin{equation}
        E(\Phi_b \mid \mathcal{B}_a) = \Phi_a
\end{equation}
when $a, b \in \mathcal{I}$ and $a \subseteq b$.  Hence $\Phi_a$, $a
\in \mathcal{I}$, defines a martingale with respect to this family of
$\sigma$-algebras.

        Martingales with more general indices like this are discussed
in \cite{ks2}.  This point of view is very natural in connection with
rearrangement of sums, and convergence of sums in the generalized
sense.  Note that the arguments for estimating maximal functions as in
Section \ref{maximal functions, 3} do not work for partially-ordered
sets of indices.  The corresponding problems with pointwise
convergence have already been seen at least implicitly in Section
\ref{rademacher sums}, in the case where $I = {\bf Z}_+$, $X_i = \{1,
-1\}$, and $\mu_i(\{1\}) = \mu_i(\{-1\}) = 1/2$ for each $i \in I$.
However, if the $\sigma$-algebras $\mathcal{B}_a$ are associated to
partitions consisting of intervals in the real line, then one can use
a covering argument as in Section \ref{maximal functions}.

\section{Stopping times}
\label{stopping times}
\setcounter{equation}{0}

        Let $(X, \mathcal{A}, \mu)$ be a probability space, and let
$\mathcal{B}_1 \subseteq \mathcal{B}_2 \subseteq \cdots$ be an
increasing sequence of $\sigma$-subalgebras of $\mathcal{A}$.  A
function $\tau : X \to {\bf Z}_+$ is said to be a \emph{stopping time}
if
\begin{equation}
        \tau^{-1}(n) = \{x \in X : \tau(x) = n\} \in \mathcal{B}_n
\end{equation}
for each $n \ge 1$.  This is equivalent to the condition that
\begin{equation}
 \tau^{-1}(\{1, \ldots, n\}) = \{x \in X : \tau(x) \le n\} \in \mathcal{B}_n
\end{equation}
for each $n$, since
\begin{equation}
        \tau^{-1}(\{1, \ldots, n\}) = \bigcup_{l = 1}^n \tau^{-1}(l)
\end{equation}
and
\begin{equation}
        \tau^{-1}(n) =
         \tau^{-1}(\{1, \ldots, n\}) \backslash \tau^{-1}(\{1, \ldots, n-1\})
\end{equation}
when $n \ge 2$.  Alternatively, $\tau$ is a stopping time if
\begin{equation}
        \{x \in X : \tau(x) > n\} \in \mathcal{B}_n
\end{equation}
for each $n$, because
\begin{equation}
        \{x \in X : \tau(x) > n\} = X \backslash \tau^{-1}(\{1, \ldots, n\}).
\end{equation}
One can also allow $\tau$ to take values in ${\bf Z}_+ \cup \{+\infty\}$,
in which case
\begin{equation}
\tau^{-1}(+\infty) = X \backslash \Big(\bigcup_{n = 1}^\infty \tau^{-1}(n)\Big)
\end{equation}
is in the $\sigma$-algebra $\mathcal{B}_\infty$ generated by
$\bigcup_{n = 1}^\infty \mathcal{B}_n$.

        If $A_1, A_2, \ldots$ is a sequence of pairwise-disjoint
subsets of $X$ with $A_n \in \mathcal{B}_n$ for each $n$, then there
is a unique stopping time $\tau$ on $X$ such that $\tau^{-1}(n) = A_n$
for each $n$.  More precisely, $\tau(x) < +\infty$ for every $x \in X$
if and only if $\bigcup_{n = 1}^\infty A_n = X$.  Similarly, if $E_1
\subseteq E_2 \subseteq \cdots$ is an increasing sequence of subsets
of $X$ with $E_n \in \mathcal{B}_n$ for each $n$, then there is a
unique stopping time $\tau$ on $X$ such that
\begin{equation}
\label{{x in X : tau(x) le n} = E_n}
        \{x \in X : \tau(x) \le n\} = E_n
\end{equation}
for each $n$.  Of course, this corresponds to taking $A_1 = E_1$ and
$A_n = E_n \backslash E_{n - 1}$ when $n \ge 2$ in the previous
statement.  As before, $\tau(x) < +\infty$ for every $x \in X$ if and
only if $\bigcup_{n = 1}^\infty E_n = X$.

        If $\tau$, $\tau'$ are stopping times on $X$, then $\max(\tau,
\tau')$ and $\min(\tau, \tau')$ are stopping times too, because
\begin{eqnarray}
\lefteqn{\{x \in X : \max(\tau(x), \tau'(x)) \le n\} =}  \\
 & & \{x \in X : \tau(x) \le n\} \cap \{x \in X : \tau'(x) \le n\} \nonumber
\end{eqnarray}
and
\begin{eqnarray}
\lefteqn{\{x \in X : \min(\tau(x), \tau'(x)) \le n\} =} \\
 & & \{x \in X : \tau(x) \le n\} \cup \{x \in X : \tau'(x) \le n\}. \nonumber
\end{eqnarray}
In particular,
\begin{equation}
\label{tau_N(x) = min(tau(x), N)}
        \tau_N(x) = \min(\tau(x), N)
\end{equation}
is a stopping time on $X$ when $\tau$ is a stopping time and $N$ is a
positive integer.

        Suppose that $\{f_n\}_{n = 1}^\infty$ is a martingale on $X$
with respect to this filtration, and let
\begin{equation}
\label{f^*(x) = sup_{n ge 1} |f_n(x)|, stopping times}
        f^*(x) = \sup_{n \ge 1} |f_n(x)|
\end{equation}
be the corresponding maximal function.  Let $t > 0$ be given, and
remember that $f^*(x) > t$ if and only if $|f_n(x)| > t$ for some $n$.
If $f^*(x) > t$, then let $\tau(x)$ be the smallest positive integer such that
\begin{equation}
        |f_{\tau(x)}(x)| > t,
\end{equation}
and put $\tau(x) = +\infty$ when $f^*(x) \le t$.  Thus $\tau(x) = n$
exactly when $|f_n(x)| > t$ and $|f_l(x)| \le t$ for $l < n$.  This
implies that $\tau^{-1}(n) \in \mathcal{B}_n$ for each $n$, because
$f_l$ is measurable with respect to $\mathcal{B}_l \subseteq
\mathcal{B}_n$ when $l \le n$.

        Let $\tau$ be a stopping time on $X$ such that $\tau(x) <
+\infty$ for every $x \in X$, and let $\mathcal{B}_\tau$ be the
collection of subsets $A$ of $X$ such that
\begin{equation}
        A \cap \tau^{-1}(n) \in \mathcal{B}_n
\end{equation}
for each $n$.  It is easy to see that this is a $\sigma$-algebra,
because $\mathcal{B}_n$ is a $\sigma$-algebra for each $n$, and that
$\mathcal{B}_\tau \subseteq \mathcal{B}_\infty$.  If $N$ is a positive
integer and $\tau(x) \le N$ for each $x \in X$, then
\begin{equation}
        \mathcal{B}_\tau \subseteq \mathcal{B}_N.
\end{equation}
More precisely, if $\tau$ is any finite stopping time, $A \in
\mathcal{B}_\tau$, and $\tau(x) \le N$ for every $x \in A$, then $A
\in \mathcal{B}_N$.  If $\tau'$ is another stopping time such that
\begin{equation}
        \tau(x) \le \tau'(x) < +\infty
\end{equation}
for every $x \in X$, then $\mathcal{B}_\tau \subseteq \mathcal{B}_{\tau'}$.

        Let $\{f_n\}_{n = 1}^\infty$ be a martingale on $X$ with
respect to this filtration, and let $\tau$ be a finite stopping time
on $X$.  If $f_\tau$ is the function on $X$ defined by
\begin{equation}
        f_\tau(x) = f_{\tau(x)}(x),
\end{equation}
then $f_\tau$ is measurable with respect to $\mathcal{B}_\tau$,
because $f_n$ is measurable with respect to $\mathcal{B}_n$ for each
$n$.  Let us check that
\begin{equation}
 \int_{\tau^{-1}(\{1, \ldots, N\})} |f_\tau| \, d\mu \le \int_X |f_N| \, d\mu
\end{equation}
for each positive integer $N$.  By the definition of $f_\tau$,
\begin{equation}
        \int_{\tau^{-1}(\{1, \ldots, N\})} |f_\tau| \, d\mu
         = \sum_{n = 1}^N \int_{\tau^{-1}(n)} |f_n| \, d\mu.
\end{equation}
Hence
\begin{eqnarray}
        \int_{\tau^{-1}(\{1, \ldots, N\})} |f_\tau| \, d\mu
         & \le & \sum_{n = 1}^N \int_{\tau^{-1}(n)} |f_N| \, d\mu \\
          & = & \int_{\tau^{-1}(\{1, \ldots, N\})} |f_N| \, d\mu, \nonumber
\end{eqnarray}
because $|f_n| \le E(|f_N| \mid \mathcal{B}_n)$ when $n \le N$.

        If $\tau(x) \le N$ for every $x \in X$, then it follows that
$f_\tau$ is integrable on $X$.  Let us verify that
\begin{equation}
        f_\tau = E(f_N \mid \mathcal{B}_\tau),
\end{equation}
remembering that $\mathcal{B}_\tau \subseteq \mathcal{B}_N$ in this case.
To see this, it suffices to show that
\begin{equation}
        \int_A f_\tau \, d\mu = \int_A f_N \, d\mu
\end{equation}
for every $A \in \mathcal{B}_\tau$.  Under these conditions,
\begin{eqnarray}
        \int_A f_\tau \, d\mu
            & = & \sum_{n = 1}^N \int_{A \cap \tau^{-1}(n)} f_n \, d\mu \\
            & = & \sum_{n = 1}^N \int_{A \cap \tau^{-1}(n)} f_N \, d\mu
              = \int_A f_N \, d\mu, \nonumber
\end{eqnarray}
because $A \cap \tau^{-1}(n) \in \mathcal{B}_n$ and $f_n = E(f_N \mid
\mathcal{B}_n)$ when $n \le N$.

        Similarly, if the $f_n$'s have bounded $L^1$ norms and $\tau$
is any finite stopping time on $X$, then we get that
\begin{equation}
        \int_{\tau^{-1}(\{1, \ldots, N\})} |f_\tau| \, d\mu
         \le \sup_{n \ge 1} \int_X |f_n| \, d\mu
\end{equation}
for every positive integer $N$.  This implies that $f_\tau$ is
integrable on $X$, and that
\begin{equation}
        \int_X |f_\tau| \, d\mu \le \sup_{n \ge 1} \int_X |f_n| \, d\mu.
\end{equation}
In particular, this holds when there is an $f \in L^1(X, \mathcal{A})$
such that $f_n = E(f \mid \mathcal{B}_n)$ for each $n$.  In this case,
one can check that
\begin{equation}
        f_\tau = E(f \mid \mathcal{B}_\tau).
\end{equation}
As before, one can show that
\begin{equation}
        \int_A f_\tau \, d\mu = \int_A f \, d\mu
\end{equation}
when $A \in \mathcal{B}_\tau$, by expressing $A$ as the union of $A
\cap \tau^{-1}(n)$, $n \ge 1$, and using the fact that $f_\tau = f_n =
E(f \mid \mathcal{B}_n)$ on $A \cap \tau^{-1}(n) \in \mathcal{B}_n$.

        Now let $\tau$ be a stopping time on $X$ that takes values in
${\bf Z}_+ \cup \{+\infty\}$, so that $\tau_N = \min(\tau, N)$ is a
finite stopping time on $X$ for each $N$.  Let $\{f_n\}_{n =
1}^\infty$ be a martingale on $X$ with respect to this filtration, and
note that $f_{\tau_N}$ is integrable on $X$ for each $N$, since
$\tau_N$ is bounded.  Let us check that
\begin{equation}
        f_{\tau_N} = E(f_{\tau_{N + 1}} \mid \mathcal{B}_N)
\end{equation}
for each $N$, so that $\{f_{\tau_N}\}_{N = 1}^\infty$ is a martingale as well.
As usual, we would like to show that
\begin{equation}
        \int_A f_{\tau_N} \, d\mu = \int_A f_{\tau_{N + 1}} \, d\mu
\end{equation}
when $A \in \mathcal{B}_N$.  Consider
\begin{equation}
        A_1 = \{x \in A : \tau(x) \le N\}
\end{equation}
and
\begin{equation}
        A_2 = \{x \in A : \tau > N\}.
\end{equation}
Thus $A_1 \cup A_2 = A$, $A_1 \cap A_2 = \emptyset$, and
\begin{equation}
        A_1, A_2 \in \mathcal{B}_N,
\end{equation}
since $\tau$ is a stopping time.  If $x \in A_1$, then $\tau_N(x) =
\tau_{N + 1}(x) = \tau(x)$, and hence $f_{\tau_N}(x) = f_{\tau_{N +
1}}(x) = f_\tau(x)$.  This implies that
\begin{equation}
        \int_{A_1} f_{\tau_N} \, d\mu = \int_{A_1} f_{\tau_{N + 1}} \, d\mu.
\end{equation}
If $x \in A_2$, then $\tau_N(x) = N$, $\tau_{N + 1}(x) = N + 1$, and
so $f_{\tau_N}(x) = f_N(x)$, $f_{\tau_{N + 1}}(x) = f_{N + 1}(x)$.
It follows that
\begin{equation}
        \int_{A_2} f_{\tau_N} \, d\mu = \int_{A_2} f_N \, d\mu
 = \int_{A_2} f_{N + 1} \, d\mu = \int_{A_2} f_{\tau_{N + 1}} \, d\mu,
\end{equation}
because $A_2 \in \mathcal{B}_N$ and $f_N = E(f_{N + 1} \mid
\mathcal{B}_N)$, as desired.

        If $\tau(x) < \infty$ for every $x \in X$, then
$\{f_{\tau_N}\}_{N = 1}^\infty$ converges to $f_\tau$ pointwise on
$X$, because $f_{\tau_N}(x) = f_\tau(x)$ when $N \ge \tau(x)$.  If the
$f_n$'s have bounded $L^1$ norms, then the $f_{\tau_N}$'s also have
bounded $L^1$ norms, and $f_\tau$ is integrable.  A necessary and
sufficient condition for $\{f_{\tau_N}\}_{N = 1}^\infty$ to converge
to $f_\tau$ in the $L^1$ norm is that
\begin{equation}
\label{int_{{x in X : tau(x) > N}} to 0 as N to infty}
        \int_{\{x \in X : \tau(x) > N\}} |f_{\tau_N}(x)| \, d\mu(x)
         = \int_{\{x \in X : \tau(x) > N\}} |f_N(x)| \, d\mu(x) \to 0
\end{equation}
as $N \to \infty$.  This holds automatically when the $f_n$'s are
uniformly integrable, and otherwise depends on both the $f_n$'s and
$\tau$.

\section{Ultrametrics}
\label{ultrametrics}
\setcounter{equation}{0}

        A metric $d(x, y)$ on a set $M$ is said to be an \emph{ultrametric} if
\begin{equation}
        d(x, z) \le \max(d(x, y), d(y, z))
\end{equation}
for every $x, y, z \in M$.  If $X_1, X_2, \ldots$ is a sequence of
nonempty sets, and $r_1, r_2, \ldots$ is a decreasing sequence of
positive real numbers that converges to $0$, then one can define an
ultrametric on the Cartesian product $X = \prod_{j = 1}^\infty X_j$ as
follows.  Each element $x$ of $X$ is a sequence $\{x_j\}_{j =
1}^\infty$ with $x_j \in X_j$ for every $j$, and we put $d(x, x) = 0$, and
\begin{equation}
        d(x, y) = r_l
\end{equation}
when $x \ne y$ and $l$ is the smallest positive integer such that $x_l
\ne y_l$.  It is easy to see that this is an ultrametric on $X$, and
that the corresponding topology is the product topology associated to
the discrete topology on $X_j$ for each $j$.

        If $d(x, y)$ is an ultrametric on a set $M$, $p, q \in M$, and
$t \ge r > 0$, then either
\begin{equation}
B(p, r) \subseteq B(q, t) \quad\hbox{or}\quad B(p, r) \cap B(q, t) = \emptyset.
\end{equation}
More precisely, the first alternative holds when $d(p, q) < t$, and
the second alternative holds when $d(p, q) \ge t$.  Using this, one
can check that open balls are closed subsets of ultrametric spaces.
There is an analogous dichotomy for closed balls, which implies that
closed balls are open subsets of ultrametric spaces.  It follows that
ultrametric spaces are totally disconnected, in the sense that they do
not contain connected subsets with more than one element.

        Another consequence of the previous dichotomy is that
\begin{equation}
        B(p, r) = B(q, r)
\end{equation}
when $d(p, q) < r$.  Thus every element of an open ball in $M$ can be
used as a center of that ball.  The collection of open balls in $M$
with the same radius $r$ forms a partition of $M$, because any two
such balls are either the same or disjoint as subsets of $M$.  If $t
\ge r$, then the partition of $M$ into open balls of radius $r$ is a
refinement of the partition of $M$ into open balls of radius $t$,
since every ball of radius $r$ is contained in a ball of radius $t$.

        The geometry of an ultrametric space is very similar to a
probability space with an increasing sequence of $\sigma$-subalgebras
of the $\sigma$-algebra of measurable sets.  In particular, one can
consider $\sigma$-subalgebras of the Borel sets in an ultrametric
space corresponding to partitions by balls of a given radius.  One can
also deal directly with Hardy--Littlewood type maximal functions,
using the nesting properties of balls to reduce of covering of a set
by balls of bounded radius to a disjoint union of balls that are
maximal elements of the covering.  Of course, there are more
complicated covering arguments for Euclidean spaces and other metric
spaces, including the basic property of intervals in the real line
mentioned in Section \ref{maximal functions}.  These can also be used
to estimate maximal functions, and so on.

\part{Vector-valued functions}

\section{Some randomized sums}
\label{randomized sums}
\setcounter{equation}{0}

        Let $(X, \mathcal{A}, \mu)$ be a probability space, and let
$\phi_1, \ldots, \phi_n$ be bounded real or complex-valued measurable
functions on $X$, with
\begin{equation}
        \|\phi_j\|_\infty \le C
\end{equation}
for some $C \ge 0$ and $j = 1, \ldots, n$.  Also let $\{1, -1\}^n$ be
the set of sequences $\epsilon = \{\epsilon_j\}_{j = 1}^n$ of length
$n$ with $\epsilon_j = 1$ or $-1$ for each $j$.  If $2 \le p <
\infty$, then there is a positive real number $C(p)$ such that
\begin{equation}
\label{2^{-n} sum_{epsilon in {1, -1}^n} int_X ...}
 2^{-n} \sum_{\epsilon \in \{1, -1\}^n} \int_X \biggl|\sum_{j = 1}^n \epsilon_j
                                        \, a_j \, \phi_j(x)\biggr|^p \, d\mu(x)
     \le C(p) \, \Big(\sum_{j = 1}^n |a_j|^2\Big)^{p/2}
\end{equation}
for all $a_1, \ldots, a_n \in {\bf R}$ or ${\bf C}$, as appropriate.
Of course, the left side is the same as
\begin{equation}
 \int_X 2^{-n} \sum_{\epsilon \in \{1, -1\}^n} \biggl|\sum_{j = 1}^n \epsilon_j
                                       \, a_j \, \phi_j(x)\biggr|^p \, d\mu(x).
\end{equation}
As in Section \ref{L^p estimates},
\begin{equation}
\label{2^{-n} sum_{epsilon in {1, -1}^n} ...}
  2^{-n} \sum_{\epsilon \in \{1, -1\}^n} \biggl|\sum_{j = 1}^n \epsilon_j
                                                \, a_j \, \phi_j(x)\biggr|^p
            \le C'(p) \, \Big(\sum_{j = 1}^n |a_j \, \phi_j(x)|^2 \Big)^{p/2}
\end{equation}
for some $C'(p) > 0$ and all $a_1, \ldots, a_n \in {\bf R}$ or ${\bf
C}$ and $x \in X$.  This implies (\ref{2^{-n} sum_{epsilon in {1,
-1}^n} int_X ...}), by integrating in $x$ and using the uniform
boundedness of the $\phi_j$'s.  More precisely, $C(p)$ depends only on
$C$ and $p$, and not on $a_1, \ldots, a_n$ or $n$.

        If $p = 2$, then we have that
\begin{equation}
 2^{-n} \sum_{\epsilon \in \{1, -1\}^n} \biggl|\sum_{j = 1}^n \epsilon_j
  \, a_j \, \phi_j(x)\biggr|^2 = \sum_{j = 1}^n |a_j \, \phi_j(x)|^2.
\end{equation}
This implies that
\begin{equation}
\label{2^{-n} sum_{epsilon in {1, -1}^n} int_X ... = sum_{j = 1}^n |a_j|^2}
 2^{-n} \sum_{\epsilon \in \{1, -1\}^n} \int_X \biggl|\sum_{j = 1}^n \epsilon_j
  \, a_j \, \phi_j\biggr|^2 \, d\mu(x) = \sum_{j = 1}^n |a_j|^2
\end{equation}
when $\|\phi_j\|_2 = 1$ for each $j$.  Otherwise, if $\|\phi_j\|_2 \ge c$
for some $c > 0$ and each $j$, then we get that
\begin{equation}
 2^{-n} \sum_{\epsilon \in \{1, -1\}^n} \int_X \biggl|\sum_{j = 1}^n \epsilon_j
                                           \, a_j \, \phi_j\biggr|^2 \, d\mu(x)
                                            \ge c^2 \, \sum_{j = 1}^n |a_j|^2.
\end{equation}
Note that
\begin{equation}
\label{(2^{-n} sum_{epsilon in {1, -1}} int_X ...)^{1/p}}
 \Big(2^{-n} \sum_{\epsilon \in \{1, -1\}} \int_X \biggl|\sum_{j = 1}^n
                  \epsilon_j \, a_j \, \phi_j(x)\biggr|^p \, d\mu(x)\Big)^{1/p}
\end{equation}
is monotone increasing in $p$, by Jensen's inequality.  This is the
same as the $L^p$ norm of $\sum_{j = 1}^n \epsilon_j \, a_j \,
\phi_j(x)$ as a function of $(x, \epsilon)$ on $X \times \{1, -1\}^n$,
with respect to the product of $\mu$ on $X$ and $2^{-n}$ times
counting measure on $\{1, -1\}^n$.

        Under these conditions, if $0 < p < 2$, then there is a $C(p)
> 0$ such that
\begin{equation}
\label{C(p)^{-1} (sum_{j = 1}^n |a_j|^2)^{p/2} le ...}
        C(p)^{-1} \, \Big(\sum_{j = 1}^n |a_j|^2\Big)^{p/2} \le
 2^{-n} \sum_{\epsilon \in \{1, -1\}^n} \int_X \biggl|\sum_{j = 1}^n \epsilon_j
                                           \, a_j \, \phi_j\biggr|^p \, d\mu(x)
\end{equation}
for all $a_1, \ldots, a_n \in {\bf R}$ or ${\bf C}$.  This can be
derived from the previous estimates and H\"older's inequality, as in
Section \ref{L^p estimates}.  More precisely, H\"older's inequality
can be used to estimate the $L^2$ norm of $\sum_{j = 1}^n \epsilon_j
\, a_j \, \phi_j(x)$ on $X \times \{1, -1\}^n$ in terms of its $L^p$
and $L^4$ norms, as before.  Under the present conditions, the $L^2$
norm is bounded from below by a constant multiple of $\Big(\sum_{j =
1}^n |a_j|^2\Big)^{1/2}$, and the $L^4$ norm is bounded from above by
a multiple of the same expression, which leads to a lower bound for
the $L^p$ norm as in (\ref{C(p)^{-1} (sum_{j = 1}^n |a_j|^2)^{p/2} le
...}).  As usual, the constant $C(p)$ in (\ref{C(p)^{-1} (sum_{j =
1}^n |a_j|^2)^{p/2} le ...}) depends on $c$, $C$, and $p$, and not on
$a_1, \ldots, a_n$ or $n$.

\section{Randomized sums, 2}
\label{randomized sums, 2}
\setcounter{equation}{0}

        Let $(X, \mathcal{A}, \mu)$ be a probability space again, and
let $\phi_1, \ldots, \phi_n$ be orthonormal functions in $L^2(X)$.
As usual, this implies that
\begin{equation}
        \int_X \biggl|\sum_{j = 1}^n \alpha_j \, \phi_j(x)\biggr|^2 \, d\mu(x)
         = \sum_{j = 1}^n |\alpha_j|^2
\end{equation}
for all $\alpha_1, \ldots, \alpha_n \in {\bf R}$ or ${\bf C}$, as
appropriate.  Hence
\begin{equation}
\label{int_X ... = sum_{j = 1}^n |alpha_j|^2}
 \int_X \biggl|\sum_{j = 1}^n \epsilon_j \, \alpha_j \, \phi_j(x)\biggr|^2 \,
  d\mu(x) = \sum_{j = 1}^n |\alpha_j|^2
\end{equation}
for every $\epsilon \in \{1, -1\}^n$.  In particular, the average of
the left side of (\ref{int_X ... = sum_{j = 1}^n |alpha_j|^2}) over
$\epsilon \in \{1, -1\}^n$ has the same value, as in (\ref{2^{-n}
sum_{epsilon in {1, -1}^n} int_X ... = sum_{j = 1}^n |a_j|^2}).

        Suppose that $\phi_1, \phi_2, \ldots$ is an orthonormal basis
for $L^2(X)$, and that the $\phi_j$'s are uniformly bounded on $X$, as
in the previous section.  Thus every function in $L^2(X)$ can be
approximated in the $L^2$ norm by a finite sum of the form
\begin{equation}
        \sum_{j = 1}^n \alpha_j \, \phi_j(x).
\end{equation}
Moreover, the average of the $L^p$ norms of
\begin{equation}
\label{sum_{j = 1}^n epsilon_j alpha_j phi_j(x)}
        \sum_{j = 1}^n \epsilon_j \, \alpha_j \, \phi_j(x)
\end{equation}
over $\epsilon \in \{1, -1\}^n$ is bounded by a constant multiple of
the $L^2$ norm for every $p < \infty$, as before.  However, this does
not mean that the $L^p$ norm of (\ref{sum_{j = 1}^n epsilon_j alpha_j
phi_j(x)}) is bounded by a multiple of the $L^2$ norm for every
$\epsilon \in \{1, -1\}$, or even for only $\epsilon = (1, \ldots,
1)$.  If we start with a function in $L^2(X)$ which is not in $L^p(X)$
for some $p > 2$, then the $L^p$ norms of its approximations are
necessarily unbounded.  Note that Fourier series and Walsh functions
are examples of this type of situation.  Lacunary series and
Rademacher functions correspond to subsets of these bases for which
the $L^p$ norms are bounded by constant multiples of the $L^2$ norms
when $2 < p < \infty$.

\section[\ The unit square]{The unit square}
\label{unit square}
\setcounter{equation}{0}

        Let $X = [0, 1) \times [0, 1)$ be the version of the unit
square associated to dyadic intervals, equipped with $2$-dimensional
Lebesgue measure.  If $I, L \subseteq [0, 1)$ are dyadic intervals
with the same length $2^{-j}$, then their Cartesian product $I \times
L$ is a dyadic square in $X$ with side length $2^{-j}$ and area $2^{-2
j}$.  There are $2^{2 j}$ dyadic squares in $X$ with side length
$2^{-j}$, they are pairwise disjoint, and their union is equal to $X$.
Let $\mathcal{A}_j$ be the collection of subsets of $X$ which can be
expressed as unions of dyadic squares with side length $2^{-j}$,
including the empty set.  This is the same as the $\sigma$-algebra of
subsets of $X$ generated by the partition $\mathcal{P}_j$ of $X$ into
dyadic squares of side length $2^{-j}$, as in Section
\ref{partitions}.  Note that $\mathcal{A}_j$ is a $\sigma$-subalgebra
of the $\sigma$ algebra of Borel subsets of $X$, and that
$\mathcal{A}_j \subseteq \mathcal{A}_{j + 1}$ for each $j$.  As usual,
a function on $X$ is measurable with respect to $\mathcal{A}_j$ if and
only if it is constant on dyadic squares with side length $2^{-j}$.

        Let $f_j(x, y)$ be the function on $X$ defined by
\begin{equation}
        f_j(x, y) = 2^j
\end{equation}
when $x$, $y$ are contained in the same dyadic interval of length $2^{-j}$, and
\begin{equation}
        f_j(x, y) = 0
\end{equation}
when $x$, $y$ are contained in distinct dyadic intervals of length $2^{-j}$.
In particular,
\begin{equation}
        \int_{I \times I} f_j(x, y) \, dx dy = 2^{-j}
\end{equation}
for each dyadic interval $I$ of length $2^{-j}$.  Summing over $I$, we
get that
\begin{equation}
        \int_{[0, 1) \times [0, 1)} f_j(x, y) \, dx dy = 1
\end{equation}
for each $j$, because there are $2^j$ dyadic intervals of length
$2^{-j}$.  Clearly $f_j(x, y)$ is measurable with respect to
$\mathcal{A}_j$ for each $j$.  It is easy to see that
\begin{equation}
\label{f_j = E(f_{j + 1} mid mathcal{A}_j)}
        f_j = E(f_{j + 1} \mid \mathcal{A}_j)
\end{equation}
for each $j$, so that $\{f_j\}_j$ is a martingale with respect to the
$\mathcal{A}_j$'s.

        Let $\nu$ be the Borel measure on $X$ defined by
\begin{equation}
        \nu(A) = |\{x \in [0, 1) : (x, x) \in A\}|,
\end{equation}
where $|E|$ denotes the Lebesgue measure of $E \subseteq [0, 1)$.
Alternatively, if
\begin{equation}
        \Delta = \{(x, x) : x \in [0, 1)\}
\end{equation}
is the diagonal in $X$, then
\begin{equation}
        \nu(A) = |\pi(A \cap \Delta)|,
\end{equation}
where $\pi(x, x) = x$ is the natural projection of $\Delta$ onto $[0, 1)$.
Of course, the restriction of $\nu$ to $\mathcal{A}_j$ is absolutely
continuous with respect to the restriction of $2$-dimensional Lebesgue
measure to $\mathcal{A}_j$ for each $j$.  One can also think of $f_j$
as the conditional expectation of $\nu$ with respect to $\mathcal{A}_j$,
as in Section \ref{other measures}.

        If $x, y \in [0, 1)$ and $x \ne y$, then $f_j(x, y) = 0$ for
all sufficiently large $j$.  In particular, $\{f_j(x, y)\}_j$
converges to $0$ almost everywhere on $X$.  Basically, $\{f_j\}_j$
converges to $\nu$ in a suitable weak sense.

        Now let $\mathcal{B}_j$ be the collection of subsets of $X$
that can be expressed as the union of sets of the form $I \times
A(I)$, where $I$ runs through the dyadic subintervals of $[0, 1)$ of
length $2^{-j}$, and $A(I)$ is a Borel set in $[0, 1)$ for each such
$I$.  Equivalently, $A \in \mathcal{B}_j$ if for each dyadic interval
$I \subseteq [0, 1)$ with $|I| = 2^{-j}$ there is a Borel set $A(I)
\subseteq [0, 1)$ such that
\begin{equation}
\label{A cap (I times [0, 1)) = I times A(I)}
        A \cap (I \times [0, 1)) = I \times A(I).
\end{equation}
Thus $\mathcal{B}_j$ is a $\sigma$-subalgebra of the $\sigma$-algebra
of Borel sets in $X$, $\mathcal{A}_j \subseteq \mathcal{B}_j$, and
$\mathcal{B}_j \subseteq \mathcal{B}_{j + 1}$ for each $j$.  A
function $f(x, y)$ on $X$ is measurable with respect to
$\mathcal{B}_j$ if and only if it is constant in $x$ on each dyadic
interval $I$ of length $2^{-j}$ and Borel measurable in $y$.

        In particular, $f_j(x, y)$ is measurable with respect to
$\mathcal{B}_j$ for each $j$.  One can also check that
\begin{equation}
\label{E(f_{j + 1} mid mathcal{B}_j) = f_j}
        E(f_{j + 1} \mid \mathcal{B}_j) = f_j
\end{equation}
for each $j$, so that $\{f_j\}_j$ is a martingale with respect to the
$\mathcal{B}_j$'s as well.  The main point is that
\begin{equation}
\label{int_{I times A} f_j(x, y) dx dy = |A cap I|}
        \int_{I \times A} f_j(x, y) \, dx dy = |A \cap I|
\end{equation}
for each dyadic interval $I$ of length $2^{-j}$ and Borel set $A
\subseteq [0, 1)$.  Similarly, if $I_1$, $I_2$ are the dyadic
intervals of length $2^{-j - 1}$ such that $I = I_1 \cup I_2$, then
\begin{eqnarray}
\lefteqn{\int_{I \times A} f_{j + 1}(x, y) \, dx dy} \\
          & = & \int_{I_1 \times A} f_{j + 1}(x, y) \, dx dy +
            \int_{I_2 \times A} f_{j + 1}(x, y) \, dx dy \nonumber \\
          & = & |A \cap I_1| + |A \cap I_2| = |A \cap I|. \nonumber
\end{eqnarray}
This implies (\ref{E(f_{j + 1} mid mathcal{B}_j) = f_j}), which can
also be seen by viewing $f_j$ as the conditional expectation of $\nu$
with respect to $\mathcal{B}_j$, by (\ref{int_{I times A} f_j(x, y) dx
dy = |A cap I|}).

        Let $\Phi_j$ be the function on $[0, 1)$ with values in
$L^1([0, 1))$ defined by
\begin{equation}
        \Phi_j(x)(y) = f_j(x, y).
\end{equation}
This may be considered as a martingale on $[0, 1)$ with values in
$L^1([0, 1))$, with respect to the usual filtration associated to
dyadic intervals of length $2^{-j}$.  Note that the $L^1$ norm of
$\Phi_j(x)$ is equal to $1$ for each $x$ and $j$, but
$\{\Phi_j(x)\}_j$ does not converge in $L^1([0, 1))$ for any $x \in
[0, 1)$.  If we identify integrable functions on $[0, 1)$ with
absolutely continuous Borel measures on $[0, 1]$, which determine
bounded linear functionals on the space of continuous functions on
$[0, 1]$ with respect to the supremum norm, then $\{\Phi_j(x)\}_j$
converges in the weak$^*$ topology to the Dirac mass at $x$.

\section[\ Partitions and products]{Partitions and products}
\label{partitions, products}
\setcounter{equation}{0}

        Let $(X_1, \mathcal{A}_1, \mu_1)$, $(X_2, \mathcal{A}_2,
\mu_2)$ be probability spaces, and let $X = X_1 \times X_2$ be their
Cartesian product, with the product probability measure $\mu = \mu_1
\times \mu_2$.  Suppose that $\mathcal{P}_1$, $\mathcal{P}_2$ are
partitions of $X_1$, $X_2$ into finitely or countably many measurable
sets, respectively, as in Section \ref{partitions}.  The corresponding
product partition $\mathcal{P}_{1, 2}$ of $X$ consists of all products
$A_1 \times A_2$, with $A_1 \in \mathcal{P}_1$ and $A_2 \in
\mathcal{P}_2$.  It is easy to see that this is a partition of $X$
into finitely or countably many measurable sets, and that the
$\sigma$-algebra generated by $\mathcal{P}_{1, 2}$ is the same as the
one associated to the $\sigma$-algebras generated by $\mathcal{P}_1$,
$\mathcal{P}_2$ in the product space.  A function $f(x_1, x_2)$ on $X$
is measurable with respect to this $\sigma$-algebra if and only if it
is constant on $A_1 \times A_2$ for each $A_1 \in \mathcal{P}_1$ and
$A_2 \in \mathcal{P}_2$.

        Now let $\mathcal{P}_1$ be a partition of $X_1$ into finitely
or countably many measurable sets, and let $\mathcal{B}_2$ be a
$\sigma$-subalgebra of $\mathcal{A}_2$.  This leads to a
$\sigma$-subalgebra $\mathcal{B}_{1, 2}$ of the $\sigma$-algebra of
measurable subsets of $X$ associated to the $\sigma$-algebra generated
by $\mathcal{P}_1$ and $\mathcal{B}_2$ in the product space.  As in
the special case described in the previous section, $\mathcal{B}_{1,
2}$ consists of the sets $A \subseteq X$ such that for each $A_1 \in
\mathcal{P}_1$ there is an $A_2 \in \mathcal{B}_2$ such that
\begin{equation}
        A \cap (A_1 \times X_2) = A_1 \times A_2.
\end{equation}
Equivalently, $A \in \mathcal{B}_{1, 2}$ if $A$ can be expressed as a
union of sets of the form $A_1 \times A_2$, where $A_1$ runs through
the elements of $\mathcal{P}_1$, and $A_2 \in \mathcal{B}_2$ for each
$A_1 \in \mathcal{P}_1$.  Thus a function $f(x_1, x_2)$ on $X$ is
measurable with respect to $\mathcal{B}_{1, 2}$ if it is constant in
$x_1$ on each $A_1 \in \mathcal{P}_1$, and measurable in $x_2$ with
respect to $\mathcal{B}_2$ for each $x_1 \in X_1$.

        As in Section \ref{partitions}, it will be convenient to ask
that $\mu_1(A_1) > 0$ for each $A_1 \in \mathcal{P}_1$.  If
$\mathcal{B}_2 = \mathcal{A}_2$ and $f$ is an integrable function on
$X$, then the conditional expectation of $f$ with respect to
$\mathcal{B}_{1,2}$ is given by
\begin{equation}
        E(f \mid \mathcal{B}_{1, 2})(x_1, x_2)
         = \frac{1}{\mu(A_1)} \int_{A_1} f(t, x_2) \, d\mu_1(t)
\end{equation}
when $x_1 \in A_1 \in \mathcal{P}_1$.  This can be seen as a
combination of the conditional expectations associated to partitions
and product spaces, as in Sections \ref{conditional expectation} and
\ref{partitions}.  If $\mathcal{B}_2$ is any $\sigma$-subalgebra of
$\mathcal{A}_2$, then $E(f \mid \mathcal{B}_{1, 2})$ can be obtained
by first averaging $f(x_1, x_2)$ over $x_1 \in A_1$ for each $A_1 \in
\mathcal{P}_1$, as before, and then taking the conditional expectation
of the resulting functions of $x_2$ with respect to $\mathcal{B}_2$.
In this case, $\mathcal{B}_{1, 2}$ is a $\sigma$-subalgebra of the
$\sigma$-algebra associated to $\mathcal{P}_1$ and $A_2$.

\section[\ Partitions and vectors]{Partitions and vectors}
\label{partitions, vectors}
\setcounter{equation}{0}

        Let $(X, \mathcal{A}, \mu)$ be a probability space, and let
$\mathcal{P}$ be a partition of $X$ into finitely or countably many
measurable sets, as in Section \ref{partitions}.  As usual, it will be
convenient to ask that $\mu(A) > 0$ for each $A \in \mathcal{P}$.
Also let $\mathcal{B}(\mathcal{P})$ be the $\sigma$-subalgebra of
$\mathcal{A}$ generated by $\mathcal{P}$, consisting of unions of
elements of $\mathcal{P}$, including the empty set.  Thus a function
on $X$ is measurable with respect to $\mathcal{B}(\mathcal{P})$ if and
only if it is constant on the elements of $\mathcal{P}$.

        Let $V$ be a real or complex vector space with a norm $\|v\|$,
and let $f(x)$ be a $V$-valued function on $X$ that is constant on the
elements of $\mathcal{P}$.  In particular, $\|f(x)\|$ is a nonnegative
real-valued function on $X$ that is constant on the elements of
$\mathcal{P}$.  If $f(A)$ denotes the value of $f$ on $A \in \mathcal{P}$, then
\begin{equation}
\label{int_X ||f(x)|| d mu(x) = sum_{A in mathcal{P}} ||f(A)|| mu(A)}
 \int_X \|f(x)\| \, d\mu(x) = \sum_{A \in \mathcal{P}} \|f(A)\| \, \mu(A).
\end{equation}
More precisely, if $\mathcal{P}$ is a partition of $X$ into finitely
many sets, then the sum on the right is a finite sum, and $\|f(x)\|$
is automatically integrable on $X$.  If $\mathcal{P}$ consists of
infinitely many measurable subsets of $X$, then the sum on the right
is interpreted as the supremum of the corresponding sums over finite
subsets of $\mathcal{P}$, which may be infinite.

        If $\mathcal{P}$ has only finitely many elements, then we can put
\begin{equation}
\label{int_X f(x) d mu(x) = sum_{A in mathcal{P}} f(A) mu(A)}
        \int_X f(x) \, d\mu(x) = \sum_{A \in \mathcal{P}} f(A) \, \mu(A).
\end{equation}
This also makes sense when $\mathcal{P}$ has infinitely many elements,
$\|f(x)\|$ is integrable on $X$, and $V$ is complete.  In this case,
the sum on the right side of (\ref{int_X ||f(x)|| d mu(x) = sum_{A in
mathcal{P}} ||f(A)|| mu(A)}) is finite, and the sum on the right side
of (\ref{int_X f(x) d mu(x) = sum_{A in mathcal{P}} f(A) mu(A)})
converges in the generalized sense, as in Section \ref{summable functions, 2}.
In both cases,
\begin{equation}
\label{||int_X f(x) d mu(x)|| le int_X ||f(x)|| d mu(x)}
        \biggl\|\int_X f(x) \, d\mu(x)\biggr\| \le \int_X \|f(x)\| \, d\mu(x).
\end{equation}
Similarly, if $B \in \mathcal{B}(\mathcal{P})$, then we would like to
put
\begin{equation}
\label{int_B f(x) d mu(x) = sum_{A in mathcal{P} atop A subseteq B} f(A) mu(A)}
        \int_B f(x) \, d\mu(x)
           = \sum_{A \in \mathcal{P} \atop A \subseteq B} f(A) \, \mu(A).
\end{equation}
As before, this makes sense when $B$ is the union of finitely many
elements of $\mathcal{P}$, and when $B$ contains infinitely many
elements of $\mathcal{P}$, $\|f(x)\|$ is integrable, and $V$ is
complete.  We also have the analogue of (\ref{||int_X f(x) d mu(x)||
le int_X ||f(x)|| d mu(x)}) with $X = B$.

        Using the Bochner integral, one can integrate much more
complicated vector-valued functions.  We shall restrict our attention
here to sums over partitions for the sake of simplicity.

\section[\ Vector-valued martingales]{Vector-valued martingales}
\label{vector-valued martingales}
\setcounter{equation}{0}

        Let $(X, \mathcal{A}, \mu)$ be a probability space, and
suppose that $\mathcal{P}_1, \mathcal{P}_2, \ldots$ is a sequence of
partitions of $X$ into finitely or countably many measurable subsets
such that $\mathcal{P}_{j + 1}$ is a refinement of $\mathcal{P}_j$ for
each $j$.  This means that each $B \in \mathcal{P}_j$ is the union of
the $A \in \mathcal{P}_{j + 1}$ such that $A \subseteq B$.  If
$\mathcal{B}_j = \mathcal{B}(\mathcal{P}_j)$ is the $\sigma$-algebra
generated by $\mathcal{P}_j$, then it follows that $\mathcal{B}_j
\subseteq \mathcal{B}_{j + 1}$ for each $j$.  As usual, it is
convenient to ask that $\mu(A) > 0$ for each $A \in \mathcal{P}_j$.

        Let $V$ be a real or complex vector space with a norm $\|v\|$,
and let $f_l$ is a $V$-valued function on $X$ that is constant on
elements of $\mathcal{P}_l$.  We would like to define the conditional
expectation of $f_l$ with respect to $\mathcal{B}_j$ for $j < l$ by
\begin{equation}
\label{E(f_l mid mathcal{B}_j)(x) = ...}
 E(f_l \mid \mathcal{B}_j)(x) = \frac{1}{\mu(B)}\int_B f_l \, d\mu =
\sum_{A \in \mathcal{P}_l \atop A \subseteq B} f_l(A) \, \frac{\mu(A)}{\mu(B)}
\end{equation}
when $x \in B \in \mathcal{P}_j$, where $f_l(A)$ denotes the value of
$f_l$ on $A \in \mathcal{P}_l$, as in the previous section.  This
makes sense when each $B \in \mathcal{P}_j$ is the union of finitely
many $A \in \mathcal{P}_l$, and when $\|f_l\|$ is integrable and $V$
is complete.  In both cases, it is easy to see that
\begin{equation}
\label{||E(f_l mid mathcal{B}_j)|| le E(||f_l|| mid mathcal{B}_j)}
        \|E(f_l \mid \mathcal{B}_j)\| \le E(\|f_l\| \mid \mathcal{B}_j).
\end{equation}
If $j < k < l$, then one can also check that
\begin{equation}
\label{E(E(f_l mid mathcal{B}_k) mid mathcal{B}_j) = E(f_l mid mathcal{B}_j)}
 E(E(f_l \mid \mathcal{B}_k) \mid \mathcal{B}_j) = E(f_l \mid \mathcal{B}_j),
\end{equation}
under these conditions, just as in the context of real or
complex-valued functions.

        Now let $\{f_j\}_{j = 1}^\infty$ be a sequence of $V$-valued
functions on $X$ such that $f_j$ is constant on the elements of
$\mathcal{P}_j$ for each $j$.  As usual, $\{f_j\}_{j = 1}^\infty$
is said to be a martingale with respect to this filtration if
\begin{equation}
\label{f_j = E(f_l mid mathcal{B}_j)}
        f_j = E(f_l \mid \mathcal{B}_j)
\end{equation}
for each $j \le l$.  More precisely, this makes sense when each
element of $\mathcal{P}_j$ is the union of finitely many elements of
$\mathcal{P}_l$, and when each $\|f_l\|$ is integrable and $V$ is
complete.  Note that (\ref{f_j = E(f_l mid mathcal{B}_j)}) holds for
all $j \le l$ when it holds for $l = j + 1$, because of (\ref{E(E(f_l
mid mathcal{B}_k) mid mathcal{B}_j) = E(f_l mid mathcal{B}_j)}).

        Of course, the simplest type of situation occurs when
$\mathcal{P}_j$ consists of only finitely many measurable subsets of
$X$ for each $j$.  All of the sums involved in the conditional
expectations are then finite sums, and the functions $\|f_l\|$ are
automatically bounded.

\section[\ $L^1$-Valued martingales]{$L^1$-Valued martingales}
\label{L^1-valued martingales}
\setcounter{equation}{0}

        Let us continue with the same notations and hypotheses as in
the previous section.  As in Section \ref{unit square}, we can get an
example of a $V$-valued martingale on $X$ with $V = L^1(X,
\mathcal{A})$ by taking
\begin{equation}
        f_l(A) = \mu(A)^{-1} \, {\bf 1}_A
\end{equation}
for each $A \in \mathcal{P}_l$.  Here ${\bf 1}_A$ denotes the
indicator function associated to $A$ on $X$, equal to $1$ on $A$ and
$0$ on $X \backslash A$, as usual.  Thus $\|f_l(x)\|_1 = 1$ for every
$x \in X$ and $l \ge 1$, and it is easy to check that (\ref{f_j =
E(f_l mid mathcal{B}_j)}) holds.

        Now let $(Y, \mathcal{B}, \nu)$ be a $\sigma$-finite measure
space, and let us consider functions on $X$ with values in $V =
L^1(Y)$.  If $f_l(x)$ is an $L^1(Y)$-valued function on $X$ that is
constant on the elements of $\mathcal{P}_l$, then
\begin{equation}
        F_l(x, y) = f_l(x)(y)
\end{equation}
defines a function on $X \times Y$ that is constant in $x$ on each
element of $\mathcal{P}_l$ and measurable in $y$ for each $x \in X$.
If $\|f_l(x)\|_{L^1(Y)}$ is integrable on $X$, then $F_l(x, y)$ is
integrable on $X \times Y$, and
\begin{eqnarray}
        \int_X \|f_l(x)\|_{L^1(Y)} \, d\mu(x)
 & = & \int_X \Big(\int_Y |F_l(x, y)| \, d\nu(y)\Big) \, d\mu(x) \\
 & = & \int_{X \times Y} |F_l(x, y)| \, d(\mu \times \nu)(x, y). \nonumber
\end{eqnarray}
Conversely, if $F_l(x, y)$ is an integrable function on $X \times Y$
that is constant in $x$ on each element of $\mathcal{P}_l$, then we
get an $L^1(Y)$-valued function $f_l(x)$ on $X$ that is constant on
each element of $\mathcal{P}_l$ and for which $\|f_l(x)\|_{L^1(Y)}$ is
integrable on $X$.

        Let $\widehat{\mathcal{B}}_l$ be the $\sigma$-algebra of
subsets of $X \times Y$ that corresponds to $\mathcal{B}_l =
\mathcal{B}(\mathcal{P}_l)$ on $X$ and $\mathcal{B}$ on $Y$ in the
product space.  As in Section \ref{partitions, products}, a set
$\widehat{A} \subseteq X \times Y$ is in $\widehat{\mathcal{B}}_l$ if
and only if for each $A \in \mathcal{P}_l$ there is a $B \in
\mathcal{B}$ such that
\begin{equation}
        \widehat{A} \cap (A \times Y) = A \times B.
\end{equation}
Equivalently, $\widehat{A} \in \widehat{\mathcal{B}}_l$ if it can be
expressed as the union of sets of the form $A \times B(A)$, where $A$
runs through the elements of $\mathcal{P}_l$, and $B(A) \in
\mathcal{B}$ for each $A \in \mathcal{P}_l$.  In the context of the
preceding paragraph, the functions $F_l(x, y)$ are measurable with
respect to $\widehat{\mathcal{B}}_l$.

        Suppose that $\{f_l\}_{l = 1}^\infty$ is a sequence of
$L^1(Y)$-valued functions on $X$ such that $f_l(x)$ is constant on
each element of $\mathcal{P}_l$ and $\|f_l(x)\|_{L^1(Y)}$ is
integrable on $X$ for each $l$.  This corresponds exactly to a
sequence $\{F_l\}_{l = 1}^\infty$ of integrable functions on $X \times
Y$ such that $F_l(x, y)$ is measurable with respect to
$\widehat{\mathcal{B}}_l$ for each $l$, as in the previous paragraphs.
If $Y$ is a probability space, then $X \times Y$ is also a probability
space, and it is easy to see that $\{f_l\}_{l = 1}^\infty$ is an
$L^1(Y)$-valued martingale on $X$ with respect to the
$\mathcal{B}_l$'s if and only if $\{F_l\}_{l = 1}^\infty$ is a
martingale on $X \times Y$ with respect to the $\widehat{\mathcal{B}}_l$'s.
This basically works as well when $Y$ is $\sigma$-finite, by extending
the relevant definitions in a natural way.

\section[\ Pointwise convergence]{Pointwise convergence}
\label{pointwise convergence}
\setcounter{equation}{0}

        Let us continue with the same notation and hypotheses as in
Section \ref{vector-valued martingales}, with the additional condition
that $V$ be complete.  Suppose that $\{f_j\}_{j = 1}^\infty$ is a
sequence of $V$-valued functions on $X$ such that $f_j(x)$ is constant
on each element of $\mathcal{P}_j$, $\|f_j(x)\|$ is integrable on $X$
for each $j$, and $\{f_j\}_{j = 1}^\infty$ is a martingale with
respect to $\mathcal{B}_j = \mathcal{B}(\mathcal{P}_j)$.  If
\begin{equation}
        f_n^*(x) = \max_{1 \le j \le n} \|f_j(x)\|
\end{equation}
is the usual maximal function and
\begin{equation}
        A_n(t) = \{x \in X : f_n^*(x) > t\}
\end{equation}
for each $t > 0$, then
\begin{equation}
\label{t mu(A_n(t)) le int_X ||f_n(x)|| d mu(x)}
        t \, \mu(A_n(t)) \le \int_X \|f_n(x)\| \, d\mu(x)
\end{equation}
for every $t > 0$ and $n \ge 1$.  This can be shown in the standard
way.  In particular, one can use the fact that $\{\|f_j\|\}_{j =
1}^\infty$ is a submartingale, because of (\ref{||E(f_l mid
mathcal{B}_j)|| le E(||f_l|| mid mathcal{B}_j)}).

        Suppose now that $\|f_n(x)\|$ has uniformly bounded $L^1$
norm, and put
\begin{equation}
        f^*(x) = \sup_{j \ge 1} \|f_j(x)\|.
\end{equation}
If
\begin{equation}
        A(t) = \{x \in X : f^*(x) > t\}
\end{equation}
for each $t > 0$, then
\begin{equation}
        A(t) = \bigcup_{n = 1}^\infty A_n(t),
\end{equation}
and of course $A_n(t) \subseteq A_{n + 1}(t)$.  It follows that
\begin{equation}
\label{t mu(A(t)) le sup_{n ge 1} int_X ||f_n(x)|| d mu(x)}
        t \, \mu(A(t)) \le \sup_{n \ge 1} \int_X \|f_n(x)\| \, d\mu(x)
\end{equation}
for each $t > 0$, by taking the limit as $n \to \infty$ in (\ref{t
mu(A_n(t)) le int_X ||f_n(x)|| d mu(x)}).

        As in Section \ref{convergence almost everywhere}, we can also
consider $\{f_j - f_l\}_{j = l}^\infty$ as a $V$-valued martingale on
$X$ with respect to the $\mathcal{B}_j$'s with $j \ge l$.  If
\begin{equation}
 B_l(t) = \bigg\{x \in X : \sup_{j \ge l} \|f_j(x) - f_l(x)\| > t\bigg\},
\end{equation}
then we get that
\begin{equation}
 t \, \mu(B_l(t)) \le \sup_{j \ge l} \int_X \|f_j(x) - f_l(x)\| \, d\mu(x)
\end{equation}
for each $t > 0$ and $l \ge 1$.  Hence
\begin{equation}
\label{t mu(bigcap_{l = 1}^infty B_l(t)) le ...}
        t \, \mu\Big(\bigcap_{l = 1}^\infty B_l(t)\Big)
 \le \lim_{l \to \infty} \sup_{j \ge l} \int_X \|f_j(x) - f_l(x)\| \, d\mu(x)
\end{equation}
for each $t > 0$.

        Suppose that
\begin{equation}
\label{lim_{l to infty} sup_{j ge l} int_X ||f_j(x) - f_l(x)|| d mu(x) = 0}
 \lim_{l \to \infty} \sup_{j \ge l} \int_X \|f_j(x) - f_l(x)\| \, d\mu(x) = 0,
\end{equation}
which means that $\{f_j\}_{j = 1}^\infty$ is a Cauchy sequence with
respect to the $L^1$ norm for $V$-valued functions on $X$.  This
together with (\ref{t mu(bigcap_{l = 1}^infty B_l(t)) le ...}) implies
that
\begin{equation}
        \mu\Big(\bigcap_{l = 1}^\infty B_l(t)\Big) = 0
\end{equation}
for every $t > 0$.  Of course,
\begin{eqnarray}
\lefteqn{X \backslash \Big(\bigcap_{l = 1}^\infty B_l(t)\Big) =} \\
             & & \{x \in X : \sup_{j \ge l} \|f_j(x) - f_l(x)\| \le t
                       \hbox{ for some } l \in {\bf Z}_+\}, \nonumber
\end{eqnarray}
and it follows that
\begin{equation}
        \lim_{l \to \infty} \sup_{j \ge l} \|f_j(x) - f_l(x)\| = 0
\end{equation}
for almost every $x \in X$, by taking $t = 1/n$ for $n \in {\bf Z}_+$.
This shows that $\{f_j(x)\}_{j = 1}^\infty$ is a Cauchy sequence in
$V$ for almost every $x \in X$, and hence that $\{f_j(x)\}_{j =
1}^\infty$ converges for almost every $x \in X$, because $V$ is
complete.  Thus this criterion for convergence almost everywhere works
as well in the vector-valued case as for real or complex-valued
functions.

\section[\ Another scenario]{Another scenario}
\label{another scenario}
\setcounter{equation}{0}

        Let $(X_1, \mathcal{A}_1, \mu_1), (X_2, \mathcal{A}_2, \mu_2),
\ldots$ be a sequence of probability spaces, and let $X = \prod_{j =
1}^\infty X_j$ be their Cartesian product, with the product measure
$\mu$.  As usual, let $\mathcal{B}_n$ be the $\sigma$-subalgebra of
the $\sigma$-algebra of measurable subsets of $X$ of the form
\begin{equation}
\label{A times prod_{j = n + 1}^infty X_j}
        A \times \prod_{j = n + 1}^\infty X_j,
\end{equation}
where $A$ is a measurable subset of $\prod_{j = 1}^n X_j$.  If each
$X_j$ has only finitely or countably many elements, and every subset
of $X_j$ is measurable, then $\mathcal{B}_n$ consists of the sets of
the form (\ref{A times prod_{j = n + 1}^infty X_j}), where $A$ is any
subset of $\prod_{j = 1}^n X_j$.  In this case, $\mathcal{B}_n$ is the
$\sigma$-algebra generated by the partition $\mathcal{P}_n$ of subsets
of $X$ of the form (\ref{A times prod_{j = n + 1}^infty X_j}), where
$A \subseteq \prod_{j = 1}^n X_j$ has exactly one element.

        Let $a_1(x_1), a_2(x_2), \ldots$ be a sequence of integrable
real or complex-valued functions on $X_1, X_2, \ldots$ such that
\begin{equation}
\label{int_{X_j} a_j(x_j) d mu(x_j) = 0}
        \int_{X_j} a_j(x_j) \, d\mu(x_j) = 0
\end{equation}
for each $j$.  Also let $V$ be a real or complex vector space with a
norm $\|v\|$, and let $v_1, v_2, \ldots$ be a sequence of elements of
$V$.  Under these conditions, it is natural to consider
\begin{equation}
\label{f_n(x) = sum_{j = 1}^n a_j(x_j) v_j}
        f_n(x) = \sum_{j = 1}^n a_j(x_j) \, v_j
\end{equation}
as a $V$-valued martingale on $X$ with respect to the
$\mathcal{B}_n$'s.  In this case, it is very easy to understand the
meaning of the vector-valued integrals, because of the special form of
the functions.  This is also consistent with the discussion in Section
\ref{vector-valued martingales} when the $X_j$'s have only finitely or
countably many elements, and all of their subsets are measurable.

        By construction, $f_n(x)$ takes values in a linear subspace of
$V$ with dimension less than or equal to $n$ for each $n \in {\bf
Z}_+$.  Thus one can identify $f_n$ with a function on $X$ with values
in ${\bf R}^n$ or ${\bf C}^n$ whose components are measurable.  One
can also check that $\|f_n(x)\|$ is measurable as a nonnegative
real-valued function on $X$, using the fact that any norm on ${\bf
R}^n$ or ${\bf C}^n$ is bounded by a constant multiple of the standard
norm, and hence is continuous with respect to the standard topology.
Moreover, $\{\|f_n(x)\|\}_{n = 1}^\infty$ is a submartingale with
respect to the $\mathcal{B}_n$'s, basically because the norm of the
integral of a $V$-valued function is less than or equal to the
integral of the norm of the function.

        Suppose that $\|f_n(x)\|$ has uniformly bounded $L^1$ norm, and let
\begin{equation}
        f^*(x) = \sup_{n \ge 1} \|f_n(x)\|
\end{equation}
be the corresponding maximal function.  As in the previous section,
\begin{equation}
 t \, \mu(\{x \in X : f^*(x) > t\}) \le \sup_{n \ge 1} \int_X \|f_n\| \, d\mu
\end{equation}
for every $t > 0$.  This permits one to show that
\begin{equation}
        \lim_{n \to \infty} \sup_{l \ge n} \|f_l(x) - f_n(x)\| = 0
\end{equation}
for almost every $x \in X$ when
\begin{equation}
        \lim_{n \to \infty} \sup_{l \ge n} \int_X \|f_l - f_n\| \, d\mu = 0,
\end{equation}
as before.  Hence $\{f_n(x)\}_{n = 1}^\infty$ is a Cauchy sequence in
$V$ for almost every $x$ in $X$ under these conditions.  If $V$ is
complete, then it follows that $\{f_n(x)\}_{n = 1}^\infty$ converges
for almost every $x \in X$.

\section[\ Hilbert space martingales]{Hilbert space martingales}
\label{hilbert space martingales}
\setcounter{equation}{0}

        Let $(X, \mathcal{A}, \mu)$ be a probability space, and
suppose that $\mathcal{B}_1 \subseteq \mathcal{B}_2 \subseteq \cdots$
is an increasing sequence of $\sigma$-subalgebra of $\mathcal{A}$ as
in Section \ref{vector-valued martingales} or the preceding section.
Also let $(V, \langle v, w \rangle)$ be a real or complex Hilbert
space, and let $\{f_j\}_{j = 1}^\infty$ be a $V$-valued martingale
with respect to the $\mathcal{B}_j$'s such that $\|f_j(x)\| \in
L^2(X)$ for each $j$.

        As in Section \ref{martingales}, one can check that
\begin{equation}
        \int_X \langle f_j(x), f_{l + 1}(x) - f_l(x)\rangle \, d\mu(x) = 0
\end{equation}
for each $j \le l$.   If $j < l$, then we get that
\begin{equation}
 \int_X \langle f_{j + 1}(x) - f_j(x), f_{l + 1}(x) - f_l(x)\rangle \, d\mu(x)
                                                                      = 0.
\end{equation}
Using the identity $f_n = f_1 + \sum_{j = 1}^{n - 1} (f_{j + 1} - f_j)$,
it follows that
\begin{eqnarray}
\label{int_X ||f_n(x)||^2 d mu(x) = ...}
\lefteqn{\int_X \|f_n(x)\|^2 \, d\mu(x) =} \\
 & &  \int_X \|f_1(x)\|^2 \, d\mu(x) +
  \sum_{j = 1}^{n - 1} \int_X \|f_{j + 1}(x) - f_j(x)\|^2 \, d\mu(x) \nonumber
\end{eqnarray}
for each $n$.  Similarly,
\begin{equation}
\label{int_X ||f_n(x) - f_l(x)||^2 d mu(x) = ...}
        \int_X \|f_n(x) - f_l(x)\|^2 \, d\mu(x) =
         \sum_{j = l}^{n - 1} \int_X \|f_{j + 1}(x) - f_j(x)\|^2 \, d\mu(x)
\end{equation}
when $n > l$.

        If $\|f_n(x)\|$ has bounded $L^2$ norm, then (\ref{int_X
||f_n(x)||^2 d mu(x) = ...}) implies that
\begin{equation}
 \sum_{j = 1}^\infty \int_X \|f_{j + 1}(x) - f_j(x)\|^2 \, d\mu(x) < \infty.
\end{equation}
Under these conditions,
\begin{equation}
 \lim_{l \to \infty} \sum_{j = l}^\infty \int_X \|f_{j + 1}(x) - f_j(x)\|^2
                                                               \, d\mu(x) = 0,
\end{equation}
and hence
\begin{equation}
 \lim_{l \to \infty} \sup_{n > l} \int_X \|f_n(x) - f_l(x)\|^2 \, d\mu(x) = 0.
\end{equation}
In particular, $\{f_j(x)\}_{j = 1}^\infty$ converges in $V$ for almost
every $x \in X$, as in the previous sections.

\section[\ Nonnegative submartingales]{Nonnegative submartingales}
\label{nonnegative submartingales}
\setcounter{equation}{0}

        Let $(X, \mathcal{A}, \mu)$ be a probability space, and let
$\mathcal{B}_1 \subseteq \mathcal{B}_2 \subseteq \cdots$ be an
increasing sequence of $\sigma$-subalgebras of $\mathcal{A}$.  Also
let $\{f_j\}_{j = 1}^\infty$ be a submartingale with respect to this
filtration such that $f_j \ge 0$ for each $j$.  This includes the case
of the norm of a vector-valued martingale, as before.  If
\begin{equation}
        f_n^*(x) = \max_{1 \le j \le n} f_j(x)
\end{equation}
and
\begin{equation}
        A_n(t) = \{x \in X : f_n^*(x) > t\},
\end{equation}
then
\begin{equation}
        t \, \mu(A_n(t)) \le \int_{A_n(t)} f_n \, d\mu \le \int_X f_n \, d\mu
\end{equation}
for each $t > 0$ and $n \ge 1$, as shown previously.  If
\begin{equation}
        f^*(x) = \sup_{j \ge 1} f_j(x)
\end{equation}
and
\begin{equation}
        A(t) = \{x \in X : f^*(x) > t\},
\end{equation}
then
\begin{equation}
        A(t) = \bigcup_{n = 1}^\infty A_n(t)
\end{equation}
and
\begin{equation}
        t \, \mu(A(t)) \le \sup_{n \ge 1} \int_X f_n \, d\mu
\end{equation}
for each $t > 0$ when the $L^1$ norms of the $f_n$'s are bounded.

        By hypothesis,
\begin{equation}
\label{0 le f_j le E(f_n mid mathcal{B}_j)}
        0 \le f_j \le E(f_n \mid \mathcal{B}_j)
\end{equation}
when $j \le n$, and of course $E(f_n \mid \mathcal{B}_j)$ is a
martingale in $j$ for each $n$.  If $f_n \in L^p(X)$, $1 < p <
\infty$, then
\begin{equation}
 \int_X \Big(\max_{1 \le j \le n} E(f_n \mid \mathcal{B}_j)\Big)^p \, d\mu
                      \le \frac{p \, 2^{p - 1}}{p - 1} \int_X f_n^p \, d\mu,
\end{equation}
as in Section \ref{maximal functions, 4}.  Hence
\begin{equation}
\label{int_X (f_n^*)^p d mu le frac{p 2^{p - 1}}{p - 1} int_X f_n^p d mu}
\int_X (f_n^*)^p \, d\mu \le \frac{p \, 2^{p - 1}}{p - 1} \int_X f_n^p \, d\mu.
\end{equation}
If the $L^p$ norm of $f_n$ is uniformly bounded in $n$, then the
monotone convergence theorem implies that $f^* \in L^p$, with
\begin{equation}
\label{int_X (f^*)^p d mu le frac{p 2^{p - 1}}{p - 1} sup_n int_X f_n^p d mu}
        \int_X (f^*)^p \, d\mu
         \le \frac{p \, 2^{p - 1}}{p - 1} \sup_{n \ge 1} \int_X f_n^p \, d\mu.
\end{equation}
Thus one gets the same $L^p$ estimates for nonnegative submartingales
as for martingales.

\section[\ $L^p$-Valued martingales]{$L^p$-Valued martingales}
\label{L^p-valued martingales}
\setcounter{equation}{0}

        As in Section \ref{L^1-valued martingales}, we can look at
$L^p$-valued martingales in terms of functions on a product space.
Let $(X, \mathcal{A}, \mu)$ be a probability space, and let
$\mathcal{P}_1, \mathcal{P}_2, \ldots$ be a sequence of partitions of
$X$ into finitely or countably many measurable subsets with positive
measure such that $\mathcal{P}_{j + 1}$ is a refinement of
$\mathcal{P}_j$ for each $j$.  Also let $(Y, \mathcal{B}, \nu)$ be a
$\sigma$-finite measure space, and fix $p$, $1 < p < \infty$.

        If $f_l(x)$ is an $L^p(Y)$-valued function on $X$ that is
constant on the elements of $\mathcal{P}_l$, then
\begin{equation}
        F_l(x, y) = f_l(x)(y)
\end{equation}
is a function on $X \times Y$ that is constant in $x$ on each element
of $\mathcal{P}_l$ and measurable in $y$ for each $x \in X$.  If
$\|f_l(x)\|_{L^p(Y)} \in L^p(X)$, then $F_l(x, y) \in L^p(X \times Y)$, and
\begin{eqnarray}
        \int_X \|f_l(x)\|_{L^p(Y)}^p \, d\mu(x)
         & = & \int_X \Big(\int_Y |F_l(x, y)|^p \, d\nu(y)\Big) \, d\mu(x) \\
 & = & \int_{X \times Y} |F_l(x, y)|^p \, d(\mu \times \nu)(x, y). \nonumber
\end{eqnarray}
Conversely, if $F_l(x, y) \in L^p(X \times Y)$ is constant in $x$ on
each element of $\mathcal{P}_l$, then we get an $L^p(Y)$-valued
function $f_l(x)$ on $X$ that is constant on each element of
$\mathcal{P}_l$ and for which $\|f_l(x)\|_{L^p(Y)} \in L^p(X)$.  If
$\widehat{B}_l$ is the $\sigma$-algebra of subsets of $X \times Y$
that corresponds to $\mathcal{B}_l = \mathcal{B}(\mathcal{P}_l)$ on
$X$ and $\mathcal{B}$ on $Y$ as before, then $F_l(x, y)$ is measurable
with respect to $\widehat{\mathcal{B}}_l$.

        Now let $\{f_l\}_{l = 1}^\infty$ be a sequence of
$L^p(Y)$-valued functions on $X$ such that $f_l(x)$ is constant on
each element of $\mathcal{P}_l$ and $\|f_l(x)\|_{L^p(Y)} \in L^p(X)$
for each $l$.  This corresponds exactly to a sequence of functions
$\{F_l\}_{l = 1}^\infty$ in $L^p(X \times Y)$ such that $F_l(x, y)$ is
measurable with respect to $\widehat{\mathcal{B}}_l$ for each $l$, as
in the preceding paragraph.  If $\{f_l\}_{l = 1}^\infty$ is an
$L^p(Y)$-valued martingale on $X$ with respect to the
$\mathcal{B}_l$'s and $Y$ is a probability space, then $X \times Y$ is
also a probability space. and $\{F_l\}_{l = 1}^\infty$ is a martingale
on $X \times Y$ with respect to the $\widehat{\mathcal{B}}_l$'s.  If
the $L^p(X)$ norm of $\|f_l(x)\|_{L^p(Y)}$ is bounded, then the $L^p(X
\times Y)$ norm of $F_l(x, y)$ is bounded, and hence $\{F_l\}_{l =
1}^\infty$ converges in $L^p(X \times Y)$.  In particular, $\{F_l\}_{l
= 1}^\infty$ is a Cauchy sequence in $L^p(X \times Y)$, which implies that
\begin{equation}
\label{lim_{l to infty} sup_{j ge l} int_X ||f_j - f_l||_{L^p(Y)}^p d mu = 0}
 \lim_{l \to \infty} \sup_{j \ge l} \int_X \|f_j(x) - f_l(x)\|_{L^p(Y)}^p
                                                              \, d\mu(x) = 0.
\end{equation}

        Of course, the same conclusion holds when $0 < \nu(Y) <
\infty$, by dividing by $\nu(Y)$ to get a probability space.
Otherwise, let $\rho$ be a strictly positive measurable function on
$Y$ such that
\begin{equation}
        \int_Y \rho(y) \, d\nu(y) = 1,
\end{equation}
which is possible because $(Y, \mathcal{B}, \nu)$ is supposed to be
$\sigma$-finite.  Thus
\begin{equation}
        \nu_\rho(B) = \int_B \rho(y) \, d\nu(y)
\end{equation}
is a probability measure on $(Y, \mathcal{B})$.  If $\phi(y) \in
L^p(Y, \nu)$, then
\begin{equation}
        \phi_\rho(y) = \phi(y) \, \rho(y)^{-1/p} \in L^p(Y, \nu_\rho),
\end{equation}
and
\begin{equation}
 \int_Y |\phi_\rho(y)|^p \, d\nu_\rho(y) = \int_Y |\phi(y)|^p \, d\nu(y).
\end{equation}
Using this, one can check that (\ref{lim_{l to infty} sup_{j ge l}
int_X ||f_j - f_l||_{L^p(Y)}^p d mu = 0}) holds for any
$\sigma$-finite measure space $(Y, \mathcal{B}, \nu)$, by reducing to
the probability space $(Y, \mathcal{B}, \nu_\rho)$.

\section[\ Another criterion]{Another criterion}
\label{another criterion}
\setcounter{equation}{0}

        Let $(X, \mathcal{A}, \mu)$ be a probability space, and let
$\mathcal{B}_1 \subseteq \mathcal{B}_2 \subseteq \cdots$ be an
increasing sequence of $\sigma$-subalgebras of $\mathcal{A}$ as in
Section \ref{vector-valued martingales} or \ref{another scenario}.
Also let $V$ be a real or complex Banach space with a norm $\|v\|$,
and let $\{f_j\}_{j = 1}^\infty$ be a $V$-valued martingale on $X$
with respect to the $\mathcal{B}_j$'s such that $\|f_j(x)\| \in
L^1(X)$ for each $j$.  Suppose that for each $\epsilon > 0$ there is a
$V$-valued martingale $\{g_j\}_{j = 1}^\infty$ on $X$ with respect to
the $\mathcal{B}_j$'s such that
\begin{equation}
\label{int_X ||f_j(x) - g_j(x)|| d mu(x) le epsilon}
        \int_X \|f_j(x) - g_j(x)\| \, d\mu(x) \le \epsilon
\end{equation}
for each $j$, and $\{g_j(x)\}_{j = 1}^\infty$ converges in $V$ for
almost every $x \in X$.  Let us check that $\{f_j(x)\}_{j = 1}^\infty$
converges in $V$ for almost every $x \in X$ under these conditions.

        Of course, it suffices to show that
\begin{equation}
\label{lim_{l to infty} sup_{j ge l} ||f_j(x) - f_l(x)|| = 0}
        \lim_{l \to \infty} \sup_{j \ge l} \|f_j(x) - f_l(x)\| = 0
\end{equation}
for almost every $x \in X$, so that $\{f_j(x)\}_{j = 1}^\infty$ is a
Cauchy sequence in $V$ for almost every $x \in X$.  Put $h_j = f_j -
g_j$, so that $\{h_j\}_{j = 1}^\infty$ is also a $V$-valued martingale
on $X$ with respect to the $\mathcal{B}_j$'s.  Observe that
\begin{eqnarray}
        \lim_{l \to \infty} \sup_{j \ge l} \|f_j(x) - f_l(x)\|
         & \le & \lim_{l \to \infty} \sup_{j \ge l} \|g_j(x) - g_l(x)\| \\
 & &    + \lim_{l \to \infty} \sup_{j \ge l} \|h_j(x) - h_l(x)\| \nonumber
\end{eqnarray}
for every $x \in X$.  This implies that
\begin{equation}
        \lim_{l \to \infty} \sup_{j \ge l} \|f_j(x) - f_l(x)\|
         \le \lim_{l \to \infty} \sup_{j \ge l} \|h_j(x) - h_l(x)\|
\end{equation}
for almost every $x \in X$, because $\{g_j(x)\}_{j = 1}^\infty$ is a
Cauchy sequence in $V$ for almost every $x \in X$.  Hence
\begin{equation}
\label{lim_{l to infty} sup_{j ge l} ||f_j(x) - f_l(x)|| le 2 h^*(x)}
        \lim_{l \to \infty} \sup_{j \ge l} \|f_j(x) - f_l(x)\|
         \le 2 \, \sup_{j \ge 1} \|h_j(x)\| = 2 \, h^*(x)
\end{equation}
for almost every $x \in X$.

        By the usual maximal function estimate,
\begin{equation}
\label{t mu({x in X : h^*(x) > t}) le epsilon}
        t \, \mu(\{x \in X : h^*(x) > t\})
         \le \sup_{j \ge 1} \int_X \|h_j(x)\| \, d\mu(x) \le \epsilon
\end{equation}
for every $t > 0$.  If
\begin{equation}
        E(t) = \bigg\{x \in X : \lim_{l \to \infty} \sup_{j \ge l}
                                         \|f_j(x) - f_l(x)\| > 2 \, t \bigg\},
\end{equation}
then
\begin{equation}
        \mu(E(t)) \le \mu(\{x \in X : h^*(x) > t\}),
\end{equation}
by (\ref{lim_{l to infty} sup_{j ge l} ||f_j(x) - f_l(x)|| le 2 h^*(x)}),
and so
\begin{equation}
        t \, \mu(E(t)) \le \epsilon
\end{equation}
for every $\epsilon, t > 0$.  Because $E(t)$ does not depend on
$\epsilon$, we may conclude that $\mu(E(t)) = 0$ for every $t > 0$.
This implies that (\ref{lim_{l to infty} sup_{j ge l} ||f_j(x) -
f_l(x)|| = 0}) holds for almost every $x \in X$, as desired.

        Note that this criterion is satisfied when
\begin{equation}
 \lim_{l \to \infty} \sup_{j \ge l} \int_X \|f_j(x) - f_l(x)\| \, d\mu(x) = 0.
\end{equation}
To see this, one can take $g_j(x)$ to be of the form $f_{\min(j,
N)}(x)$ for large positive integers $N$.  This converges as $j \to
\infty$ for each fixed $N$ trivially, and (\ref{int_X ||f_j(x) -
g_j(x)|| d mu(x) le epsilon}) holds for sufficiently large $N$ by
hypothesis.

\section[\ $\ell^1$-Valued martingales]{$\ell^1$-Valued martingales}
\label{ell^1-valued martingales}
\setcounter{equation}{0}

        As before, let $(X, \mathcal{A}, \mu)$ be a probability space,
and let $\mathcal{P}_1, \mathcal{P}_2, \ldots$ be a sequence of
partitions of $X$ into finitely or countably many measurable subsets
with positive measure such that $\mathcal{P}_{j + 1}$ is a refinement
of $\mathcal{P}_j$ for each $j$.  Suppose that $\{f_j\}_{j =
1}^\infty$ is a sequence of functions on $X$ with values in $\ell^1 =
\ell^1({\bf Z}_+)$.  Thus for each $x \in X$ and $j \ge 1$ we get a
summable sequence $\{f_{j, k}(x)\}_{k = 1}^\infty$ of real or complex
numbers.  Of course, $f_j(x)$ is constant on each element of
$\mathcal{P}_j$ if and only if $f_{j, k}(x)$ is constant on each
element of $\mathcal{P}_j$ for each $k \ge 1$.  If
$\|f_j(x)\|_{\ell^1}$ is integrable on $X$, then $f_{j, k}(x)$ is
integrable on $X$ for each $k$, and
\begin{equation}
\label{int_X ||f_j(x)||_{ell^1} d mu(x)}
        \quad \int_X \|f_j(x)\|_{\ell^1} \, d\mu(x) =
        \int_X \sum_{k = 1}^\infty |f_{j, k}(x)| \, d\mu(x)
          = \sum_{k = 1}^\infty \int_X |f_{j, k}(x)| \, d\mu(x).
\end{equation}

        Suppose now that $\{f_j\}_{j = 1}^\infty$ is an
$\ell^1$-valued martingale on $X$ with respect to $\mathcal{B}_j =
\mathcal{B}(\mathcal{P}_j)$.  This implies that $\{f_{j, k}\}_{j =
1}^\infty$ is a martingale on $X$ with respect to the
$\mathcal{B}_j$'s for each $k$.  In particular,
\begin{equation}
\label{int_X |f_{j, k}(x)| d mu(x) le int_X |f_{j + 1, k}(x)| d mu(x)}
  \int_X |f_{j, k}(x)| \, d\mu(x) \le \int_X |f_{j + 1, k}(x)| \, d\mu(x)
\end{equation}
for each $j, k \ge 1$, and hence
\begin{equation}
\label{int_X ||f_j||_{ell^1} d mu le int_X ||f_{j + 1}||_{ell^1} d mu}
        \int_X \|f_j(x)\|_{\ell^1} \, d\mu(x)
         \le \int_X \|f_{j + 1}(x)\|_{\ell^1} \, d\mu(x)
\end{equation}
for each $j$.

        Suppose also that the $L^1(X)$ norm of $\|f_j(x)\|_{\ell^1}$
is bounded.  Because of monotonicity,
\begin{equation}
        \sup_{j \ge 1} \int_X \|f_j(x)\|_{\ell^1} \, d\mu(x)
          = \lim_{j \to \infty} \int_X \|f_j(x)\|_{\ell^1} \, d\mu(x),
\end{equation}
and similarly
\begin{equation}
        \sup_{j \ge 1} \int_X |f_{j, k}(x)| \, d\mu(x)
          = \lim_{j \to \infty} \int_X |f_{j, k}(x)| \, d\mu(x)
\end{equation}
for each $k$.  The monotone convergence theorem for sums implies that
\begin{equation}
\lim_{j \to \infty} \int_X \|f_j(x)\|_{\ell^1} \, d\mu(x) = \sum_{k = 1}^\infty            \Big(\lim_{j \to \infty} \int_X |f_{j, k}(x)| \, d\mu(x)\Big).
\end{equation}
Therefore
\begin{equation}
          \sup_{j \ge 1} \int_X \|f_j(x)\|_{\ell^1} \, d\mu(x) =
 \sum_{k = 1}^\infty \Big(\sup_{j \ge 1} \int_X |f_{j, k}(x)| \, d\mu(x)\Big).
\end{equation}

        Let $N$ be a large positive integer, and put
\begin{eqnarray}
 g_{j, k}(x) = f_{j, k}(x), \, h_{j, k}(x) = 0 & \hbox{when} & k \le N, \\
 g_{j, k}(x) = 0, \, h_{j, k}(x) = f_{j, k}(x) & \hbox{when} & k > N.
\end{eqnarray}
If $g_j(x) = \{g_{j, k}(x)\}_{k = 1}^\infty$, $h_j(x) = \{h_{j,
k}(x)\}_{k = 1}^\infty$, then $g_j(x), h_j(x) \in \ell^1$ and
\begin{equation}
        f_j(x) = g_j(x) + h_j(x).
\end{equation}
Note that $\{g_j(x)\}_{j = 1}^\infty$ converges for almost every $x
\in X$, as a consequence of the convergence almost everywhere of real
or complex martingales with bounded $L^1$ norm.  One can also check
that $\|h_j(x)\|_{\ell^1}$ has small $L^1(X)$ norm, uniformly in $j$,
and for sufficiently large $N$, by the discussion in the preceding
paragraph.  Thus $\{f_j\}_{j = 1}^\infty$ satisfies the criterion
described in the previous section, and it follows that $\{f_j(x)\}_{j
= 1}^\infty$ converges in $\ell^1$ for almost every $x \in X$.

\section[\ Differentiability of paths]{Differentiability of paths}
\label{differentiability of paths}
\setcounter{equation}{0}

        Let $(V, \|v\|)$ be a real or complex Banach space, and let $f
: [a, b] \to V$ be a path of finite length.  Suppose that for each
$\epsilon > 0$ there is a path $g : [a, b] \to V$ of finite length
such that the length of $f - g$ on $[a, b]$ is less than or equal to
$\epsilon$ and $g$ is differentiable almost everywhere on $[a, b]$.
We would like to show that $f$ is also differentiable almost
everywhere on $[a, b]$.

        If $a \le x \le b$ and $r > 0$, then let $\delta_r(f)(x)$
be the set of difference quotients
\begin{equation}
\label{frac{f(x) - f(y)}{x - y}}
        \frac{f(x) - f(y)}{x - y},
\end{equation}
where $a \le y \le b$ and $0 < |x - y| < r$.  One can check that $f$
is differentiable at $x$ if and only if
\begin{equation}
\label{lim_{r to 0} diam delta_r(f)(x) = 0}
        \lim_{r \to 0} \diam \delta_r(f)(x) = 0,
\end{equation}
using the completeness of $V$ for the ``if'' part.  Put $h = f - g$,
and observe that
\begin{equation}
\label{diam delta_r(f)(x) le diam delta_r(g)(x) + diam delta_r(h)(x)}
        \diam \delta_r(f)(x) \le \diam \delta_r(g)(x) + \diam \delta_r(h)(x)
\end{equation}
for every $x \in [a, b]$ and $r > 0$.

        By hypothesis,
\begin{equation}
        \lim_{r \to 0} \diam \delta_r(g)(x) = 0
\end{equation}
for almost every $x \in [a, b]$.  Hence
\begin{equation}
        \lim_{r \to 0} \delta_r(f)(x) \le \sup_{r > 0} \, \diam \delta_r(h)(x)
\end{equation}
for almost every $x \in [a, b]$.

        Using maximal functions as in Section \ref{maximal functions,
2}, we get that for each $t > 0$ there is an open set $E_t(h) \subseteq
{\bf R}$ such that
\begin{equation}
        \|h(x) - h(y)\| \le t \, \|x - y|
\end{equation}
when $E_t(h)$ does not contain the interval connecting $x, y \in [a, b]$,
and
\begin{equation}
        |E_t(h)| \le 2 \, \epsilon \, t^{-1}.
\end{equation}
Here $|E_t(h)|$ denotes the Lebesgue measure of $E_t(h)$, as usual.
Thus
\begin{equation}
        \sup_{r > 0} \, \delta_r(h)(x) \le 2 \, t
\end{equation}
for every $x \in [a, b] \backslash E_t(h)$.  It follows that
\begin{equation}
        \lim_{r \to 0} \delta_r(f)(x) \le 2 \, t
\end{equation}
for almost every $x \in [a, b] \backslash E_t(h)$.  Using these
estimates for every $\epsilon, t > 0$, we get that (\ref{lim_{r to 0}
diam delta_r(f)(x) = 0}) holds for almost every $x \in [a, b]$, as
desired.

\section[\ Paths in $\ell^1$]{Paths in $\ell^1$}
\label{paths in ell^1}
\setcounter{equation}{0}

        Let $f : [a, b] \to \ell^1 = \ell^1({\bf Z}_+)$ be a path of
finite length.  Thus $f(x) = \{f_j(x)\}_{j = 1}^\infty$, where each
$f_j$ is a real or complex-valued function on $[a, b]$ of bounded
variation.  More precisely, let $l$ be a positive integer, and let
$\mathcal{P}_1, \ldots, \mathcal{P}_l$ be partitions of $[a, b]$, as
in Section \ref{lengths of paths}.  Also let $\mathcal{P}$ be a
partition of $[a, b]$ that is a common refinement of $\mathcal{P}_1,
\ldots, \mathcal{P}_l$.  If $\Lambda_a^b(f, \mathcal{P})$ denotes the
approximation to the length of $f$ associated to $\mathcal{P}$, and
similarly for the $f_j$'s and $\mathcal{P}_j$'s, then
\begin{equation}
\label{sum_{j = 1}^l Lambda_a^b(f_j, mathcal{P}_j) le ...}
        \sum_{j = 1}^l \Lambda_a^b(f_j, \mathcal{P}_j)
         \le \sum_{j = 1}^l \Lambda_a^b(f_j, \mathcal{P})
          \le \Lambda_a^b(f, \mathcal{P}).
\end{equation}
Hence
\begin{equation}
        \sum_{j = 1}^l \Lambda_a^b(f_j, \mathcal{P}_j) \le \Lambda_a^b(f),
\end{equation}
where $\Lambda_a^b(f)$ denotes the length of $f$ on $[a, b]$.  This
implies that
\begin{equation}
        \sum_{j = 1}^l \Lambda_a^b(f_j) \le \Lambda(f).
\end{equation}
because $\mathcal{P}_1, \ldots, \mathcal{P}_l$ are arbitrary
partitions of $[a, b]$.  Therefore
\begin{equation}
        \sum_{j = 1}^\infty \Lambda_a^b(f_j) \le \Lambda(f),
\end{equation}
because $l \ge 1$ is arbitrary.

        Similarly, if $\mathcal{P}$ is any partition of $[a, b]$, then
\begin{equation}
\label{Lambda_a^b(f, mathcal{P}) le ...}
        \Lambda_a^b(f, \mathcal{P})
         = \sum_{j = 1}^\infty \Lambda_a^b(f_j, \mathcal{P}_j)
          \le \sum_{j = 1}^\infty \Lambda_a^b(f_j).
\end{equation}
This implies that
\begin{equation}
        \Lambda_a^b(f) \le \sum_{j = 1}^\infty \Lambda_a^b(f_j).
\end{equation}
It follows that
\begin{equation}
        \Lambda_a^b(f) = \sum_{j = 1}^\infty \Lambda_a^b(f_j),
\end{equation}
by the remarks in the preceding paragraph.

        Let $N$ be a large positive integer, and put
\begin{eqnarray}
        g_j(x) = f_j(x), \, h_j(x) = 0 & \hbox{when} & j \le N, \\
        g_j(x) = 0, \, h_j(x) = f_j(x) & \hbox{when} & j > N. \nonumber
\end{eqnarray}
If $g(x) = \{g_j(x)\}_{j = 1}^\infty$, $h(x) = \{h_j(x)\}_{j =
1}^\infty$, then $g(x), h(x) \in \ell^1$ for every $x$ in $[a, b]$, and
\begin{equation}
        f(x) = g(x) + h(x).
\end{equation}
Observe that $f, g : [a, b] \to \ell^1$ have finite length, and that
the length of $h$ on $[a, b]$ tends to $0$ as $N \to \infty$.  We also
know that $g$ is differentiable almost everywhere on $[a, b]$, by the
corresponding results for real or complex-valued functions.  It
follows that $f$ is also differentiable almost everywhere on $[a, b]$,
as in the previous section.

\section[\ $L^p$-Valued functions]{$L^p$-Valued functions}
\label{L^p-valued functions}
\setcounter{equation}{0}

        Let $(Y, \mathcal{B}, \nu)$ be a $\sigma$-finite measure
space, and consider ${\bf R} \times Y$, equipped with the product
measure corresponding to Lebesgue measure on the real line.  A
function $F(x, y) \in L^p({\bf R} \times Y)$, $1 \le p < \infty$, may
be considered as representing an $L^p$ function on ${\bf R}$ with
values in $L^p(Y)$.  Put
\begin{equation}
        F_p(x) = \Big(\int_Y |F(x, y)|^p \, d\nu(y)\Big)^{1/p},
\end{equation}
which is the $L^p(Y)$ norm of $F(x, y)$ in $y$.  By Fubini's theorem,
\begin{equation}
        \Big(\int_{{\bf R} \times Y} |F(x, y)|^p \, dx \, d\nu(y)\Big)^{1/p}
         = \Big(\int_{\bf R} F_p(x)^p \, dx\Big)^{1/p}.
\end{equation}
Thus the $L^p({\bf R} \times Y)$ norm of $F(x, y)$ is the
same as starting with the $L^p(Y)$ norm of $F(x, y)$ in $y$, and then
taking the $L^p({\bf R})$ norm of the result in $x$.

        Put
\begin{equation}
        \ L(F)(x) = \limsup_{r \to 0} \frac{1}{2 r} \int_{x - r}^{x + r}
             \Big(\int_Y |F(t, y) - F(x, y)|^p \, d\nu(y)\Big)^{1/p} \, dt.
\end{equation}
As in Section \ref{lebesgue's theorem}, we would like to say that
\begin{equation}
        L(F)(x) = 0
\end{equation}
for almost every $x \in {\bf R}$.  As usual, there are two main
ingredients in the proof.  The first is that this condition holds for
a dense class of functions $F(x, y)$ in $L^p(x, y)$.  In this case,
one can use finite linear combinations of functions of the form $f(x)
\, g(y)$, where $f(x) \in L^p({\bf R})$ and $g(y) \in L^p(Y)$.  If one
also takes $f(x)$ to be continuous, then the limit is equal to $0$ at
every $x \in {\bf R}$.  If $Y$ is a locally compact Hausdorff
topological space and $\nu$ is a Borel measure on $Y$ with suitable
regularity properties, then one can use continuous functions on ${\bf
R} \times Y$ with compact support as the dense class.  Again the limit
is equal to $0$ for every $x \in {\bf R}$ in this situation.  The
second main ingredient is an estimate for an appropriate maximal
function, which reduces here to the Hardy--Littlewood maximal function
of $F_p(x) \in L^p({\bf R})$.

\section[\ Continuous $L^p$-valued functions]{Continuous $L^p$-valued functions}
\label{continuous L^p-valued functions}
\setcounter{equation}{0}

        Let $(Y, \mathcal{B}, \nu)$ be a measure space, and let $f$ be
a continuous function on the real line with values in $L^p(Y)$, $1 \le
p < \infty$.  If $g \in L^p(Y)$, then
\begin{equation}
        \{y \in Y : |g(y)| \ge 1/n\}
\end{equation}
has finite measure for each $n \in {\bf Z}_+$, and hence
\begin{equation}
        \{y \in Y : g(y) \ne 0\}
\end{equation}
is $\sigma$-finite.  Applying this to $f(r)$ for each rational number
$r$, we get that there is a $\sigma$-finite measurable set $Y_0
\subseteq Y$ such that $f(r) = 0$ on $Y \backslash Y_0$ for every $r
\in {\bf Q}$.  This implies that $f(r) = 0$ almost everywhere on $Y
\backslash Y_0$ for every $r \in {\bf R}$, because $f : {\bf R} \to
L^p(Y)$ is continuous.  Thus we may as well suppose that $Y$ is
$\sigma$-finite.

        Let us now restrict our attention to the case where $f$ has
compact support on ${\bf R}$.  More precisely, let $I = [a, b]$ be a
closed interval in the real line such that $f(x) = 0$ when $x \in {\bf
R} \backslash [a, b]$.  Consider the product ${\bf R} \times L^p(Y)$
with the product measure associated to Lebesgue measure on ${\bf R}$,
as in the preceding section.  We would like to check that there is an
$F(x, y) \in L^p({\bf R} \times Y)$ such that
\begin{equation}
\label{f(x) = F(x, cdot)}
        f(x) = F(x, \cdot)
\end{equation}
as elements of $L^p(Y)$ for almost every $x \in {\bf R}$.  In this case,
\begin{eqnarray}
\label{int_{R times Y} |F(x, y)|^p dx d nu(y) = int_R ||f(x)||_{L^p(Y)}^p dx}
        \int_{{\bf R} \times Y} |F(x, y)|^p \, dx \, d\nu(y)
         & = & \int_{\bf R} \Big(\int_Y |F(x, y)|^p \, d\nu(y)\Big) \, dx \\
         & = & \int_{\bf R} \|f(x)\|_{L^p(Y)}^p \, dx.  \nonumber
\end{eqnarray}
In particular, the $L^p({\bf R} \times Y)$ of $F(x, y)$ would be
bounded by a constant multiple of the supremum norm of
$\|f(x)\|_{L^p(Y)}$ on $I$.  If $F(x, y), \widetilde{F}(x, y) \in
L^p({\bf R} \times Y)$ both satisfy (\ref{f(x) = F(x, cdot)}) for
almost every $x \in {\bf R}$, then it follows that $F(x, y) =
\widetilde{F}(x, y)$ for almost every $(x, y) \in {\bf R} \times Y$.

        If $f(x) = \phi(x) \, g$ for some real or complex-valued
function $\phi(x)$ with compact support on ${\bf R}$ and some $g \in
L^p(Y)$, then we can simply take $F(x, y) = \phi(x) \, g(y)$.
Similarly, if $f(x)$ is a finite linear combination of $L^p(Y)$-valued
functions on ${\bf R}$ of this form, then it is easy to get $F(x, y)$.
Otherwise, one can approximate $f(x)$ by a sequence $\{f_j(x)\}_{j =
1}^\infty$ of $L^p(Y)$-valued functions of this type with respect to
the supremum norm of $\|f(x)\|_{L^p(Y)}$ on $I$.  By construction,
$f_j(x)$ corresponds to a function $F_j(x, y)$ in $L^p({\bf R} \times
Y)$ for each $j$.  Moreover, $\{F_j(x, y)\}_{j = 1}^\infty$ is a
Cauchy sequence in $L^p({\bf R} \times Y)$, because of (\ref{int_{R
times Y} |F(x, y)|^p dx d nu(y) = int_R ||f(x)||_{L^p(Y)}^p dx}).
Hence $\{F_j(x, y)\}_{j = 1}^\infty$ converges to a function $F(x, y)$
in $L^p({\bf R} \times Y)$.  It is not too difficult to verify that
this function $F(x, y)$ satisfies (\ref{f(x) = F(x, cdot)}), as
desired.

\section[\ Lipschitz $L^p$-valued functions]{Lipschitz $L^p$-valued functions}
\label{lipschitz L^p-valued functions}
\setcounter{equation}{0}

        Let $(Y, \mathcal{B}, \nu)$ be a $\sigma$-finite measure
space, and suppose that $f : {\bf R} \to L^p(Y)$ is a Lipschitz
mapping for some $1 < p < \infty$.  It will be convenient to ask also
at first that $f(x)$ have compact support in ${\bf R}$, which is to
say that there is a closed interval $[a, b]$ in the real line such
that $f(x) = 0$ when $x \in {\bf R} \backslash [a, b]$.  Let $F(x, y)$
be the function in $L^p({\bf R} \times Y)$ that corresponds to $f(x)$
as in the previous section.  Because $f(x)$ has compact support, the
ordinary Lipschitz condition implies an integrated Lipschitz condition
of the form
\begin{equation}
 \Big(\int_{\bf R} \|f(x + h) - f(x)\|_{L^p(Y)}^p \, dx\Big)^{1/p} \le C \, |h|
\end{equation}
for some $C \ge 0$ and every $h \in {\bf R}$.  This implies that
\begin{equation}
 \Big(\int_{{\bf R} \times Y} |F(x + h, y) - F(x, h)|^p \, dx \,
                                         d\nu(y)\Big)^{1/p} \le C \, |h|.
\end{equation}

        Let $q$ be the exponent conjugate to $p$, so that $1/p + 1/q =
1$.  If $h \in {\bf R}$, $h \ne 0$, and $\Phi(x, y) \in L^q({\bf R}
\times Y)$, then put
\begin{equation}
\label{lambda_h(Phi) = ...}
 \lambda_h(\Phi) = \int_{{\bf R} \times Y} \frac{F(x + h, y) - F(x, y)}{h} \,
                                              \Phi(x, y) \, dx \, d\nu(y).
\end{equation}
This defines a bounded linear functional on $L^q({\bf R} \times Y)$,
with dual norm less than or equal to $C$, by H\"older's inequality.
We also have that
\begin{equation}
\label{lambda_h(Phi) = ..., 2}
        \lambda_h(\Phi) = - \int_{{\bf R} \times Y} F(x, y) \,
                \frac{\Phi(x, y) - \Phi(x - h, y)}{h} \, dx \, d\nu(y),
\end{equation}
using the change oe variables $x \mapsto x - h$.  If
\begin{equation}
        \Phi(x, y) = \phi(x) \, \psi(y),
\end{equation}
where $\phi(x)$ is a continuously-differentiable real or
complex-valued function on the real line with compact support and
$\psi(y) \in L^q(Y)$, then we get that
\begin{equation}
        \lim_{h \to 0} \lambda_h(\Phi) = - \int_{{\bf R} \times Y} F(x, y) \,
                                     \phi'(x) \, \psi(y) \, dx \, d\nu(y).
\end{equation}
Similarly,
\begin{equation}
\label{lim_{h to 0} lambda_h(Phi)}
        \lim_{h \to 0} \lambda_h(\Phi)
\end{equation}
exists when $\Phi(x, y)$ is a finite linear combination of functions
of this form.  As in Section \ref{uniform boundedness, 4}, it follows
that (\ref{lim_{h to 0} lambda_h(Phi)}) exists for all $\Phi(x, y) \in
L^q({\bf R} \times Y)$, since it exists for a dense linear subspace of
$L^q({\bf R} \times Y)$, and since the dual norms of $\lambda_h$,
$h \in {\bf R} \backslash \{0\}$, are bounded.

        Thus (\ref{lim_{h to 0} lambda_h(Phi)}) defines a bounded
linear functional on $L^q({\bf R} \times Y)$.  By the Riesz
representation theorem, there is a function $G(x, y)$ in $L^p({\bf R}
\times Y)$ such that
\begin{equation}
\label{lim_{h to 0} lambda_h(Phi) = int_{R times Y} G(x,y) Phi(x,y) dx dnu(y)}
        \lim_{h \to 0} \lambda_h(\Phi) = \int_{{\bf R} \times Y} G(x, y) \,
                                              \Phi(x, y) \, dx \, d\nu(y)
\end{equation}
for every $\Phi(x, y) \in L^q({\bf R} \times Y)$.  In particular,
\begin{eqnarray}
\label{int_{R times Y} F(x, y) phi'(x) psi(y) dx d nu(y) = ...}
\lefteqn{\int_{{\bf R} \times Y} F(x, y) \, \phi'(x) \, \psi(y) \, dx
                                                           \, d\nu(y) =} \\
 & & - \int_{{\bf R} \times Y} G(x, y) \, \phi(x) \, \psi(y) \, dx \, d\nu(y)
                                                             \nonumber
\end{eqnarray}
when $\phi(x)$ is a real or complex-valued continuously-differentiable
function on the real line with compact support and $\psi(y) \in L^q(Y)$.
Put
\begin{equation}
        f_\psi(x) = \int_Y f(x)(y) \, \psi(y) \, d\nu(y)
\end{equation}
for each $x \in {\bf R}$.  More precisely, $f(x) \in L^p(Y)$ for every
$x \in {\bf R}$, and $f_\psi(x)$ is the integral of the product of
this function with $\psi \in L^q(Y)$ over $Y$.  Thus $f_\psi(x)$ is a
Lipschitz function on ${\bf R}$ with compact support for every $\psi
\in L^q(Y)$, because $f : {\bf R} \to L^p(Y)$ is a Lipschitz mapping
with compact support.  Using (\ref{int_{R times Y} F(x, y) phi'(x)
psi(y) dx d nu(y) = ...}), we get that
\begin{equation}
        \int_{\bf R} f_\psi(x) \, \phi'(x) \, dx 
 = - \int_{{\bf R} \times Y} G(x, y) \, \phi(x) \, \psi(y) \, dx \, d\nu(y)
\end{equation}
for every $\phi(x)$, $\psi(y)$ as before.  This implies that
\begin{equation}
\label{f_psi'(x) = int_Y G(x, y) psi(y) d nu(y)}
        f_\psi'(x) = \int_Y G(x, y) \, \psi(y) \, d\nu(y)
\end{equation}
for every $\psi \in L^q(Y)$ in the sense of distributions, as in
Section \ref{L^p derivatives}.  Hence
\begin{equation}
 f_\psi(t) - f_\psi(r) = \int_r^t \int_Y G(x, y) \, \psi(y) \, d\nu(y) \, dx
\end{equation}
for every $r, t \in {\bf R}$ with $r < t$ and $\psi \in L^q(Y)$.
It follows that
\begin{equation}
        f(t) - f(r) = \int_r^t G(x, \cdot) \, dx
\end{equation}
when $r < t$, where both sides of the equation are elements of $L^p(Y)$.

        Now that we have this expression for differences of the values
of $f$, one can use the analogue of Lebesgue's theorem in this context
to conclude that $f$ is differentiable almost everywhere as an
$L^p(Y)$-valued function on the real line.  This works as well for
Lipschitz mappings from the real line into $L^p(Y)$ that may not have
compact support, since the problem is local.  This also works for
paths of finite length in $L^p(Y)$, $1 < p < \infty$, because of the
approximation arguments in Sections \ref{maximal functions, 2} and
\ref{vector-valued functions}.

\section[\ More duality]{More duality}
\label{more duality}
\setcounter{equation}{0}

        Let $(Y, \mathcal{B}, \nu)$ be a $\sigma$-finite measure
space, and let $f$ be a continuous function from the real line into
$L^p(Y)$, $1 < p < \infty$.  Suppose also that $f$ has compact support
in ${\bf R}$, and let $q$ be the exponent conjugate to $p$, so that
$1/p + 1/q = 1$.  We would like to define a bounded linear functional
on $L^q({\bf R} \times Y)$ directly by
\begin{equation}
\label{Lambda(Phi) = int_{bf R} (int_Y f(x)(y) Phi(x, y) d nu(y)) dx}
 \Lambda(\Phi) = \int_{\bf R} \Big(\int_Y f(x)(y) \, \Phi(x, y) \, d\nu(y)\Big)
                                                                    \, dx.
\end{equation}
Because of H\"older's inequality,
\begin{equation}
        \quad  \biggl|\int_Y f(x)(y) \, \Phi(x, y) \, d\nu(y)\biggr|
 \le \|f(x)\|_{L^p(Y)} \, \Big(\int_Y |\Phi(x, y)|^q \, d\nu(y)\Big)^{1/q}
\end{equation}
and
\begin{eqnarray}
\lefteqn{\int_{\bf R} \|f(x)\|_{L^p(Y)} \, \Big(\int_Y |\Phi(x, y)|^q \,
                                                d\nu(y)\Big)^{1/q} \, dx} \\
 & \le & \Big(\int_{\bf R} \|f(x)\|_{L^p(Y)}^p \, dx\Big)^{1/p} \, 
  \Big(\int_{{\bf R} \times Y} |\Phi(x, y)|^q \, dx \, d\nu(y)\Big)^{1/q}.
                                                                   \nonumber
\end{eqnarray}
However, one should be a bit careful about the measurability of
\begin{equation}
        \int_Y f(x)(y) \, \Phi(x, y) \, d\nu(y)
\end{equation}
as a function of $x$.  If $\Phi(x, y) = \phi(x) \, \psi(y)$ for some
$\phi(x) \in L^q({\bf R})$, $\psi(y) \in L^q(Y)$, then this reduces to
\begin{equation}
        \phi(x) \int_Y f(x)(y) \, \psi(y) \, d\nu(y).
\end{equation}
The continuity of $f : {\bf R} \to L^p(Y)$ implies that
\begin{equation}
        \int_Y f(x)(y) \, \psi(y) \, d\nu(y)
\end{equation}
is continuous in $x$, and so there is no problem in this case.
Because linear combinations of functions of this type are dense in in
$L^q({\bf R} \times Y)$, one can use this to extend $\Lambda(\Phi)$ to
all $\Phi \in L^q({\bf R} \times Y)$.  Similarly,
\begin{equation}
 \lambda_h(\Phi) = \int_{\bf R} \Big(\int_Y \frac{f(x + h)(y) - f(x)(y)}{h}
                                         \, \Phi(x, y) \, d\nu(y)\Big) \, dx
\end{equation}
can be defined more directly as a bounded linear functional on
$L^q({\bf R} \times Y)$ for each $h \in {\bf R} \backslash \{0\}$.
Equivalently,
\begin{equation}
        \lambda_h(\Phi) = - \int_{\bf R} \Big(\int_Y f(x)(h) \,
          \frac{\Phi(x, y) - \Phi(x - h, y)}{h} \, d\nu(y)\Big) \, dx.
\end{equation}
If $\Phi(x, y) = \phi(x) \, \psi(y)$, where now $\phi(x)$ is a
continuously-differentiable function on ${\bf R}$ with compact support
and $\psi(y) \in L^q(Y)$, then it follows that
\begin{equation}
        \lim_{h \to 0} \lambda_h(\Phi) =
 - \int_{\bf R} \phi'(x) \Big(\int_Y f(x)(y) \, \psi(y) \, d\nu(y)\Big) \, dx.
\end{equation}
At this point, one can continue as in the preceding section when $f :
{\bf R} \to L^p(Y)$ is Lipschitz.

\section[\ $\ell^p$-Valued functions]{$\ell^p$-Valued functions}
\label{ell^p-valued functions}
\setcounter{equation}{0}

        Of course, the arguments in the previous sections can be
simplified when the functions take values in $\ell^p = \ell^p({\bf
Z}_+)$, and there are some commonalities with $p = 1$.  Suppose that
$f(x) = \{f_j(x)\}_{j = 1}^\infty$ is a Lipschitz function on the real
line with values in $\ell^p$, $1 \le p \le \infty$.  In particular,
$f_j(x)$ is a Lipschitz function on ${\bf R}$ for each $j$, and hence
is differentiable almost everywhere.  Using the Lipschitz condition
for $f : {\bf R} \to \ell^p$, one can check that $\{f_j'(x)\}_{j =
1}^\infty \in \ell^p$ for every $x \in {\bf R}$ such that $f_j'(x)$
exists for each $j$, with $\ell^p$ norm bounded by the Lipschitz
constant for $f$.  We also have that
\begin{equation}
        f_j(t) - f_j(r) = \int_r^t f_j'(x) \, dx
\end{equation}
for every $r, t \in {\bf R}$ with $r < t$.  If $p < \infty$, then one
can use this to show that $f$ is differentiable almost everywhere on
${\bf R}$ as a mapping into $\ell^p$, with derivative given by
$\{f_j'(x)\}_{j = 1}^\infty$.  As usual, it is convenient to restrict
one's attention initially to functions $f$ with compact support, so
that $\|\{f_j'(x)\}\|_{\ell^p} \in L^p({\bf R})$.  As in the $p = 1$
case, one can approximate $f$ by functions with only finitely many
nonzero components, for which differentiability almost everywhere is
already known.  One can then use maximal function estimates to show
that the errors are small most of the time.

        Note that a Lipschitz mapping from the real line into a
separable Hilbert space is differentiable almost everywhere, by the $p
= 2$ case.  This can be extended to paths of finite length in a
separable Hilbert space, because of the approximation arguments in
Sections \ref{maximal functions, 2} and \ref{vector-valued functions}.
As in Section \ref{continuity conditions}, any path of finite length
is continuous at all but finitely or countably many elements of its
domain, and hence is contained in a separable subspace of the range.
This implies that a path of finite length in any Hilbert space is
differentiable almost everywhere, because it is contained in a
separable Hilbert subspace.

\section[\ Products and $\sigma$-subalgebras]{Products and $\sigma$-subalgebras}
\label{products, sigma-subalgebras}
\setcounter{equation}{0}

        Let $(X_1, \mathcal{A}_1, \mu_1)$, $(X_2, \mathcal{A}_2,
\mu_2)$ be probability spaces, and let $X = X_1 \times X_2$ be their
Cartesian product, with the product measure $\mu = \mu_1 \times
\mu_2$.  Also let $\mathcal{B}_1$, $\mathcal{B}_2$ be
$\sigma$-subalgebras of $\mathcal{A}_1$, $\mathcal{A}_2$,
respectively, and let $\mathcal{B}$ be the corresponding
$\sigma$-subalgebra of the $\sigma$-algebra of measurable subsets of
$X$.  If $\phi_1(x_1) \in L^1(X_1)$, $\phi_2(x_2) \in L^1(X_2)$, then
$\phi(x_1, x_2) = \phi_1(x_1) \, \phi_2(x_2) \in L^1(X)$, and we would
like to check that
\begin{equation}
 E_X(\phi \mid \mathcal{B}) = E_{X_1}(\phi_1 \mid \mathcal{B}_1) \,
                                  E_{X_2}(\phi_2 \mid \mathcal{B}_2),
\end{equation}
where the subscripts of $E$ are included to indicate the spaces on
which the conditional expectations are taken.  Both sides of the
equation are measurable with respect to $\mathcal{B}$, and so it
suffices to verify that
\begin{equation}
        \int_B E_X(\phi \mid \mathcal{B}) \, d\mu
              = \int_B E_{X_1}(\phi_1 \mid \mathcal{B}_1) \,
                            E_{X_2}(\phi_2 \mid \mathcal{B}_2) \, d\mu
\end{equation}
for every $B \in \mathcal{B}$.  This reduces to
\begin{equation}
        \int_B \phi \, d\mu
              = \int_B E_{X_1}(\phi_1 \mid \mathcal{B}_1) \,
                            E_{X_2}(\phi_2 \mid \mathcal{B}_2) \, d\mu,
\end{equation}
by the definition of the conditional expectation.  If $B = B_1 \times
B_2$ with $B_1 \in \mathcal{B}_1$, $B_2 \in \mathcal{B}_2$, then both
sides of this equation are equal to
\begin{equation}
 \Big(\int_{B_1} \phi_1 \, d\mu_1\Big) \Big(\int_{B_2} \phi_2 \, d\mu_2\Big),
\end{equation}
using the definition of the conditional expectation again.  This
implies that the previous equation holds when $B$ is the union of
finitely many pairwise-disjoint sets of the form $B_1 \times B_2$,
with $B_1 \in \mathcal{B}_1$ and $B_2 \in \mathcal{B}_2$.  The
analogous statement for any $B \in \mathcal{B}$ follows by
approximation.  If $\mathcal{B}_1$ or $\mathcal{B}_2$ is generated by
a partition of $X_1$ or $X_2$ into finitely or countably many
measurable sets, then every $B \in \mathcal{B}$ can be expressed as
the union of finitely or countably many disjoint sets of the form $B_1
\times B_2$, with $B_1 \in \mathcal{B}_1$ and $B_2 \in \mathcal{B}_2$,
as in Section \ref{partitions, products}, and the approximation is
much simpler.

\section[\ $\sigma$-Subalgebras and vectors]{$\sigma$-Subalgebras and vectors}
\label{sigma-subalgebras, vectors}
\setcounter{equation}{0}

        Let $(X, \mathcal{A}, \mu)$ be a probability space, and let
$\mathcal{B}$ be a $\sigma$-subalgebra of $\mathcal{A}$.  Also let $V$
be a finite-dimensional real or complex vector space with a norm,
which can be identified with ${\bf R}^n$ or ${\bf C}^n$ for some $n$
using a basis.  Thus a $V$-valued function $f(x)$ on $X$ corresponds
an $n$-tuple $(f_1(x), \ldots, f_n(x))$ of real or complex-valued
functions on $X$.  Such a function is considered to be integrable when
its components $f_1(x), \ldots, f_n(x)$ are integrable, in which case
the integral is defined by integrating the components separately.
Similarly, the conditional expectation of a $V$-valued function $f$ on
$X$ may be defined by applying the conditional expectation to the
components of $f$.

        Let $\lambda$ be a linear functional on $V$, so that
$\lambda(v)$ can be expressed by a linear combination of the
components of $v$.  If $f(x)$ is an integrable $V$-valued function on
$X$, then $\lambda(f(x))$ is an integrable real or complex-valued
function on $X$, and
\begin{equation}
\label{lambda(int_X f(x) d mu(x)) = int_X lambda(f(x)) d mu(x)}
 \lambda\Big(\int_X f(x) \, d\mu(x)\Big) = \int_X \lambda(f(x)) \, d\mu(x).
\end{equation}
If $\|v\|$ is a norm on $V$, and $\|\lambda\|_*$ is the corresponding
dual norm on $V^*$, then it follows that
\begin{eqnarray}
        \biggl|\lambda\Big(\int_X f(x) \, d\mu(x)\Big)\biggr|
         & \le & \int_X |\lambda(f(x))| \, d\mu(x) \\
         & \le & \|\lambda\|_* \, \int_X \|f(x)\| \, d\mu(x). \nonumber
\end{eqnarray}
This implies that
\begin{equation}
        \biggl\|\int_X f(x) \, d\mu(x)\biggr\| \le \int_X \|f(x)\| \, d\mu(x),
\end{equation}
by the Hahn--Banach theorem.  The same conclusion could also be
obtained by approximating the integral by finite sums.

        Similarly,
\begin{equation}
\label{lambda(E(f mid mathcal{B})) = E(lambda(f) mid mathcal{B})}
 \lambda(E(f \mid \mathcal{B})) = E(\lambda(f) \mid \mathcal{B}),
\end{equation}
and hence
\begin{equation}
\label{|lambda(E(f mid mathcal{B}))| le ||lambda||_* E(||f|| mid mathcal{B})}
 |\lambda(E(f \mid \mathcal{B}))| \le E(|\lambda(f)| \mid \mathcal{B})
                                \le \|\lambda\|_* \, E(\|f\| \mid \mathcal{B}).
\end{equation}
This implies that
\begin{equation}
\label{||E(f mid mathcal{B})|| le E(||f|| mid mathcal{B})}
        \|E(f \mid \mathcal{B})\| \le E(\|f\| \mid \mathcal{B}).
\end{equation}
More precisely, if (\ref{|lambda(E(f mid mathcal{B}))| le ||lambda||_*
E(||f|| mid mathcal{B})}) holds at some point $x \in X$ for every
linear functional $\lambda$ on $V$, then (\ref{||E(f mid mathcal{B})||
le E(||f|| mid mathcal{B})}) also holds at $x$, by the Hahn--Banach
theorem.  This works as well when (\ref{|lambda(E(f mid mathcal{B}))| le
||lambda||_* E(||f|| mid mathcal{B})}) holds for every $\lambda$ in a
dense subset of
\begin{equation}
\label{lambda in V^* : ||lambda||_* = 1}
        \{\lambda \in V^* : \|\lambda\|_* = 1\}.
\end{equation}
Because $V$ and hence $V^*$ are finite-dimensional, there is a
countable dense set in (\ref{lambda in V^* : ||lambda||_* = 1}).  If
(\ref{|lambda(E(f mid mathcal{B}))| le ||lambda||_* E(||f|| mid
mathcal{B})}) holds almost everywhere on $X$ for each $\lambda \in
V^*$, then it holds simultaneously for a countable set of $\lambda$'s
almost everywhere on $X$.  This implies that (\ref{||E(f mid
mathcal{B})|| le E(||f|| mid mathcal{B})}) holds almost everywhere on
$X$, as desired.

\section[\ Martingales and products]{Martingales and products}
\label{martingales, products}
\setcounter{equation}{0}

        Let $(X, \mathcal{A}, \mu)$, $(Y, \mathcal{B}, \nu)$ be
probability spaces, and suppose that their Cartesian product $X \times
Y$ is equipped with the product probability measure $\mu \times \nu$.
Also let $\mathcal{B}_1 \subseteq \mathcal{B}_2 \subseteq \cdots$ be
an increasing sequence of $\sigma$-subalgebras of $\mathcal{A}$, and
let $\widehat{\mathcal{B}}_j$ be the $\sigma$-subalgebra of the
$\sigma$-algebra of measurable subsets of $X \times Y$ that corresponds
to $\mathcal{B}_j$ on $X$ and $\mathcal{B}$ on $Y$ in the product space.
As before, a function $F(x, y) \in L^p(X \times Y)$, $1 \le p < \infty$,
may be considered as representing an $L^p$ function on $X$ with values in
$L^p(Y)$, and thus a martingale on $X \times Y$ with respect to
$\widehat{\mathcal{B}}_j$ may be considered as representing a type of
vector-valued martingale on $X$ with respect to $\mathcal{B}_j$.

        If $F(x, y) \in L^1(X \times Y)$, then put
\begin{equation}
        I_Y(F)(x) = \int_Y F(x, y) \, d\nu(y).
\end{equation}
Let us check that
\begin{equation}
        I_Y(E_{X \times Y}(F \mid \widehat{\mathcal{B}}_j))
           = E_X(I_Y(F) \mid \mathcal{B}_j)
\end{equation}
for each $j$, where the subscripts of $E$ indicate the spaces on which
the conditional expectations are taken.  Both sides of the equation
are measurable functions on $X$ with respect to $\mathcal{B}_j$,
and so it is enough to show that
\begin{equation}
\label{int_A I_Y(E_{X times Y}(F mid widehat{mathcal{B}}_j) d mu = ...}
 \int_A I_Y(E_{X \times Y}(F \mid \widehat{\mathcal{B}}_j) \, d\mu
 = \int_A E_X(I_Y(F) \mid \mathcal{B}_j) \, d\mu
\end{equation}
for every $A \in \mathcal{B}_j$.  Of course,
\begin{eqnarray}
 \int_A I_Y(E_{X \times Y}(F \mid \widehat{\mathcal{B}}_j) \, d\mu
 & = & \int_{A \times Y} E_{X \times Y}(F \mid \widehat{\mathcal{B}}_j)
                                                     \, d(\mu \times \nu) \\
 & = & \int_{A \times Y} F \, d(\mu \times \nu), \nonumber
\end{eqnarray}
because $A \times Y \in \widehat{\mathcal{B}}_j$.  Similarly,
\begin{equation}
 \int_A E_X(I_Y(F) \mid \mathcal{B}_j) \, d\mu = \int_A I_Y(F) \, d\mu
   = \int_{A \times Y} F \, d(\mu \times \nu).
\end{equation}

        Let us say that a measurable function $F(x, y) \in L^1(X
\times Y)$ is \emph{nice} if there are finitely many pairwise-disjoint
measurable subsets $B_1, \ldots, B_n$ of $Y$ with positive measure
such that $F(x, y)$ is constant in $y$ on $B_k$ for $k = 1, \ldots,
n$.  If $\phi_k(x) = F(x, y)$ when $y \in B_k$, then $\phi_k(x) \in
L^1(X)$ for each $k$, and
\begin{equation}
\label{F(x, y) = sum_{k = 1}^n phi_k(x) {bf 1}_{B_k}(y)}
        F(x, y) = \sum_{k = 1}^n \phi_k(x) \, {\bf 1}_{B_k}(y),
\end{equation}
where ${\bf 1}_{B_k}(y)$ is the indicator function associated to $B_k$
on $Y$, equal to $1$ when $y \in B_k$ and to $0$ when $y \in Y
\backslash B_k$.  In this case,
\begin{equation}
\label{E_{X times Y}(F mid widehat{mathcal{B}}_j)(x, y) = ...}
        E_{X \times Y}(F \mid \widehat{\mathcal{B}}_j)(x, y)
 = \sum_{k = 1}^n E_X(\phi_k \mid \mathcal{B}_j)(x) \, {\bf 1}_{B_k}(y),
\end{equation}
as in Section \ref{products, sigma-subalgebras}.  In effect, $F$
corresponds to a function on $X$ with values in an $n$-dimensional
vector space under these conditions.

        Suppose that $F(x, y) \in L^p(X \times Y)$, $1 \le p <
\infty$, and put
\begin{equation}
        N_p(F)(x) = \Big(\int_Y |F(x, y)|^p \, d\nu(y)\Big)^{1/p}.
\end{equation}
Thus $N_p(F) \in L^p(X)$, and
\begin{equation}
        \quad  \Big(\int_X N_p(F)(x)^p \, d\mu(x)\Big)^{1/p} =
 \Big(\int_{X \times Y} |F(x, y)|^p \, d(\mu \times \nu)(x, y)\Big)^{1/p}.
\end{equation}
We would like to check that
\begin{equation}
\label{N_p(E_{X x Y}(F | widehat{mathcal{B}}_j)) le E_X(N_p(F) | mathcal{B}_j)}
        N_p(E_{X \times Y}(F \mid \widehat{\mathcal{B}}_j))(x)
         \le E_X(N_p(F) \mid \mathcal{B}_j)(x)
\end{equation}
for almost every $x \in X$ and each $j \ge 1$.  If $p = 1$, then
$N_1(F) = I_Y(|F|)$, and
\begin{eqnarray}
        I_Y(|E_{X \times Y}(F \mid \widehat{\mathcal{B}}_j)|)
          & \le & I_Y(E_{X \times Y}(|F| \mid \widehat{\mathcal{B}}_j)) \\
           & = & E_X(I_Y(|F|) \mid \mathcal{B}_j). \nonumber
\end{eqnarray}
If $p > 1$ and $F$ is nice, then (\ref{N_p(E_{X x Y}(F |
widehat{mathcal{B}}_j)) le E_X(N_p(F) | mathcal{B}_j)}) follows from
the discussion in the preceding section.  More precisely, one can take
$V$ to be the $n$-dimensional vector space spanned by ${\bf 1}_{B_1},
\ldots, {\bf 1}_{B_n}$, equipped with the $L^p(Y)$ norm.  Otherwise,
one can approximate $F$ by nice functions in $L^p(X \times Y)$.

        If $\{F_j\}_{j = 1}^\infty$ is a martingale on $X \times Y$
with respect to the $\widehat{\mathcal{B}}_j$'s such that $F_j \in
L^p(X \times Y)$ for each $j$, $1 \le p < \infty$, then it follows
that $\{N_p(F_j)\}_{j = 1}^\infty$ is a submartingale on $X$ with
respect to the $\mathcal{B}_j$'s.  This leads to the same type of
maximal function estimates as before.  If $\{F_j\}_{j = 1}^\infty$
converges in $L^p(X \times Y)$, then one may conclude that $\{F_j(x,
\cdot)\}_{j = 1}^\infty$ converges in $L^p(Y)$ for almost every $x \in
X$.  In particular, this holds when $1 < p < \infty$ and the norm of
$F_j(x, y)$ in $L^p(X \times Y)$ is uniformly bounded in $j$.

        Instead of (\ref{N_p(E_{X x Y}(F | widehat{mathcal{B}}_j)) le
E_X(N_p(F) | mathcal{B}_j)}), it is easier to show that
\begin{equation}
        I_Y(|E_{X \times Y}(F \mid \widehat{\mathcal{B}}_j)|^p)
              \le E_X(I_Y(|F|^p) \mid \mathcal{B}_j)
\end{equation}
almost everywhere on $X$.  As in the $p = 1$ case, one has that
\begin{eqnarray}
        I_Y(|E_{X \times Y}(F \mid \widehat{\mathcal{B}}_j)|^p)
         & \le & I_Y(E_{X \times Y}(|F|^p \mid \widehat{\mathcal{B}}_j)) \\
          & = & E_X(I_Y(|F|^p) \mid \mathcal{B}_j) \nonumber
\end{eqnarray}
when $p > 1$.  If $\{F_j\}_{j = 1}^\infty$ is a martingale on $X
\times Y$ with respect to the $\widehat{\mathcal{B}}_j$'s such that
$F_j \in L^p(X \times Y)$ for each $j$, then this implies the less
precise statement that $I_Y(|F_j|^p) = N_p(F_j)^p$ is a submartingale
on $X$ with respect to the $\mathcal{B}_j$'s.  One can still get some
maximal function estimates from this, which are adequate for the same
conclusions about pointwise convergence.

        If $Y'$ is a $\sigma$-finite measure space, then one can
choose a positive weight on $Y'$ to get a probability measure, as in
Section \ref{L^p-valued martingales}.  This permits one to identify
$L^p(Y')$ with $L^p(Y)$ for a probability space $Y$, as before.  Thus
martingales on $X$ with values in $L^p(Y')$ can be identified with
martingales on $X$ with values in $L^p(Y)$, to which the discussion in
this section applies.

\section[\ $\ell^p$-Valued martingales]{$\ell^p$-Valued martingales}
\label{ell^p-valued martingales}
\setcounter{equation}{0}

        Let $(X, \mathcal{A}, \mu)$ be a probability space, and let
$\{f_l(x)\}_{l = 1}^\infty$ be a sequence of real or complex-valued functions
on $X$ such that $f_l(x) \in L^p(X)$ for each $l$, $1 \le p < \infty$, and
\begin{equation}
\label{sum_{l = 1}^infty int_X |f_l(x)|^p d mu(x) < infty}
        \sum_{l = 1}^\infty \int_X |f_l(x)|^p \, d\mu(x) < \infty.
\end{equation}
This is the same as
\begin{equation}
\label{int_X sum_{l = 1}^infty |f_l(x)|^p d mu(x) < infty}
        \int_X \sum_{l = 1}^\infty |f_l(x)|^p \, d\mu(x) < \infty,
\end{equation}
which implies that $\{f_l(x)\}_{l = 1}^\infty \in \ell^p = \ell^p({\bf
Z}_+)$ for almost every $x \in X$.  One can also think of
$\{f_l(x)\}_{l = 1}^\infty$ as an element of $L^p(X \times {\bf
Z}_+)$, where $X \times {\bf Z}_+$ is equipped with the product
measure associated to counting measure on ${\bf Z}_+$.

        If $\mathcal{B}$ is a $\sigma$-subalgebra of $\mathcal{A}$,
then of course one can take the conditional expectation $E(f_l \mid
\mathcal{B})$ of $f_l$ for each $l$, and
\begin{equation}
        \int_X |E(f_l \mid \mathcal{B})|^p \, d\mu(x)
 \le \int_X E(|f_l|^p \mid \mathcal{B}) \, d\mu(x) = \int_X |f_l|^p \, d\mu.
\end{equation}
Hence
\begin{equation}
        \sum_{l = 1}^\infty \int_X |E(f_l \mid \mathcal{B})|^p \, d\mu 
         \le \sum_{l = 1}^\infty \int_X |f_l|^p \, d\mu.
\end{equation}
This is another way to look at conditional expectation of
$\ell^p$-valued functions, which is consistent with the earlier
discussions.

        More precisely,
\begin{equation}
        |E(f_l \mid \mathcal{B})|^p \le E(|f_l|^p \mid \mathcal{B})
\end{equation}
almost everywhere on $X$ for each $l$, and so
\begin{equation}
        \sum_{l = 1}^\infty |E(f_l \mid \mathcal{B})|^p
          \le \sum_{l = 1}^\infty E(|f_l|^p \mid \mathcal{B})
            = E\Big(\sum_{l = 1}^\infty |f_l|^p \mid \mathcal{B}\Big)
\end{equation}
almost everywhere on $X$.  As in Section \ref{sigma-subalgebras,
vectors},
\begin{equation}
        \Big(\sum_{l = 1}^n |E(f_l \mid \mathcal{B})|^p\Big)^{1/p}
 \le E\Big(\Big(\sum_{l = 1}^n |f_l|^p\Big)^{1/p} \mid \mathcal{B}\Big)
\end{equation}
almost everywhere on $X$ for each $n \in {\bf Z}_+$.  This implies
that
\begin{equation}
        \Big(\sum_{l = 1}^n |E(f_l \mid \mathcal{B})|^p\Big)^{1/p}
 \le E\Big(\Big(\sum_{l = 1}^\infty |f_l|^p\Big)^{1/p} \mid \mathcal{B}\Big)
\end{equation}
almost everywhere on $X$ for each $n$, and thus
\begin{equation}
        \Big(\sum_{l = 1}^\infty |E(f_l \mid \mathcal{B})|^p\Big)^{1/p}
 \le E\Big(\Big(\sum_{l = 1}^\infty |f_l|^p\Big)^{1/p} \mid \mathcal{B}\Big).
\end{equation}

        As in Section \ref{L^p-valued martingales}, one can choose a
positive weight on ${\bf Z}_+$ to identify $\ell^p$ with $L^p(Y)$,
where $Y$ is a probability space.  Thus the estimates in the preceding
paragraph can be seen as a special case of those in the previous
section, with simplifications from the discreteness of $Y$.  As
before, one can get submartingales from the norms of $\ell^p$-valued
martingales, and then maximal function estimates for these.  In
particular, it follows that an $\ell^1$-valued martingale with bounded
$L^1$ norm converges almost everywhere, as in Section
\ref{ell^1-valued martingales}.

\section[\ Approximation in product spaces]{Approximation in product spaces}
\label{approximations in product spaces}
\setcounter{equation}{0}

        Let $(X_1, \mathcal{A}_1, \mu_1)$, $(X_2, \mathcal{A}_2,
\mu_2)$ be measure spaces with $\mu_1(X_1), \mu_2(X_2) < \infty$, and
consider their Cartesian product $X_1 \times X_2$.  The
$\sigma$-algebra $\mathcal{A}$ of measurable subsets of $X$ is defined
as the smallest $\sigma$-algebra of subsets of $X$ that contains the
measurable rectangles $A_1 \times A_2$, $A_1 \in \mathcal{A}_1$, $A_2
\in \mathcal{A}_2$.  Note that the intersection of two measurable
rectangles in $X$ is also a measurable rectangle, and that the
complement of a measurable rectangle is the union of three
pairwise-disjoint measurable rectangles, since
\begin{eqnarray}
\lefteqn{(X_1 \times X_2) \backslash (A_1 \times A_2) =} \\
 & & ((X_1 \backslash A_1) \times A_2) \cup (A_1 \times (X_2 \backslash A_2))
          \cup ((X_1 \backslash A_1) \times (X_2 \backslash A_2)). \nonumber
\end{eqnarray}
Let $\mathcal{E}$ be the collection of subsets of $X$ that can be
expressed as the union of finitely many pairwise-disjoint measurable
rectangles.  This is an algebra of subsets of $X$, by the previous
observations.  Also let $\mu = \mu_1 \times \mu_2$ be the product
measure associated to $\mu_1$, $\mu_2$ on $\mathcal{A}$.  If
\begin{equation}
        d(A, B) = \mu(A \bigtriangleup B)
\end{equation}
is the corresponding semimetric on $\mathcal{A}$ as in Section
\ref{distances between sets}, then $\mathcal{E}$ is dense in
$\mathcal{A}$ with respect to $d(A, B)$.  Depending on the way that
the product measure is defined, this may be obvious from the
construction.  At any rate, this follows from the discussion in
Section \ref{distances between sets}, which implies that the closure
$\overline{\mathcal{E}}$ of $\mathcal{E}$ in $\mathcal{A}$ is a
$\sigma$-subalgebra of $\mathcal{A}$ that contains $\mathcal{E}$.  One
could also use the characterization of $\mathcal{A}$ as the smallest
monotone class of subsets of $X$ that contains $\mathcal{E}$.  If
$X_1$, $X_2$ are $\sigma$-finite and $A \subseteq X$ is a measurable
set with $\mu(A) < \infty$, then one can first approximate $A$ by
subsets of products of measurable sets with finite measure, and then
continue as before to approximate $A$ by elements of $\mathcal{E}$.
Using these approximations, one can check that nice functions are
dense in $L^p(X)$ when $1 \le p < \infty$, as in Section
\ref{martingales, products}.  Of course, these statements are much
simpler when $X_1$ or $X_2$ has only finitely or countably many
elements and all of its subsets are measurable, or when
$\mathcal{A}_1$ or $\mathcal{A}_2$ is generated by a partition of the
corresponding space into finitely or countably many subsets.

\section[\ Mixed norms]{Mixed norms}
\label{mixed norms}
\setcounter{equation}{0}

        Let $(X, \mathcal{A}, \mu)$, $(Y, \mathcal{B}, \nu)$ be
probability spaces, and let their Cartesian product $X \times Y$ be
equipped with the product measure $\mu \times \nu$, as usual.
Consider the space of real or complex-valued measurable functions
$F(x, y)$ on $X \times Y$ such that
\begin{equation}
\label{int_X (int_Y |F(x, y)|^p d nu(y)^{1/p} d mu(x)}
 \int_X \Big(\int_Y |F(x, y)|^p \, d\nu(y)\Big)^{1/p} \, d\mu(x)
\end{equation}
is finite, where $1 \le p < \infty$.  It is easy to see that this is a
vector space, and that (\ref{int_X (int_Y |F(x, y)|^p d nu(y)^{1/p} d
mu(x)}) becomes a norm on this space when we identify functions that
are equal almost everywhere.  If $F(x, y) \in L^p(X \times Y)$, then
$F(x, y)$ is in this space, and
\begin{eqnarray}
\label{int_X(int_Y |F(x, y)|^p d nu(y))^{1/p} d mu(x) le ...}
\lefteqn{\int_X\Big(\int_Y |F(x, y)|^p \, d\nu(y)\Big)^{1/p} \, d\mu(x)} \\
         & \le & \Big(\int_{X \times Y} |F(x, y)|^p \,
                           d(\mu \times \nu)(x, y)\Big)^{1/p}, \nonumber
\end{eqnarray}
by Fubini's theorem and Jensen's inequality.  Similarly,
if $F(x, y)$ is in this space, then $F(x, y) \in L^1(X \times Y)$, and
\begin{eqnarray}
\label{int_{X times Y} |F(x, y)| d(mu times nu)(x, y) le ...}
\lefteqn{\int_{X \times Y} |F(x, y)| \, d(\mu \times \nu)(x, y)} \\
 & \le & \int_X\Big(\int_Y |F(x, y)|^p \, d\nu(y)\Big)^{1/p} \, d\mu(x),
                                                                 \nonumber
\end{eqnarray}
again by Fubini's theorem and Jensen's inequality.

        Suppose that $F(x, y)$ is in this space, and put
\begin{equation}
        N_p(F)(x) = \Big(\int_Y |F(x, y)|^p \, d\nu(y)\Big)^{1/p},
\end{equation}
as in Section \ref{martingales, products}.  Thus (\ref{int_X (int_Y
|F(x, y)|^p d nu(y)^{1/p} d mu(x)}) is the same as the $L^1(X)$ norm
of $N_p(F)(x)$.  If $L \ge 0$, then define $F_L(x, y)$ on $X \times Y$ by
\begin{eqnarray}
        F_L(x, y) & = & F(x, y) \hbox{ when } N_p(F)(x) \le L, \\
                  & = & 0 \qquad\quad\hbox{when } N_p(F)(x) > L.\nonumber
\end{eqnarray}
In particular, $N_p(F_L)(x) = N_p(F)(x)$ when $N_p(F)(x) \le L$, and
$N_p(F_L)(x) = 0$ when $N_p(F)(x) > L$.  It follows that $F_L(x, y)
\in L^p(X \times Y)$ for each $L$, and that $F_L(x, y)$ converges to
$F(x, y)$ with respect to the norm (\ref{int_X (int_Y |F(x, y)|^p d
nu(y)^{1/p} d mu(x)}) as $L \to \infty$, so that $L^p(X \times Y)$
is a dense linear subspace of this space.

        Let $\{F_j(x, y)\}_{j = 1}^\infty$ be a sequence of measurable
functions on $X \times Y$.  By Fatou's lemma,
\begin{equation}
\label{int_Y liminf_{j to infty} |F_j(x, y)|^p d nu(y) le ...}
        \int_Y \liminf_{j \to \infty} |F_j(x, y)|^p \, d\nu(y)
         \le \liminf_{j \to \infty} \int_Y |F_j(x, y)|^p \, d\nu(y)
\end{equation}
for every $x \in X$.  Equivalently,
\begin{eqnarray}
\label{(int_Y (liminf_{j to infty} |F_j(x, y)|)^p d nu(y))^{1/p} le ...}
\lefteqn{\Big(\int_Y \Big(\liminf_{j \to \infty} |F_j(x, y)|\Big)^p \,
                                                      d\nu(y)\Big)^{1/p}} \\
        & \le & \liminf_{j \to \infty} \Big(\int_Y |F_j(x, y)|^p \,
                                                d\nu(y)\Big)^{1/p}. \nonumber
\end{eqnarray}
for each $x \in X$.  Applying Fatou's lemma a second time, we get that
\begin{eqnarray}
\lefteqn{\int_X\Big(\int_Y\Big(\liminf_{j \to \infty} |F_j(x, y)|\Big)^p \,
                                           d\nu(y)\Big)^{1/p} \, d\mu(x)} \\
        & \le & \liminf_{j \to \infty} \int_X\Big(\int_Y |F_j(x, y)|^p \,
                                     d\nu(y)\Big)^{1/p} \, d\mu(x). \nonumber
\end{eqnarray}

        Suppose that $\{F_j(x, y)\}_{j = 1}^\infty$ converges almost
everywhere to $F(x, y)$ on $X \times Y$.  It follows that for almost
every $x \in X$, $\{F_j(x, y)\}_{j = 1}^\infty$ converges to $F(x, y)$
for almost every $y \in Y$.  Hence
\begin{equation}
\label{int_Y |F(x, y)|^p d nu(y) le ...}
        \int_Y |F(x, y)|^p \, d\nu(y)
         \le \liminf_{j \to \infty} \int_Y |F_j(x, y)|^p \, d\nu(y)
\end{equation}
for almost every $x \in X$.  This implies that
\begin{eqnarray}
\lefteqn{\int_X\Big(\int_Y |F(x, y)|^p \, d\nu(y)\Big)^{1/p} \, d\mu(x)} \\
 & \le & \liminf_{j \to \infty} \int_X\Big(\int_Y |F_j(x, y)|^p \,
                               d\nu(y)\Big)^{1/p} \, d\mu(x), \nonumber
\end{eqnarray}
as before.

\section[\ Mixed-norm martingales]{Mixed-norm martingales}
\label{mixed-norm martingales}
\setcounter{equation}{0}

        Let $(X, \mathcal{A}, \mu)$, $(Y, \mathcal{B}, \nu)$ be
probability spaces, and let $X \times Y$ be equipped with the product
measure $\mu \times \nu$.  Also let $\mathcal{B}_1 \subseteq
\mathcal{B}_2 \subseteq \cdots$ be an increasing sequence of
$\sigma$-subalgebras of $\mathcal{A}$, and let $\widehat{\mathcal{B}}_j$
be the $\sigma$-algebra of subsets of $X \times Y$ that corresponds to
$\mathcal{B}_j$ on $X$ and $\mathcal{B}$ on $Y$ in the product space.
Suppose that $F(x, y)$ is a measurable function on $X \times Y$ for which
(\ref{int_X (int_Y |F(x, y)|^p d nu(y)^{1/p} d mu(x)}) is finite, $1 \le p < 
\infty$.  In particular, $F(x, y) \in L^1(X \times Y)$, and so
\begin{equation}
        F_j = E(F \mid \widehat{\mathcal{B}}_j)
\end{equation}
defines a martingale on $X \times Y$ with respect to $\widehat{\mathcal{B}}_j$.

        As in Section \ref{martingales, products},
\begin{equation}
\label{N_p(F_j) le E_X(N_p(F) mid mathcal{B}_j)}
        N_p(F_j) \le E_X(N_p(F) \mid \mathcal{B}_j)
\end{equation}
almost everywhere on $X$ for each $j \ge 1$, where the subscript $X$
of $E$ indicates that the conditional expectation is taken on $X$.
More precisely, this is the same as (\ref{N_p(E_{X x Y}(F |
widehat{mathcal{B}}_j)) le E_X(N_p(F) | mathcal{B}_j)}) when $F(x, y)
\in L^p(X \times Y)$, and otherwise we can approximate $F(x, y)$ by
elements of $L^p(X \times Y)$ with respect to the norm (\ref{int_X
(int_Y |F(x, y)|^p d nu(y)^{1/p} d mu(x)}), as in the previous section.
Integrating (\ref{N_p(F_j) le E_X(N_p(F) mid mathcal{B}_j)}) over $X$,
we get that
\begin{equation}
 \int_X N_p(F_j) \, d\mu \le \int_X E_X(N_p(F) \mid \mathcal{B}_j) \, d\mu
                              = \int_X N_p(F) \, d\mu
\end{equation}
for each $j$.  Thus the norm of $F_j$ with respect to (\ref{int_X
(int_Y |F(x, y)|^p d nu(y)^{1/p} d mu(x)}) is less than or equal to
(\ref{int_X (int_Y |F(x, y)|^p d nu(y)^{1/p} d mu(x)}) for each $j$.

        One can also check that $\{F_j\}_{j = 1}^\infty$ converges to
$F$ with respect to the norm (\ref{int_X (int_Y |F(x, y)|^p d
nu(y)^{1/p} d mu(x)}).  If $F \in L^p(X \times Y)$, then $\{F_j\}_{j =
1}^\infty$ converges to $F$ with respect to the $L^p$ norm, and hence
with respect to (\ref{int_X (int_Y |F(x, y)|^p d nu(y)^{1/p} d
mu(x)}).  Otherwise, one can approximate $F$ by elements of $L^p(X
\times Y)$, using the uniform bound for the norm of $F_j$ in the
previous paragraph.

        If we apply (\ref{N_p(F_j) le E_X(N_p(F) mid mathcal{B}_j)})
to $F_{j + 1}$ instead of $F_j$, then we get that
\begin{equation}
\label{N_p(F_j) le E_X(N_p(F_{j + 1}) mid mathcal{B}_j)}
        N_p(F_j) \le E_X(N_p(F_{j + 1}) \mid \mathcal{B}_j)
\end{equation}
almost everywhere on $X$ for each $j \ge 1$.  Thus $\{N_p(F_j)\}_{j =
1}^\infty$ is a submartingale on $X$ with respect to the
$\mathcal{B}_j$'s, which leads to maximal function estimates as
before.  Using convergence of $\{F_j\}_{j = 1}^\infty$ to $F$ with
respect to the norm (\ref{int_X (int_Y |F(x, y)|^p d nu(y)^{1/p} d
mu(x)}), one can show that $\{F_j(x, \cdot)\}_{j = 1}^\infty$
converges to $F(x, \cdot)$ in $L^p(Y)$ for almost every $x \in X$.
This is basically the same as in the previous situations, once we have
the same ingredients as before.

\section[\ Mixed-norm convergence]{Mixed-norm convergence}
\label{mixed-norm convergence}
\setcounter{equation}{0}

        Let $(X, \mathcal{A}, \mu)$, $(Y, \mathcal{B}, \nu)$ be
probability spaces, and let $X \times Y$ be equipped with $\mu \times
\nu$, as usual.  Also let $\mathcal{B}_1 \subseteq \mathcal{B}_2
\subseteq \cdots$ be an increasing sequence of $\sigma$-subalgebras of
$\mathcal{A}$, and let $\widehat{\mathcal{B}}_j$ be the
$\sigma$-algebra of subsets of $X \times Y$ that corresponds to
$\mathcal{B}_j$ on $X$ and $\mathcal{B}$ on $Y$ in the product space.
Suppose that $\{F_j\}_{j = 1}^\infty$ is a martingale on $X \times Y$
with respect to the $\widehat{\mathcal{B}}_j$'s whose norms as in
(\ref{int_X (int_Y |F(x, y)|^p d nu(y)^{1/p} d mu(x)}) are uniformly
bounded for some $p > 1$.  Equivalently,
\begin{equation}
        \int_X N_p(F_j) \, d\mu \le C
\end{equation}
for some $C \ge 0$ and each $j$.  Note that $\{N_p(F_j)\}_{j =
1}^\infty$ is a submartingale on $X$ with respect to the
$\mathcal{B}_j$'s, as in (\ref{N_p(F_j) le E_X(N_p(F_{j + 1}) mid
mathcal{B}_j)}).

        Suppose in addition that $\{N_p(F_j)\}_{j = 1}^\infty$ is
uniformly integrable on $X$, as in Section \ref{uniform
integrability}, and let us check that $\{F_j\}_{j = 1}^\infty$ is
uniformly integrable on $X \times Y$.  Let $\epsilon > 0$ be given,
and choose $\delta > 0$ such that
\begin{equation}
        \int_A N_p(F_j) \, d\mu < \frac{\epsilon}{2}
\end{equation}
for every measurable set $A \subseteq X$ with $\mu(A) < \delta$ and
each $j$.  If
\begin{equation}
        A_{j, L} = \{x \in X : N_p(F_j)(x) > L\},
\end{equation}
then
\begin{equation}
        \mu(A_{j, L}) < L^{-1} \, C
\end{equation}
for each $j$, $L$, by Tchebychev's inequality.  Hence $\mu(A_{j, L}) <
\delta$ for each $j$ when $L$ is sufficiently large, which implies that
\begin{equation}
        \int_{A_{j, L} \times Y} |F_j| \, d(\mu \times \nu)
         \le \int_{A_{j, L}} N_p(F_j) \, d\mu < \frac{\epsilon}{2}
\end{equation}
for each $j$ when $L$ is sufficiently large.  On the complement of
$A_{j, L} \times Y$, we have that
\begin{equation}
 \int_{(X \backslash A_{j, L}) \times Y} |F_j|^p \, d(\mu \times \nu)
            = \int_{X \backslash A_{j, L}} N_p(F_j)^p \, d\mu \le L^p
\end{equation}
for each $j$, $L$, by the definition of $A_{j, L}$.  Let $q$ be the
exponent conjugate to $p$, so that $1/p + 1/q = 1$.  If $B \subseteq X
\times Y$ is measurable, then
\begin{equation}
 \int_B |F_j| \, d(\mu \times \nu) \le ((\mu \times \nu)(B))^{1/q} \,
                    \Big(\int_B |F_j|^p \, d(\mu \times \nu)\Big)^{1/p},
\end{equation}
by H\"older's inequality.  If $B \subseteq (X \backslash A_{j, L})
\times Y$, then it follows that
\begin{equation}
 \int_B |F_j| \, d(\mu \times \nu) \le L \, ((\mu \times \nu)(B))^{1/q}.
\end{equation}
In order to show that $\{F_j\}_{j = 1}^\infty$ is uniformly
integrable, one can combine this with the earlier estimate for the
integral of $|F_j|$ over $A_{j, L} \times Y$ when $L$ is sufficiently
large.

        If $\{F_j\}_{j = 1}^\infty$ is uniformly integrable on $X
\times Y$, then $\{F_j\}_{j = 1}^\infty$ converges in $L^1(X \times
Y)$ to a function $F$, and $F_j = E(F \mid \widehat{\mathcal{B}}_j)$
for each $j$.  Moreover, $\{F_j\}_{j = 1}^\infty$ converges to $F$
almost everywhere on $X \times Y$, which implies that the norm of $F$
with respect to (\ref{int_X (int_Y |F(x, y)|^p d nu(y)^{1/p} d mu(x)})
is also finite, as in Section \ref{mixed norms}.  Thus we are back in
the situation of the preceding section.  This implies that $\{F_j\}_{j
= 1}^\infty$ also converges to $F$ with respect to the norm
(\ref{int_X (int_Y |F(x, y)|^p d nu(y)^{1/p} d mu(x)}), and that
$\{F_j(x, \cdot)\}_{j = 1}^\infty$ converges to $F(x, \cdot)$ in
$L^p(Y)$ for almost every $x \in X$.

        Suppose now that $\{N_p(F_j)\}_{j = 1}^\infty$ is still
bounded in $L^1(X)$, but may not be uniformly integrable.  Because
$\{N_p(F_j)\}_{j = 1}^\infty$ is a submartingale with respect to the
$\mathcal{B}_j$'s, the corresponding maximal function can be estimated
in the usual way.  In this case, $\{F_j\}_{j = 1}^\infty$ can be
approximated by martingales $\{G_j\}_{j = 1}^\infty$ on $X \times Y$
such that $\{N_p(G_j)\}_{j = 1}^\infty$ is uniformly integrable, as in
Section \ref{convergence almost everywhere}.  More precisely, the
approximation basically takes place in the $x$ variable.  This permits
one to show that $\{F_j(x, \cdot)\}_{j = 1}^\infty$ converges in
$L^p(Y)$ for almost every $x \in X$, as before.

\section[\ The $\ell^p$ version]{The $\ell^p$ version}
\label{ell^p version}
\setcounter{equation}{0}

        Let $(X, \mathcal{A}, \mu)$ be a probability space, and let $1
\le p < \infty$ be given.  If $\{f_l(x)\}_{l = 1}^\infty$ is a
sequence of real or complex-valued measurable functions on $X$ such that
\begin{equation}
\label{int_X (sum_{l = 1}^infty |f_l(x)|^p)^{1/p} d mu(x)}
        \int_X \Big(\sum_{l = 1}^\infty |f_l(x)|^p\Big)^{1/p} \, d\mu(x)
\end{equation}
is finite, then
\begin{equation}
        \sum_{l = 1}^\infty |f_l(x)|^p < \infty
\end{equation}
for almost every $x \in X$.  It is easy to see that the space of
sequences of functions on $X$ of this type is a vector space, and that
(\ref{int_X (sum_{l = 1}^infty |f_l(x)|^p)^{1/p} d mu(x)}) defines a
norm on this vector space when we identify functions that are equal
almost everywhere on $X$.  We can also use a weight on the set of
positive integers to identify $\ell^p$ with $L^p(Y)$ for a probability
space $Y$, so that this expression is the same as (\ref{int_X (int_Y
|F(x, y)|^p d nu(y)^{1/p} d mu(x)}).

        If $\{f_l(x)\}_{l = 1}^\infty$ is a sequence of functions in
$L^p(X)$ such that
\begin{equation}
        \sum_{l = 1}^\infty \int_X |f_l(x)|^p \, d\mu(x)
 = \int_X \sum_{l = 1}^\infty |f_l(x)|^p \, d\mu(x) < \infty,
\end{equation}
then (\ref{int_X (sum_{l = 1}^infty |f_l(x)|^p)^{1/p} d mu(x)}) is
also finite, because
\begin{equation}
 \int_X \Big(\sum_{l = 1}^\infty |f_l(x)|^p\Big)^{1/p} \, d\mu(x)
  \le \Big(\int_X \sum_{l = 1}^\infty |f_l(x)|^p \, d\mu(x)\Big)^{1/p},
\end{equation}
by Jensen's inequality.  These sequences of functions are dense among
those for which (\ref{int_X (sum_{l = 1}^infty |f_l(x)|^p)^{1/p} d
mu(x)}) is finite, with respect to the norm (\ref{int_X (sum_{l =
1}^infty |f_l(x)|^p)^{1/p} d mu(x)}), for the same reasons as in
Section \ref{mixed norms}.  Of course, these two conditions on
sequences of functions on $X$ are the same when $p = 1$.
Alternatively, if $\{f_l(x)\}_{l = 1}^\infty$ is a sequence of
functions on $X$ for which (\ref{int_X (sum_{l = 1}^infty
|f_l(x)|^p)^{1/p} d mu(x)}) is finite, then
\begin{equation}
\label{lim_{n to infty} int_X (sum_{l = n}^infty |f_l(x)|^p)^{1/p} d mu(x) = 0}
 \lim_{n \to \infty} \int_X \Big(\sum_{l = n}^\infty |f_l(x)|^p\Big)^{1/p} \,
                                                                  d\mu(x) = 0,
\end{equation}
by the dominated convergence theorem.  This implies that
$\{f_l(x)\}_{l = 1}^\infty$ can be approximated by sequences of
functions for which all but finitely many terms are equal to $0$ with
respect to the norm (\ref{int_X (sum_{l = 1}^infty |f_l(x)|^p)^{1/p} d
mu(x)}).

        Let $f(x) = \{f_l(x)\}_{l = 1}^\infty$ be a sequence of
measurable functions on $X$ for which (\ref{int_X (sum_{l = 1}^infty
|f_l(x)|^p)^{1/p} d mu(x)}) is finite, and let $\mathcal{B}$ be a
$\sigma$-subalgebra of $\mathcal{A}$.  As in Sections
\ref{sigma-subalgebras, vectors} and \ref{ell^p-valued martingales},
\begin{equation}
        \Big(\sum_{l = 1}^n |E(f_l \mid \mathcal{B})|^p\Big)^{1/p}
 \le E\Big(\Big(\sum_{l = 1}^n |f_l|^p\Big)^{1/p} \mid \mathcal{B}\Big)
\end{equation}
almost everywhere on $X$ for each $n$, and hence
\begin{equation}
        \Big(\sum_{l = 1}^\infty |E(f_l \mid \mathcal{B})|^p\Big)^{1/p}
 \le E\Big(\Big(\sum_{l = 1}^\infty |f_l|^p\Big)^{1/p} \mid \mathcal{B}\Big)
\end{equation}
almost everywhere on $X$.  In particular,
\begin{eqnarray}
 \int_X\Big(\sum_{l = 1}^\infty |E(f_l \mid \mathcal{B})|^p\Big)^{1/p} \, d\mu
 & \le & \int_X E\Big(\Big(\sum_{l = 1}^\infty |f_l|^p\Big)^{1/p} \mid
                                           \mathcal{B}\Big) \, d\mu \\
 & = & \int_X \Big(\sum_{l = 1}^\infty |f_l|^p\Big)^{1/p} \, d\mu. \nonumber
\end{eqnarray}

        Now let $\mathcal{B}_1 \subseteq \mathcal{B}_2 \subseteq
\cdots$ be an increasing sequence of $\sigma$-subalgebras of
$\mathcal{A}$.  Also let $f_j(x) = \{f_{j, l}(x)\}_{l = 1}^\infty$ be
a sequence of measurable functions with respect to $\mathcal{B}_j$ for
which (\ref{int_X (sum_{l = 1}^infty |f_l(x)|^p)^{1/p} d mu(x)}) is
finite for each $j$, and put
\begin{equation}
 A_j = \int_X \Big(\sum_{l = 1}^\infty |f_{j, l}(x)|^p\Big)^{1/p} \, d\mu(x).
\end{equation}
Suppose that $\{f_{j, l}\}_{j = 1}^\infty$ is a martingale with
respect to this filtration for each $l$, so that $f_{j, l} = E(f_{j +
1, l} \mid \mathcal{B}_j)$ for each $j, l \ge 1$.  Thus
\begin{equation}
\label{||f_j(x)||_p = (sum_{l = 1}^infty |f_{j, l}(x)|^p)^{1/p}}
        \|f_j(x)\|_p = \Big(\sum_{l = 1}^\infty |f_{j, l}(x)|^p\Big)^{1/p}
\end{equation}
is a submartingale with respect to this filtration, as in the previous
paragraph.  This implies that $\{A_j\}_{j = 1}^\infty$ is monotone
increasing, as usual.  Similarly, if
\begin{equation}
 A_{j, n} = \int_X \Big(\sum_{l = 1}^n |f_{j, l}(x)|^p\Big)^{1/p} \, d\mu(x),
\end{equation}
then $A_{j, n} \le A_{j + 1, n}$ for each $j, n \ge 1$.  Note that
\begin{equation}
        \lim_{n \to \infty} A_{j, n} = A_j
\end{equation}
for each $j$, by the dominated convergence theorem.

        Suppose that the $A_j$'s are bounded, and put
\begin{equation}
        A = \sup_{j \ge 1} A_j.
\end{equation}
Let $\delta > 0$ be given, and choose $j_0$ such that
\begin{equation}
        A_{j_0} > A - \delta.
\end{equation}
Because $A_{j_0, n} \to A_{j_0}$ as $n \to \infty$, we can choose
$n_0$ so that
\begin{equation}
        A_{j_0, n_0} > A - \delta.
\end{equation}
If $j \ge j_0$, then monotonicity implies that
\begin{equation}
        A_{j, n_0} > A - \delta.
\end{equation}

        Let us pause a moment to record some elementary inequalities
that will be helpful later.  If $a, b > 0$, then
\begin{equation}
        (a + b)^{1/p} \ge a^{1/p} + p^{-1} \, (a + b)^{(1/p) - 1} \, b.
\end{equation}
This follows from calculus, because
\begin{equation}
        \frac{d}{dt} (a + t)^{1/p} = p^{-1} \, (a + t)^{(1/p) - 1}
\end{equation}
is minimized on $[0, b]$ at $t = b$.  Remember that $0 < 1/p \le 1$,
because $1 \le p < \infty$.  If $b \ge \epsilon \, a$ for some
$\epsilon > 0$, then $a + b \le (\epsilon^{-1} + 1) \, b$, and so
\begin{equation}
(a + b)^{1/p} \ge a^{1/p} + p^{-1}\, (\epsilon^{-1} + 1)^{(1/p) - 1}\, b^{1/p}.
\end{equation}
This implies that
\begin{equation}
\label{b^{1/p} le ...}
        b^{1/p} \le \epsilon^{1/p} \, a^{1/p}
           + p \, (\epsilon^{-1} + 1)^{1 - (1/p)} \, ((a + b)^{1/p} - a^{1/p}),
\end{equation}
for every $\epsilon > 0$.  More precisely, $b^{1/p}$ is less than or
equal to the second term on the right when $b \ge \epsilon \, a$, by
the previous inequality, and otherwise $b^{1/p}$ is less than or equal
to the first term on the right, because $b < \epsilon \, a$.  Note
that (\ref{b^{1/p} le ...}) also holds when $a = 0$ or $b = 0$.

        Let us apply (\ref{b^{1/p} le ...}) to
\begin{equation}
        a = \sum_{l = 1}^n |f_{j, l}(x)|^p, \quad
         b = \sum_{l = n + 1}^\infty |f_{j, l}(x)|^p,
\end{equation}
using also the fact that $a^{1/p} \le (a + b)^{1/p} = \|f_j(x)\|_p$.
This implies that
\begin{eqnarray}
          & & \Big(\sum_{l = n + 1}^\infty |f_{j, l}(x)|^p\Big)^{1/p} \\
        & \le & \epsilon^{1/p} \, \|f_j(x)\|_p +
            p \, (\epsilon^{-1} + 1)^{1 - (1/p)} \, \Big(\|f_j(x)\|_p
             - \Big(\sum_{l = 1}^n |f_{j, l}(x)|^p\Big)^{1/p}\Big). \nonumber
\end{eqnarray}
Integrating over $X$, we get that
\begin{eqnarray}
\lefteqn{\int_X \Big(\sum_{l = n + 1}^\infty |f_{j, l}(x)|^p\Big)^{1/p}
                                                             \, d\mu(x)} \\
 & \le & \epsilon^{1/p} \, A + p \, (\epsilon^{-1} + 1)^{1 - (1/p)} \,
                                                   (A - A_{j, n}), \nonumber
\end{eqnarray}
using also the fact that $A_j \le A$ for each $j$, by the definition
of $A$.  Taking $n = n_0$, we get that
\begin{equation}
\label{int_X (sum_{l = n_0 + 1}^infty |f_{j, l}(x)|^p)^{1/p} d mu(x) le ...}
 \int_X \Big(\sum_{l = n_0 + 1}^\infty |f_{j, l}(x)|^p\Big)^{1/p} \, d\mu(x)
 \le \epsilon^{1/p} \, A + p \, (\epsilon^{-1} + 1)^{1 - (1/p)} \, \delta
\end{equation}
when $j \ge j_0$.

        If $\eta > 0$ is given, then we can first choose $\epsilon >
0$ so that $\epsilon^{1/p} \, A < \eta/2$, and then choose $\delta$
depending on $\epsilon$ such that $p \, (\epsilon^{-1} + 1)^{1 -
(1/p)} \, \delta < \eta/2$.  If $j_0, n_0$ are as before, then
(\ref{int_X (sum_{l = n_0 + 1}^infty |f_{j, l}(x)|^p)^{1/p} d mu(x) le
...}) implies that
\begin{equation}
\label{int_X (sum_{l = n_0 + 1}^infty |f_{j, l}(x)|^p)^{1/p} d mu(x) < eta}
 \int_X \Big(\sum_{l = n_0 + 1}^\infty |f_{j, l}(x)|^p\Big)^{1/p} \, d\mu(x)
         < \eta
\end{equation}
when $j \ge j_0$.  The integral on the left side of (\ref{int_X
(sum_{l = n_0 + 1}^infty |f_{j, l}(x)|^p)^{1/p} d mu(x) < eta}) is
actually monotone increasing in $j$, for the usual submartingale
reasons, which implies that (\ref{int_X (sum_{l = n_0 + 1}^infty
|f_{j, l}(x)|^p)^{1/p} d mu(x) < eta}) holds for every $j$.  Put
$g_{j, l}(x) = f_{j, l}(x)$ when $l \le n_0$ and $g_{j, l}(x) = 0$
when $l > n_0$, so that $\{g_{j, l}\}_{j = 1}^\infty$ is a martingale
with respect to the $\mathcal{B}_j$'s for each $l$, and $f_j(x) =
\{f_{j, l}(x)\}_{l = 1}^\infty$ is approximated by $g_j(x) = \{g_{j,
l}(x)\}_{l = 1}^\infty$ uniformly in $j$ with respect to the norm
(\ref{int_X (sum_{l = 1}^infty |f_l(x)|^p)^{1/p} d mu(x)}), by
(\ref{int_X (sum_{l = n_0 + 1}^infty |f_{j, l}(x)|^p)^{1/p} d mu(x) <
eta}).  Using this approximation and maximal function estimates for
$\|f_j(x) - g_j(x)\|_p$, one can show that $\{f_j(x)\}_{j = 1}^\infty$
converges in $\ell^p$ for almost every $x \in X$, as in Sections
\ref{another criterion} and \ref{ell^1-valued martingales}.

\section[\ The doubling condition]{The doubling condition}
\label{doubling condition}
\setcounter{equation}{0}

        Let $(X, \mathcal{A}, \mu)$ be a probability space, and let
$\mathcal{P}_0, \mathcal{P}_1, \mathcal{P}_2, \ldots$ be a sequence of
partitions of $X$ into finitely many measurable sets of positive measure such
that $\mathcal{P}_{j + 1}$ is a refinement of $\mathcal{P}_j$ for each $j$
and $\mathcal{P}_0$ is the trivial partition consisting of only $X$ itself.
We say that the $\mathcal{P}_j$'s satisfy a doubling condition if there is
a $C \ge 1$ such that
\begin{equation}
\label{mu(A) le C mu(B)}
        \mu(A) \le C \, \mu(B)
\end{equation}
when $A \in \mathcal{P}_j$, $B \in \mathcal{P}_{j + 1}$, and $B \subseteq A$.
This implies that for each $A \in \mathcal{P}_j$ there are less than or
equal to $C$ sets $B \in \mathcal{P}_{j + 1}$ such that $B \subseteq A$.
In particular, this implies that $\mathcal{P}_j$ has less than or equal to
$C^j$ elements for each $j$.  If $X = [0, 1)$ is equipped with Lebesgue measure
and $\mathcal{P}_j$ consists of the dyadic subintervals of $[0, 1)$ with
length $2^{-j}$, then (\ref{mu(A) le C mu(B)}) holds with $C = 2$.

        Let $\mathcal{B}_j = \mathcal{B}(\mathcal{P}_j)$ be the
$\sigma$-subalgebra of $\mathcal{A}$ generated by $\mathcal{P}_j$, as
in Section \ref{partitions}.  Thus $\mathcal{B}_j \subseteq
\mathcal{B}_{j + 1}$ for each $j$, since $\mathcal{P}_{j + 1}$ is
supposed to be a refinement of $\mathcal{P}_j$.  If $f_{j + 1}(x)$ is
a nonnegative real-valued function on $X$ which is measurable with
respect to $\mathcal{B}_{j + 1}$ for some $j \ge 0$, then
\begin{equation}
\label{f_{j + 1} le C E(f_{j + 1} mid mathcal{B}_j)}
        f_{j + 1} \le C \, E(f_{j + 1} \mid \mathcal{B}_j).
\end{equation}
If $\{f_j\}_{j = 0}^\infty$ is a martingale with respect to this
filtration consisting of nonnegative real-valued functions, then
\begin{equation}
\label{f_{j + 1} le C f_j}
        f_{j + 1} \le C \, f_j
\end{equation}
for each $j$.

        Suppose now that $\{\phi_j\}_{j = 0}^\infty$ is a
submartingale on $X$ with respect to this filtration consisting of
nonnegative real-valued functions, and put
\begin{equation}
\label{psi_j = E(phi_{j + 1} mid mathcal{B}_j)}
        \psi_j = E(\phi_{j + 1} \mid \mathcal{B}_j).
\end{equation}
Thus
\begin{equation}
\label{phi_j le psi_j}
        \phi_j \le \psi_j
\end{equation}
for each $j \ge 0$, because $\{\phi_j\}_{j = 1}^\infty$ is a
submartingale.  Hence
\begin{equation}
\label{psi_j = ... le E(psi_{j + 1} mid mathcal{B}_j)}
        \psi_j = E(\phi_{j + 1} \mid \mathcal{B}_j)
                \le E(\psi_{j + 1} \mid \mathcal{B}_j),
\end{equation}
which implies that $\{\psi_j\}_{j = 0}^\infty$ is also a submartingale.
The doubling condition implies that
\begin{equation}
\label{phi_{j + 1} le C psi_j}
        \phi_{j + 1} \le C \, \psi_j
\end{equation}
for each $j \ge 0$, as in (\ref{f_{j + 1} le C E(f_{j + 1} mid
mathcal{B}_j)}).  If the $\phi_j$'s have bounded $L^p$ norm for some
$p \ge 1$, then the $\psi_j$'s have bounded $L^p$ norm as well, and
with the same bound.

        Let $V$ be a real or complex vector space with a norm $\|v\|$,
and let $\{f_j(x)\}_{j = 0}^\infty$ be a $V$-valued martingale on $X$
with respect to the $\mathcal{B}_j$'s, as in Section
\ref{vector-valued martingales}.  Thus $\phi_j(x) = \|f_j(x)\|$ is a
nonnegative real-valued submartingale on $X$, and $\psi_j(x)$ can be
defined as in the previous paragraph.  Note that $f_0(x)$ is constant
on $X$, and let $t \ge \|f_0(x)\|$ be given.  Put $\tau(x) = \infty$
when $\psi_j(x) \le t$ for each $j \ge 0$, and otherwise let $\tau(x)$
be the smallest nonnegative integer $l$ such that $\psi_l(x) > t$.
This is a stopping time, as in Section \ref{stopping times}.  If
$\tau_n(x) = \min(\tau(x), n)$, then
\begin{equation}
\label{g_n(x) = f_{tau_n(x)}(x)}
        g_n(x) = f_{\tau_n(x)}(x)
\end{equation}
is also a $V$-valued martingale on $X$, as before.  This is basically
the same as the approximation to $\{f_j(x)\}_{j = 0}^\infty$ described
in Section \ref{convergence almost everywhere}, except that we use the
maximal function associated to $\psi_j$ instead of $\phi_j$.  
By construction,
\begin{equation}
\label{||g_n(x)|| = ||f_n(x)|| le psi_n(x) le t}
        \|g_n(x)\| = \|f_n(x)\| \le \psi_n(x) \le t
\end{equation}
when $n < \tau(x)$, and
\begin{equation}
\label{||g_n(x)|| = ||f_{tau(x)}(x)|| le C psi_{tau(x) - 1}(x) le C t}
 \|g_n(x)\| = \|f_{\tau(x)}(x)\| \le C \, \psi_{\tau(x) - 1}(x) \le C \, t
\end{equation}
when $0 < \tau(x) \le n$, because of the doubling condition. 
It follows that
\begin{equation}
\label{||g_n(x)|| le C t}
        \|g_n(x)\| \le C \, t
\end{equation}
for every $x \in X$ and $n \ge 0$, since $\|g_n(x)\| = \|f_0(x)\| \le
t$ when $\tau(x) = 0$, by hypothesis.  This is analogous to
(\ref{|g_n| le max (|h|, t)}), with the integrable function $h(x)$
replaced by $C \, t$.  If $\phi_j(x) = \|f_j(x)\|$ has bounded $L^1$
norm, so that $\psi_j(x)$ has bounded $L^1$ norm too, then the measure
of the set where $\tau(x) < \infty$ can be estimated as before.  Of
course, $g_n(x) = f_n(x)$ for every $n \ge 0$ when $\tau(x) = \infty$.
If every uniformly bounded $V$-valued martingale on $X$ converges
almost everywhere, then every $V$-valued martingale $\{f_j(x)\}_{j =
1}^\infty$ such that $\|f_j(x)\|$ has bounded $L^1$ norm also
converges almost everywhere, as in Section \ref{convergence almost
everywhere}.

\section[\ Paths and martingales]{Paths and martingales}
\label{paths, martingales}
\setcounter{equation}{0}

        Let $V$ be a real or complex vector space with a norm $\|v\|$,
and let $F$ be a $V$-valued function on $[0, 1]$.  If $[a, b)$ is a
dyadic subinterval of $[0, 1)$ of length $b - a = 2^{-j}$, then put
\begin{equation}
\label{f_j(x) = 2^j (F(b) - F(a))}
        f_j(x) = 2^j \, (F(b) - F(a))
\end{equation}
for every $x \in [a, b)$.  This defines $f_j(x)$ as a $V$-valued
function on $[0, 1)$ which is constant on the dyadic intervals of
length $2^{-j}$.  It is easy to see that the $f_j$'s form a $V$-valued
martingale on $[0, 1)$ with respect to Lebesgue measure and the
$\sigma$-subalgebras of measurable sets generated by the partitions of
$[0, 1)$ by dyadic intervals of length $2^{-j}$, as in Section
\ref{vector-valued martingales}.  Note that
\begin{equation}
        \int_0^1 \|f_j(x)\| \, dx
         = \sum_{l = 0}^{2^j - 1} \|F((l + 1) \, 2^{-j}) - F(l \, 2^{-j})\|.
\end{equation}
If $F : [0, 1] \to V$ has finite length $\Lambda$, then
\begin{equation}
        \int_0^1 \|f_j(x)\| \, dx \le \Lambda
\end{equation}
for each $j$.  If $F$ is Lipschitz, then the $f_j$'s are uniformly
bounded.  If $F$ is differentiable at $x$, then
\begin{equation}
        \lim_{j \to \infty} f_j(x) = F'(x).
\end{equation}

        If $V = L^1([0, 1])$ and $F(x)$ is the indicator function of
$[0, x]$, then $F$ is a Lipschitz function on $[0, 1]$ with values in
$L^1([0, 1])$, as in Section \ref{vector-valued functions}.  The
corresponding martingale is the same as the one described in Section
\ref{unit square}.

        Now let $V = L^\infty({\bf R})$, and let $\phi$ be a real or
complex-valued Lipschitz function on the real line.  Also let $\phi_x$
be the translate of $\phi$ by $x$, so that $\phi_x(y) = \phi(y - x)$.
If $\phi$ is bounded, then $F(x) = \phi_x$ defines a Lipschitz mapping
from ${\bf R}$ into $L^\infty({\bf R})$, as in Section
\ref{vector-valued functions}.  Otherwise, $F(x) = \phi_x - \phi$
defines a Lipschitz mapping from ${\bf R}$ into $L^\infty({\bf R})$,
using only the hypothesis that $\phi$ is Lipschitz on ${\bf R}$.  The
restriction of $F(x)$ to $x \in [0, 1]$ defines a martingale
$\{f_j\}_j$ with values in $L^\infty({\bf R})$ as before.  If $\phi$
is continuously-differentiable with uniformly continuous derivative,
then $F$ is differentiable at every $x \in {\bf R}$ as an
$L^\infty({\bf R})$-valued function on ${\bf R}$.  In this case, the
derivative of $F$ at $x$ corresponds to $-1$ times the derivative of
$\phi$ translated by $x$.  If $x \in [0, 1)$, then it is easy to see
that $\{f_j(x)\}_{j = 1}^\infty$ converges to the same limit in
$L^\infty({\bf R})$.  Conversely, if $\{f_j(x)\}_{j = 1}^\infty$
converges in $L^\infty({\bf R})$ for any $x \in [0, 1)$, then one can
show that $\phi$ is continuously differentiable with uniformly
continuous derivative.  This is analogous to the fact that $\phi$ is
continuously-differentiable with uniformly continuous derivative when
$F$ is differentiable at a single point, but slightly more
complicated, since we are only using ``dyadic'' difference quotients
of $F$.  If $\{f_j(x)\}_{j = 1}^\infty$ converges in $L^\infty({\bf
R})$ for some $x \in [0, 1)$, then the limit determines a bounded
uniformly continuous function $\psi$ on ${\bf R}$, because $f_j(x)$
corresponds to a bounded Lipschitz function on ${\bf R}$ for each $j$
that converges uniformly on ${\bf R}$ as $j \to \infty$.  One can
check that $\psi = -\phi_x'$ where $\phi_x$ is differentiable, and
then use the fact that Lipschitz functions are differentiable almost
everywhere and can be represented by integrals of their derivatives to
show that $\phi_x$ is continuously differentiable with derivative $-\psi$.
Alternatively, one can argue that $\phi_x' = -\psi$ in the sense of
distributions, and hence that $\phi_x$ is continuously differentiable
with derivative $-\psi$.

        Of course, one can just as well take $V$ to be the space
$C_b({\bf R})$ of bounded continuous functions on the real line with
the supremum norm here, which can be identified with a closed linear
subspace of $L^\infty({\bf R})$.  There is also a simple way to embed
$C_b({\bf R})$ linearly and isometrically into $\ell^\infty$, by
restricting a bounded continuous function on the real line to the
rationals, and then enumerating the latter by a sequence to get
bounded sequences of real or complex numbers.  If $\phi$ has compact
support, then one can view to restriction of $F(x)$ to $x \in [0, 1]$
as a Lipschitz mapping into the space of continuous functions on a
sufficiently large closed interval in the real line.

\section[\ $L^\infty$ Norms]{$L^\infty$ Norms}
\label{L^infty norms}
\setcounter{equation}{0}

        Let $(Y, \mathcal{B}, \nu)$ be a probability space, and
suppose that $g \in L^\infty(Y)$, so that
\begin{equation}
\label{||g||_p = (int_Y |g(y)|^p d nu(y))^{1/p} le ||g||_infty}
        \|g\|_p = \Big(\int_Y |g(y)|^p \, d\nu(y)\Big)^{1/p} \le \|g\|_\infty
\end{equation}
for each $p < \infty$.  Of course, $\|g\|_p$ is monotone increasing in
$p$, by Jensen's inequality, and it is well known and not difficult to
show that
\begin{equation}
\label{lim_{p to infty} ||g||_p = ||g||_infty}
        \lim_{p \to \infty} \|g\|_p = \|g\|_\infty.
\end{equation}
Similarly, if $g$ is a measurable function on $Y$ that is not
essentially bounded, and if $g \in L^p(Y)$ for each $p < \infty$, then
$\|g\|_p \to \infty$ as $p \to \infty$.

        Let $(X, \mathcal{A}, \nu)$ be another probability space, and
consider their Cartesian product $X \times Y$, equipped with the
product measure $\mu \times \nu$.  If $F(x, y)$ is a measurable
function on $X \times Y$, then
\begin{equation}
\label{N_infty(F)(x) = lim_{p to infty} N_p(F)(x) = ...}
        N_\infty(F)(x) = \lim_{p \to \infty} N_p(F)(x)
         = \lim_{p \to \infty} \Big(\int_Y |F(x, y)|^p \, d\nu(y)\Big)^{1/p}
\end{equation}
is a convenient way to express the norm of $F(x, y)$ as a function of
$y$ in $L^\infty(Y)$ for each $x \in X$.  More precisely, it is often
helpful to restrict $p$ to be a positive integer here, so that
$N_\infty(f)(x)$ is expressed as the limit of a monotone increasing
sequence of functions.  This makes it easy to derive properties of
$N_\infty(F)(x)$ like those for $N_p(F)(x)$ when $p < \infty$
discussed earlier.

        If $f(x) = \{f_l(x)\}_{l = 1}^\infty$ is a sequence of real or
complex-valued measurable functions on $X$, then the $\ell^\infty$
norm of $f(x)$ can be expressed as
\begin{equation}
\label{||f(x)||_infty = ... = lim_{n to infty} max_{1 le l le n} |f_l(x)|}
        \|f(x)\|_\infty = \sup_{l \ge 1} |f_l(x)|
                        = \lim_{n \to \infty} \max_{1 \le l \le n} |f_l(x)|,
\end{equation}
which implies that $\|f(x)\|_\infty$ is measurable on $X$.  If
$\|f(x)\|_\infty$ is integrable on $X$, then it is very easy to see
that the $\ell^\infty$ norm of the conditional expectation of the
$f_l(x)$'s with respect to some $\sigma$-subalgebra of $\mathcal{A}$
is bounded by the conditional expectation of $\|f(x)\|_\infty$.  One
can simply use the fact that $|f_l(x)| \le \|f(x)\|_\infty$ for each
$l$ to get that the conditional expectation of $f_l(x)$ is bounded by
the conditional expectation of $\|f(x)\|_\infty$ for each $l$, and
then take the supremum over $l$.

\section[\ Paths and measures]{Paths and measures}
\label{paths, measures}
\setcounter{equation}{0}

        Let $(V, \|v\|)$ be a real or complex Banach space, and let $F
: [a, b] \to V$ be a path of finite length.  As in Section
\ref{continuity conditions}, the one-sided limit $F(x+) = \lim_{y \to
x+} F(y)$ exists for every $x \in [a, b)$, and similarly $F(x-) =
\lim_{y \to x-} F(y)$ exists for every $x \in (a, b]$.  We can extend
$F$ to the whole real line by putting $F(x) = F(a)$ when $x < a$ and
$F(x) = F(b)$ when $x > b$, so that $F(a-) = F(a)$ and $F(b+) = F(b)$.

        As in Section \ref{functions, measures}, we can put
\begin{equation}
        \nu((r, t)) = F(t-) - F(r+)
\end{equation}
when $a \le r < t \le b$, and
\begin{equation}
        \nu([r, t]) = F(t+) - F(r-)
\end{equation}
when $a \le r \le t \le b$.  Similarly, we can put
\begin{equation}
        \nu([r, t)) = F(t-) - F(t+), \quad \nu((r, t]) = F(t+) - F(r-)
\end{equation}
when $a \le r < t \le b$.  This determines a finitely-additive
$V$-valued measure on the algebra $\mathcal{E}$ of subsets of $[a, b]$
that can be expressed as the union of finitely many intervals, where
the intervals may be open, closed, or half-open and half-closed.  Of
course, this is a bit simpler when $F$ is continuous.

        Let $\alpha(x)$ be the length of the restriction of $F$ to
$[a, x]$ when $a \le x \le b$.  This can be extended to all $x \in
{\bf R}$ by setting $\alpha(x) = 0$ when $x < a$ and $\alpha(x) =
\alpha(b)$ when $x > b$.  Thus $\alpha(x)$ is a monotone increasing
function on ${\bf R}$, which determines a nonnegative Borel measure
$\mu$ on ${\bf R}$ as in Section \ref{functions, measures}.  It is
easy to see that
\begin{equation}
\label{||nu(A)|| le mu(A)}
        \|\nu(A)\| \le \mu(A)
\end{equation}
for every $A \in \mathcal{E}$, because
\begin{equation}
        \|F(t) - F(r)\| \le \alpha(t) - \alpha(r)
\end{equation}
when $r \le t$.  Note that $\alpha(t) - \alpha(r)$ is the same as the
length of the restriction of $F$ to $[r, t]$ when $r \le t$, as in
Section \ref{lengths of paths}.

        Let $\mathcal{B}$ be the $\sigma$-algebra of Borel subsets of
$[a, b]$.  Thus $\mathcal{E} \subseteq \mathcal{B}$, and $\mathcal{B}$
is the smallest $\sigma$-algebra of subsets of $[a, b]$ that contains
$\mathcal{E}$.  If $d(A, B) = \mu(A \bigtriangleup B)$ is the distance
between $A, B \in \mathcal{B}$ associated to $\mu$ as in Section
\ref{distances between sets}, then it follows that the closure of
$\mathcal{E}$ in $\mathcal{B}$ with respect to $d(A, B)$ is equal to
$\mathcal{B}$.  This can also be seen more directly from the
construction of $\mu$.

        If $A, B \in \mathcal{E}$, then
\begin{eqnarray}
 \quad  \nu(A) - \nu(B) & = & (\nu(A \backslash B) + \nu(A \cap B))
                            - (\nu(B \backslash A) - \nu(A \cap B))   \\
                 & = & \nu(A \backslash B) - \nu(B \backslash A), \nonumber
\end{eqnarray}
and hence
\begin{eqnarray}
\|\nu(A) - \nu(B)\| & \le & \|\nu(A \backslash B)\| + \|\nu(B \backslash A)\|\\
         & \le & \mu(A \backslash B) + \mu(B \backslash A) = d(A, B),\nonumber
\end{eqnarray}
by (\ref{||nu(A)|| le mu(A)}).  This permits $\nu$ to be extended to a
$V$-valued function on $\mathcal{B}$, using uniform continuity and
completeness.  More precisely, if $A \in \mathcal{B}$, then there is a
sequence $\{A_j\}_{j = 1}^\infty$ of elements of $\mathcal{E}$ that
converges to $A$ with respect to $d(A, B)$.  This implies that
$\{\nu(A_j)\}_{j = 1}^\infty$ is a Cauchy sequence in $V$, because of
the uniform continuity of $\nu$ with respect to $d(\cdot,\cdot)$ just
established.  It follows that $\{\nu(A_j)\}_{j = 1}^\infty$ converges
in $V$, because $V$ is complete, and $\nu(A)$ is defined to be the
limit of this sequence.  One can also check that this does not depend
on the particular sequence $\{A_j\}_{j = 1}^\infty$ converging to $A$,
using the uniform continuity of $\nu$ with respect to $d(\cdot,
\cdot)$ again.  Note that this extension satisfies
\begin{equation}
\label{||nu(A) - nu(B)|| le d(A, B)}
        \|\nu(A) - \nu(B)\| \le d(A, B)
\end{equation}
for every $A, B \in \mathcal{B}$, since this holds when $A, B \in
\mathcal{E}$ and is preserved under limits.  In particular,
(\ref{||nu(A)|| le mu(A)}) holds for every $A \in \mathcal{B}$.

        Let $A, B \in \mathcal{B}$ be given, and let $\{A_j\}_{j =
1}^\infty$, $\{B_j\}_{j = 1}^\infty$ be sequences of elements of
$\mathcal{B}$ that converge to $A$, $B$ with respect to $d(\cdot,
\cdot)$, respectively.  This implies that $\{A_j \cap B_j\}_{j =
1}^\infty$ converges to $A \cap B$, and that $\{A_j \cup B_j\}_{j =
1}^\infty$ converges to $A \cup B$, as in Section \ref{distances
between sets}.  Of course,
\begin{equation}
        \nu(A_j) + \nu(B_j) = \nu(A_j \cap B_j) + \nu(A_j \cup B_j)
\end{equation}
for each $j$, because $\nu$ is finitely additive on $\mathcal{E}$.
Taking the limit as $j \to \infty$, we get that
\begin{equation}
        \nu(A) + \nu(B) = \nu(A \cap B) + \nu(A \cup B),
\end{equation}
because of (\ref{||nu(A) - nu(B)|| le d(A, B)}).  This shows that
$\nu$ is finitely additive on $\mathcal{B}$.

        If $E_1, E_2, \cdots$ is a sequence of elements of
$\mathcal{B}$ that are pairwise-disjoint, then
\begin{equation}
        \sum_{l = 1}^\infty \|\nu(E_l)\| \le \sum_{l = 1}^\infty \mu(E_l)
         = \mu\Big(\bigcup_{l = 1}^\infty E_l\Big),
\end{equation}
since (\ref{||nu(A)|| le mu(A)}) holds for every $A \in \mathcal{B}$.
Moreover,
\begin{equation}
 \biggl\|\nu\Big(\bigcup_{l = n + 1}^\infty E_l\Big)\biggr\| \le
  \mu\Big(\bigcup_{l = n + 1}^\infty E_l\Big) \to 0 \hbox{ as } n \to \infty.
\end{equation}
This implies that
\begin{equation}
        \sum_{l = 1}^\infty \nu(E_l) = \nu\Big(\bigcup_{l = 1}^\infty E_l\Big),
\end{equation}
because we already know that $\nu$ is finitely additive on $\mathcal{B}$.

\section[\ Paths and integrals]{Paths and integrals}
\label{paths, integrals}
\setcounter{equation}{0}

        Let $(V, \|v\|)$ be a real or complex Banach space, and let $F
: [a, b] \to V$ be a path of finite length, as in the previous
section.  Also let $\phi$ be a continuous real or complex-valued
function on $[a, b]$, as appropriate.  Suppose that $\mathcal{P} =
\{t_j\}_{j = 0}^n$ is a partition of $[a, b]$, and that $t_{j - 1} \le
r_j \le t_j$ for $j = 1, \ldots, n$, and consider
\begin{equation}
\label{sum_{j = 1}^n phi(r_j) (F(t_j) - F(t_{j - 1}))}
        \sum_{j = 1}^n \phi(r_j) \, (F(t_j) - F(t_{j - 1})).
\end{equation}
This is an approximation to the Riemann--Stieltjes integral of $\phi$
with respect to $F$, whose existence and basic properties will be
discussed now.  Basically, this is very similar to the
Riemann--Stieltjes integral of a continuous function with respect to a
real or complex-valued function of bounded variation on $[a, b]$.

        If $t_{j - 1} \le r'_j \le t_j$ is another collection of
intermediate points, then the difference of the corresponding sums can
be expressed as
\begin{eqnarray}
\lefteqn{\sum_{j = 1}^n \phi(r_j) \, (F(t_j) - F(t_{j - 1}))
          - \sum_{j = 1}^n \phi(r'_j) \, (F(t_j) - F(t_{j - 1}))} \\
 & = & \sum_{j = 1}^n (\phi(r_j) - \phi(r'_j)) \, (F(t_j) - F(t_{j - 1})).
                                                                \nonumber
\end{eqnarray}
Of course, $\phi$ is uniformly continuous on $[a, b]$, since it is
continuous and $[a, b]$ is compact.  Thus for each $\epsilon > 0$
there is a $\delta > 0$ such that
\begin{equation}
\label{|phi(r) - phi(r')| le epsilon}
        |\phi(r) - \phi(r')| \le \epsilon
\end{equation}
when $r, r' \in [a, b]$ and $|r - r'| < \delta$.  In particular,
\begin{eqnarray}
\lefteqn{\biggl\|\sum_{j = 1}^n \phi(r_j) \, (F(t_j) - F(t_{j - 1}))
         - \sum_{j = 1}^n \phi(r'_j) \, (F(t_j) - F(t_{j - 1}))\biggr\|} \\
 & \le & \sum_{j = 1}^n |\phi(r_j) - \phi(r'_j)| \, \|F(t_j) - F(t_{j - 1})\|
    \le \epsilon \, \Lambda_a^b                       \nonumber
\end{eqnarray}
when the mesh size $\max_{1 \le j \le n} (t_j - t_{j - 1})$ of
$\mathcal{P}$ is strictly less than $\delta$, where $\Lambda_a^b$
denotes the length of $F$ on $[a, b]$.

        If $\mathcal{P}$, $\widetilde{\mathcal{P}}$ are two partitions
of $[a, b]$ with sufficiently small mesh size, then one can check that
the difference between the corresponding sums (\ref{sum_{j = 1}^n
phi(r_j) (F(t_j) - F(t_{j - 1}))}) is also small.  As usual, it is
helpful to let $\widehat{\mathcal{P}}$ be a common refinement of
$\mathcal{P}$ and $\widetilde{\mathcal{P}}$, and to look at the
differences between the sums corresponding to $\mathcal{P}$,
$\widetilde{\mathcal{P}}$ and the sum corresponding to
$\widehat{\mathcal{P}}$.  These differences can be estimated in much
the same way as in the previous paragraph, using the uniform
continuity of $\phi$.  If $\mathcal{P}_1, \mathcal{P}_2, \ldots$ is a
sequence of partitions of $[a, b]$ whose mesh sizes are converging to
$0$, then the corresponding sums form a Cauchy sequence in $V$, and
hence converges, by completeness of $V$.  The limit does not depend on
the particular sequence of partitions, because the difference between
the sums associated to partitions with small mesh size is small, as
before.

        The Riemann-Stieltjes integral
\begin{equation}
        \int_a^b \phi \, dF
\end{equation}
of $\phi$ with respect to $F$ is the limit of the sums (\ref{sum_{j =
1}^n phi(r_j) (F(t_j) - F(t_{j - 1}))}) described in the previous
paragraph.  Observe that
\begin{equation}
        \biggl\|\sum_{j = 1}^n \phi(r_j) \, (F(t_j) - F(t_{j - 1}))\biggr\|
         \le \Big(\sup_{a \le r \le b} |\phi(r)|\Big) \, \Lambda_a^b
\end{equation}
for every partition $\mathcal{P}$ of $[a, b]$, and hence
\begin{equation}
\label{||int_a^b phi dF|| le (sup_{a le r le b} |phi(r)|) Lambda_a^b}
        \biggl\|\int_a^b \phi \, dF\biggr\|
         \le \Big(\sup_{a \le r \le b} |\phi(r)|\Big) \, \Lambda_a^b.
\end{equation}
If $\alpha(x)$ is the length of the restriction of $F$ to $[a, x]$ for
each $x \in [a, b]$, then one can improve this to get that
\begin{equation}
\label{||int_a^b phi dF|| le int_a^b |phi| d alpha}
        \biggl\|\int_a^b \phi \, dF\biggr\| \le \int_a^b |\phi| \, d\alpha,
\end{equation}
where the right side is a classical Riemann-Stieltjes integral.  This
is a more localized version of (\ref{||int_a^b phi dF|| le (sup_{a le
r le b} |phi(r)|) Lambda_a^b}), which can be derived using the
analogue of (\ref{||int_a^b phi dF|| le (sup_{a le r le b} |phi(r)|)
Lambda_a^b}) on small subintervals of $[a, b]$.  As in Section
\ref{functions, measures}, the Riemann--Stieltjes integral of a
continuous function on $[a, b]$ with respect to $\alpha$ can be
extended to the Lebesgue-Stieltjes integral with respect to a positive
Borel measure $\mu_\alpha$ on $[a, b]$.  As usual, continuous
functions on $[a, b]$ form a dense linear subspace of
$L^1(\mu_\alpha)$.  Using (\ref{||int_a^b phi dF|| le int_a^b |phi| d
alpha}), the Riemann--Stieltjes integral of $\phi$ with respect to $F$
can be extended to $\phi \in L^1(\mu_\alpha)$.  More precisely, if
$\phi$ is an integrable function on $[a, b]$ with respect to
$\mu_\alpha$, then there is a sequence $\{\phi_j\}_{j = 1}^\infty$ of
continuous functions on $[a, b]$ which converge to $\phi$ in
$L^1(\mu_\alpha)$.  Because of (\ref{||int_a^b phi dF|| le int_a^b
|phi| d alpha}), the corresponding sequence of Riemann--Stieltjes
integrals of the $\phi_j$'s with respect to $F$ form a Cauchy sequence
in $V$, and therefore converges, by completeness.  One can also check
that the limit depends only on $\phi$, and not on the particular
sequence of continuous approximations $\{\phi_j\}_{j = 1}^\infty$.
Hence the Lebesgue--Stieltjes integral of $\phi$ with respect to $F$
may be defined as this limit in $V$.  Of course, this is very similar
to the argument in the previous section.

\section[\ Integrating vector measures]{Integrating vector measures}
\label{integrating vector measures}
\setcounter{equation}{0}

        Let $(X, \mathcal{A})$ be a measurable space, and let $(V,
\|v\|)$ be a real or complex Banach space.  Also let $\mu$ be a
$V$-valued function on $\mathcal{A}$ such that for any sequence $A_1,
A_2, \ldots$ of pairwise-disjoint measurable subsets of $X$,
\begin{equation}
\label{sum_{j = 1}^infty ||mu(A_j)||}
        \sum_{j = 1}^\infty \|\mu(A_j)\|
\end{equation}
converges, and
\begin{equation}
\label{sum_{j = 1}^infty mu(A_j) = mu(bigcup_{j = 1}^infty A_j)}
        \sum_{j = 1}^\infty \mu(A_j) = \mu\Big(\bigcup_{j = 1}^\infty A_j\Big).
\end{equation}
As in Section \ref{vector-valued measures}, there is a nonnegative real-valued
measure $\|\mu\|$ on $X$ associated to $\mu$ such that
\begin{equation}
\label{||mu(A)|| le ||mu||(A)}
        \|\mu(A)\| \le \|\mu\|(A)
\end{equation}
for each $A \in \mathcal{A}$, and $\|\mu\|(X) < \infty$.

        Suppose that $f(x)$ is a real or complex-valued measurable
simple function on $X$, as appropriate.  This means that there are
finitely many pairwise-disjoint measurable subsets $A_1, \ldots, A_n$
of $X$ and real or complex numbers $\alpha_1, \ldots, \alpha_n$ such that
\begin{equation}
\label{f(x) = sum_{j = 1}^n alpha_j {bf 1}_{A_j}(x)}
        f(x) = \sum_{j = 1}^n \alpha_j \, {\bf 1}_{A_j}(x).
\end{equation}
Here ${\bf 1}_A(x)$ is the indicator function associated to $A
\subseteq X$ on $X$, equal to $1$ when $x \in A$ and to $0$ when $x
\in X \backslash A$.  The integral of $f$ with respect to $\mu$ is given by
\begin{equation}
\label{int_X f d mu = sum_{j = 1}^n alpha_j mu(A_j)}
        \int_X f \, d\mu = \sum_{j = 1}^n \alpha_j \, \mu(A_j),
\end{equation}
and satisfies
\begin{equation}
\label{||int_X f d mu|| le ... = int_X |f| d||mu||}
 \biggl\|\int_X f \, d\mu\biggr\| \le \sum_{j = 1}^n |\alpha_j| \, \|\mu(A_j)\|
                                   = \int_X |f| \, d\|\mu\|.
\end{equation}
More precisely, (\ref{int_X f d mu = sum_{j = 1}^n alpha_j mu(A_j)})
does not depend on the particular representation (\ref{f(x) = sum_{j =
1}^n alpha_j {bf 1}_{A_j}(x)}) of $f$, and it also works when the
$A_j$'s are not pairwise disjoint.

        Let $f(x)$ be an integrable real or complex-valued function on
$X$ with respect to $\|\mu\|$, as appropriate, and let $\{f_l\}_{l =
1}^\infty$ be a sequence of measurable simple functions on $X$ that
converge to $f$ in $L^1(X, \|\mu\|)$.  Using (\ref{||int_X f d mu|| le
... = int_X |f| d||mu||}), one can check that
\begin{equation}
\label{{int_X f_l d mu}_{l = 1}^infty}
        \bigg\{\int_X f_l \, d\mu\bigg\}_{l = 1}^\infty
\end{equation}
is a Cauchy sequence in $V$, and hence converges, by completeness.
The integral of $f$ with respect to $\mu$ can be defined by
\begin{equation}
\label{int_X f d mu = lim_{l to infty} int_X f_l d mu}
        \int_X f \, d\mu = \lim_{l \to \infty} \int_X f_l \, d\mu.
\end{equation}
As usual, one can also check that this does not depend on the sequence
$\{f_l\}_{l = 1}^\infty$ of simple functions converging to $f$, and that
\begin{equation}
\label{||int_X f d mu|| le int_X |f| d||mu||}
        \biggl\|\int_X f \, d\mu\biggr\| \le \int_X |f| \, d\|\mu\|,
\end{equation}
by (\ref{||int_X f d mu|| le ... = int_X |f| d||mu||}).

        If $\lambda$ is a bounded linear functional on $V$, then
\begin{equation}
        \mu_\lambda(A) = \lambda(\mu(A))
\end{equation}
defines a real or complex measure on $X$, as appropriate.  Note that
\begin{equation}
        |\mu_\lambda(A)| = |\lambda(\mu(A))| \le \|\lambda\|_* \, \|\mu(A)\|
                                              \le \|\lambda\|_* \, \|\mu\|(A),
\end{equation}
and hence
\begin{equation}
\label{|mu_lambda|(A) le ||lambda||_* ||mu||(A)}
        |\mu_\lambda|(A) \le \|\lambda\|_* \, \|\mu\|(A)
\end{equation}
for every $A \in \mathcal{A}$.  If $f$ is a measurable simple function
on $X$, then it is easy to see that
\begin{equation}
\label{lambda(int_X f d mu) = int_X f d mu_lambda}
        \lambda\Big(\int_X f \, d\mu\Big) = \int_X f \, d\mu_\lambda
\end{equation}
for every $\lambda \in V^*$.  This also works when $f \in L^1(X,
\|\mu\|)$, by approximating $f$ by simple functions, as in the
previous paragraph.  The integral of $f$ with respect to $\mu$ is
uniquely determined by this property, because of the Hahn--Banach
theorem.

\section[\ Measures and orthogonality]{Measures and orthogonality}
\label{measures, orthogonality}
\setcounter{equation}{0}

        Let $(X, \mathcal{A})$ be a measurable space, and let $(V,
\langle v, w \rangle)$ be a real or complex Hilbert space.  Also let
$\nu(A)$ be a finitely-additive $V$-valued measure on $(X,
\mathcal{A})$ such that
\begin{equation}
\label{langle nu(A), nu(B)rangle = 0}
        \langle \nu(A), \nu(B)\rangle = 0
\end{equation}
whenever $A$, $B$ are disjoint measurable subsets of $X$.  In
particular,
\begin{equation}
\label{||nu(A cup B)||^2 = ||nu(A)||^2 + ||nu(B)||^2}
        \|\nu(A \cup B)\|^2 = \|\nu(A)\|^2 + \|\nu(B)\|^2
\end{equation}
when $A$, $B$ are disjoint.  It follows that
\begin{equation}
\label{sum_{j = 1}^n ||nu(A_j)||^2 + ... =  ||nu(bigcup_{j = 1}^infty A_j)||^2}
        \sum_{j = 1}^n \|\nu(A_j)\|^2
              + \biggl\|\nu\Big(\bigcup_{j = n + 1}^\infty A_j\Big)\biggr\|^2
         =  \biggl\|\nu\Big(\bigcup_{j = 1}^\infty A_j\Big)\biggr\|^2
\end{equation}
for any sequence $A_1, A_2, \ldots$ of pairwise-disjoint measurable
subsets of $X$ and $n \ge 1$, and hence
\begin{equation}
\label{sum_{j = 1}^n ||nu(A_j)||^2 le ||nu(bigcup_{j = 1}^infty A_j)||^2}
        \sum_{j = 1}^n \|\nu(A_j)\|^2
         \le \biggl\|\nu\Big(\bigcup_{j = 1}^\infty A_j\Big)\biggr\|^2.
\end{equation}
Thus
\begin{equation}
        \sum_{j = 1}^\infty \|\nu(A_j)\|^2
         \le \biggl\|\nu\Big(\bigcup_{j = 1}^\infty A_j\Big)\biggr\|^2,
\end{equation}
which implies that $\sum_{j = 1}^\infty \nu(A_j)$ converges in $V$
when $A_1, A_2, \ldots$ are disjoint.  In this case, we ask also that
\begin{equation}
\label{sum_{j = 1}^infty nu(A_j) = nu(bigcup_{j = 1}^infty A_j)}
        \sum_{j = 1}^\infty \nu(A_j) = \nu\Big(\bigcup_{j = 1}^\infty A_j\Big),
\end{equation}
which implies that
\begin{equation}
\label{sum_{j = 1}^infty ||nu(A_j)||^2 = ||nu(bigcup_{j = 1}^infty A_j)||^2}
        \sum_{j = 1}^\infty \|\nu(A_j)\|^2
         = \biggl\|\nu\Big(\bigcup_{j = 1}^\infty A_j\Big)\biggr\|^2.
\end{equation}
This shows that $\|\nu(A)\|^2$ is a nonnegative real-valued measure on
$X$ under these conditions, which may be denoted $\|\nu\|^2$.

        As a basic example of this type of situation, let $\mu$ be a
nonnegative real-valued measure on $X$, and consider $V = L^2(X,
\mu)$, with the standard integral inner product.  Let $g \in L^2(X,
\mu)$ be given, and let $\nu_g$ be the $L^2(X, \mu)$-valued function
on $\mathcal{A}$ defined by
\begin{equation}
\label{nu_g(A) = g {bf 1}_A}
        \nu_g(A) = g \, {\bf 1}_A.
\end{equation}
Equivalently, $\nu_g(A)$ is the function on $X$ equal to $g$ on $A$
and to $0$ on $X \backslash A$ for each measurable set $A \subseteq
X$.  In particular,
\begin{equation}
\label{||nu_g(A)||^2 = int_A |g|^2 d mu}
        \|\nu_g(A)\|^2 = \int_A |g|^2 \, d\mu.
\end{equation}
It is easy to see that $\nu_g$ satisfies all of the conditions
described in the previous paragraph.

        Let $V$ be any Hilbert space again, and let $\nu$ be a
$V$-valued function on $\mathcal{A}$ that satisfies the same
conditions as before.  Let $A_1, \ldots, A_n$ be finitely many
pairwise-disjoint measurable subsets of $X$, and let $\alpha_1,
\ldots, \alpha_n$ be real or complex numbers, as appropriate.  If $f =
\sum_{j = 1}^n \alpha_j {\bf 1}_{A_j}$ is the corresponding simple
function, then its integral with respect to $\nu$ is given by
\begin{equation}
\label{int_X f d nu = sum_{j = 1}^n alpha_j nu(A_j)}
        \int_X f \, d\nu = \sum_{j = 1}^n \alpha_j \, \nu(A_j).
\end{equation}
In this case,
\begin{equation}
\label{||int_X f dnu||^2 = ...  = int_X |f|^2 d||nu||^2}
        \biggl\|\int_X f \, d\nu\biggr\|^2
         = \sum_{j = 1}^n |\alpha_j|^2 \, \|\nu(A_j)\|^2
         = \int_X |f|^2 \, d\|\nu\|^2.
\end{equation}
Using standard arguments based on continuity and completeness, the
integral of $f$ with respect to $\nu$ can be extended to an isometric
linear mapping from $L^2(X, \|\nu\|^2)$ into $V$.

        Suppose that $V = L^2(X, \mu)$ for some nonnegative
real-valued measure $\mu$ on $X$, and that $\nu = \nu_g$ for some $g
\in L^2(X, \mu)$.  If $f$ is a measurable simple function on $X$, then
it is easy to see that
\begin{equation}
\label{int_X f d nu_g = f g}
        \int_X f \, d\nu_g = f \, g
\end{equation}
as an element of $L^2(X, \mu)$, and that
\begin{equation}
\label{int_X |f|^2 d||nu_g||^2 = int_X |f|^2 |g|^2 d mu}
        \int_X |f|^2 \, d\|\nu_g\|^2 = \int_X |f|^2 \, |g|^2 \, d\mu.
\end{equation}
If $f \in L^2(X, \|\nu_g\|^2)$, then $f \, g \in L^2(X, \mu)$, and the
same statements hold.

\section[\ Paths and orthogonality]{Paths and orthogonality}
\label{paths, orthogonality}
\setcounter{equation}{0}

        Let $(V, \langle v, w \rangle)$ be a real or complex Hilbert
space, and let $p(t)$ be a $V$-valued function on a closed interval
$[a, b]$ in the real line.  Suppose that
\begin{equation}
\label{langle p(t_2) - p(t_1), p(t_3) - p(t_2) rangle = 0}
        \langle p(t_2) - p(t_1), p(t_3) - p(t_2) \rangle = 0
\end{equation}
whenever $a \le t_1 \le t_2 \le t_3 \le b$, which implies that
\begin{equation}
\label{langle p(t_2) - p(t_1), p(t_4) - p(t_3) rangle = 0}
        \langle p(t_2) - p(t_1), p(t_4) - p(t_3) \rangle = 0
\end{equation}
when $t_3 \le t_4 \le b$ too.  More precisely, (\ref{langle p(t_2) -
p(t_1), p(t_3) - p(t_2) rangle = 0}) also holds with $t_3$ replaced by
$t_4$ in this case, and (\ref{langle p(t_2) - p(t_1), p(t_4) - p(t_3)
rangle = 0}) follows by expressing $p(t_4) - p(t_3)$ as the difference
of $p(t_4) - p(t_2)$ and $p(t_3) - p(t_2)$.  If we put
\begin{equation}
\label{alpha(t) = ||p(t) - p(a)||^2}
        \alpha(t) = \|p(t) - p(a)\|^2
\end{equation}
for $a \le t \le b$, then
\begin{eqnarray}
\label{alpha(t) = ... = alpha(r) + ||p(t) - p(r)||^2 ge alpha(r)}
        \alpha(t) & = & \|p(r) - p(a)\|^2 + \|p(t) - p(r)\|^2        \\
                  & = & \alpha(r) + \|p(t) - p(r)\|^2 \ge \alpha(r) \nonumber
\end{eqnarray}
when $a \le r \le t \le b$, so that $\alpha(t)$ is monotone increasing
on $[a, b]$.  One can show that the one-sided limit $p(t+)$ exists
when $a \le t < b$, and similarly that $p(t-)$ exists when $a < t \le
b$, in analogy with Section \ref{continuity conditions}.  Note that
$p(t)$ is continuous at the same points where $\alpha(t)$ is
continuous, because of (\ref{alpha(t) = ... = alpha(r) + ||p(t) -
p(r)||^2 ge alpha(r)}).  It is convenient to extend $p(t)$ to the
whole real line, by putting $p(t) = p(a)$ when $t < a$ and $p(t) =
p(b)$ when $t > b$, so that $p(a-) = p(a)$ and $p(b+) = p(b)$ are
defined as well.  We can extend $\alpha(t)$ to ${\bf R}$ in the same
way, so that $\alpha(t) = 0$ when $t < a$ and $\alpha(t) = \alpha(b)$
when $t > b$.

        As in Sections \ref{functions, measures} and \ref{paths, measures}, put
\begin{equation}
        \nu((r, t)) = p(t-) - p(t+)
\end{equation}
and
\begin{equation}
        \nu([r, t)) = p(t-) - p(r-), \quad \nu((r, t]) = p(t+) - p(r+)
\end{equation}
when $a \le r < t \le b$, and
\begin{equation}
        \nu([r, t]) = p(t+) - p(r-)
\end{equation}
when $a \le r \le t \le b$.  This determines a finitely-additive
$V$-valued measure on the algebra $\mathcal{E}$ of subsets of $[a, b]$
that can be expressed as the union of finitely many intervals, where
the intervals may be open, closed, or half-open and half-closed.
By hypothesis,
\begin{equation}
        \langle \nu(I), \nu(I') \rangle = 0
\end{equation}
for every pair $I$, $I'$ of disjoint subintervals of $[a, b]$.
If $\mu_\alpha$ is the nonnegative Borel measure associated to
$\alpha(t)$ as in Section \ref{functions, measures}, then
\begin{equation}
\label{||nu(A)||^2 = mu_alpha(A)}
        \|\nu(A)\|^2 = \mu_\alpha(A)
\end{equation}
for every subinterval $A$ of $[a, b]$.  This also works when $A \in
\mathcal{E}$, because $A$ is then the union of finitely many
pairwise-disjoint subintervals $I_1, \ldots I_n$ of $[a, b]$, and
$\nu(I_1), \ldots, \nu(I_n)$ are orthogonal to each other in $V$.

        Let
\begin{equation}
        f(t) = \sum_{j = 1}^n c_j \, {\bf 1}_{I_j}(t)
\end{equation}
be a step function on $[a, b]$, where $I_1, \ldots, I_n$ are
pairwise-disjoint subintervals of $[a, b]$, and $c_1, \ldots, c_n$ are
real or complex numbers, as appropriate.  The integral of $f$ with
respect to $\nu$ can be defined by
\begin{equation}
       \int_a^b f \, d\nu = \sum_{j = 1}^n c_j \, \nu(I_j).
\end{equation}
In this case,
\begin{equation}
\biggl\|\int_a^b f \, d\nu\biggr\|^2 = \sum_{j = 1}^n |c_j|^2\, \|\nu(I_j)\|^2,
\end{equation}
because $\nu(I_1), \ldots, \nu(I_n)$ are orthogonal to each other in $V$.
Hence
\begin{equation}
        \biggl\|\int_a^b f \, d\nu\biggr\|^2 = \int_a^b |f|^2 \, d\mu_\alpha,
\end{equation}
as in (\ref{||nu(A)||^2 = mu_alpha(A)}).  Thus the integral of $f$
with respect to $\nu$ defines a linear isometry from the subspace of
$L^2([a, b], \mu_\alpha)$ consisting of step functions into $V$.  This
can be extended to a linear isometry from $L^2([a, b], \mu_\alpha)$
into $V$, by standard arguments of continuity and completeness.  In
particular, $\nu$ can be extended to a $V$-valued Borel measure on
$[a, b]$ as in the previous section, by applying this extension to
indicator functions of measurable subsets of $[a, b]$.

        If $\mu$ is a finite nonnegative Borel measure on $[a, b]$,
then $p(t) = {\bf 1}_{[a, t]}$ defines a mapping from $[a, b]$ into
$L^2([a, b], \mu)$ that satisfies the conditions mentioned at the
beginning of the section.  One could also use the indicator function
associated to $(a, t)$, $[a, t)$, or $(a, t]$, and the corresponding
differences of one-sided limits of $p$ would be the same.  Note that
these indicator functions are already the same in $L^2([a, b], \mu)$
when $\mu(\{x\}) = 0$ for each $x \in [a, b]$, in which case $p$ is
continuous.  One can check that $\mu_\alpha = \mu$ in this situation,
and that the embedding described in the preceding paragraph reduces to
the identity mapping on $L^2([a, b], \mu)$.

\section[\ Minkowski's integral inequality]{Minkowski's integral inequality}
\label{minkowski's integral inequality}
\setcounter{equation}{0}

        Let $(X, \mathcal{A}, \mu)$, $(Y, \mathcal{B}, \nu)$ be
measure spaces, with finite or $\sigma$-finite measure.  If $F(x, y)$
is a nonnegative measurable function on the Cartesian product $X
\times Y$ and $1 \le p < \infty$, then \emph{Minkowksi's integral
inequality} states that
\begin{eqnarray}
\label{(int_Y (int_X F dmu)^p dnu)^{1/p} le int_X (int_Y F^p dnu)^{1/p} dmu}
\lefteqn{\Big(\int_Y \Big(\int_X F(x, y) \, d\mu(x)\Big)^p \,
                                                 d\nu(y)\Big)^{1/p}} \\
& \le & \int_X \Big(\int_Y F(x, y)^p \, d\nu(y)\Big)^{1/p} \, d\mu(x).\nonumber
\end{eqnarray}
This is an integrated version of the triangle inequality for the $L^p$
norm, which is also known as Minkowski's inequality.  Note that one
has equality in (\ref{(int_Y (int_X F dmu)^p dnu)^{1/p} le int_X
(int_Y F^p dnu)^{1/p} dmu}) when $p = 1$, by Fubini's theorem.  We
have basically encountered versions of this already in connection with
conditional expectation, and we would like to mention a couple of
other approaches now.

        Let $A_1, \ldots, A_n$ be finitely many pairwise-disjoint
measurable subsets of $X$ whose union is equal to $X$.  If $F(x, y)$
is constant in $x$ on each $A_j$, then (\ref{(int_Y (int_X F dmu)^p
dnu)^{1/p} le int_X (int_Y F^p dnu)^{1/p} dmu}) reduces to the
ordinary Minkowski inequality for finite sums.  Otherwise, one can get
(\ref{(int_Y (int_X F dmu)^p dnu)^{1/p} le int_X (int_Y F^p dnu)^{1/p}
dmu}) by approximating $F(x, y)$ by functions of this type.  This is
analogous to the earlier discussion of ``nice functions'' on $X \times
Y$, but with the roles of $X$ and $Y$ exchanged.  A key point is that
measurable subsets of $X \times Y$ with finite measure can be
approximated by finite unions of measurable rectangles, as in Section
\ref{approximations in product spaces}.

        Alternatively, put
\begin{equation}
        N_p(F)(x) = \Big(\int_Y F(x, y)^p \, d\nu(y)\Big)^{1/p},
\end{equation}
as before.  If $\mu$ is a probability measure on $X$, then
\begin{equation}
        \Big(\int_X F(x, y) \, d\mu(x)\Big)^p \le \int_X F(x, y)^p \, d\mu(x)
\end{equation}
for each $y \in Y$, by Jensen's inequality.  Hence
\begin{eqnarray}
        \int_Y\Big(\int_X F(x, y) \, d\mu(x)\Big)^p \, d\nu(y)
         & \le & \int_Y \int_X F(x, y)^p \, d\mu(x) \, d\nu(y) \\
          & = & \int_X N_p(F)(x)^p \, d\mu(y), \nonumber
\end{eqnarray}
by Fubini's theorem.  If $N_p(F)(x) \le 1$ for $\mu$-almost every $x
\in X$, then it follows that
\begin{equation}
        \int_Y\Big(\int_X F(x, y) \, d\mu(x)\Big)^{1/p} \, d\nu(y) \le 1.
\end{equation}
This may be considered as a special case of (\ref{(int_Y (int_X F
dmu)^p dnu)^{1/p} le int_X (int_Y F^p dnu)^{1/p} dmu}), and the
general case may be derived from it using homogeneity, as follows.  If
the right side of (\ref{(int_Y (int_X F dmu)^p dnu)^{1/p} le int_X
(int_Y F^p dnu)^{1/p} dmu}) is equal to $0$, then $F(x, y) = 0$ almost
everywhere on $X \times Y$, the left side of (\ref{(int_Y (int_X F
dmu)^p dnu)^{1/p} le int_X (int_Y F^p dnu)^{1/p} dmu}) is also equal
to $0$, and there is nothing to do.  There is also nothing to do when
the right side of (\ref{(int_Y (int_X F dmu)^p dnu)^{1/p} le int_X
(int_Y F^p dnu)^{1/p} dmu}) is $+\infty$.  Thus we may suppose that
the right side of (\ref{(int_Y (int_X F dmu)^p dnu)^{1/p} le int_X
(int_Y F^p dnu)^{1/p} dmu}) is positive and finite, and we can even
take it to be equal to $1$, by multiplying $F$ by a positive constant.
We may also suppose that $N_p(F)(x) > 0$ for every $x \in X$, because
the $x \in X$ for which $N_p(F)(x) = 0$ do not play a role in
(\ref{(int_Y (int_X F dmu)^p dnu)^{1/p} le int_X (int_Y F^p dnu)^{1/p}
dmu}).  If we put
\begin{equation}
        F'(x, y) = N_p(x)^{-1} \, F(x, y),
\end{equation}
then $N_p(F')(x) = 1$ for every $x \in X$ automatically.  Similary, if we put
\begin{equation}
        \mu'(A) = \int_A N_p(F)(x) \, d\mu(x),
\end{equation}
then $\mu'$ is a probability measure on $X$, because the right side of
(\ref{(int_Y (int_X F dmu)^p dnu)^{1/p} le int_X (int_Y F^p dnu)^{1/p}
dmu}) is supposed to be equal to $1$.  The special case of Minkowski's
integral inequality under consideration implies that
\begin{equation}
        \int_Y\Big(\int_X F'(x, y) \, d\mu'(x)\Big)^p \, d\nu(y) \le 1.
\end{equation}
This implies that the left side of (\ref{(int_Y (int_X F dmu)^p
dnu)^{1/p} le int_X (int_Y F^p dnu)^{1/p} dmu}) is less than or equal
to $1$, as desired.

        Let $N_\infty(F)(x)$ be the essential supremum of $F(x, y)$
over $y \in Y$.  The $p = \infty$ version of (\ref{(int_Y (int_X F
dmu)^p dnu)^{1/p} le int_X (int_Y F^p dnu)^{1/p} dmu}) states that
the essential supremum of
\begin{equation}
\label{int_X F(x, y) d mu(y)}
        \int_X F(x, y) \, d\mu(y)
\end{equation}
over $y \in Y$ is less than or equal to
\begin{equation}
\label{int_X N_infty(F)(x) d mu(x)}
        \int_X N_\infty(F)(x) \, d\mu(x).
\end{equation}
If $Y$ is a probability space, then this can be obtained from
(\ref{(int_Y (int_X F dmu)^p dnu)^{1/p} le int_X (int_Y F^p dnu)^{1/p}
dmu}) by taking the limit as $p \to \infty$ with $p \in {\bf Z}_+$, as
in Section \ref{L^infty norms}.  Otherwise, one can reduce to the case
of probability spaces by approximating $Y$ by subsets of finite
measure, or using a positive weight on $Y$ with integral $1$.
Alternatively, if $N_\infty(F)(x) \le 1$ for almost every $x \in X$,
then $F(x, y) \le 1$ for almost every $(x, y) \in X \times Y$, by
Fubini's theorem.  If $\mu$ is a probability measure on $X$, then it
follows that (\ref{int_X F(x, y) d mu(y)}) is less than or equal to
$1$ for almost every $y \in Y$.  As in the previous paragraph, this
may be considered as a special case of the desired estimate, and the
general case can be derived from it in the same way as before.

\section[\ Spaces of measures]{Spaces of measures}
\label{spaces of measures}
\setcounter{equation}{0}

        Let $(X, \mathcal{A})$ be a measurable space, and let $(V,
\|v\|)$ be a real or complex Banach space.  Consider the space
$\mathcal{M}(X, V)$ of $V$-valued functions $\mu$ on $\mathcal{A}$ such that
\begin{equation}
        \sum_{j = 1}^\infty \|\mu(A_j)\| < \infty
\end{equation}
and
\begin{equation}
        \sum_{j = 1}^\infty \mu(A_j) = \mu\Big(\bigcup_{j = 1}^\infty A_j\Big)
\end{equation}
for every sequence $A_1, A_2, \ldots$ of pairwise-disjoint measurable
subsets of $X$.  As usual, the first condition already implies that
$\sum_{j = 1}^\infty \mu(A_j)$ converges in $V$.  The second condition
is equivalent to asking that $\mu$ be finitely additive and have the
continuity property that
\begin{equation}
\label{lim_{n to infty} mu(bigcup_{j = 1}^n A_j)= mu(bigcup_{j = 1}^infty A_j)}
        \lim_{n \to \infty} \mu\Big(\bigcup_{j = 1}^n A_j\Big)
         = \mu\Big(\bigcup_{j = 1}^\infty A_j\Big),
\end{equation}
just as for real or complex measures.

        Note that $\mathcal{M}(X, V)$ is a vector space over the real
or complex numbers, as appropriate.  If $\mu \in \mathcal{M}(X, V)$,
then $p(A) = \|\mu(A)\|$ satisfies the conditions in Section
\ref{uniform boundedness, 3}, and $\|\mu\|(A) = p^*(A)$ is a finite
nonnegative measure on $X$, as in Section \ref{vector-valued
measures}.  By construction,
\begin{equation}
        \|\mu(A)\| \le \|\mu\|(A)
\end{equation}
for every measurable set $A \subseteq X$, and $\|\mu\|(A)$ is the
smallest nonnegative measure on $X$ with this property, as in Section
\ref{uniform boundedness, 3}.  It is easy to check that $\|\mu\|(X)$
defines a norm on $\mathcal{M}(X, V)$.

        Suppose that $\mu_1, \mu_2, \ldots$ is a sequence of elements
of $\mathcal{M}(X, V)$ which is a Cauchy sequence with respect to this
norm.  Thus for each $\epsilon > 0$ there is an $L \ge 1$ such that
\begin{equation}
        \|\mu_l - \mu_n\|(X) < \epsilon
\end{equation}
for every $l, n \ge L$.  Of course,
\begin{equation}
 \|\mu_l(A) - \mu_n(A)\| \le \|\mu_l - \mu_n\|(A) \le \|\mu_l - \mu_n\|(X)
\end{equation}
for every measurable set $A \subseteq X$ and $l, n \ge 1$, which
implies that $\{\mu_l(A)\}_{l = 1}^\infty$ is a Cauchy sequence in $V$
for every $A \in \mathcal{A}$.  Let $\mu(A)$ be the limit of this
sequence in $V$, which converges by completeness.  Note that
$\{\mu_l(A)\}_{l = 1}^\infty$ actually converges to $\mu(A)$ uniformly
on $\mathcal{A}$, because the Cauchy condition holds uniformly over $A
\in \mathcal{A}$.

        If $A_1, A_2, \ldots$ is a sequence of pairwise-disjoint
measurable subsets of $X$, then
\begin{equation}
 \sum_{j = 1}^\infty \|\mu_l(A_j)\| \le \sum_{j = 1}^\infty \|\mu_l\|(A_j)
  = \|\mu_l\|\Big(\bigcup_{j = 1}^\infty A_j\Big) \le \|\mu_l\|(X)
\end{equation}
for each $l$.  In the limit as $l \to \infty$, we get that
\begin{equation}
\label{sum_{j = 1}^infty ||mu(A_j)|| le sup_{l ge 1} ||mu_l||(X)}
        \sum_{j = 1}^\infty \|\mu(A_j)\| \le \sup_{l \ge 1} \|\mu_l\|(X).
\end{equation}
The right side is finite because $\{\mu_l\}_{l = 1}^\infty$ is a
Cauchy sequence, and hence is bounded.  It is easy to see that
$\mu(A)$ is finitely additive, since $\mu_l(A)$ is finitely additive
for each $l$.  The continuity condition (\ref{lim_{n to infty}
mu(bigcup_{j = 1}^n A_j)= mu(bigcup_{j = 1}^infty A_j)}) can also be
derived from the corresponding property of the $\mu_l$'s, using the
fact that $\{\mu_l(A)\}_{l = 1}^\infty$ converges to $\mu(A)$
uniformly on $\mathcal{A}$.  Similarly, if $A_1, A_2, \ldots$ is a sequence of
pairwise-disjoint measurable subsets of $X$ whose union is equal to $X$, then
\begin{equation}
        \sum_{j = 1}^\infty \|\mu_l(A_j) - \mu_n(A_j)\|
 \le \sum_{j = 1}^\infty \|\mu_l(A_j) - \mu_n\|(A_j) = \|\mu_l - \mu_n\|(X)
\end{equation}
for each $l, n \ge 1$, as before.  This implies that
\begin{equation}
        \sum_{j = 1}^\infty \|\mu(A_j) - \mu_n(A_j)\|
         \le \sup_{l \ge n} \|\mu_l - \mu_n\|(X)
\end{equation}
for each $n \ge 1$, by taking the limit as $l \to \infty$, as in
(\ref{sum_{j = 1}^infty ||mu(A_j)|| le sup_{l ge 1} ||mu_l||(X)}).  It
follows that
\begin{equation}
        \|\mu - \mu_n\|(X) \le \sup_{l \ge n} \|\mu_l - \mu_n\|(X)
\end{equation}
for each $n$, by taking the supremum over all such partitions
$\{A_j\}_{j = 1}^\infty$ of $X$.  This shows that $\mu \in
\mathcal{M}(X, V)$ and that $\{\mu_n\}_{n = 1}^\infty$ converges to
$\mu$ with respect to the norm $\|\mu\|(X)$, and hence that
$\mathcal{M}(X, V)$ is complete.

\section[\ Products and measures]{Products and measures}
\label{products, measures}
\setcounter{equation}{0}

        Let $(X, \mathcal{A}, \mu)$, $(Y, \mathcal{B}, \nu)$ be finite
or $\sigma$-finite measure spaces, and let $F(x, y)$ be a measurable
function on $X \times Y$.  As usual, we put
\begin{equation}
\label{def of N_p(F)(x)}
        N_p(F)(x) = \Big(\int_Y |F(x, y)|^p \, d\nu(y)\Big)^{1/p}
\end{equation}
when $1 \le p < \infty$, and we let $N_\infty(F)(x)$ be the essential
supremum of $|F(x, y)|$ over $y \in Y$.  Suppose that
\begin{equation}
        \int_X N_p(F)(x) \, d\mu(x) < \infty
\end{equation}
for some $p$, $1 \le p \le \infty$, and put
\begin{equation}
        \phi(A)(y) = \int_A F(x, y) \, d\mu(x)
\end{equation}
for each measurable set $A \subseteq X$.  This defines $\phi(A)$
as a measurable function on $Y$ which is in $L^p(Y)$ and satisfies
\begin{equation}
\label{||phi(A)||_{L^p(Y)} le int_A N_p(F)(x) d mu(x)}
        \|\phi(A)\|_{L^p(Y)} \le \int_A N_p(F)(x) \, d\mu(x),
\end{equation}
by Minkowski's integral inequality.  If $A_1, A_2, \ldots$ is a
sequence of pairwise-disjoint measurable subsets of $X$, then
\begin{eqnarray}
        \sum_{j = 1}^\infty \|\phi(A_j)\|_{L^p(Y)}
         & \le & \sum_{j = 1}^\infty \int_{A_j} N_p(F)(x) \, d\mu(x) \\
& = & \int_{\bigcup_{j = 1}^\infty A_j} N_p(F)(x) \, d\mu(x) < \infty.\nonumber
\end{eqnarray}
Thus $\sum_{j = 1}^\infty \phi(A_j)$ converges in $L^p(Y)$, and it is
easy to see that
\begin{equation}
 \sum_{j = 1}^\infty \phi(A_j) = \phi\Big(\bigcup_{j = 1}^\infty A_j\Big).
\end{equation}
Hence $\phi \in \mathcal{M}(X, L^p(Y))$.  If $\|\phi\|(A)$ is as in the
previous section, then
\begin{equation}
\label{||phi||(A) le int_A N_p(F)(x) d mu(x)}
        \|\phi\|(A) \le \int_A N_p(F)(x) \, d\mu(x),
\end{equation}
because of (\ref{||phi(A)||_{L^p(Y)} le int_A N_p(F)(x) d mu(x)}).

\section[\ $L^p$-Valued measures]{$L^p$-Valued measures}
\label{L^p-valued measures}
\setcounter{equation}{0}

        Let $(X, \mathcal{A})$ be a measurable space, and let $(Y,
\mathcal{B}, \nu)$ be a $\sigma$-finite measure space.  Suppose that
$\mu \in \mathcal{M}(X, L^p(Y))$ for some $p$, $1 < p \le \infty$.
Thus $\|\mu\|$ is a finite nonnegative real measure on $X$, and we can
consider the product measure $\|\mu\| \times \nu$ on $X \times Y$.  We
would like to represent $\mu$ by a function on $X \times Y$, as in the
previous section.

        Let $A_1, \ldots, A_n$ be finitely many pairwise-disjoint
measurable subsets of $X$ such that $\bigcup_{j = 1}^n A_j = X$, and
let $g_1(y), \ldots, g_n(y)$ be elements of $L^q(Y)$, where $1 \le q <
\infty$ is the exponent conjugate to $p$, $1/p + 1/q = 1$. 
Put $G(x, y) = g_j(y)$ when $x \in A_j$, and
\begin{equation}
\label{L(G) = sum_{j = 1}^n int_Y mu(A_j)(y) g_j(y) d nu(y)}
        L(G) = \sum_{j = 1}^n \int_Y \mu(A_j)(y) \, g_j(y) \, d\nu(y).
\end{equation}
By H\"older's inequality,
\begin{eqnarray}
\lefteqn{\biggl|\int_Y \mu(A_j)(y) \, g_j(y) \, d\nu(y)\biggr|} \\
         & &  \le \|\mu(A_j)\|_{L^p(Y)} \, \|g_j\|_{L^q(Y)} 
                     \le \|\mu\|(A_j) \, \|g_j\|_{L^q(Y)} \nonumber
\end{eqnarray}
for each $j$.  This implies that
\begin{equation}
\label{|L(G)| le int_X N_q(G)(x) d||mu||}
        |L(G)| \le \int_X N_q(G)(x) \, d\|\mu\|,
\end{equation}
where $N_q(G)(x)$ denotes the $L^q(Y)$ norm of $G(x, y)$ as a function
of $y$, as usual.  In particular,
\begin{equation}
\label{|L(G)| le ||mu||(X)^{1/p} ||G||_{L^q(X times Y, ||mu|| times nu)}}
 |L(G)| \le \|\mu\|(X)^{1/p} \, \|G\|_{L^q(X \times Y, \|\mu\| \times \nu)}.
\end{equation}

        It is easy to see that (\ref{L(G) = sum_{j = 1}^n int_Y
mu(A_j)(y) g_j(y) d nu(y)}) does not depend on the particular
representation of $G(x, y)$ in the preceding paragraph, because $\mu$
is finitely additive.  One can also check that the collection of these
functions $G(x, y)$ forms a linear subspace of $L^q(X \times Y,
\|\mu\| \times \nu)$, and that $L(G)$ defines a linear functional on
this subspace.  The main point is that any two partitions of $X$ into
finitely many measurable sets has a common refinement, and so any two
functions of this type can be represented in this way using the same
partition of $X$.  This subspace is also dense in $L^q(X \times Y,
\|\mu\| \times \nu)$, because $q < \infty$.  We also know from
(\ref{|L(G)| le ||mu||(X)^{1/p} ||G||_{L^q(X times Y, ||mu|| times
nu)}}) that $L(G)$ is a bounded linear functional on this subspace,
with respect to the $L^q$ norm, and hence has a unique extension to a
bounded linear functional on $L^q(X \times Y, \|\mu\| \times \nu)$.

        The Riesz representation theorem implies that there is a unique
element $F(x, y)$ of $L^p(X \times Y, \|\mu\| \times \nu)$ such that
\begin{equation}
\label{L(G) = int_{X times Y} F(x, y) G(x, y) d||mu||(x) d nu(y)}
 L(G) = \int_{X \times Y} F(x, y) \, G(x, y) \, d\|\mu\|(x) \, d\nu(y)
\end{equation}
for every $G \in L^q(X \times Y, \|\mu\| \times \nu)$.  If $A$ is a
measurable subset of $X$ and $g(y) \in L^q(Y)$, then we can apply this
to $G(x, y) = {\bf 1}_A(x) \, g(y)$, to get that
\begin{equation}
        \int_Y \mu(A)(y) \, g(y) \, d\nu(y)
         = \int_Y \Big(\int_A F(x, y) \, d\mu(x)\Big) \, g(y) \, d\nu(y).
\end{equation}
It follows that
\begin{equation}
\label{mu(A)(y) = int_A F(x, y) d mu(x)}
        \mu(A)(y) = \int_A F(x, y) \, d\mu(x)
\end{equation}
as elements of $L^p(Y)$ for every measurable set $A \subseteq Y$, as
in the previous section.  Moreover,
\begin{equation}
\label{||F||_{L^p(X times Y, ||mu|| times nu)} le ||mu||(X)^{1/p}}
        \|F\|_{L^p(X \times Y, \|\mu\| \times \nu)} \le \|\mu\|(X)^{1/p},
\end{equation}
because of (\ref{|L(G)| le ||mu||(X)^{1/p} ||G||_{L^q(X times Y,
||mu|| times nu)}}).  If $p = \infty$, then this say that the
$L^\infty$ norm of $F(x, y)$ is less than or equal to $1$ on $X \times
Y$.  Otherwise, if $p < \infty$, and if $A$ is a measurable subset of
$X$, then (\ref{|L(G)| le int_X N_q(G)(x) d||mu||}) implies that
\begin{equation}
\label{|L(G)| le ||mu||(A)^{1/p} ||G||_{L^q(A times Y, ||mu|| times nu)}}
 |L(G)| \le \|\mu\|(A)^{1/p} \, \|G\|_{L^q(A \times Y, \|\mu\| \times \nu)}
\end{equation}
when $G(x, y) = 0$ for every $x \in X \backslash A$.  Hence
\begin{equation}
        \Big(\int_A \int_Y |F(x, y)|^p \, d\mu(x) \, d\nu(y)\Big)^{1/p}
         \le \|\mu\|(A)^{1/p},
\end{equation}
or equivalently,
\begin{equation}
        \int_A N_p(F)(x)^p \, d\mu(x) \le \|\mu\|(A).
\end{equation}
This shows that $N_p(F)(x) \le 1$ almost everywhere on $X$ with
respect to $\|\mu\|$.

\section[\ $\ell^1$-Valued measures]{$\ell^1$-Valued measures}
\label{ell^1-valued measures}
\setcounter{equation}{0}

        Let $(X, \mathcal{A})$ be a measurable space, and let $\mu_1,
\mu_2, \ldots$ be a sequence of real or complex-valued measures on $X$
such that $\sum_{j = 1}^\infty |\mu_j|(X) < \infty$.  This implies that
\begin{equation}
        \sum_{j = 1}^\infty |\mu_j(A)| \le \sum_{j = 1}^\infty |\mu_j|(A)
         \le \sum_{j = 1}^\infty |\mu_j|(X) < \infty
\end{equation}
for every measurable set $A \subseteq X$, which means that $\mu(A) =
\{\mu_j(A)\}_{j = 1}^\infty \in \ell^1$ for each $A \in \mathcal{A}$.
Put $\rho(A) = \sum_{j = 1}^\infty |\mu_j|(A)$, so that $\rho$ is a finite
nonnegative real measure on $X$ by hypothesis, and
\begin{equation}
\label{||mu(A)||_1 = sum_{j = 1}^infty |mu_j(A)| le rho(A)}
        \|\mu(A)\|_1 = \sum_{j = 1}^\infty |\mu_j(A)| \le \rho(A)
\end{equation}
for each $A \in \mathcal{A}$.  Using this, one can check that $\mu \in
\mathcal{M}(X, \ell^1)$, and that $\|\mu\|(A) \le \rho(A)$ for each $A
\in \mathcal{A}$.

        This construction is actually equivalent to the one in Section
\ref{products, measures}, with $p = 1$ and $Y = {\bf Z}_+$, equipped
with counting measure.  This is because $\mu_j$ is absolutely
continuous with respect to $\rho$ for each $j$, and hence can be
expressed in terms of an integrable function $f_j$ with respect to
$\rho$, as in the Radon--Nikodym theorem.  The $L^1$ norm of $f_j$
with respect to $\rho$ is equal to $|\mu_j|(X)$ for each $j$, and is
summable over $j$.  Thus the sequence of $f_j$'s can be identified
with an integrable function on $X \times {\bf Z}_+$, using $\rho$ as
the measure on $X$.

        Conversely, suppose that $\mu \in \mathcal{M}(X, \ell^1)$.  Thus
$\mu(A) = \{\mu_j(A)\}_{j = 1}^\infty$ for some real or complex-valued
functions $\mu_j$ on $\mathcal{A}$, as appropriate.  It is easy to see
that $\mu_j$ is a real or complex measure on $X$ for each $j$, because
of the corresponding properties of $\mu$.  A key point now is that
\begin{equation}
        \sum_{j = 1}^\infty |\mu_j|(A) \le \|\mu\|(A)
\end{equation}
for every $A \in \mathcal{A}$.  Of course, it suffices to show that
\begin{equation}
        \sum_{j = 1}^n |\mu_j|(A) \le \|\mu\|(A)
\end{equation}
for every $A \in \mathcal{A}$ and $n \ge 1$.  Remember that
$|\mu_j|(A) = p_j^*(A)$ is defined as in Section \ref{uniform
boundedness, 3}, using $p_j(A) = |\mu_j(A)|$.  More precisely,
$p_j^*(A)$ can be defined as the supremum of sums of $p_j$ over
partitions of $A$ into finitely many measurable subsets.  If we use
the same partition of $A$ for each $j$, then the desired estimate
would follow from the definition of $\|\mu\|(A)$ as $p^*(A)$ with
$p(A) = \|\mu(A)\|_{\ell^1}$.  If instead we have different partitions
of $A$ for $j = 1, \ldots, n$, then we can use a common refinement of
them to reduce to the case of a single partition of $A$.

        Suppose now that $\mu \in \mathcal{M}(X, \ell^p)$, $1 \le p
\le \infty$.  As in the preceding paragraph, $\mu(A) = \{\mu_j(A)\}_{j
= 1}^\infty$, where each $\mu_j$ is a real or complex measure on $X$.
It is easy to see that $\mu_j$ is absolutely continuous with respect
to $\|\mu\|$ for each $j$, and so can be expressed in terms of an
integrable function with respect to $\|\mu\|$, by the Radon--Nikodym
theorem.  If $p = 1$, then the $L^1$ norms of these functions are
summable, as before.  If $p > 1$, then we are back in the situation of
the previous section, with $Y = {\bf Z}_+$ equipped with counting
measure.

\section[\ Finite sums]{Finite sums}
\label{sums}
\setcounter{equation}{0}

        Let $(X, \mathcal{A})$ be a measurable space, and let $(V,
\|v\|)$ be a real or complex Banach space.  Suppose that $\mu_1,
\ldots, \mu_n$ are finitely many real or complex measures on $X$, as
appropriate, and that $v_1, \ldots, v_n$ are vectors in $V$.  It is
easy to see that
\begin{equation}
\label{mu(A) = sum_{j = 1}^n mu_j(A) v_j}
        \mu(A) = \sum_{j = 1}^n \mu_j(A) \, v_j
\end{equation}
defines an element of $\mathcal{M}(X, V)$.  Of course,
\begin{equation}
        \|\mu(A)\| \le \sum_{j = 1}^n |\mu_j(A)| \, \|v_j\|
                    \le \sum_{j = 1}^n |\mu_j|(A) \, \|v_j\|
\end{equation}
for each $A \in \mathcal{A}$, which implies that
\begin{equation}
        \|\mu\|(A) \le \sum_{j = 1}^n |\mu_j|(A) \, \|v_j\|.
\end{equation}

        Let $\rho$ be a finite nonnegative real measure on $X$ such
that $\mu_j$ is absolutely continuous with respect to $\rho$ for each
$j$.  One can take
\begin{equation}
        \rho = \sum_{j = 1}^n |\mu_j|,
\end{equation}
for instance.  By the Radon--Nikodym theorem, there are integrable
functions $f_1, \ldots, f_n$ on $X$ with respect to $\rho$ such that
\begin{equation}
\label{mu_j(A) = int_A f_j d rho}
        \mu_j(A) = \int_A f_j \, d\rho
\end{equation}
for each $A \in \mathcal{A}$ and $j = 1, \ldots, n$.  If $f(x) =
\sum_{j = 1}^n f_j(x) \, v_j$, then
\begin{equation}
        \|\mu(A)\| = \biggl\|\sum_{j = 1}^n v_j \int_A f_j \, d\rho\biggr\|
                   \le \int_A \|f\| \, d\rho
\end{equation}
for each $A \in \mathcal{A}$, as in Section \ref{sigma-subalgebras, vectors}.
This implies that
\begin{equation}
        \|\mu\|(A) \le \int_A \|f\| \, d\rho
\end{equation}
for each $A \in \mathcal{A}$.

        More precisely,
\begin{equation}
        \|\mu\|(A) = \int_A \|f\| \, d\rho
\end{equation}
for each $A \in \mathcal{A}$ under these conditions.  To see this,
remember that
\begin{equation}
        \sum_{k = 1}^l \|\mu(A_k)\| \le \|\mu\|(A)
\end{equation}
when $A_1, \ldots, A_l$ are pairwise-disjoint measurable sets whose
union is $A$, by definition of $\|\mu\|(A)$.  In order to show that
\begin{equation}
\label{int_A ||f|| d rho le ||mu||(A)}
        \int_A \|f\| \, d\rho \le \|\mu\|(A),
\end{equation}
one can choose measurable sets $A_k$ on which the $f_j$'s are
approximately constant.

        Let us now start with a measure $\mu \in \mathcal{M}(X, V)$
that takes values in a finite-dimensional linear subspace of $V$.  If
$v_1, \ldots, v_n$ is a basis for this linear subspace, then there are
unique real or complex measures $\mu_1, \ldots, \mu_n$ on $X$ for
which $\mu$ can be expressed as in (\ref{mu(A) = sum_{j = 1}^n mu_j(A)
v_j}).  Because any two norms on a finite-dimensional real or complex
vector space are equivalent,
\begin{equation}
 \bigg\|\sum_{j = 1}^n t_j \, v_j\biggr\| \ge c \, \max_{1 \le j \le n} |t_j|
\end{equation}
for some $c > 0$ and every $t_1, \ldots, t_n \in {\bf R}$ or ${\bf C}$,
as appropriate.  This implies that
\begin{equation}
\label{c max_{1 le j le n} |mu_j(A)| le ||mu(A)|| le ||mu||(A)}
        c \, \max_{1 \le j \le n} |\mu_j(A)| \le \|\mu(A)\| \le \|\mu\|(A)
\end{equation}
for each $A \in \mathcal{A}$, and hence that $\mu_j$ is absolutely
continuous with respect to $\|\mu\|$ for each $j$.  Thus we can take
$\rho = \|\mu\|$ in the previous paragraphs, and it follows that the
corresponding function $f$ satisfies $\|f(x)\| = 1$ for almost every
$x \in X$ with respect to $\|\mu\|$.

\section[\ Approximations]{Approximations}
\label{approximations}
\setcounter{equation}{0}

        Let $(X, \mathcal{A})$ be a measurable space, and let $(V,
\|v\|)$ be a real or complex Banach space.  Suppose that $\mu_1,
\mu_2, \ldots$ is a sequence of elements of $\mathcal{M}(X, V)$ such
that $\mu_j$ takes values in a finite-dimensional linear subspace
$V_j$ of $V$ for each $j$.  Suppose also that $\{\mu_j\}_{j =
1}^\infty$ converges to $\mu \in \mathcal{M}(X, V)$ with respect to
the total variation norm, so that
\begin{equation}
\label{lim_{j to infty} ||mu_j - mu||(X) = 0}
        \lim_{j \to \infty} \|\mu_j - \mu\|(X) = 0.
\end{equation}
Let $\rho$ be a finite nonnegative real measure on $X$ such that
$\|\mu_j\|$ is absolutely continuous with respect to $\rho$ for each
$j$, such as
\begin{equation}
        \rho(A) = \sum_{j = 1}^\infty a_j \, \|\mu_j\|(A)
\end{equation}
for some $a_j > 0$ with $\sum_{j = 1}^\infty a_j \, \|\mu_j\|(X) < \infty$.
Thus each $\mu_j$ can be expressed as
\begin{equation}
        \mu_j(A) = \int_A f_j \, d\rho
\end{equation}
for some $V_j$-valued integrable function $f_j$ on $X$ with respect to
$\rho$, by applying the Radon--Nikodym theorem to the components of
$\mu(A)$ with respect to a basis for $V_j$ as in the previous section.
More precisely, each $f_j$ is the sum of finitely many real or
complex-valued integrable functions on $X$ with respect to $\rho$
times basis vectors of $V_j$, and the integral of $f_j$ over $A$ is
the sum of the integrals of the components of $f_j$ over $A$ times the
corresponding basis vectors of $V_j$.  We also have that
\begin{equation}
\label{int_X ||f_j - f_l|| d rho = ||mu_j - mu_l||(X)}
\int_X \|f_j - f_l\| \, d\rho = \|\mu_j - \mu_l\|(X) \to 0
\end{equation}
as $j, l \to \infty$, because of (\ref{lim_{j to infty} ||mu_j -
mu||(X) = 0}).

\section[\ Uniform convexity]{Uniform convexity}
\label{uniform convexity}
\setcounter{equation}{0}

        Let $V$ be a vector space with a norm $\|v\|$.  It will be
convenient to take $V$ to be a real vector space here, but complex
vector spaces can also be considered as real vector spaces, and so
everything in this section works as well in that case.  We say that
$V$ is \emph{uniformly convex} if for every $\epsilon > 0$ there is a
$\delta > 0$ such that
\begin{equation}
        v, w \in V, \ \|v\| = \|w\| = 1, \hbox{ and } 
          \biggl\|\frac{v + w}{2}\biggr\| > 1 - \delta
\end{equation}
imply that
\begin{equation}
        \|v - w\| < \epsilon.
\end{equation}
It is easy to see that inner product spaces are uniformly convex,
because of the parallelogram law.  It is well known that real and
complex $L^p$ spaces are uniformly convex when $1 < p < \infty$.

        Suppose that $v, w \in V$, $\|v\|, \|w\| \le 1$, and
\begin{equation}
        \biggl\|\frac{v + w}{2}\biggr\| > 1 - \delta_1
\end{equation}
for some $\delta_1 \in (0, 1/2)$.  In particular,
\begin{equation}
        \frac{\|v\| + \|w\|}{2} > 1 - \delta_1 > \frac{1}{2},
\end{equation}
and so $\|v\|, \|w\| > 0$.  If $v' = v/\|v\|$, $w' = w/\|w\|$, then
\begin{equation}
        \|v' - v\| = (\|v\|^{-1} - 1) \, \|v\| = 1 - \|v\|,
\end{equation}
and similarly for $w$.  Thus
\begin{equation}
 \frac{\|v' - v\| + \|w' - w\|}{2} = 1 - \frac{\|v\| + \|w\|}{2} < \delta_1,
\end{equation}
which implies that
\begin{eqnarray}
 \biggl\|\frac{v + w}{2}\biggr\| & \le & \biggl\|\frac{v' + w'}{2}\biggr\|
                                     + \frac{\|v - v'\| + \|w - w'\|}{2} \\
         & < & \biggl\|\frac{v' + w'}{2}\biggr\| + \delta_1 \nonumber
\end{eqnarray}
and
\begin{equation}
        \biggl\|\frac{v' + w'}{2}\biggr\| > 1 - 2 \, \delta_1.
\end{equation}
If $\delta_1$ is sufficiently small, then
\begin{equation}
        \|v' - w'\| < \epsilon/2,
\end{equation}
by uniform convexity.  If also $\delta_1 \le \epsilon/4$, then
\begin{equation}
        \|v - w\| \le \|v' - w'\| + \|v - v'\| + \|w - w'\|
             < \epsilon/2 + 2 \, \delta_1 \le \epsilon.
\end{equation}
This shows that uniform convexity implies the analogous condition
in which $\|v\|, \|w\| \le 1$.

        Suppose that $v_1, \ldots, v_n \in V$, $\|v_j\| \le 1$ for $j
= 1, \ldots, n$, $t_1, \ldots, t_n$ are nonnegative real numbers, and
that $\sum_{j = 1}^n t_j = 1$.  Let $\epsilon > 0$ be given, and put
\begin{equation}
        a = \sum_{j = 1}^n t_j \, v_j.
\end{equation}
Thus $\|a\| \le 1$, and we would like to show that there is an $\eta >
0$ such that $\|a\| > 1 - \eta$ implies that
\begin{equation}
\label{sum_{j = 1}^n t_j ||v_j - a|| < epsilon}
        \sum_{j = 1}^n t_j \|v_j - a\| < \epsilon,
\end{equation}
where $\eta$ does not depend on $n$, the $v_j$'s, or the $t_j$'s.  Let
$\lambda$ be a bounded linear functional on $V$ such that
$\|\lambda\|_* = 1$ and $\lambda(a) = \|a\|$, the existence of which
follows from the Hahn--Banach theorem, as usual.  Hence
\begin{equation}
        \sum_{j = 1}^n t_j \, \lambda(v_j) = \lambda(a) = \|a\| > 1 - \eta,
\end{equation}
which implies that
\begin{equation}
\label{sum_{j = 1}^n t_j (1 - lambda_j(v_j)) < eta}
        \sum_{j = 1}^n t_j \, (1 - \lambda_j(v_j)) < \eta.
\end{equation}
Note that $1 - \lambda(v_j) \ge 0$ for each $j$, because
$|\lambda(v_j)| \le 1$.  In addition,
\begin{equation}
 \biggl\|\frac{v_j + a}{2}\biggr\| \ge \lambda\Big(\frac{v_j + a}{2}\Big)
                                      = \frac{\lambda(v_j) + \|a\|}{2}.
\end{equation}

        Let $\delta_2$ be associated to $\epsilon/2$ as in the second
version of uniform convexity.  If $\lambda(v_j) > 1 - \delta_2$ and
$\eta \le \delta_2$, then
\begin{equation}
 \|(v_j + a)/2\| > \frac{(1 - \delta_2) + (1 - \eta)}{2} \ge 1 - \delta_2,
\end{equation}
and so
\begin{equation}
        \|v_j - a\| < \epsilon/2.
\end{equation}
Let $I_1$ be the set of $j = 1, \ldots, n$ such that $\lambda(v_j) > 1
- \delta_2$, and let $I_2$ be the set of $j = 1, \ldots, n$ such that
$\lambda(v_j) \le 1 - \delta_2$.  If $\eta \le \delta_2$, then
\begin{equation}
        \sum_{j \in I_1} t_j \, \|v_j - a\| < \epsilon/2,
\end{equation}
by the preceding computation.  Of course, $\|v_j - a\| \le \|v_j\| +
\|a\| \le 2$ for each $j$, and so
\begin{equation}
 \sum_{j \in I_2} t_j \|v_j - a\| \le 2 \sum_{j \in I_2} t_j.
\end{equation}
Using (\ref{sum_{j = 1}^n t_j (1 - lambda_j(v_j)) < eta}), we get that
\begin{equation}
        \sum_{j \in I_2} t_j \, \delta_2
         \le \sum_{j \in I_2} t_j \, (1 - \lambda(v_j)) < \eta,
\end{equation}
which implies that
\begin{equation}
        \sum_{j \in I_2} t_j \, \|v_j - a\| \le 2 \sum_{j \in I_2} t_j
                                               < 2 \, \delta_2^{-1} \, \eta.
\end{equation}
Therefore
\begin{eqnarray}
 \sum_{j = 1}^n t_j \, \|v_j - a\| & = & \sum_{j \in I_1} t_j \, \|v_j - a\| 
                                      + \sum_{j \in I_2} t_j \, \|v_j - a\| \\
        & < & \epsilon/2 + 2 \, \delta_2^{-1} \, \eta \le \epsilon \nonumber
\end{eqnarray}
when $\eta$ is sufficiently small, as desired.

\section[\ Uniform convexity and measures]{Uniform convexity and measures}
\label{uniform convexity, measures}
\setcounter{equation}{0}

        Let $(X, \mathcal{A})$ be a measurable space, and let $(V,
\|v\|)$ be a uniformly convex Banach space.  Also let $\epsilon > 0$
be given, and let $\eta$ be as in the previous section.  Suppose that
$\mu \in \mathcal{M}(X, V)$ satisfies
\begin{equation}
        \|\mu(X)\| > (1 - \eta) \, \|\mu\|(X).
\end{equation}
Let $\mu_0 \in \mathcal{M}(X, V)$ be defined by
\begin{equation}
\label{mu_0(A) = frac{mu(X)}{||mu||(X)} ||mu||(A)}
        \mu_0(A) = \frac{\mu(X)}{\|\mu\|(X)} \, \|\mu\|(A),
\end{equation}
so that $\mu_0$ is the vector $\mu(X)/\|\mu\|(X)$ times the
nonnegative real measure $\|\mu\|$ on $X$.  We would like to show that
\begin{equation}
        \|\mu - \mu_0\|(X) \le \epsilon \, \|\mu\|(X)
\end{equation}
under these conditions.

        We may as well suppose also that $\|\mu\|(X) = 1$, since
otherwise we can divide $\mu$ by $\|\mu\|(X) > 0$.  Let $A_1, \ldots,
A_n$ be finitely many pairwise disjoint measurable subsets of $X$ such
that $X = \bigcup_{j = 1}^n A_j$, and let us check that
\begin{equation}
\label{sum_{j = 1}^n ||mu(A_j) - mu_0(A_j)|| < epsilon}
        \sum_{j = 1}^n \|\mu(A_j) - \mu_0(A_j)\| < \epsilon.
\end{equation}
If $\|\mu\|(A_j) = 0$ for some $j$, then $\mu(A_j) = \mu_0(A_j) = 0$,
and we can absorb $A_j$ into one of the other $A_l$'s without
affecting the sum.  Thus we may as well ask that $\|\mu\|(A_j) > 0$
for each $j$ too.  If we put
\begin{equation}
\label{t_j = ||mu||(A_j), v_j = frac{mu(A_j)}{||mu||(A_j)}}
 t_j = \|\mu\|(A_j) \quad\hbox{and}\quad  v_j = \frac{\mu(A_j)}{\|\mu\|(A_j)},
\end{equation}
then $\sum_{j = 1}^n t_j = 1$ and $\|v_j\| \le 1$ for each $j$,
because $\|\mu(A_j)\| \le \|\mu\|(A_j)$.  Also,
\begin{equation}
        \sum_{j = 1}^n t_j \, v_j = \sum_{j = 1}^n \mu(A_j) = \mu(X),
\end{equation}
and
\begin{equation}
\sum_{j = 1}^n \|\mu(A_j) - \mu_0(A_j)\| = \sum_{j = 1}^n t_j \|v_j - \mu(X)\|.
\end{equation}
Thus (\ref{sum_{j = 1}^n ||mu(A_j) - mu_0(A_j)|| < epsilon}) reduces
to (\ref{sum_{j = 1}^n t_j ||v_j - a|| < epsilon}), with $a = \mu(X)$.

        Now let $\mu$ be any element of $\mathcal{M}(X, V)$, and let
$\theta$ be a small positive real number.  By the definition of
$\|\mu\|(X)$, there are finitely many pairwise-disjoint measurable
sets $X_1, \ldots, X_r$ such that $X = \bigcup_{l = 1}^r X_l$ and
\begin{equation}
\label{||mu||(X) < sum_{l = 1}^r ||mu(X_l)| + theta}
        \|\mu\|(X) < \sum_{l = 1}^r \|\mu(X_l)\| + \theta.
\end{equation}
Of course, $\|\mu\|(X) = \sum_{l = 1}^r \|\mu\|(X_l)$, and so
\begin{equation}
\label{sum_{l = 1}^r (||mu||(X_l) - ||mu(X_l)||) < theta}
        \sum_{l = 1}^r (\|\mu\|(X_l) - \|\mu(X_l)\|) < \theta.
\end{equation}
Each term in the sum is nonnegative, since $\|\mu(X_l)\| \le \|\mu\|(X_l)$.
If $L_2$ is the set of $l = 1, \ldots, r$ such that
\begin{equation}
\label{||mu(X_l)|| le (1 - eta) ||mu||(X_l)}
        \|\mu(X_l)\| \le (1 - \eta) \, \|\mu\|(X_l),
\end{equation}
where $\eta > 0$ is as before, then it follows that
\begin{equation}
\label{eta sum_{l in L_2} ||mu||(X_l) le ... < theta}
        \eta \sum_{l \in L_2} \|\mu\|(X_l)
         \le \sum_{l \in L_2} (\|\mu\|(X_l) - \|\mu(X_l)\|) < \theta.
\end{equation}

        Let $L_1$ be the set of $l = 1, \ldots, r$ such that
$\|\mu(X_l)\| > (1 - \eta) \, \|\mu\|(X_l)$, and for each $l \in L_1$,
let $\mu_l \in \mathcal{M}(X, V)$ be defined by
\begin{equation}
\label{mu_l(A) = frac{mu(X_l)}{||mu||(X_l)} ||mu||(A cap X_l)}
        \mu_l(A) = \frac{\mu(X_l)}{\|\mu\|(X_l)} \, \|\mu\|(A \cap X_l).
\end{equation}
This is analogous to (\ref{mu_0(A) = frac{mu(X)}{||mu||(X)}
||mu||(A)}), applied to the restriction of $\mu$ to $X_l$,
and it follows from the earlier discussion that
\begin{equation}
        \|\mu - \mu_l\|(X_l) \le \epsilon \, \|\mu\|(X_l)
\end{equation}
for each $l \in L_1$.  Combining this with the earlier estimate
(\ref{eta sum_{l in L_2} ||mu||(X_l) le ... < theta}) for $L_2$, we get that
\begin{equation}
        \biggl\|\mu - \sum_{l \in L_1} \mu_l\biggr\|
          < \epsilon \, \sum_{l \in L_1} \|\mu\|(X_l) + \eta^{-1} \, \theta
             \le \epsilon \, \|\mu\|(X) + \eta^{-1} \, \theta.
\end{equation}
Remember that $\eta$ depends on $\epsilon$, while $\theta$ can be
chosen independently of $\epsilon$, $\eta$.  Thus the right side can
be made arbitrarily small, by choosing $\epsilon$ and then $\theta$
appropriately.

\section[\ Uniform convexity and paths]{Uniform convexity and paths}
\label{uniform convexity, paths}
\setcounter{equation}{0}

        Let $(V, \|v\|)$ be a uniformly convex Banach space, and let
$f : [x, y] \to V$ be a path of finite length $\Lambda_x^y$.  Also let
$\epsilon > 0$ be given, and let $\eta = \eta(\epsilon)$ be as in
Section \ref{uniform convexity}.  Suppose that
\begin{equation}
        \|f(x) - f(y)\| > (1 - \eta) \, \Lambda_x^y.
\end{equation}
Put
\begin{equation}
\label{f_0(z) = frac{f(y) - f(x)}{Lambda_x^y} Lambda_x^z}
        f_0(z) = \frac{f(y) - f(x)}{\Lambda_x^y} \, \Lambda_x^z,
\end{equation}
where $\Lambda_x^z$ is the length of $f$ on $[x, z]$, $x \le z
\le y$.  We would like to show that
\begin{equation}
\label{the length of f - f_0 on [x, y] is le epsilon Lambda_x^y}
 \hbox{the length of $f - f_0$ on $[x, y]$ is } \le \epsilon \, \Lambda_x^y.
\end{equation}
This is basically the same as the argument for measures in the
previous section.  As before, we may as well suppose that $\Lambda_x^y
= 1$, since otherwise we can divide $f$ by $\Lambda_x^y$.

        If $\{r_j\}_{j = 0}^n$ is any partition of $[x, y]$, then we
would like to show that
\begin{eqnarray}
\label{sum_{j = 1}^n ||(f(r_j) - f_0(r_j)) - (f(r_{j - 1}) - f_0(r_{j - 1}))||}
\lefteqn{\sum_{j = 1}^n \|(f(r_j) - f_0(r_j)) -
                                      (f(r_{j - 1}) - f_0(r_{j - 1}))\|} \\
 & = & \sum_{j = 1}^n \|(f(r_j) - f(r_{j - 1})) - (f_0(r_j) - f_0(r_{j - 1}))\|
         < \epsilon. \nonumber
\end{eqnarray}
We may as well ask that the length $\Lambda_{r_{j - 1}}^{r_j}$ of $f$
on $[r_{j - 1}, r_j]$ be positive for each $j = 1, \ldots, n$, since
otherwise $f$, $f_0$ are constant on $[r_{j - 1}, r_j]$, and $r_j$ or
$r_{j - 1}$ could be removed from the partition without affecting the
sum.  Put
\begin{equation}
        t_j = \Lambda_{r_{j - 1}}^{r_j} \quad\hbox{and}\quad
          v_j = \frac{f(r_j) - f(r_{j - 1})}{\Lambda_{r_{j - 1}}^{r^j}},
\end{equation}
so that $\sum_{j = 1}^n t_j = 1$ and $\|v_j\| \le 1$ for each $j$,
because $\|f(r_j) - f(r_{j - 1})\| \le \Lambda_{r_{j - 1}}^{r_j}$.
Moreover,
\begin{equation}
        \sum_{j = 1}^n t_j \, v_j = \sum_{j = 1}^n (f(r_j) - f(r_{j - 1}))
                                   = f(y) - f(x)
\end{equation}
and
\begin{eqnarray}
\lefteqn{\sum_{j = 1}^n \|(f(r_j) - f(r_{j - 1}))
                                       - (f_0(r_j) - f_0(r_{j - 1}))\|} \\
         & = & \sum_{j = 1}^n t_j \, \|v_j - (f(y) - f(x))\|. \nonumber
\end{eqnarray}
Thus (\ref{sum_{j = 1}^n ||(f(r_j) - f_0(r_j)) - (f(r_{j - 1}) -
f_0(r_{j - 1}))||}) follows from (\ref{sum_{j = 1}^n t_j ||v_j - a|| <
epsilon}), with $a = f(y) - f(x)$.

        Now let $f : [a, b] \to V$ be a path of finite length
$\Lambda_a^b$, and let $\theta$ be a small positive real number.  
By the definition of $\Lambda_a^b$, there is a partition 
$\{x_l\}_{l = 0}^r$ of $[a, b]$ such that
\begin{equation}
        \Lambda_a^b < \sum_{j = 1}^r \|f(x_l) - f(x_{l - 1})\| + \theta.
\end{equation}
This implies that
\begin{equation}
\label{sum_{l = 1}^r (Lambda_{x_{l - 1}}^{x_l} - ||f(x_l) - f(x_{l - 1})||)}
        \sum_{l = 1}^r (\Lambda_{x_{l - 1}}^{x_l} - \|f(x_l) - f(x_{l - 1})\|)
            < \theta,
\end{equation}
because $\Lambda_a^b = \sum_{l = 1}^r \Lambda_{x_{l - 1}}^{x_l}$.
Note that the terms in the sum are nonnegative, since $\|f(x_l) -
f(x_{l - 1})\| \le \Lambda_{x_{l - 1}}^{x_l}$.  If $L_2$ is the set of
$l = 1, \ldots, r$ such that
\begin{equation}
        \|f(x_l) - f(x_{l - 1})\| \le (1 - \eta) \, \Lambda_{x_{l - 1}}^{x_l},
\end{equation}
where $\eta > 0$ is as before, then
\begin{equation}
\label{eta sum_{l in L_2} Lambda_{x_{l - 1}}^{x_l} le ... < theta}
        \eta \, \sum_{l \in L_2} \Lambda_{x_{l - 1}}^{x_l}
  \le \sum_{l \in L_2} (\Lambda_{x_{l - 1}}^{x_l} - \|f(x_l) - f(x_{l - 1})\|)
             < \theta.
\end{equation}

        Let $L_1$ be the set of $l = 1, \ldots, r$ such that
\begin{equation}
        \|f(x_l) - f(x_{l - 1})\| > (1 - \eta) \, \Lambda_{x_l}^{x_{l - 1}}.
\end{equation}
If $l \in L_2$, the define $f_l : [a, b] \to V$ by
\begin{equation}
        f_l(z) = \frac{f(x_l) - f(x_{l - 1})}{\Lambda_{x_{l - 1}}^{x_l}} 
                                                    \, \Lambda_{x_{l - 1}}^z
\end{equation}
when $x_{l - 1} \le z \le x_l$, and put $f(z) = 0$ when $z \le x_{l -
1}$, $f(z) = f(x_l) - f(x_{l - 1})$ when $z \ge x_l$.  This is the
same as (\ref{f_0(z) = frac{f(y) - f(x)}{Lambda_x^y} Lambda_x^z}) on
$[x_{l - 1}, x_l]$ with $x = x_{l -1}$, $y = x_l$.  As in (\ref{the
length of f - f_0 on [x, y] is le epsilon Lambda_x^y}), the length of
$f - f_l$ on $[x_{l - 1}, x_l]$ is less than or equal to $\epsilon \,
\Lambda_{x_{l - 1}}^{x_l}$.  Combining this with (\ref{eta sum_{l in
L_2} Lambda_{x_{l - 1}}^{x_l} le ... < theta}), we get that the length
of $f - \sum_{j \in L_1} f_j$ on $[a, b]$ is less than or equal to
\begin{equation}
 \sum_{l \in L_1} \epsilon \, \Lambda_{x_{l - 1}}^{x_l} + \eta^{-1} \, \theta
                            \le \epsilon \, \Lambda_a^b + \eta^{-1} \, \theta.
\end{equation}
This uses the fact that the length of a path on $[a, b]$ is the sums
of the lengths of its restrictions to the intervals $[x_{l - 1},
x_l]$, $1 \le l \le r$.  If $l \in L_1$, then $f_j$ is constant on
$[x_{l - 1}, x_l]$ when $j \ne l$, by construction, and so the length
of $f - \sum_{j \in L_1} f_j$ is the same as the length of $f - f_l$
on this interval.  Similarly, if $l \in L_2$, then $f_j$ is constant
on $[x_{l - 1}, x_l]$ for each $j \in L_1$, and the length of $f -
\sum_{j \in L_1} f_j$ is the same as the length of $f$ on this
interval.  It follows from this estimate that the length of $f -
\sum_{j \in L_1} f_j$ can be made arbitrarily small, first by choosing
$\epsilon$ to be very small, and then choosing $\theta$ to be
sufficiently small, depending on $\eta$, which also depends on
$\epsilon$.

\section[\ Uniform convexity and martingales]{Uniform convexity and martingales}\label{uniform convexity, martingales}
\setcounter{equation}{0}

        Let $(V, \|v\|)$ be a uniformly convex Banach space.  Also let
$\epsilon > 0$ be given, and let $\eta > 0$ be as in Section
\ref{uniform convexity}.  We may as well ask that $\eta \le \epsilon$
too, which is practically unavoidable anyway.

        Let $(X, \mathcal{A}, \mu)$ be a probability space, and let
$\mathcal{P}_1, \mathcal{P}_2, \ldots$ be a sequence of partitions of
$X$ into finitely or countably many pairwise disjoint measurable
subsets of positive measure such that $\mathcal{P}_{j + 1}$ is a
refinement of $\mathcal{P}_j$ for each $j$.  As usual, the arguments
that follows are a bit simpler when each $\mathcal{P}_j$ has only
finitely many elements, but countable partitions and other situations
can be accommodated as well.  Let $\mathcal{B}_j =
\mathcal{B}(\mathcal{P}_j)$ be the $\sigma$-algebra of measurable
subsets of $X$ generated by $\mathcal{P}_j$, as in Section
\ref{partitions}, so that $\mathcal{B}_j \subseteq \mathcal{B}_{j +
1}$ for each $j$.

        We would like to consider $V$-valued martingales on $X$ with
respect to this filtration, as in Section \ref{vector-valued
martingales}.  Remember that a $V$-valued function $f_j$ on $X$ is
measurable with respect to $\mathcal{B}_j$ if and only if it is
constant on the elements of $\mathcal{P}_j$.  Suppose that we have a
sequence $\{f_j\}_{j = 1}^\infty$ of $V$-valued functions on $X$ such
that $f_j$ is measurable with respect to $\mathcal{B}_j$ for each $j$
and $\|f_j\|$ has bounded $L^1$ norm.  Suppose also that $\{f_j\}_{j =
1}^\infty$ is a martingale with respect to the $\mathcal{B}_j$'s, so
that the value of $f_j$ on $B \in \mathcal{P}_j$ is equal to the
average of the values of $f_{j + 1}$ on the sets $A \in \mathcal{P}_{j
+ 1}$ with $A \subseteq B$.

        Under these conditions, $\{\|f_j\|\}_{j = 1}^\infty$ is a
submartingale on $X$ with respect to the $\mathcal{B}_j$'s.  In
particular, the $L^1$ norm of $\|f_j\|$ is monotone increasing in $j$,
and so
\begin{equation}
        \lim_{j \to \infty} \int_X \|f_j\| \, d\mu 
            = \sup_{j \ge 1} \int_X \|f_j\| \, d\mu.
\end{equation}
Let $\theta$ be a small positive real number, and suppose that
\begin{equation}
\label{sup_{n ge 1} int_X ||f_n|| d mu < int_X ||f_j|| d mu + theta}
 \sup_{n \ge 1} \int_X \|f_n\| \, d\mu < \int_X \|f_j\| \, d\mu + \theta.
\end{equation}
Note that
\begin{equation}
\label{int_B ||f_{n + 1}|| d mu ge int_B ||f_n|| d mu}
        \int_B \|f_{n + 1}\| \, d\mu \ge \int_B \|f_n\| \, d\mu
\end{equation}
when $B \in \mathcal{P}_j$ and $n \ge j$, because $\{\|f_n\|\}_{n =
1}^\infty$ is a submartingale.  Moreover,
\begin{eqnarray}
\label{lim_{n to infty} int_X ||f_n|| d mu = ...}
        \lim_{n \to \infty} \int_X \|f_n\| \, d\mu
 & = & \lim_{n \to \infty} \sum_{B \in \mathcal{P}_j} \int_B \|f_n\| \, d\mu
                                                                          \\
 & = & \sum_{B \in \mathcal{P}_j} \lim_{n \to \infty} \int_B \|f_n\| \, d\mu.
                                                               \nonumber
\end{eqnarray}
This is obvious when $\mathcal{P}_j$ has only finitely many elements,
and otherwise one can use the monotone convergence theorem for sums.
It follows that
\begin{equation}
\label{sum_{B in mathcal{P}_j} ... < theta}
 \sum_{B \in \mathcal{P}_j} \Big(\lim_{n \to \infty} \int_B \|f_n\| \, d\mu
                                      - \int_B \|f_j\| \, d\mu\Big) < \theta,
\end{equation}
where each term in the sum is nonnegative.

        Let $\mathcal{P}_j'$ be the set of $B \in \mathcal{P}_j$ such that
\begin{equation}
        \int_B \|f_j\| \, d\mu > (1 - \eta) \,
                        \lim_{n \to \infty} \int_B \|f_n\| \, d\mu.
\end{equation}
Thus $\mathcal{P}_j'' = \mathcal{P}_j \backslash \mathcal{P}_j'$
consists of $B \in \mathcal{P}_j$ such that
\begin{equation}
        \int_B \|f_j\| \, d\mu \le (1 - \eta) \,
                        \lim_{n \to \infty} \int_B \|f_n\| \, d\mu,
\end{equation}
and satisfies
\begin{equation}
\label{eta sum_{B in P_j''} lim_{n to infty} int_B ||f_n|| d mu < theta}
        \eta \sum_{B \in \mathcal{P}_j''}
         \lim_{n \to \infty} \int_B \|f_n\| \, d\mu < \theta,
\end{equation}
by (\ref{sum_{B in mathcal{P}_j} ... < theta}).

        Let $f_n(A)$ be the value of $f_n$ on $A \in \mathcal{P}_n$,
as in Section \ref{vector-valued martingales}.  Thus
\begin{equation}
        f_j(B) = \sum_{A \in \mathcal{P}_n \atop A \subseteq B} f_n(A) \,
                                                    \frac{\mu(A)}{\mu(B)}
\end{equation}
when $B \in \mathcal{P}_j$ and $n \ge j$, because $\{f_n\}_{n =
1}^\infty$ is a martingale.  In addition,
\begin{equation}
        \int_B \|f_j\| \, d\mu = \|f_j(B)\| \, \mu(B)
\end{equation}
and
\begin{equation}
 \int_B \|f_n\| \, d\mu = \sum_{A \in \mathcal{P}_n \atop A \subseteq B}
                                                     \|f_n(A)\| \, \mu(A).
\end{equation}
Note that
\begin{equation}
        \int_B \|f_n\| \, d\mu \ge \int_B \|f_j\| \, d\mu > 0
\end{equation}
when $B \in \mathcal{P}_j'$ and $n \ge j$, and put
\begin{equation}
        t_n(A) = \|f_n(A)\| \, \mu(A) \, \Big(\int_B \|f_n\| \, d\mu\Big)^{-1}
\end{equation}
for each $A \in \mathcal{P}_n$ with $A \subseteq B$, so that
\begin{equation}
        \sum_{A \in \mathcal{P}_n \atop A \subseteq B} t_n(A) = 1,
\end{equation}
by construction.  Also put $v_n(A) = f_n(A) / \|f_n(A)\|$ when $A \in
\mathcal{P}_n$ and $f_n(A) \ne 0$, and $v_n(A) = 0$ when $f_n(A) = 0$,
so that
\begin{eqnarray}
 \sum_{A \in \mathcal{P}_n \atop A \subseteq B} v_n(A) \, t_n(A)
 & = & \sum_{A \in \mathcal{P}_n \atop A \subseteq B} f_n(A) \, \mu(A) \, 
                                      \Big(\int_B \|f_n\| \, d\mu\Big)^{-1} \\
 & = & f_j(B) \, \mu(B) \, \Big(\int_B \|f_n\| \, d\mu\Big)^{-1}. \nonumber
\end{eqnarray}
It follows that
\begin{equation}
\label{||sum_{A in mathcal{P}_n atop A subseteq B} v_n(A) t_n(A)|| > 1 - eta}
\biggl\|\sum_{A \in \mathcal{P}_n \atop A \subseteq B} v_n(A) \, t_n(A)\biggr\|
          > 1 - \eta
\end{equation}
when $B \in \mathcal{P}_j'$ and $n \ge j$.

        This is exactly the situation discussed in Section
\ref{uniform convexity}, except that the sum in (\ref{||sum_{A in
mathcal{P}_n atop A subseteq B} v_n(A) t_n(A)|| > 1 - eta}) may have
infinitely many terms, which can be handled in the same way as before.
If
\begin{equation}
 a_{j, n}(B) = f_j(B) \, \mu(B) \, \Big(\int_B \|f_n\| \, d\mu)\Big)^{-1},
\end{equation}
then we get that
\begin{equation}
\label{sum_{A in mathcal{P}_n atop A subseteq B} ... < epsilon}
        \sum_{A \in \mathcal{P}_n \atop A \subseteq B} 
                            \|v_n(A) - a_{j, n}(B)\| \, t_n(A) < \epsilon
\end{equation}
when $B \in \mathcal{P}_j'$ and $n \ge j$.  Put $a_j(B) =
f_j(B)/\|f_j(B)\|$, which is the same as $a_{j, j}(B)$, and observe that
\begin{equation}
 a_{j, n}(B) = a_j(B) \, \frac{\int_B \|f_j\| \, d\mu}{\int_B \|f_n\| \, d\mu}.
\end{equation}
This implies that
\begin{equation}
        \|a_j(B) - a_{j, n}(B)\| < \eta
\end{equation}
when $B \in \mathcal{P}_j'$ and $n \ge j$.  Combining this with
(\ref{sum_{A in mathcal{P}_n atop A subseteq B} ... < epsilon}), we get that
\begin{equation}
 \sum_{A \in \mathcal{P}_n \atop A \subseteq B} \|v_n(A) - a_j(B)\| \, t_n(A)
                                          < \epsilon + \eta \le 2 \, \epsilon
\end{equation}
when $B \in \mathcal{P}_j'$ and $n \ge j$.

        Equivalently,
\begin{equation}
 \sum_{A \in \mathcal{P}_n \atop A \subseteq B} \|v_n(A) - a_j(B)\| \,
           \|f_n(A)\| \, \mu(A) < 2 \, \epsilon \, \int_B \|f_n\| \, d\mu
\end{equation}
when $B \in \mathcal{P}_j'$ and $n \ge j$, which reduces to
\begin{equation}
 \sum_{A \in \mathcal{P}_n \atop A \subseteq B} \bigl\|f_n(A) - a_j(B) \,
         \|f_n(A)\|\bigr\| \, \mu(A) < 2 \, \epsilon \, \int_B \|f_n\| \, d\mu,
\end{equation}
using the definition of $v_n(A)$.  The sum on the left can be
expressed as in integral, so that
\begin{equation}
        \int_B \bigl\|f_n - a_j(B) \, \|f_n\|\bigr\| \, d\mu
          < 2 \, \epsilon \, \int_B \|f_n\| \, d\mu
\end{equation}
when $B \in \mathcal{P}_j'$ and $n \ge j$.  Put $a_j(B) = 0$ when $B
\in \mathcal{P}_j''$, and let $a_j(x)$ be the $V$-valued function on
$X$ equal to $a_j(B)$ when $x \in B \in \mathcal{P}$.  Summing the
previous estimate over $B \in \mathcal{P}_j'$, and using (\ref{eta
sum_{B in P_j''} lim_{n to infty} int_B ||f_n|| d mu < theta}) for $B
\in \mathcal{P}_j''$, we get that
\begin{equation}
\label{int_X ||f_n - a_j ||f_n|| || d mu < ...}
        \int_X \bigl\|f_n - a_j \, \|f_n\|\bigr\| \, d\mu
          < 2 \, \epsilon \, \int_X \|f_n\| \, d\mu + \eta^{-1} \, \theta
\end{equation}
when $n \ge j$.

        As usual, the right side of (\ref{int_X ||f_n - a_j ||f_n|| ||
d mu < ...}) can be made arbitrarily small, by first choosing
$\epsilon$ to be as small as one likes, and then choosing $\theta$
depending on $\eta$, which depends on $\epsilon$.  This works
uniformly over $n \ge j$, because the $L^1$ norm of $\|f_n\|$ is
bounded, by hypothesis.  Because $\{\|f_n\|\}_{n = 1}^\infty$ is a
submartingale on $X$ with bounded integral, there is a real-valued
martingale $\{g_n\}_{n = 1}^\infty$ on $X$ such that $\|f_n\| \le g_n$ and
\begin{equation}
        \int_X g_n \, d\mu = \lim_{l \to \infty} \int_X \|f_l\| \, d\mu,
\end{equation}
for each $n$, as in Section \ref{submartingales}.  Of course, the
integral of $g_n$ over $X$ is independent of $n$, because of the
martingale condition.  In particular,
\begin{equation}
\label{int_X (g_n - ||f_n||) d mu < theta}
        \int_X (g_n - \|f_n\|) \, d\mu < \theta
\end{equation}
when $n \ge j$, by (\ref{sup_{n ge 1} int_X ||f_n|| d mu < int_X
||f_j|| d mu + theta}) and the monotonicity of the integral of $\|f_n\|$.
Using (\ref{int_X ||f_n - a_j ||f_n|| || d mu < ...}), we get that
\begin{equation}
\label{int_X ||f_n - a_j g_n|| d mu < ...}
        \int_X \|f_n - a_j \, g_n\| \, d\mu 
 < 2 \, \epsilon \, \int_X \|f_n\| \, d\mu + (\eta^{-1} + 1) \, \theta
\end{equation}
when $n \ge j$, since $\|a_j(x)\| \le 1$ for every $x \in X$, by construction.
Note that $\{a_j \, g_n\}_{n = j}^\infty$ is a $V$-valued martingale
on $X$, because $\{g_n\}_{n = 1}^\infty$ is a martingale on $X$ and
$a_j$ is constant on the elements of $\mathcal{P}_j$.

\section[\ Strict convexity]{Strict convexity}
\label{strict convexity}
\setcounter{equation}{0}

        Let $V$ be a real vector space with a norm $\|v\|$.  As
before, a complex vector space is automatically a real vector space
too, and so everything in this section can be used in that case as
well.  The closed unit ball
\begin{equation}
        B_1 = \{v \in V : \|v\| \le 1\}
\end{equation}
in $V$ is said to be \emph{strictly convex} if for every $v, w \in
B_1$ with $v \ne w$ and every $t \in {\bf R}$ with $0 < t < 1$ we have
that
\begin{equation}
\label{||t v + (1 - t) w|| < 1}
        \|t \, v + (1 - t) \, w\| < 1.
\end{equation}
Of course, (\ref{||t v + (1 - t) w|| < 1}) holds automatically when
$\|v\| < 1$ or $\|w\| < 1$, and so it suffices to check this when
$\|v\| = \|w\| = 1$.  One can show that the unit ball in an inner
product space is strictly convex by determining when equality occurs
in the Cauchy--Schwarz inequality.  The unit ball in an $L^p$ space is
strictly convex when $1 < p < \infty$, because of the strict convexity
of the function $|x|^p$ on the real line.  This is similar to the
proof of the convexity of the unit ball in $L^p$ using the convexity
of $|x|^p$, as in Section \ref{p-summable functions}.  Note that $B_1$
is strictly convex when $V$ is uniformly convex.

        Let $\lambda$ be a nonzero bounded linear functional on $V$,
and suppose that $v$, $w$ are vectors in $V$ such that $\|v\| = \|w\|
= 1$ and $\lambda(v) = \lambda(w) = \|\lambda\|_*$.  Thus
\begin{equation}
 \lambda(t \, v + (1 - t) \, w) = t \, \lambda(v) + (1 - t) \, \lambda(w)
                                 = \|\lambda\|_*
\end{equation}
when $0 < t < 1$, and hence
\begin{equation}
        \|\lambda\|_* = |\lambda(t \, v + (1 - t) \, w)|
                          \le \|\lambda\|_* \, \|t \, v + (1 - t) \, w\|,
\end{equation}
which implies that
\begin{equation}
        \|t \, v + (1 - t) \, w\| \ge 1.
\end{equation}
By the triangle inequality, $\|t \, v + (1 - t) \, w\| \le 1$ when $0
< t < 1$, and so
\begin{equation}
        \|t \, v + (1 - t) \, w\| = 1.
\end{equation}
If $B_1$ is strictly convex, then it follows that $v = w$ under these
conditions.  Conversely, let us check that this property characterizes
strict convexity of $B_1$.

        Suppose that $v, w \in V$, $\|v\| = \|w\| = 1$, $0 < t < 1$,
and that $a = t \, v + (1 - t) \, w$ satisfies $\|a\| = 1$.  As usual,
there is a bounded linear functional $\lambda$ on $V$ such that
$\lambda(a) = \|\lambda\|_* = 1$, because of the Hahn--Banach theorem.
This implies that $|\lambda(v)|, |\lambda(w)| \le 1$ and
\begin{equation}
 1 = \lambda(t \, v + (1 - t) \, w) = t \, \lambda(v) + (1 - t) \, \lambda(w),
\end{equation}
so that $\lambda(v) = \lambda(w) = 1$.  If we have the uniqueness
property described in the previous paragraph, then we get that $v =
w$, which means that $B_1$ is strictly convex.

        If $V$ is not uniformly convex, then there is an $\epsilon >
0$ and sequences of vectors $\{v_j\}_{j = 1}^\infty$, $\{w_j\}_{j =
1}^\infty$ in $V$ such that $\|v_j\| = \|w_j\| = 1$ and $\|v_j - w_j\|
\ge \epsilon$ for each $j$, and
\begin{equation}
        \lim_{j \to \infty} \biggl\|\frac{v_j + w_j}{2}\biggr\| = 1.
\end{equation}
If $V$ has finite dimension $n$, then there is a one-to-one linear
mapping from ${\bf R}^n$ onto $V$.  This mapping is also a
homeomorphism with respect to the standard topology on ${\bf R}^n$ and
the topology on $V$ determined by the metric associated to the norm.
In particular, closed and bounded subsets of $V$ are compact in this
case.  Thus we may suppose in addition that $\{v_j\}_{j = 1}^\infty$,
$\{w_j\}_{j = 1}^\infty$ converge to some vectors $v, w \in V$,
respectively, by passing to subsequences.  By hypothesis, $\|v\| =
\|w\| = 1$, $\|v - w\| \ge \epsilon > 0$, and $\|(v + w)/2\| = 1$,
which is impossible when $B_1$ is strictly convex.  This shows that
$V$ is uniformly convex when $V$ is finite-dimensional and $B_1$ is
strictly convex.

        Suppose that $B_1$ is strictly convex, and that
\begin{equation}
        \|v + w\| = \|v\| + \|w\|
\end{equation}
for some $v, w \in V$ with $v, w \ne 0$.  If
\begin{equation}
        v' = \frac{v}{\|v\|}, \ w' = \frac{w}{\|w\|}, \hbox{ and } 
         t = \frac{\|v\|}{\|v\| + \|w\|},
\end{equation}
then $1 - t = \|w\|/(\|v\| + \|w\|)$ and
\begin{equation}
        t \, v' + (1 - t) \, w' = \frac{v + w}{\|v\| + \|w\|}.
\end{equation}
This has norm $1$ by hypothesis, so that $v' = w'$ by strict
convexity.  Equivalently, $w = r \, v$, where $r = \|w\| / \|v\|$.

        Let $(X, \mathcal{A})$ be a measurable space, and suppose that
$\mu \in \mathcal{M}(X, V)$.  If $A \subseteq X$ is measurable, then
$\mu(X) = \mu(A) + \mu(X \backslash A)$, which implies that
\begin{eqnarray}
\label{||mu(X)|| le ||mu(A)|| + ||mu(X backslash A)|| le ... = ||mu||(X)}
        \|\mu(X)\| & \le & \|\mu(A)\| + \|\mu(X \backslash A)\| \\
 & \le & \|\mu\|(A) + \|\mu\|(X \backslash A) = \|\mu\|(X). \nonumber
\end{eqnarray}
If $\|\mu(X)\| = \|\mu\|(X)$, then it follows that
\begin{equation}
        \|\mu(X)\| = \|\mu(A)\| + \|\mu(X \backslash A)\|
\end{equation}
and
\begin{equation}
        \|\mu(A)\| = \|\mu\|(A)
\end{equation}
for every measurable set $A \subseteq X$.  If $\|\mu\|(X) > 0$ and
$B_1$ is strictly convex, then one can argue as in the preceding
paragraph to get that
\begin{equation}
        \mu(A) = \mu(X) \, \frac{\|\mu\|(A)}{\|\mu\|(X)}
\end{equation}
for every measurable set $A \subseteq X$.

        Suppose now that $f : [a, b] \to V$ is a path of finite
length, and let $\Lambda_x^y$ be the length of the restriction of $f$
to $[x, y] \subseteq [a, b]$.  Thus
\begin{eqnarray}
\label{||f(b) - f(a)|| le ... = Lambda_a^b}
        \|f(b) - f(a)\| & \le & \|f(x) - f(a)\| + \|f(b) - f(x)\| \\
          & \le & \Lambda_a^x + \Lambda_x^b = \Lambda_a^b \nonumber
\end{eqnarray}
when $a \le x \le b$.  If $\|f(b) - f(a)\| = \Lambda_a^b$, then it follows that
\begin{equation}
        \|f(b) - f(a)\| = \|f(x) - f(a)\| + \|f(b) - f(x)\|
\end{equation}
and
\begin{equation}
        \|f(x) - f(a)\| = \Lambda_a^x
\end{equation}
when $a \le x \le b$.  If $\Lambda_a^b > 0$ and $B_1$ is strictly
convex, then one can argue as before to get that
\begin{equation}
\label{f(x) - f(a) = (f(b) - f(a)) frac{Lambda_a^x}{Lambda_a^b}}
        f(x) - f(a) = (f(b) - f(a)) \, \frac{\Lambda_a^x}{\Lambda_a^b}
\end{equation}
when $a \le x \le b$.

\section[\ Minimizing distances]{Minimizing distances}
\label{minimizing distances}
\setcounter{equation}{0}

        Let $(V, \|v\|)$ be a uniformly convex Banach space, and let
$E$ be a nonempty closed convex set in $V$.  Also let $v \in V$ be
given, and let $\rho$ be the distance from $v$ to $E$,
\begin{equation}
        \rho = \inf \{\|v - w\| : w \in E\}.
\end{equation}
Let $\{w_j\}_{j = 1}^\infty$ be a sequence of elements of $E$ such that
\begin{equation}
\label{lim_{j to infty} ||v - w_j|| = rho}
        \lim_{j \to \infty} \|v - w_j\| = \rho.
\end{equation}
Because $E$ is convex, $(w_j + w_l)/2 \in E$ for every $j, l \ge 1$, and so
\begin{equation}
\label{||v - frac{w_j + w_l}{2}|| ge rho}
        \biggl\|v - \frac{w_j + w_l}{2}\biggr\| \ge \rho.
\end{equation}
Suppose that $v \not\in E$, so that $\rho > 0$, and put
\begin{equation}
        u_j = \frac{v - w_j}{\|v - w_j\|}
\end{equation}
for each $j$.  Thus $\|u_j\| = 1$ for each $j$, and hence $\|(u_j +
u_l)/2\| \le 1$ for every $j, l \ge 1$, by the triangle inequality.
Using (\ref{lim_{j to infty} ||v - w_j|| = rho}) and (\ref{||v -
frac{w_j + w_l}{2}|| ge rho}), it is easy to see that
\begin{equation}
        \lim_{j, l \to \infty} \biggl\|\frac{u_j + u_l}{2}\biggr\| = 1.
\end{equation}
This implies that
\begin{equation}
        \lim_{j, l \to \infty} \|u_j - u_l\| = 0,
\end{equation}
because of uniform convexity.  Using (\ref{lim_{j to infty} ||v -
w_j|| = rho}) again, it is easy to check that
\begin{equation}
 \|w_j - w_l\| = \|(v - w_j) - (v - w_l)\| \to 0 \hbox{ as } j, l \to \infty.
\end{equation}
This shows that $\{w_j\}_{j = 1}^\infty$ is a Cauchy sequence, which
therefore converges to some $w \in V$.  We also have that $w \in E$,
because $E$ is closed.  Of course, $\|v - w\| = \rho$, so that $w$
minimizes the distance to $v$ from elements of $E$.

        Suppose that $w'$ is another element of $E$ such that $\|v -
w'\| = \rho$.  If $0 < t < 1$, then $t \, w + (1 - t) \, w' \in E$,
because $E$ is convex, and so 
\begin{equation}
        \|v - (t \, w + (1 - t) \, w')\| \ge \rho.
\end{equation}
Moreover,
\begin{equation}
 \|v - (t \, w + (1 - t) \, w')\| \le t \, \|v - w\| + (1 - t) \, \|v - w'\|
                                       = \rho,
\end{equation}
which implies that
\begin{equation}
        \|v - (t \, w + (1 - t) \, w')\| = \rho.
\end{equation}
Put $u = \rho^{-1} \, (v - w)$, $u' = \rho^{-1} \, (v - w')$, so that
$\|u\| = \|u'\| = 1$ and
\begin{equation}
\label{||t u + (1 - t) u'|| = 1}
        \|t \, u + (1 - t) \, u'\| = 1
\end{equation}
when $0 < t < 1$.  Strict convexity of the closed unit ball in $V$
implies that $u = u'$, which is the same as saying that $w = w'$.

        Let $\lambda$ be a nonzero bounded linear functional on $V$,
and let $E$ be the set of $w \in V$ such that $\lambda(w) =
\|\lambda\|_*$.  This is a closed affine subspace of $V$, which is
convex in particular.  The distance $\rho$ from $E$ to $0$ is the same
as the infimum of $\|w\|$ over $w \in E$, which is equal to $1$ in
this case, by the definition of the dual norm of $\lambda$.  The
arguments in the previous paragraphs imply that there is a unique $w
\in E$ such that $\|w\| = 1$.  This shows that the supremum is
attained in the definition of the dual norm of a bounded linear
functional on a uniformly convex Banach space.

\section[\ Another approximation argument]{Another approximation argument}
\label{another approximation argument}
\setcounter{equation}{0}

        Let $V_1$ be a real vector space with a norm $\|v\|$.  As
usual, everything in this section can also be applied to complex
vector spaces, since they are real vector spaces too.  Suppose that
$V_1$ is uniformly convex, so that for each $\epsilon > 0$ there is a
$\delta(\epsilon) > 0$ such that for every $v, w \in V_1$ with
$\|v\| = \|w\| = 1$ and
\begin{equation}
        \biggl\|\frac{v + w}{2}\biggr\| > 1 - \delta(\epsilon)
\end{equation}
we have that $\|v - w\| < \epsilon$, as in Section \ref{uniform
convexity}.  Although uniform convexity follows from strict convexity
of the unit ball in finite dimensions, as in Section \ref{strict
convexity}, the estimates in this section will only depend on
$\delta(\epsilon)$, and not on the particular norm $\|v\|$, or the
dimension of $V_1$.  Hence these estimates hold uniformly over all
finite-dimensional subspaces of a uniformly convex Banach space, for
instance.

        Let $(X, \mathcal{A}, \mu)$ be a probability space, and let
$\mathcal{B}$ be a $\sigma$-subalgebra of $\mathcal{A}$.  As in
Section \ref{sigma-subalgebras, vectors}, it is easy to deal with
integrals of $V_1$-valued functions on $X$, by integrating the
components of these functions with respect to a basis for $V_1$.
Similarly, the conditional expectation of a $V_1$-valued function on
$X$ with respect to $\mathcal{B}$ can be defined by taking the
conditional expectation of the components of the function with respect
to a basis.  It is easy to see that this does not depend on the choice
of a basis for $V_1$, using the linearity of integration and
conditional expectation.

        Let $f$ be an integrable $V_1$-valued function on $X$ with
respect to $\mu$, which means that the components of $f$ with respect
to a basis are integrable real-valued functions.  Also let
$f_\mathcal{B} = E(f \mid \mathcal{B})$ be the conditional expectation
of $f$ with respect to $\mathcal{B}$, as usual.  Remember that
\begin{equation}
\label{||f_mathcal{B}|| le E(||f|| mid mathcal{B})}
        \|f_\mathcal{B}\| \le E(\|f\| \mid \mathcal{B})
\end{equation}
almost everywhere on $X$, as in Section \ref{sigma-subalgebras,
vectors}, so that
\begin{equation}
        \int_X \|f_\mathcal{B}\| \, d\mu \le \int_X \|f\| \, d\mu,
\end{equation}
in particular.  Let $\theta$ be a small positive real number, and suppose that
\begin{equation}
\label{int_X ||f|| d mu < int_X ||f_mathcal{B}|| d mu + theta}
        \int_X \|f\| \, d\mu < \int_X \|f_\mathcal{B}\| \, d\mu + \theta.
\end{equation}
This implies that
\begin{equation}
\label{int_X (E(||f|| mid mathcal{B}) - ||f_mathcal{B}||) d mu < theta}
 \int_X (E(\|f\| \mid \mathcal{B}) - \|f_\mathcal{B}\|) \, d\mu < \theta,
\end{equation}
because the integrals of $\|f\|$ and $E(\|f\| \mid \mathcal{B})$ over
$X$ are the same, since $X \in \mathcal{B}$.

        Let $\eta$ be another small positive real number, and put
\begin{eqnarray}
        X_1 & = & \{x \in X : \|f_\mathcal{B}(x)\| > (1 - \eta) \, 
                               E(\|f\| \mid \mathcal{B})(x)\}, \\
        X_2 & = & \{x \in X : \|f_\mathcal{B}(x)\| \le (1 - \eta) \,
                               E(\|f\| \mid \mathcal{B})(x)\}.
\end{eqnarray}
Thus $X_1, X_2 \in \mathcal{B}$, because $\|f_\mathcal{B}\|$, $E(\|f\|
\mid \mathcal{B})$ are measurable with respect to $\mathcal{B}$.  Note
that
\begin{eqnarray}
\label{eta int_{X_2} ||f|| d mu le ... < theta}
 \eta \int_{X_2} \|f\| \, d\mu 
            & = & \int_{X_2} \eta \, E(\|f\| \mid \mathcal{B}) \, d\mu  \\
           & \le & \int_X (E(\|f\| \mid \mathcal{B}) - \|f_\mathcal{B}\|)
                                              \, d\mu < \theta, \nonumber
\end{eqnarray}
where we use the fact that $X_2 \in \mathcal{B}$ in the first step,
and (\ref{||f_mathcal{B}|| le E(||f|| mid mathcal{B})}) and the
definition of $X_2$ in the second step.

        In order to see what happens on $X_1$, it will be convenient
to use linear functionals on $V_1$.  Of course, every linear
functional on $V_1$ is bounded, because $V_1$ has finite dimension,
and the dual $V_1^*$ of $V_1$ has finite dimension equal to the
dimension of $V_1$.  In particular, there is a sequence of linear
functionals $\{\lambda_j\}_{j = 1}^\infty$ on $V_1$ such that
$\|\lambda_j\|_* = 1$ for each $j$ and the $\lambda_j$'s are dense in
the set of $\lambda \in V_1^*$ with $\|\lambda\|_* = 1$.  As usual,
for each $v \in V_1$ there is a $\lambda \in V_1^*$ such that
$\|\lambda\|_* = 1$ and $\lambda(v) = \|v\|$, because of the
Hahn--Banach theorem.  This implies that
\begin{equation}
\label{||v|| = sup_{j ge 1} lambda_j(v)}
        \|v\| = \sup_{j \ge 1} \lambda_j(v)
\end{equation}
for each $v \in V_1$, by approximating $\lambda$ by $\lambda_j$'s,
and using the fact that $\|\lambda_j\|_* = 1$ for each $j$.

        Put
\begin{equation}
        A_j = \{x \in X : \lambda_j(f_\mathcal{B}(x)) > (1 - \eta) \,
                                           E(\|f\| \mid \mathcal{B})(x)\}
\end{equation}
for each $j \ge 1$, so that $A_j \subseteq X_1$ and $A_j \in
\mathcal{B}$ for each $j$, and
\begin{equation}
        \bigcup_{j = 1}^\infty A_j = X_1,
\end{equation}
by (\ref{||v|| = sup_{j ge 1} lambda_j(v)}).  It is better to have
disjoint sets, and so we let $B_1 = A_1$ and $B_n = A_n \backslash
\Big(\bigcup_{j = 1}^{n - 1} A_j\Big)$ when $n \ge 2$.  Thus $B_n
\subseteq A_n \subseteq X_1$ and $B_n \in \mathcal{B}$ for each $n$,
$B_l \cap B_n = \emptyset$ when $l < n$, and
\begin{equation}
        \bigcup_{n = 1}^\infty B_n = \bigcup_{j = 1}^\infty A_j = X_1,
\end{equation}
as before.  Note that $\lambda \circ f_\mathcal{B} = E(\lambda \circ f
\mid \mathcal{B})$ for each linear functional $\lambda$ on $V_1$.
This implies that
\begin{equation}
        \int_{B_n} \lambda_n \circ f_\mathcal{B} \, d\mu
         = \int_{B_n} \lambda_n \circ f \, d\mu,
\end{equation}
since $B_n \in \mathcal{B}$, while
\begin{equation}
 \int_{B_n} E(\|f\| \mid \mathcal{B}) \, d\mu = \int_{B_n} \|f\| \, d\mu.
\end{equation}
Because $B_n \subseteq A_n$,
\begin{equation}
        \int_{B_n} \lambda_n \circ f_\mathcal{B} \, d\mu
          > (1 - \eta) \int_{B_n} E(\|f\| \mid \mathcal{B}) \, d\mu
\end{equation}
when $\mu(B_n) > 0$, and hence
\begin{equation}
 \int_{B_n} \lambda_n \circ f \, d\mu > (1 - \eta) \int_{B_n} \|f\| \, d\mu.
\end{equation}
Equivalently,
\begin{equation}
\label{int_{B_n} (||f|| - lambda_n circ f) d mu < eta int_{B_n} ||f|| d mu}
 \int_{B_n} (\|f\| - \lambda_n \circ f) \, d\mu < \eta \int_{B_n} \|f\| \, d\mu
\end{equation}
when $\mu(B_n) > 0$, where the integrand on the left is nonnegative,
since $\|\lambda_n\|_* = 1$.

        Let $\epsilon > 0$ be given, and put $\delta =
\delta(\epsilon)$.  Also put
\begin{eqnarray}
B_{n, 1} & = & \{x \in B_n : \lambda_n(f(x)) > (1 - \delta) \, \|f(x)\|\}, \\
B_{n, 2} & = & \{x \in B_n : \lambda_n(f(x)) \le (1 - \delta) \, \|f(x)\| \}.
\end{eqnarray}
Thus
\begin{equation}
\label{delta int_{B_{n, 2}} ||f|| d mu le ... < eta int_{B_n} ||f|| d mu}
        \delta \int_{B_{n, 2}} \|f\| \, d\mu
                \le \int_{B_{n, 2}} (\|f\| - \lambda_n \circ f) \, d\mu 
                  < \eta \int_{B_n} \|f\| \, d\mu
\end{equation}
when $\mu(B_n) > 0$, by (\ref{int_{B_n} (||f|| - lambda_n circ f) d mu
< eta int_{B_n} ||f|| d mu}).  As before, we shall be interested in
$\eta$'s that are small compared to $\delta$, so that the integral
of $\|f\|$ over $B_{n, 2}$ is relatively small.

        If $x \in X_1$, then $f_\mathcal{B}(x) \ne 0$, and we put
$a(x) = f_\mathcal{B}(x) / \|f_\mathcal{B}(x)\|$.  Otherwise, if $x
\in X_2$, then we put $a(x) = 0$.  If $x \in B_n \subseteq A_n
\subseteq X_1$, then
\begin{equation}
\label{lambda_n(a(x)) > 1 - eta}
        \lambda_n(a(x)) > 1 - \eta,
\end{equation}
using also (\ref{||f_mathcal{B}|| le E(||f|| mid mathcal{B})}).  If $x
\in B_{n, 1}$, then $f(x) \ne 0$, and we put $b(x) = f(x) / \|f(x)\|$.
Note that
\begin{equation}
        \lambda_n(b(x)) > 1 - \delta,
\end{equation}
by definition of $B_{n, 1}$.  Thus $\|a(x)\| = \|b(x)\| = 1$ and
\begin{equation}
        \biggl\|\frac{a(x) + b(x)}{2}\biggr\|
         \ge \lambda_n\Big(\frac{a(x) + b(x)}{2}\Big) 
           > 1 - \frac{\delta + \eta}{2} \ge 1 - \delta
\end{equation}
when $x \in B_{n, 1}$ and $\eta \le \delta$.  This implies that
$\|a(x) - b(x)\| < \epsilon$, because of uniform convexity.
Equivalently,
\begin{equation}
        \bigl\|f(x) - a(x) \, \|f(x)\| \bigr\| < \epsilon \, \|f(x)\|
\end{equation}
when $x \in B_{n, 1}$ and $\eta \le \delta$.

        It follows that
\begin{equation}
 \int_{B_n} \bigl\|f(x) - a(x) \, \|f(x)\| \bigr\| \, d\mu \le
\int_{B_{n, 1}} \epsilon \, \|f\| \, d\mu + \int_{B_{n, 2}} 2 \, \|f\| \, d\mu
\end{equation}
when $\eta \le \delta(\epsilon)$, and hence
\begin{equation}
\label{int_{B_n} ||f(x) - a(x) ||f(x)|| || d mu le ...}
 \int_{B_n} \bigl\|f(x) - a(x) \, \|f(x)\| \bigr\| \, d\mu \le
 (\epsilon + 2 \, \delta(\epsilon)^{-1} \, \eta) \int_{B_n} \|f\| \, d\mu,
\end{equation}
because of (\ref{delta int_{B_{n, 2}} ||f|| d mu le ... < eta
int_{B_n} ||f|| d mu}).  This also holds trivially when $\eta >
\delta(\epsilon)$, since the coefficient on the right would be greater
than $2$.  Summing over $n$, we get that
\begin{equation}
 \int_{X_1} \bigl\|f(x) - a(x) \, \|f(x)\| \bigr\| \, d\mu \le
 (\epsilon + 2 \, \delta(\epsilon)^{-1} \, \eta) \int_{X_1} \|f\| \, d\mu.
\end{equation}
Combining this with (\ref{eta int_{X_2} ||f|| d mu le ... < theta}), we obtain
\begin{equation}
\label{int_X ||f(x) - a(x) ||f(x)|| || d mu < ...}
 \int_X \bigl\|f(x) - a(x) \, \|f(x)\| \bigr\| \, d\mu <
 (\epsilon + 2 \, \delta(\epsilon)^{-1} \, \eta) \int_X \|f\| \, d\mu
  + \eta^{-1} \, \theta.
\end{equation}

        Alternatively, one might prefer to take $a(x) =
f_\mathcal{B}(x) / \|f_\mathcal{B}(x)\|$ for every $x$ in $X$ such
that $f_\mathcal{B}(x) \ne 0$, even when $x \in X_2$.  This would
ensure that $a(x)$ does not depend on $f$ even indirectly, through the
definition of $X_2$.  In this case, we would get that
\begin{equation}
\label{int_X ||f(x) - a(x) ||f(x)|| || d mu < ..., 2}
 \quad  \int_X \bigl\|f(x) - a(x) \, \|f(x)\| \bigr\| \, d\mu <
 (\epsilon + 2 \, \delta(\epsilon)^{-1} \, \eta) \int_X \|f\| \, d\mu
  + 2 \, \eta^{-1} \, \theta,
\end{equation}
which is to say that we would multiply $\eta^{-1} \, \theta$ by $2$ in
the previous estimate.  In both situations, $a(x)$ is measurable with
respect to $\mathcal{B}$, because $f_\mathcal{B}$ is measurable with
respect to $\mathcal{B}$ and $X_2 \in \mathcal{B}$.

\section[ Examples in $\ell^p$]{Examples in $\ell^p$}
\label{examples in ell^p}
\setcounter{equation}{0}

        Let $a_1, a_2, \ldots$ be a sequence of real or complex
numbers, and consider
\begin{equation}
        f_n(x) = \sum_{j = 1}^n a_j \, r_j(x) \, \delta_j.
\end{equation}
Here $r_1(x), r_2(x), \ldots$ are the Rademacher functions, and
$\delta_j = \{\delta_{j, l}\}_{l = 1}^\infty$ is the sequence defined
by $\delta_{j, l} = 1$ when $j = l$ and $\delta_{j, l} = 0$ when $j
\ne l$.  Thus $\{f_n\}_{n = 1}^\infty$ is a martingale on the dyadic
unit interval with respect to the usual filtration associated to
dyadic subintervals, and with values in the vector space of sequences
of real or complex numbers, as appropriate.  In particular,
$\{f_n\}_{n = 1}^\infty$ is a martingale with values in $\ell^p$ for
each $p$, $1 \le p \le \infty$.  Note that the $\ell^p$ norm of
$f_n(x)$ is equal to the $\ell^p$ norm of the finite sequence $a_1,
\ldots, a_n$ for each $x$ and $n$.  Hence the $L^1$ norm of
$\|f_n(x)\|_{\ell^p}$ is equal to the $\ell^p$ norm of $a_1, \ldots,
a_n$ for each $n$.  It follows that the $L^1$ norm of
$\|f_n(x)\|_{\ell^p}$ is uniformly bounded over $n$ if and only if
$\{a_j\}_{j = 1}^\infty$ is in $\ell^p$.  If $\{a_j\}_{j = 1}^\infty
\in \ell^p$ and $p < \infty$, then it is easy to see that $f_n(x)$
converges in $\ell^p$ as $n \to \infty$ for each $x$.  Similarly, if
$\{a_j\}_{j = 1}^\infty$ converges to $0$, then $f_n(x)$ converges in
$c_0$ equipped with the $\ell^\infty$ norm as $n \to \infty$ for each
$x$.  If $\{a_j\}_{j = 1}^\infty$ is bounded, then $f_n(x)$ is
uniformly bounded in $\ell^\infty$, but it does not converge in the
$\ell^\infty$ norm as $n \to \infty$ for any $x$ unless $\{a_j\}_{j =
1}^\infty$ converges to $0$.

\section[\ Uniform convergence]{Uniform convergence}
\label{uniform convergence}
\setcounter{equation}{0}

        Let $(V, \|v\|)$ be a real or complex Banach space, and let
$v_1, v_2, \ldots$ be a sequence of elements of $V$.  As in Section
\ref{rademacher functions}, let $X$ be the set of sequences $x =
\{x_j\}_{j = 1}^\infty$ with $x_j = 1$ or $-1$ for each $j$, which is
the same as the Cartesian product of a sequence of copies of $\{1,
-1\}$.  Consider
\begin{equation}
        f_n(x) = \sum_{j = 1}^n x_j \, v_j
\end{equation}
for each positive integer $n$ and $x \in X$.  This is basically the
same as the sequence of functions considered in the previous section
when $V = \ell^p$ and $v_j = a_j \, \delta_j$, since $r_j(x) = x_j$ is
another version of the Rademacher functions.  Let us check that
$\{f_n\}_{n = 1}^\infty$ converges uniformly on $X$ when $\sum_{j \in
{\bf Z}_+} v_j$ converges in the generalized sense, as in Section
\ref{generalized convergence, 2}.  In particular, $\{f_n\}_{n =
1}^\infty$ converges uniformly on $X$ when $\sum_{j = 1}^\infty v_j$
converges absolutely.  In this case, it is very easy to show directly
that $\{f_n\}_{n = 1}^\infty$ converges uniformly, by the same
argument as in Weierstrass' $M$-test.

        Suppose that $\sum_{j \in {\bf Z}_+} v_j$ converges in the
generalized sense, which implies that it satisfies the generalized
Cauchy criterion, as in Section \ref{generalized convergence, 2}.
This means that for each $\epsilon > 0$ there is a finite set
$A_\epsilon \subseteq {\bf Z}_+$ such that
\begin{equation}
        \biggl\|\sum_{j \in B} v_j\biggr\| < \epsilon
\end{equation}
for every finite set $B \subseteq {\bf Z}_+$ with $A_\epsilon \cap B =
\emptyset$.  Let $L_\epsilon$ be the maximum of the elements of
$A_\epsilon$, with $L_\epsilon = 0$ when $A_\epsilon = \emptyset$.  
If $n > l$, then
\begin{equation}
 f_n(x) - f_l(x) = \sum_{j = l + 1}^n x_j \, v_j
                 = \sum_{j \in B_{l, n}^+} v_j - \sum_{j \in B_{l, n}^-} v_j,
\end{equation}
where $B_{l, n}^+$, $B_{l, n}^-$ are the sets of positive integers $j$
such that $l < j \le n$ and $x_j = 1$ or $-1$, respectively. 
If $l \ge L_\epsilon$, then $B_{l, n}^+ \cap A_\epsilon = B_{l, n}^-
\cap A_\epsilon = \emptyset$, and so
\begin{equation}
 \|f_n(x) - f_l(x)\| \le \biggl\|\sum_{j \in B_{l, n}^+} v_j\biggr\|
                       + \biggl\|\sum_{j \in B_{l, n}^-} v_j\biggr\| 
                      < \epsilon + \epsilon = 2 \, \epsilon.
\end{equation}
This shows that $\{f_n\}_{n = 1}^\infty$ is a Cauchy sequence with
respect to the supremum norm on the space of $V$-valued functions on
$X$.  It follows that $\{f_n\}_{n = 1}^\infty$ converges uniformly on
$X$, because $V$ is complete.  As usual, one can observe first that
$\{f_n(x)\}_{n = 1}^\infty$ is a Cauchy sequence in $V$ for each $x
\in X$, which converges because of completeness, and then check that
$\{f_n\}_{n = 1}^\infty$ converges uniformly on $X$ to the pointwise
limit, because of the uniform version of the Cauchy condition.

        Conversely, suppose that $\{f_n\}_{n = 1}^\infty$ converges
uniformly on $X$, and hence satisfies the uniform version of the
Cauchy condition.  This means that for each $\epsilon > 0$ there is an
$N_\epsilon \ge 0$ such that
\begin{equation}
        \|f_n(x) - f_l(x)\| < \epsilon
\end{equation}
for every $n > l \ge N_\epsilon$ and $x \in X$, or equivalently
\begin{equation}
\label{||sum_{j = l + 1}^n x_j v_j|| < epsilon}
        \biggl\|\sum_{j = l + 1}^n x_j \, v_j\biggr\| < \epsilon
\end{equation}
for every $n > l \ge N_\epsilon$ and $x \in X$.  Let $B \subseteq {\bf
Z}_+$ be a nonempty finite set whose minimal element is greater than
$N_\epsilon$.  If $y, z \in X$ are defined by $y_j = 1$ for every $j$,
$z_j = 1$ when $j \in B$, and $z_j = -1$ otherwise, then
\begin{equation}
\sum_{j = N_\epsilon + 1}^n y_j \, v_j + \sum_{j = N_\epsilon + 1}^n z_j \, v_j
                                                        = 2 \sum_{j \in B} v_j
\end{equation}
when the maximal element of $B$ is less than or equal to $n$.  Hence
\begin{equation}
        2 \, \biggl\|\sum_{j \in B} v_j\biggr\| 
                 \le \biggl\|\sum_{j = N_\epsilon + 1}^n y_j \, v_j\biggr\|
                      + \biggl\|\sum_{j = N_\epsilon + 1}^n z_j \, v_j\biggr\|
                     < \epsilon + \epsilon = 2 \, \epsilon,
\end{equation}
by (\ref{||sum_{j = l + 1}^n x_j v_j|| < epsilon}).  This is the same
as saying that $\Big\|\sum_{j \in B} v_j\Big\| < \epsilon$ when $B
\subseteq {\bf Z}_+$ is a finite set disjoint from $\{1, \ldots,
N_\epsilon\}$, which implies that $\sum_{j \in {\bf Z}_+} v_j$
satisfies the generalized Cauchy criterion.  Thus $\sum_{j \in {\bf
Z}_+} v_j$ converges in the generalized sense, because $V$ is
complete.

        Actually, the same conclusion holds when $\{f_n(x)\}_{n =
1}^\infty$ converges in $V$ for every $x \in X$, which is the same as
saying that $\sum_{j = 1}^\infty x_j \, v_j$ converges for every $x
\in X$.  To see this, suppose for the sake of a contradiction that
$\sum_{j \in {\bf Z}_+} v_j$ does not satisfy the generalized Cauchy
condition.  This means that for each $\epsilon > 0$ and finite set $A
\subseteq {\bf Z}_+$ there is a finite set $B \subseteq {\bf Z}_+$
such that
\begin{equation}
\label{||sum_{j in B} v_j|| ge epsilon}
        \biggl\|\sum_{j \in B} v_j\biggr\| \ge \epsilon.
\end{equation}
By applying this repeatedly, we can get an infinite sequence $B_1,
B_2, \ldots$ of finite subsets of ${\bf Z}_+$ such that the maximal
element of $B_l$ is strictly less than the minimal element of $B_{l +
1}$ for each $l$, and (\ref{||sum_{j in B} v_j|| ge epsilon}) holds
with $B = B_l$ for each $l$.  Let $y, z \in X$ be defined by $y_j = 1$
for each $j$, $z_j = 1$ when $j \in B_l$ for some $l \ge 1$, and
$z_j = -1$ otherwise.  If $b_n$ is the maximal element of $B_n$, then
\begin{equation}
        \sum_{j = 1}^{b_n} y_j \, v_j + \sum_{j = 1}^{b_n} z_j \, v_j
                     = 2 \sum_{l = 1}^n \Big(\sum_{j \in B_l} v_j\Big).
\end{equation}
Thus the convergence of $\sum_{j = 1}^\infty y_j \, v_j$ and $\sum_{j
= 1}^\infty z_j \, v_j$ imply the convergence of
\begin{equation}
        \sum_{l = 1}^\infty \Big(\sum_{j \in B_l} v_j\Big).
\end{equation}
This implies in turn that
\begin{equation}
        \lim_{l \to \infty} \sum_{j \in B_l} v_j = 0,
\end{equation}
a contradiction.  This shows that $\sum_{j \in {\bf Z}_+} v_j$
satisfies the generalized Cauchy condition, and hence converges in the
generalized sense, because $V$ is complete.  Therefore $\sum_{j \in
{\bf Z}_+} v_j$ converges in the generalized sense if and only if
$\sum_{j = 1}^\infty x_j \, v_j$ converges for every $x \in X$, in
which case the partial sums $f_n$ converge uniformly on $X$.

\section[\ Bounded sums]{Bounded sums}
\label{bounded sums}
\setcounter{equation}{0}

        Let $V$ be a real or complex vector space with a norm $\|v\|$,
and let $\{1, -1\}^n$ be the Cartesian product of $n$ copies of $\{1,
-1\}$, consisting of all sequences $\epsilon = \{\epsilon_j\}_{j =
1}^n$ of length $n$ with $\epsilon_j = 1$ or $-1$ for each $j$.  Also
let $Z(V)$ be the collection of sequences $v_1, v_2, \ldots$ of
vectors in $V$ for which the sums $\sum_{j = 1}^n \epsilon_j \, v_j$
are uniformly bounded in $V$ over $\epsilon \in \{1, -1\}^n$ and all
positive integers $n$.  This is a vector space over the real or
complex numbers, as appropriate, with respect to termwise addition and
scalar multiplication.  If $\{v_j\}_{j = 1}^\infty \in Z(V)$, then put
\begin{equation}
 \|\{v_j\}_{j = 1}^\infty\|_{Z(V)} = \sup\Big\{\biggl\|\sum_{j = 1}^n
 \epsilon_j \, v_j\biggr\|: \epsilon \in \{1, -1\}^n, \, n \in {\bf Z}_+\Big\}.
\end{equation}
Note that $Z(V)$ is a linear subspace of the space $X(V)$ of sequences
$\{v_j\}_{j = 1}^\infty$ of vectors in $V$ with bounded partial sums
$\sum_{j = 1}^n v_j$, discussed in Section \ref{bounded partial sums},
since we can take $\epsilon_j = 1$ for each $j$.  Similarly,
\begin{equation}
        \|\{v_j\}_{j = 1}^\infty\|_{X(V)} \le \|\{v_j\}_{j = 1}^\infty\|_{Z(V)}
\end{equation}
for each $\{v_j\}_{j = 1}^\infty \in Z(V)$.  It is easy to see that
$\|\{v_j\}_{j = 1}^\infty\|_{Z(V)}$ is a norm on $Z(V)$, and in
particular that $v_j = 0$ for every $j$ when $\|\{v_j\}_{j =
1}^\infty\|_{Z(V)} = 0$.

        Let $\{v_j\}_{j = 1}^\infty$ be a sequence of vectors in $V$,
let $B$ be a finite nonempty set of positive integers, and let $n$ be
the maximal element of $B$.  If $\alpha, \beta \in \{1, -1\}^n$ are
defined by $\alpha_j = 1$ for each $j$, $\beta_j = 1$ when $j \in B$,
and $\beta_j = -1$ otherwise, then
\begin{equation}
 \sum_{j = 1}^n \alpha_j \, v_j + \sum_{j = 1}^n \beta_j \, v_j
        = 2 \sum_{j \in B} v_j,
\end{equation}
and hence
\begin{equation}
        2 \, \biggl\|\sum_{j \in B} v_j\biggr\| 
           \le \biggl\|\sum_{j = 1}^n \alpha_j \, v_j\biggr\|
                + \biggl\|\sum_{j = 1}^n \beta_j \, v_j\biggr\|.
\end{equation}
If $\{v_j\}_{j = 1}^\infty \in Z(V)$, then we get that
\begin{equation}
 \biggl\|\sum_{j \in B} v_j\biggr\| \le \|\{v_j\}_{j = 1}^\infty\|_{Z(V)},
\end{equation}
which implies that $\{v_j\}_{j = 1}^\infty$ is in the space $Y({\bf
Z}_+, V)$ discussed in Section \ref{bounded finite subsums}, and that
\begin{equation}
\label{||{v_j}_{j = 1}^infty||_{Y(Z_+, V)} le ||{v_j}_{j = 1}^infty||_{Z(V)}}
        \|\{v_j\}_{j = 1}^\infty\|_{Y({\bf Z}_+, V)}
         \le \|\{v_j\}_{j = 1}^\infty\|_{Z(V)}.
\end{equation}
Conversely, if $\{v_j\}_{j = 1}^\infty \in Y({\bf Z}_+, V)$, $n \in
{\bf Z}_+$, and $\epsilon \in \{1, -1\}^n$, then
\begin{equation}
        \sum_{j = 1}^n \epsilon_j \, v_j
          = \sum_{1 \le j \le n \atop \epsilon_j = 1} v_j
              - \sum_{1 \le j \le n \atop \epsilon_j = -1} v_j,
\end{equation}
which implies that
\begin{equation}
        \biggl\|\sum_{j = 1}^n \epsilon_j \, v_j\biggr\|
         \le \biggl\|\sum_{1 \le j \le n \atop \epsilon_j = 1} v_j\biggr\|
             + \biggl\|\sum_{1 \le j \le n \atop \epsilon_j = -1} v_j\biggr\|
          \le 2 \, \|\{v_j\}_{j = 1}^\infty\|_{Y({\bf Z}_+, V)}.
\end{equation}
Thus $\{v_j\}_{j = 1}^\infty \in Z(V)$ and
\begin{equation}
        \|\{v_j\}_{j = 1}^\infty\|_{Z(V)}
          \le 2 \, \|\{v_j\}_{j = 1}^\infty\|_{Y({\bf Z}_+, V)},
\end{equation}
which shows that $Y({\bf Z}_+, V) = Z(V)$, and that the corresponding
norms are equivalent.

        Let $Z_0(V)$ be the closure in $Z(V)$ of the collection of
sequences $\{v_j\}_{j = 1}^\infty$ of vectors in $V$ with $v_j = 0$
for all but finitely many $j$.  This is the same as the closure of
this set in $Y({\bf Z}_+, V)$, which is also the same as the
collection $Y_0({\bf Z}_+, V)$ of sequences $\{v_j\}_{j = 1}^\infty$
of elements of $V$ such that $\sum_{j \in {\bf Z}_+} v_j$ satisfies
the generalized Cauchy criterion.  If $V$ is complete, then this is
the same as the collection of sequences $\{v_j\}_{j = 1}^\infty$ of
vectors in $V$ such that $\sum_{j \in {\bf Z}_+} v_j$ converges in the
generalized sense, as usual.  This characterization of $Z_0(V)$ is
basically equivalent to the discussion in the previous section.

        Suppose that $\{v_j\}_{j = 1}^\infty$ is a sequence of vectors
in $V$ that is not in $Z(V)$.  Thus for each $N \ge 1$ there is an $n
\in {\bf Z}_+$ and an $\epsilon \in \{1, -1\}^n$ such that
\begin{equation}
        \biggl\|\sum_{j = 1}^n \epsilon_j \, v_j\biggr\| \ge N.
\end{equation}
Equivalently, for each $l, L \ge 1$ there is an $n \ge l$ and
$\epsilon_l, \ldots, \epsilon_n \in \{1, -1\}$ such that
\begin{equation}
\label{||sum_{j = l}^n epsilon_j v_j|| ge L + sum_{j = 1}^{l - 1} ||v_j||}
        \biggl\|\sum_{j = l}^n \epsilon_j \, v_j\biggr\|
         \ge L + \sum_{j = 1}^{l - 1} \|v_j\|.
\end{equation}
This follows from the previous statement by taking $N = L + 2 \sum_{j
= 1}^{l - 1} \|v_j\|$, and using the triangle inequality to get that
\begin{equation}
        \biggl\|\sum_{j = 1}^n \epsilon_j \, v_j\biggr\|
         \le \biggl\|\sum_{j = l}^n \epsilon_j \, v_j\biggr\|
               + \sum_{j = 1}^{l - 1} \|v_j\|.
\end{equation}
Applying (\ref{||sum_{j = l}^n epsilon_j v_j|| ge L + sum_{j = 1}^{l -
1} ||v_j||}) repeatedly, we get a strictly increasing sequence $n_1,
n_2, \ldots$ of positive integers and a sequence $\epsilon_1,
\epsilon_2, \ldots$ with $\epsilon_j \in \{1, -1\}$ for each $j$ such
that
\begin{equation}
        \biggl\|\sum_{j = 1}^{n_1} \epsilon_j \, v_j\biggr\| \ge 1
\end{equation}
and
\begin{equation}
        \biggl\|\sum_{j = n_k + 1}^{n_{k + 1}} \epsilon_j \, v_j\biggr\|
           \ge k + 1 + \sum_{j = 1}^{n_k} \, \|v_j\|
\end{equation}
for each $k \ge 1$.  Using the triangle inequality again, we get that
\begin{equation}
        \biggl\|\sum_{j = n_k + 1}^{n_{k + 1}} \epsilon_j \, v_j\biggr\|
         \le \biggl\|\sum_{j = 1}^{n_{k + 1}} \epsilon_j \, v_j\biggr\|
                + \sum_{j = 1}^{n_k} \|v_j\|
\end{equation}
for each $k$.  Hence
\begin{equation}
        \biggl\|\sum_{j = 1}^{n_k} \epsilon_j \, v_j\biggr\| \ge k
\end{equation}
for each $k$, so that the partial sums $\sum_{j = 1}^n \epsilon_j \,
v_j$ are not uniformly bounded over $n \in {\bf Z}_+$ even for this
single sequence $\epsilon = \{\epsilon_j\}_{j = 1}^\infty$.  
If $\{v_j\}_{j = 1}^\infty$ is a sequence of vectors in $V$ for which
the partial sums $\sum_{j = 1}^\infty \epsilon_j \, v_j$ are uniformly
bounded over $n \in {\bf Z}_+$ for each sequence $\epsilon =
\{\epsilon_j\}_{j = 1}^\infty$ of elements of $\{1, -1\}$, then
it follows that $\{v_j\}_{j = 1}^\infty \in Z(V)$.

\section[\ Bounded coefficients]{Bounded coefficients}
\label{bounded coefficients}
\setcounter{equation}{0}

        Let $E$ be a nonempty set, and let $V$ be a real or complex
vector space with a norm $\|v\|$.  Also let $f \in Y(E, V)$ be given,
as in Section \ref{bounded finite subsums}.  If $A \subseteq E$, then
we let ${\bf 1}_A(x)$ be the indicator function associated to $A$ on
$E$, equal to $1$ when $x \in A$ and to $0$ when $x \in E \backslash
A$.  Thus
\begin{equation}
        \sum_{x \in B} {\bf 1}_A(x) \, f(x) = \sum_{x \in A \cap B} f(x)
\end{equation}
for every finite set $B \subseteq E$, which implies that ${\bf 1}_A \,
f \in Y(E, V)$, and that
\begin{equation}
        \|{\bf 1}_A \, f\|_{Y(E, V)} \le \|f\|_{Y(E, V)}.
\end{equation}

        Now let $a$ be a real-valued function on $E$ such that $0 \le
a(x) \le 1$ for every $x \in E$.  Let $A_1$ be the set of $x \in E$
such that $a(x) \ge 1/2$, and put
\begin{equation}
        a_1(x) = a(x) - \frac{1}{2} \, {\bf 1}_{A_1}(x).
\end{equation}
Thus $0 \le a_1(x) \le 1/2$ for every $x \in E$, and we can repeat the
process by taking $A_2$ to be the set of $x \in E$ such that $a_1(x)
\ge 1/4$.  Continuing in this manner, we get a sequence of subsets
$A_1, A_2, \ldots$ of $E$ such that
\begin{equation}
\label{a(x) = sum_{j = 1}^infty 2^{-j} {bf 1}_{A_j}(x)}
        a(x) = \sum_{j = 1}^\infty 2^{-j} \, {\bf 1}_{A_j}(x)
\end{equation}
for each $x \in E$.  If $f \in Y(E, V)$, as before, then it follows
that $a \, f \in Y(E, V)$ too, and that
\begin{equation}
        \|a \, f\|_{Y(E, V)} \le \|f\|_{Y(E, V)}.
\end{equation}

        If $a$ is a bounded nonnegative real-valued function on $E$,
then we get that $a \, f \in Y(E, V)$, with
\begin{equation}
        \|a \, f\|_{Y(E, V)} \le \|a\|_\infty \, \|f\|_{Y(E, V)}.
\end{equation}
If $a$ is any bounded real-valued function on $E$, then we can apply
the previous remarks to the positive and negative parts of $a$, to get
that $a \, f \in Y(E, V)$ and
\begin{equation}
\label{||a f||_{Y(E, V)} le 2 ||a||_infty ||f||_{Y(E, V)}}
        \|a \, f\|_{Y(E, V)} \le 2 \, \|a\|_\infty \, \|f\|_{Y(E, V)}.
\end{equation}
If $V$ is complex and $a$ is a bounded complex-valued function on $E$,
then we can apply this to the real and imaginary parts of $a$, to get
that $a \, f \in Y(E, V)$ and
\begin{equation}
\label{||a f||_{Y(E, V)} le 4 ||a||_infty ||f||_{Y(E, V)}}
        \|a \, f\|_{Y(E, V)} \le 4 \, \|a\|_\infty \, \|f\|_{Y(E, V)}.
\end{equation}
In particular, multiplication by $a$ defines a bounded linear operator
on $Y(E, V)$ in each case.

        Of course, if $f(x) \ne 0$ for only finitely many $x \in E$,
then $a \, f$ has the same property.  This implies that $a \, f \in
Y_0(E, V)$ when $f \in Y_0(E, V)$ and $a$ is bounded, because $Y_0(E,
V)$ is the closure in $Y(E, V)$ of the linear subspace of functions on
$E$ with finite support.  Equivalently, $\sum_{x \in E} a(x) \, f(x)$
satisfies the generalized Cauchy condition when $\sum_{x \in E} f(x)$
satisfies the generalized Cauchy condition and $a$ is a bounded.  If
$V$ is complete, then it follows that $\sum_{x \in E} a(x) \, f(x)$
converges in the generalized sense when $\sum_{x \in E} f(x)$
converges in the generalized sense and $a$ is bounded.

\section[\ Another norm]{Another norm}
\label{another norm}
\setcounter{equation}{0}

        Let $E$ be a nonempty set, and let $V$ be a real or complex
vector space with a norm $\|v\|$.  Suppose that $f(x)$ is a $V$-valued
function on $E$, and consider sums of the form
\begin{equation}
        \sum_{x \in B} \beta(x) \, f(x),
\end{equation}
where $B \subseteq E$ is a nonempty finite set, and $\beta$ is a
function on $B$ with values in $\{1, -1\}$.  Of course, this is the same as
\begin{equation}
        \sum_{x \in B_+} f(x) - \sum_{x \in B_-} f(x),
\end{equation}
where $B_\pm = \{x \in B : \beta(x) = \pm 1\}$.  If $Z(E, V)$ is the
space of $V$-valued functions on $E$ for which these sums have bounded
norm, then it is easy to see that $Z(E, V)$ is the same as the space
$Y(E, V)$ discussed in Section \ref{bounded finite subsums}.  More
precisely, $Z(E, V) \subseteq Y(E, V)$ because one can take $\beta(x)
= 1$ for each $x \in B$, while $Y(E, V) \subseteq Z(E, V)$ by the
triangle inequality.  If $f \in Y(E, V) = Z(E, V)$, then put
\begin{equation}
        \|f\|_{Z(E, V)}
         = \sup_{B, \beta} \, \biggl\|\sum_{x \in B} \beta(x) \, f(x)\biggr\|,
\end{equation}
where the supremum is taken over all nonempty finite sets $B \subseteq E$
and functions $\beta : B \to \{-1, 1\}$.  Note that
\begin{equation}
        \|f\|_{Y(E, V)} \le \|f\|_{Z(E, V)} \le 2 \, \|f\|_{Y(E, V)},
\end{equation}
for the same reasons that $Y(E, V) = Z(E, V)$.

        If $E = {\bf Z}_+$, then the $Z(E, V)$ norm reduces to the
$Z(V)$ norm described in Section \ref{bounded sums}, where we identify
$V$-valued functions on ${\bf Z}_+$ with sequences whose terms are in
$V$.  Clearly
\begin{equation}
        \|f\|_{Z(V)} \le \|f\|_{Z({\bf Z}_+, V)}
\end{equation}
for each $f \in Z(V) = Y({\bf Z}_+, V)$, because the $Z(V)$
corresponds to taking $B$ to be of the form $\{1, \ldots, n\}$, $n \in
{\bf Z}_+$, in the previous paragraph.  Conversely, if $B$ is any
nonempty finite set of positive integers, and $\beta : B \to \{1,
-1\}$, then we can take $n$ to be the maximal element of $B$, and
put $\epsilon_j = \epsilon'_j = \beta(j)$ when $j \in B$, and $\epsilon_j
= 1$ and $\epsilon'_j = -1$ when $1 \le j \le n$ and $j \not\in B$.
Thus
\begin{equation}
 2 \sum_{j \in B} \beta(j) \, f(j) = \sum_{j = 1}^n \epsilon_j \, f(j)
                                    + \sum_{j = 1}^n \epsilon'_j \, f(j),
\end{equation}
and hence
\begin{equation}
 2 \, \biggl\|\sum_{j \in B} \beta(j) \, f(j)\biggr\|
        = \biggl\|\sum_{j = 1}^n \epsilon_j \, f(j)\biggr\|
   + \biggl\|\sum_{j = 1}^n \epsilon'_j \, f(j)\biggr\| \le 2 \, \|f\|_{Z(V)}.
\end{equation}
This implies that
\begin{equation}
        \|f\|_{Z({\bf Z}_+, V)} \le \|f\|_{Z(V)},
\end{equation}
by taking the supremum over $B$, $\beta$.

        Let $E$ be any nonempty set again, and let $A$, $B$ be
disjoint nonempty finite subsets of $E$.  Also let $\alpha$, $\beta$
be functions on $A$, $B$, respectively, with values in $\{1, -1\}$.
Let $\gamma$, $\gamma'$ be the functions on $C = A \cup B$ defined by
$\gamma(x) = \gamma'(x) = \alpha(x)$ when $x \in A$ and $\gamma(x) =
-\gamma'(x) = \beta(x)$ when $x \in B$.  If $f(x)$ is any $V$-valued
function on $E$, then
\begin{equation}
 2 \sum_{x \in A} \alpha(x) \, f(x) = \sum_{x \in C} \gamma(x) \, f(x)
                                       + \sum_{x \in C} \gamma'(x) \, f(x)
\end{equation}
and
\begin{equation}
 2 \sum_{x \in B} \beta(x) \, f(x) = \sum_{x \in C} \gamma(x) \, f(x)
                                      - \sum_{x \in C} \gamma'(x) \, f(x).
\end{equation}
In particular,
\begin{equation}
 2 \, \biggl\|\sum_{x \in A} \alpha(x) \, f(x)\biggr\|
        \le \biggl\|\sum_{x \in C} \gamma(x) \, f(x)\biggr\|
           + \biggl\|\sum_{x \in C} \gamma'(x) \, f(x)\biggr\|,
\end{equation}
as in the preceding paragraph.

        Suppose now that $V$ is uniformly convex, and let $\epsilon >
0$ be given.  As in Section \ref{uniform convexity}, there is a
$\delta_1 > 0$ such that $\|v - w\| < \epsilon$ whenever $v, w \in V$
satisfy $\|v\|, \|w\| \le 1$ and $\|(v + w)/2\| > \delta_1$.
Equivalently, $\|v - w\| < \epsilon \, R$ when $\|v\|, \|w\| \le R$
and $\|(v + w)/2\| > (1 - \delta_1) \, R$ for any $R > 0$, by dividing
by $R$.  Let $f \in Y(E, V)$ with $f \not\equiv 0$ be given, and let
us apply this with $R = \|f\|_{Z(E, V)}$.  By definition of
$\|f\|_{Z(E, V)}$, there is a nonempty finite set $A \subseteq E$ and
a function $\alpha : A \to \{1, -1\}$ such that
\begin{equation}
        \biggl\|\sum_{x \in A} \alpha(x) \, f(x)\biggr\|
                                      > (1 - \delta_1) \, \|f\|_{Z(E, v)}.
\end{equation}
Let $B$ be another nonempty finite subset of $E$ that is disjoint from
$A$, and let $\beta$ be a function on $B$ with values in $\{1, -1\}$.
If $C$, $\gamma$, and $\gamma'$ are as in the previous paragraph and
\begin{equation}
        v = \sum_{x \in C} \gamma(x) \, f(x),
         \quad  w = \sum_{x \in C} \gamma'(x) \, f(x),
\end{equation}
then $\|v\|, \|w\| \le \|f\|_{Z(E, V)}$, and
\begin{equation}
        \biggl\|\frac{v + w}{2}\biggr\| 
          = \biggl\|\sum_{x \in A} \alpha(x) \, f(x)\biggr\|
           > (1 - \delta_1) \, \|f\|_{Z(E, V)}.
\end{equation}
Because of uniform convexity, we get that
\begin{equation}
        2 \biggl\|\sum_{x \in B} \beta(x) \, f(x)\biggr\|
         = \|v - w\| < \epsilon \, \|f\|_{Z(E, V)}.
\end{equation}
It follows that $\sum_{x \in E} f(x)$ satisfies the generalized Cauchy
condition, and hence converges in the generalized sense when $V$ is
also complete.

\section[\ Additional properties]{Additional properties}
\label{additional properties}
\setcounter{equation}{0}

        Let $E$ be a nonempty set, and let $V$ be a real or complex
vector space with a norm $\|v\|$.  If $a : E \to \{1, -1\}$ and $f \in
Y(E, V)$, then $a \, f \in Y(E, V)$, and in fact
\begin{equation}
\label{||f||_{Z(E, V)} = sup_a ||a f||_{Y(E, V)}}
        \|f\|_{Z(E, V)} = \sup_a \|a \, f\|_{Y(E, V)},
\end{equation}
where the supremum is taken over all such mappings $a$.  In particular,
\begin{equation}
        \|b \, f\|_{Z(E, V)} = \|f\|_{Z(E, V)}
\end{equation}
for every $f \in Y(E, V)$ and $b : E \to \{1, -1\}$.  If $f \in Y(E,
V)$ and $b$ is a bounded real-valued function on $E$, then $b \, f \in
Y(E, V)$, as in Section \ref{bounded coefficients}, and
\begin{equation}
\label{||b f||_{Z(E, V)} le ||b||_infty ||f||_{Z(E, V)}}
        \|b \, f\|_{Z(E, V)} \le \|b\|_\infty \, \|f\|_{Z(E, V)}.
\end{equation}
This follows from the analogous statement for the $Y(E, V)$ norm in
Section \ref{bounded coefficients} when $b$ is nonnegative, and
otherwise one can express $b$ as the product of a nonnegative function
and a function with values in $\{1, -1\}$.

        Suppose now that $V$ is a complex vector space, and let ${\bf
T}$ be the unit circle in the complex plane, consisting of the complex
numbers $z$ with $|z| = 1$.  If $a : E \to {\bf T}$ and $f \in Y(E,
V)$, then $a \, f \in Y(E, V)$ and
\begin{equation}
\label{||a f||_{Y(E, V)} le 4 ||f||_{Y(E, V)}}
        \|a \, f\|_{Y(E, V)} \le 4 \, \|f\|_{Y(E, V)},
\end{equation}
as in Section \ref{bounded coefficients}.  Put
\begin{equation}
        \|f\|_{W(E, V)} = \sup_a \|a \, f\|_{Y(E, V)},
\end{equation}
where the supremum is taken over all mappings $a : E \to {\bf T}$.
It is easy to see that this is a norm on $Y(E, V)$, and that
\begin{equation}
        \|f\|_{Y(E, V)} \le \|f\|_{W(E, V)} \le 4 \, \|f\|_{Y(E, V)}
\end{equation}
for every $f \in Y(E, V)$.  Equivalently,
\begin{equation}
        \|f\|_{W(E, V)}
         = \sup_{B, \beta} \, \biggl\|\sum_{x \in B} \beta(x) \, f(x)\biggr\|,
\end{equation}
where the supremum is taken over all nonempty finite sets $B \subseteq
E$ and functions $\beta : B \to {\bf T}$.

        More precisely, one can also check that
\begin{equation}
        \|f\|_{Z(E, V)} \le \|f\|_{W(E, V)} \le 2 \, \|f\|_{Z(E, V)}
\end{equation}
for every $f \in Y(E, V)$.  The first inequality follows from the
definitions and the fact that $1, -1 \in {\bf T}$.  The second
inequality uses the estimate
\begin{equation}
        \|a \, f\|_{Z(E, V)} \le 2 \, \|a\|_\infty \, \|f\|_{Z(E, V)}
\end{equation}
for every bounded complex-valued function $a$ on $E$ and $f \in Y(E,
V)$.  This follows from (\ref{||b f||_{Z(E, V)} le ||b||_infty
||f||_{Z(E, V)}}) applied to the real and imaginary parts of $a$.

        By construction,
\begin{equation}
        \|b \, f\|_{W(E, V)} = \|f\|_{W(E, V)}
\end{equation}
for every $f \in Y(E, V)$ and $b : E \to {\bf T}$.  If $b$ is a
bounded complex-valued function on $E$, then
\begin{equation}
\label{||b f||_{W(E, V)} le ||b||_infty ||f||_{W(E, V)}}
        \|b \, f\|_{W(E, V)} \le \|b\|_\infty \, \|f\|_{W(E, V)}
\end{equation}
for every $f \in Y(E, V)$.  In the case where $b$ is a bounded
nonnegative real-valued function on $E$, this follows from the
corresponding statement for the $Y(E, V)$ norm in Section \ref{bounded
coefficients}.  Otherwise, one can express $b$ as the product of a
nonnegative real-valued function and a function with values in ${\bf
T}$, to get the same conclusion from the previous two cases.

\section[\ Tori]{Tori}
\label{tori}
\setcounter{equation}{0}

        Let ${\bf T}$ be the unit circle in the complex plane,
as before.  It is well known that
\begin{equation}
        \int_{\bf T} z \, |dz| = 0,
\end{equation}
where $|dz|$ denotes the element of integration with respect to arc length.
One way to see this is to compare this integral with a line integral,
\begin{equation}
        \int_{\bf T} i \, z \, |dz| = \oint_{\bf T} dz = 0,
\end{equation}
using the fact that the unit tangent vector to ${\bf T}$ at a point $z
\in {\bf T}$ corresponds to $i \, z$ with respect to the standard
orientation.  Alternatively, one can use the change of variables $z
\mapsto -z$ to get that
\begin{equation}
        \int_{\bf T} z \, |dz| = - \int_{\bf T} z \, |dz|,
\end{equation}
and hence that the integral is $0$, because arc length is not affected
by this transformation.

        Of course, ${\bf T}$ is a compact Hausdorff topological space,
and a probability space with respect to arc length measure divided by
$2 \pi$.  As usual, the $n$-dimensional torus ${\bf T}^n$ is the
Cartesian product of $n$ copies of ${\bf T}$, consisting of ordered
$n$-tuples $z = (z_1, \ldots, z_n)$ with $z_j \in {\bf T}$ for $j = 1,
\ldots, n$.  This is also a compact Hausdorff topological space for
each $n$, and a probability space with respect to the corresponding
product measure.  The coordinate functions $z_1, \ldots, z_n$ may be
considered as complex-valued independent random variables on ${\bf
T}^n$.

        Similarly, we can consider the space ${\bf T}^\infty$ of
sequences $z = \{z_j\}_{j = 1}^\infty$ such that $z_j \in {\bf T}$ for
each $j$, which is the Cartesian product of a sequence of copies of
${\bf T}$.  This is a compact Hausdorff topological space with respect
to the product topology, and a probability space with respect to the
product measure.  The coordinate functions $z_1, z_2, \ldots$ form an
infinite sequence of independent random variables on this
infinite-dimensional torus, as before.  Note that the sequences $x =
\{x_j\}_{j = 1}^\infty$ with $x_j = 1$ or $-1$ for each $j$ form a
closed set in ${\bf T}^\infty$.

        Let $(V, \|v\|)$ be a complex Banach space, and let
$\{v_j\}_{j = 1}^\infty$ be a sequence of elements of $V$.  Consider
the $V$-valued functions
\begin{equation}
\label{f_n(z) = sum_{j = 1}^n z_j v_j}
        f_n(z) = \sum_{j = 1}^n z_j \, v_j
\end{equation}
on ${\bf T}^\infty$ for each $n \ge 1$.  If $\sum_{j \in {\bf Z}_+}
v_j$ converges in the generalized sense, then $\{f_n\}_{n = 1}^\infty$
converges uniformly on ${\bf T}^\infty$.  This is similar to the
discussion in Section \ref{uniform convergence}, using also the
estimates in Section \ref{bounded coefficients}, or the $W({\bf Z}_+,
V)$ norm in the previous section, which is basically the same.  The
converse statements discussed in Section \ref{uniform convergence} are
already applicable in this situation, because $1, -1 \in {\bf T}$.

\section[\ Norms and linear functionals]{Norms and linear functionals}
\label{norms, linear functionals}
\setcounter{equation}{0}

        Let $E$ be a nonempty set, and let $V$ be a real or complex
vector space with a norm $\|v\|$.  If $f \in Y(E, V)$ and $\lambda$ is
a bounded linear functional on $V$, then $\lambda(f(x))$ is a summable
function on $E$, as in Section \ref{sums, linear functionals}.  Put
\begin{equation}
        \|f\|_{L(E, V)} = \sup \Big\{\sum_{x \in E} |\lambda(f(x))| :
                                  \lambda \in V^*, \ \|\lambda\|_* \le 1\Big\}.
\end{equation}
As in Section \ref{sums, linear functionals}, this is less than or
equal to $2 \, \|f\|_{Y(E, V)}$ in the real case, less than or equal
to $4 \, \|f\|_{Y(E, V)}$ in the complex case, and greater than or
equal to $\|f\|_{Y(E, V)}$ in both cases.  It is easy to see from the
definition that $\|f\|_{L(E, V)}$ is a norm on $Y(E, V)$, and that
\begin{equation}
        \|b \, f\|_{L(E, V)} \le \|b\|_\infty \, \|f\|_{L(E, V)}
\end{equation}
for every $f \in Y(E, V)$ and bounded real or complex-valued function
$b$ on $E$, as appropriate.  If $B \subseteq E$ is a finite set,
$\beta$ is a real or complex-valued function on $B$ such that
$|\beta(x)| = 1$ for each $x \in B$, and $\lambda \in V^*$, then
\begin{equation}
        \biggl|\lambda\Big(\sum_{x \in B} \beta(x) \, f(x)\Big)\biggr|
         = \biggl|\sum_{x \in B} \beta(x) \, \lambda(f(x))\biggr|
          \le \sum_{x \in B} |\lambda(f(x))|,
\end{equation}
with equality in the last step for suitable choices of $\beta$.  Using
this, one can check that $\|f\|_{L(E, V)}$ is equal to $\|f\|_{Z(E,
V)}$ in the real case, and is equal to $\|f\|_{W(E, V)}$ in the
complex case.

\section[\ Sums and $c_0(E)$]{Sums and $c_0(E)$}
\label{sums, c_0(E)}
\setcounter{equation}{0}

        Let $E$ be a nonempty set, and let $(V, \|v\|)$ be a real or
complex Banach space.  If $f \in Y(E, V)$ and $a$ is a bounded real or
complex-valued function on $E$, as appropriate, then $a \, f \in Y(E,
V)$ and
\begin{equation}
        \|a \, f\|_{Y(E, V)} \le 2 \, \|a\|_\infty \, \|f\|_{Y(E, V)}
\end{equation}
in the real case, and
\begin{equation}
        \|a \, f\|_{Y(E, V)} \le 4 \, \|a\|_\infty \, \|f\|_{Y(E, V)}
\end{equation}
in the complex case, as in Section \ref{bounded coefficients}.  If $a
\in c_0(E)$, then it follows that $a \, f$ is in $Y_0(E, V)$, since
$a$ can be approximated by functions with finite support in the
$\ell^\infty$ norm.  This is the same as saying that $\sum_{x \in E}
a(x) \, f(x)$ satisfies the generalized Cauchy criterion when $a \in
c_0(E)$, and hence converges in the generalized sense because $V$ is
complete.  Thus
\begin{equation}
        T_f(a) = \sum_{x \in E} a(x) \, f(x)
\end{equation}
defines a bounded linear mapping from $c_0(E)$ into $V$.  One can
check that the operator norm of $T_f$ is equal to the $Z(E, V)$ norm
of $f$ in the real case, and is equal to the $W(E, V)$ norm of $f$ in
the complex case.  Conversely, if $T$ is a bounded linear mapping from
$c_0(E)$ in $V$, then $T = T_f$ for some $f \in Y(E, V)$.  To see
this, one can take
\begin{equation}
        f(x) = T(\delta_x),
\end{equation}
where $\delta_x$ is the function on $E$ defined by $\delta_x(x) = 1$
and $\delta_x(y) = 0$ when $x \ne y$.  If $a$ is a real or
complex-valued function on $E$ with finite support, then $a$ is a
linear combination of finitely many $\delta_x$'s, and so $T(a)$ is
given by the same expression as $T_f(a)$, because of linearity.  
Using this and the boundedness of $T$, one can show that $f \in Y(E,
V)$, and more precisely that the $Z(E, V)$ norm of $f$ is less than or
equal to the operator norm of $T$ in the real case, and that the $W(E,
V)$ norm of $f$ is less than or equal to the operator norm of $T$ in
the complex case.  This implies that $T(a) = T_f(a)$ for every $a \in
c_0(E)$, because $T$, $T_f$ are bounded linear operators which agree on
the dense linear subspace of $c_0(E)$ consisting of functions $a$ with
finite support.

\section[\ Integrability]{Integrability}
\label{integrability}
\setcounter{equation}{0}

        Let $(X, \mathcal{A}, \mu)$ be a measure space, and let $(V,
\|v\|)$ be a real or complex Banach space.  As in Section
\ref{sigma-subalgebras, vectors}, it is easy to deal with integration
of functions with values in a finite-dimensional subspace of $V$.
Suppose that $\{f_j\}_{j = 1}^\infty$ is a sequence of $V$-valued
functions on $X$ such that each $f_j$ takes values in a
finite-dimensional subspace of $V$, each $f_j$ is integrable in the
sense of Section \ref{sigma-subalgebras, vectors}, and
\begin{equation}
        \lim_{j, l \to \infty} \int_X \|f_j - f_l\| \, d\mu = 0.
\end{equation}
This implies in particular that the sequence of integrals
\begin{equation}
        \int_X f_j \, d\mu
\end{equation}
is a Cauchy sequence in $V$, and hence converges in $V$, by completeness.

        A sufficient condition for this type of convergence to hold is that
\begin{equation}
\label{sum_{j = 1}^infty int_X ||f_j - f_{j + 1}|| d mu < infty}
        \sum_{j = 1}^\infty \int_X \|f_j - f_{j + 1}\| \, d\mu < \infty.
\end{equation}
This is the same as
\begin{equation}
        \int_X \sum_{j = 1}^\infty \|f_j - f_{j + 1}\| \, d\mu < \infty,
\end{equation}
which implies that
\begin{equation}
        \sum_{j = 1}^\infty \|f_j(x) - f_{j + 1}(x)\| < \infty
\end{equation}
for almost every $x \in X$.  It follows that
\begin{equation}
        \sum_{j = 1}^\infty (f_j(x) - f_{j + 1}(x))
\end{equation}
converges in $V$ for almost every $x \in X$, by completeness again.
Put
\begin{equation}
        f(x) = \lim_{j \to \infty} f_j(x),
\end{equation}
which exists for almost every $x \in X$ by the convergence of the previous sum.
Of course, any sequence of $V$-valued functions as in the preceding paragraph
has a subsequence that satisfies this summability condition, and hence 
converges almost everywhere.

        Under these conditions, put
\begin{equation}
\label{int_X f d mu = lim_{j to infty} int_X f_j d mu}
        \int_X f \, d\mu = \lim_{j \to \infty} \int_X f_j \, d\mu.
\end{equation}
This is basically the definition of the Bochner integral.  Note that
$\{\|f_j(x)\|\}_{j = 1}^\infty$ converges in $L^1(X)$ to $\|f(x)\|$,
which implies that
\begin{equation}
\label{||int_X f d mu|| le int_X ||f|| d mu}
        \biggl\|\int_X f \, d\mu\biggr\| \le \int_X \|f\| \, d\mu.
\end{equation}
Similarly, if $\lambda$ is a bounded linear functional on $V$, then
$\lambda(f_j(x))$ converges in $L^1(X)$ to $\lambda(f(x))$, and hence
\begin{equation}
        \lambda\Big(\int_X f \, d\mu\Big) = \int_X \lambda \circ f \, d\mu.
\end{equation}
This shows that the integral of $f$ does not depend on the particular
sequence of approximations.

        Remember that a function on $X$ with values in a topological
space is said to be measurable if the inverse image of every open set
in the range is measurable.  Thus the composition of a measurable
function with a continuous mapping to another topological space is
also measurable.  If $f : X \to V$ is measurable with respect to the
topology on $V$ associated to the norm, then it follows that
$\|f(x)\|$ is measurable too.  If in addition $\|f(x)\|$ is integrable
and $V$ is separable, then $f$ can be approximated by integrable
functions with values in finite-dimensional subspaces of $V$, as
before.  To see this, one can start by using the integrability of
$\|f(x)\|$ to approximate $f$ by bounded measurable $V$-valued
functions that are equal to $0$ on the complements of suitable subsets
of finite measure.  One can then use the separability of $V$ to
approximate these functions by $V$-valued simple functions.  The same
argument would work if $f$ takes values in a separable subspace of $V$
almost everywhere on $X$.

\section[\ Bounded measures]{Bounded measures}
\label{bounded measures}
\setcounter{equation}{0}

        Let $X$ be a set, let $\mathcal{A}$ be an algebra of subsets
of $X$, and let $V$ be a real or complex vector space.  A $V$-valued
function $\mu$ on $\mathcal{A}$ is said to be a finitely-additive
$V$-valued measure on $(X, \mathcal{A})$ if
\begin{equation}
        \mu(A \cup B) = \mu(A) + \mu(B)
\end{equation}
for every $A, B \in \mathcal{A}$ with $A \cap B = \emptyset$.  
If $A_1, \ldots, A_n \in \mathcal{A}$ and $t_1, \ldots, t_n \in
{\bf R}$ or ${\bf C}$, as appropriate, then
\begin{equation}
        f(x) = \sum_{j = 1}^n t_j \, {\bf 1}_{A_j}(x)
\end{equation}
is a measurable simple function on $X$, and we put
\begin{equation}
        \int_X f \, d\mu = \sum_{j = 1}^n t_j \, \mu(A_j).
\end{equation}
It is easy to see that this does not depend on the particular
representation of $f$ as a linear combination of indicator functions,
and that it defines a linear mapping from the vector space of
measurable simple functions on $X$ into $V$.  If $\lambda$ is a linear
functional on $V$, then $\mu_\lambda(A) = \lambda(\mu(A))$ is a
finitely-additive real or complex measure on $(X, \mathcal{A})$, as
appropriate, and
\begin{equation}
        \lambda\Big(\int_X f \, d\mu\Big) = \int_X f \, d\mu_\lambda.
\end{equation}

        Suppose now that $V$ is equipped with a norm $\|v\|$, and that
$\mu$ is bounded, so that
\begin{equation}
        C(\mu) = \sup \{\|\mu(A)\| : A \in \mathcal{A}\} < \infty.
\end{equation}
If $A_1, \ldots, A_n$ are finitely many pairwise-disjoint measurable
subsets of $X$ and $E_n = \{1, \ldots, n\}$, then $\mu(A_j)$ may be
considered as a $V$-valued function on $E_n$ whose $Y(E_n, V)$ norm is
less than or equal to $C(\mu)$, because of the finite additivity of
$\mu$.  As in Section \ref{bounded coefficients}, it follows that
\begin{equation}
 \biggl\|\int_X f \, d\mu\biggr\| \le k \, C(\mu) \, \sup_{x \in X} |f(x)|
\end{equation}
for every measurable simple function $f$ on $X$, where $k = 1$ when
$f$ is real-valued and nonnegative, $k = 2$ when $f$ is real-valued,
and $k = 4$ when $f$ is complex-valued.  If $\mathcal{A}$ is a
$\sigma$-algebra and $V$ is complete, then the integral can be
extended to bounded measurable real or complex-valued functions $f$ on
$X$, as appropriate, because simple functions are dense in the space
of bounded measurable functions with respect to the supremum norm.  If
$\lambda$ is a bounded linear functional on $V$, then $\mu_\lambda$ is
a bounded finitely-additive real or complex measure on $(X,
\mathcal{A})$, with
\begin{equation}
        C(\mu_\lambda) \le \|\lambda\|_* \, C(\mu),
\end{equation}
and we get the same relationship with the integral of a bounded
measurable function as for simple functions.

        In particular, this works when $\mathcal{A}$ is a
$\sigma$-algebra and $\mu$ is countably additive, in the sense that
\begin{equation}
        \sum_{j = 1}^\infty \mu(A_j) = \mu\Big(\bigcup_{j = 1}^\infty A_j\Big)
\end{equation}
for every sequence $A_1, A_2, \ldots$ of pairwise-disjoint measurable
subsets of $X$, as in Section \ref{vector-valued measures}.  More
precisely, convergence of the series on the left in $V$ is part of the
hypothesis, and we have seen that this implies that $\mu$ is bounded.
In this case, $\mu_\lambda$ is a countably-additive real or complex
measure on $(X, \mathcal{A})$ for each $\lambda \in V^*$, as before.
If $\mu$ has the additional property that $\sum_{j = 1}^\infty
\|\mu(A_j)\|$ converges for every sequence $A_1, A_2, \ldots$ of
pairwise-disjoint measurable subsets of $X$, then one can integrate
any $f \in L^1(X, \|\mu\|)$, as in Section \ref{integrating vector
measures}.  If instead $V$ is a Hilbert space and $\mu(A)$ is
orthogonal to $\mu(B)$ when $A$, $B$ are disjoint measurable subsets
of $X$, then the integral can be defined on a suitable $L^2$ space, as
in Section \ref{measures, orthogonality}.

        As another situation like this, suppose that $V = W^*$ for
some Banach space $W$, $\mathcal{A}$ is a $\sigma$-algebra, and $\mu$
is countably additive with convergence in the weak$^*$ topology on
$V$.  This implies that $\mu_w(A) = \mu(A)(w)$ is a countably-additive
real or complex measure on $(X, \mathcal{A})$ for each $w \in W$.
This is the same as the measure $\mu_\lambda$ defined before, where
$\lambda$ is the bounded linear functional on $V$ corresponding to
evaluation at $w$.  Using the uniform boundedness principle, one can
show that $\mu$ is bounded, as in Section \ref{vector-valued
measures}.  Under these conditions, the integral of a bounded
measurable function $f$ on $X$ can be defined more directly as a
bounded linear functional on $W$ by
\begin{equation}
        \Big(\int_X f \, d\mu\Big)(w) = \int_X f \, d\mu_w,
\end{equation}
which is also satisfied by the previous definition.

\section[\ Weak$^*$ measurability]{Weak$^*$ measurability}
\label{weak^* measurability}
\setcounter{equation}{0}

        Let $W$ be a real or complex vector space with a norm $\|w\|_W$
which is separable, and let $w_1, w_2, \ldots$ be a sequence of
elements of $W$ such that $\|w_j\|_W = 1$ for each $j$ and the set of
$w_j$'s is dense in the unit sphere in $W$.  This uses the fact that a
subset of a separable metric space is also separable.  Note that
\begin{equation}
        \|\lambda\|_{W^*} = \sup_{j \ge 1} |\lambda(w_j)|
\end{equation}
for every bounded linear functional $\lambda$ on $W$.  Also let $(X,
\mathcal{A})$ be a measurable space, and let $f$ be a function on $X$
with values in $W^*$.  If $f(x)(w)$ is measurable as a real or
complex-valued function on $X$ for every $w \in W$, then it follows
that $\|f(x)\|_V$ is measurable on $X$ as well.

        Let us say that $f : X \to W^*$ is weak$^*$ measurable if $f$
is measurable with respect to the weak$^*$ topology on $W^*$.  This
automatically implies that $f(x)(w)$ is measurable for each $w \in W$,
since evaluation at $w$ is a continuous function on $W^*$.
Conversely, $f$ is weak$^*$ measurable when $f(x)(w)$ is measurable
for every $w \in W$ and $W$ is separable.  To see this, one may as
well suppose that $f$ is bounded, because one can use the
measurability of $\|f(x)\|_{W^*}$ to express $X$ as the union of a
sequence of measurable sets on which $f$ is bounded.  If $B$ is a ball
in $W^*$, then the topology on $B$ induced by the weak$^*$ topology on
$W^*$ is metrizable, because $W$ is separable, as in Section
\ref{seminorms, 2}.  If $B$ is a closed ball in $W^*$, then $B$ is
also compact in the weak$^*$ topology, by the Banach--Alaoglu theorem.
Thus $B$ is compact and metrizable with respect to the topology
induced by the weak$^*$ topology, and hence is separable with respect
to this topology.  This implies that relatively open subsets of $B$ in
the weak$^*$ topology can be given in terms of countable unions of
basic open sets, which permits the weak$^*$ measurability of $f$ to be
obtained from the measurability of $f(x)(w)$ for each $w \in W$.

        Of course, $f$ is weak$^*$ measurable if $f$ is measurable
with respect to the topology on $W^*$ associated to the dual norm,
because every open set in $W^*$ with respect to the weak$^*$ topology
is also open in the norm topology.  Conversely, if $f$ is weak$^*$
measurable and $W^*$ is separable, then $f$ is measurable with respect
to the norm topology on $W^*$.  Indeed, separability of $W^*$ implies
that each open set $U \subseteq W^*$ in the norm topology is a
countable union of closed balls.  If $B$ is a closed ball in $W^*$,
then $B$ is a closed set in $W^*$ in the weak$^*$ topology by the
definition of the dual norm, and so $f^{-1}(B)$ is measurable in $X$
by weak$^*$ measurability.  It follows that $f^{-1}(U)$ is the union
of countably many measurable subsets of $X$, and hence is measurable.

        Similarly, if $V$ is a real or complex vector space with a
norm $\|v\|_V$, then we say that $f : X \to V$ is weakly measurable if
$f$ is measurable with respect to the weak topology on $V$.  If $f$ is
measurable with respect to the topology on $V$ associated to the norm,
then $f$ is weakly measurable, because every open set in $V$ with
respect to the weak topology is also an open set in the norm topology.
Conversely, if $f$ is weakly measurable and $V$ is separable, then $f$
is measurable with respect to norm topology on $V$.  As before,
separability of $V$ implies that every open set $U \subseteq V$ in the
norm topology is the countable union of closed balls.  In this case,
the fact that a closed ball $B$ in $V$ is also closed in the weak
topology uses the Hahn--Banach theorem.  If $f$ is weakly measurable,
then it follows that $f^{-1}(B)$ is a measurable set in $X$ for each
closed ball $B$ in $V$, and hence that $f^{-1}(U)$ is measurable in
$X$ for each open set $U \subseteq V$ in the norm topology.  If $V^*$
is separable, then one can argue as before that $\|f(x)\|$ is
measurable on $X$ when $\lambda(f(x))$ is measurable for each $\lambda
\in V^*$.  The same argument shows that $\|f(x) - v\|$ is measurable
on $X$ for every $v \in V$ under these conditions, so that $f^{-1}(B)$
is measurable in $X$ for each ball $B$ in $V$.  One can then use
separability of $V$ again to get that $f$ is measurable with respect
to the norm topology on $V$.

\section[\ Weak$^*$ measures]{Weak$^*$ measures}
\label{weak^* measures}
\setcounter{equation}{0}

        Let $(X, \mathcal{A})$ be a measurable space, and let $(W,
\|w\|_W)$ be a real or complex Banach space.  Let us say that a
function $\mu$ on $\mathcal{A}$ with values in the dual $W^*$ of $W$
is a \emph{weak$^*$ measure} if
\begin{equation}
        \sum_{j = 1}^\infty \mu(A_j) = \mu\Big(\bigcup_{j = 1}^\infty A_j\Big)
\end{equation}
for every sequence $A_1, A_2, \ldots$ of pairwise-disjoint measurable
subsets of $X$, where the series is supposed to converge in the
weak$^*$ topology on $W^*$.  This is equivalent to asking that $\mu$
be finitely additive, and that
\begin{equation}
        \lim_{j \to \infty} \mu(B_j) = \mu\Big(\bigcup_{j = 1}^\infty B_j\Big)
\end{equation}
in the weak$^*$ topology for every increasing sequence $B_1, B_2,
\ldots$ of measurable subsets of $X$.  This is also equivalent to the
condition that $\mu$ be finitely additive and satisfy
\begin{equation}
        \lim_{j \to \infty} \mu(C_j) = \mu\Big(\bigcap_{j = 1}^\infty C_j\Big)
\end{equation}
in the weak$^*$ topology for every decreasing sequence $C_1, C_2,
\ldots$ of measurable subsets of $X$.  This is also the same as saying
that
\begin{equation}
        \mu_w(A) = \mu(A)(w)
\end{equation}
is a countably-additive real or complex measure on $(X, \mathcal{A})$,
as appropriate, for every $w \in W$.

        Remember that convergent sequences in $W^*$ in the weak$^*$
topology are bounded with respect to the dual norm when $W$ is
complete, by the theorem of Banach and Steinhaus.  If $\mu$ is a
weak$^*$ measure on $(X, \mathcal{A})$ with values in $W^*$, then
there is a $C \ge 0$ such that
\begin{equation}
\label{||mu(A)||_{W^*} le C}
        \|\mu(A)\|_{W^*} \le C
\end{equation}
for every $A \in \mathcal{A}$, by the same arguments as in Section
\ref{vector-valued measures}.  Equivalently,
\begin{equation}
        |\mu_w(A)| \le C \, \|w\|_W
\end{equation}
for every $w \in W$ and $A \in \mathcal{A}$, which implies that
\begin{equation}
\label{|mu_w|(X) le k C ||w||_W}
        |\mu_w|(X) \le k \, C \, \|w\|_W
\end{equation}
for every $w \in W$, where $k = 2$ in the real case and $k = 4$ in the
complex case.  Thus $w \mapsto \mu_w$ defines a bounded linear mapping
from $W$ into the space of real or complex measures on $(X,
\mathcal{A})$, as appropriate, equipped with the norm associated to
the total variation.  Conversely, a bounded linear mapping from $W$
into the space of real or complex measures on $(X, \mathcal{A})$
determines a weak$^*$ measure on $(X, \mathcal{A})$ with values in
$W^*$ in this way.

        If $\nu$ is a nonnegative measure on $(X, \mathcal{A})$ and $f
\in L^1(X, \nu)$, then
\begin{equation}
        \nu_f(A) = \int_A f \, d\nu
\end{equation}
defines a real or complex measure on $(X, \mathcal{A})$, as
appropriate.  Thus a bounded linear mapping from $W$ into $L^1(X,
\nu)$ determines a weak$^*$ measure $\mu$ on $(X, \mathcal{A})$ with
values in $W^*$, as in the previous paragraph.  In this case, $\mu$ is
absolutely continuous with respect to $\nu$, in the sense that $\mu(A)
= 0$ for every measurable set $A \subseteq X$ with $\nu(A) = 0$,
because $\nu_f$ is absolutely continuous with respect to $\nu$ for
every $f \in L^1(X, \nu)$.  Of course, any weak$^*$ measure $\mu$ on
$(X, \mathcal{A})$ with values in $W^*$ is absolutely continuous with
respect to $\nu$ in this sense if and only if $\mu_w$ is absolutely
continuous with respect to $\nu$ for each $w \in W$.  If $\nu$ is
$\sigma$-finite, then the Radon--Nikodym theorem implies that every
weak$^*$ measure $\mu$ on $(X, \mathcal{A})$ that is absolutely
continuous with respect to $\nu$ corresponds to a bounded linear
mapping from $W$ into $L^1(X, \nu)$.

        If $E$ is a nonempty set, then $Y(E, W^*)$ can be identified
with the space of bounded linear mappings from $W$ into $\ell^1(E)$.
This is basically another way of looking at the discussion in Section
\ref{sums in dual spaces}.  We can also think of $\ell^1(E)$ as being
the $L^1$ space associated to counting measure on $E$, so that
elements of $\ell^1(E)$ determine real or complex measures on $E$ as
in the preceding paragraph.  More precisely, these are measures
defined on arbitrary subsets of $E$.  It follows that elements of
$Y(E, W^*)$ determine bounded linear mappings from $W$ into real or
complex measures on $E$, as appropriate, and hence weak$^*$ measures
on $E$ with values in $W^*$.

\section[\ Weak$^*$ integrability]{Weak$^*$ integrability}
\label{weak^* integrability}
\setcounter{equation}{0}

        Let $(X, \mathcal{A}, \nu)$ be a measure space, let $W$ be a
real or complex vector space with a norm $\|w\|_W$.  Also let $f$ be a
$W^*$-valued function on $X$ such that $f(x)(w)$ is measurable on $X$
for each $w \in W$.  If $W$ is separable, then it follows that
$\|f(x)\|_{W^*}$ is measurable on $X$, as in Section \ref{weak^*
measurability}.  Alternatively, if $f : X \to W^*$ is measurable with
respect to the weak$^*$ topology on $W^*$, then we get that $f(x)(w)$
is measurable for each $w \in W$ and that $\|f(x)\|_{W^*}$ is
measurable.  The latter uses the fact that closed balls in $W^*$ are
closed sets in the weak$^*$ topology, by definition of the dual norm.

        At any rate, if $\|f(x)\|_{W^*}$ is integrable with respect to
$\nu$, then $f(x)(w)$ is also integrable with respect to $\nu$ for
each $w \in W$, and
\begin{equation}
 \int_X |f(x)(w)| \, d\nu(x) \le \|w\|_W \int_X \|f(x)\|_{W^*} \, d\nu(x)
\end{equation}
for every $w \in W$.  In particular, $w \mapsto f(x)(w)$ is a bounded
linear mapping from $W$ into $L^1(X, \nu)$, which leads to a weak$^*$
measure $\mu$ on $(X, \mathcal{A})$ with values in $W^*$, as in the
previous section.  More precisely,
\begin{equation}
        \mu_w(A) = \mu(A)(w) = \int_A f(x)(w) \, d\nu(x)
\end{equation}
for every measurable set $A \subseteq X$ and $w \in W$, which implies that
\begin{equation}
\label{||mu(A)||_{W^*} le int_A ||f(x)||_{W^*} d nu}
        \|\mu(A)\|_{W^*} \le \int_A \|f(x)\|_{W^*} \, d\nu.
\end{equation}
If $A_1, A_2, \ldots$ is a sequence of pairwise-disjoint measurable
subsets of $X$, then it is easy to see that $\sum_{j = 1}^\infty
\mu(A_j)$ converges absolutely with respect to the dual norm on $W^*$,
and that the sum is equal to $\mu\big(\bigcup_{j = 1}^\infty A_j\big)$.

        Let $\|\mu\|(A)$ be the total variation measure associated to
$\mu$ as in Section \ref{vector-valued measures}.  Thus $\|\mu\|(A) =
p^*(A)$ corresponds to $p(A) = \|\mu(A)\|_{W^*}$ as in Section
\ref{uniform boundedness, 3}.  In this case,
\begin{equation}
\label{||mu||(A) le int_A ||f(x)||_{W^*} d nu(x)}
        \|\mu\|(A) \le \int_A \|f(x)\|_{W^*} \, d\nu(x)
\end{equation}
for each $A \in \mathcal{A}$, because of (\ref{||mu(A)||_{W^*} le
int_A ||f(x)||_{W^*} d nu}).  Of course,
\begin{equation}
        |\mu_w(A)| \le \|\mu(A)\|_{W^*} \, \|w\|_W \le \|\mu\|(A) \, \|w\|_W
\end{equation}
for every $A \in \mathcal{A}$ and $w \in W$, which implies that
\begin{equation}
        |\mu_w|(A) \le \|\mu\|(A) \, \|w\|_W,
\end{equation}
where $|\mu_w|$ is the total variation measure associated to $\mu_w$.
Hence
\begin{equation}
        \int_A |f(x)(w)| \, d\nu(x) \le \|\mu\|(A) \, \|w\|
\end{equation}
for every $A \in \mathcal{A}$ and $w \in W$.

        Suppose that $W$ is separable, and let $w_1, w_2, \ldots$ be a
sequence of elements of $W$ such that $\|w_j\|_W = 1$ for each $j$ and
the set of $w_j$'s is dense in the unit sphere in $W$.  If
\begin{equation}
        \phi_n(x) = \max_{1 \le j \le n} |f_j(x)(w_j)|,
\end{equation}
then $\phi_n(x)$ is measurable on $X$ for each $n$,
\begin{equation}
        \phi_n(x) \le \phi_{n + 1}(x) \le \|f(x)\|_{W^*},
\end{equation}
for each $x \in X$ and $n \ge 1$, and
\begin{equation}
 \lim_{n \to \infty} \phi_n(x) = \sup_{n \ge 1} \phi_n(x) = \|f(x)\|_{W^*}
\end{equation}
for each $x \in X$.  Let $A \subseteq X$ be a measurable set, and let
$A_1, \ldots, A_n$ be pairwise-disjoint measurable subsets of $X$ such
that $\bigcup_{j = 1}^n A_j = A$.  Observe that
\begin{equation}
        \sum_{j = 1}^n \int_{A_j} |f(x)(w_j)| \, d\nu(x)
         \le \sum_{j = 1}^n \|\mu\|(A_j) = \|\mu\|(A).
\end{equation}
This implies that
\begin{equation}
        \int_A \phi_n(x) \, d\nu(x) \le \|\mu\|(A)
\end{equation}
for each $n$.  Using the monotone convergence theorem, we get that
\begin{equation}
        \int_A \|f(x)\|_{W^*} \, d\nu(x) \le \|\mu\|(A).
\end{equation}
It follows that
\begin{equation}
        \|\mu\|(A) = \int_A \|f(x)\|_{W^*} \, d\nu(x)
\end{equation}
for every $A \in \mathcal{A}$ when $W$ is separable.

\end{document}